\documentclass[11pt,letterpaper]{article}
\usepackage{ETfinal}

\begin{document}

\title{Proximal Estimation and Inference}

\makeatletter
\renewcommand{\@fnsymbol}[1]{\ensuremath{\ifcase#1\or *\or **\or ***\or ****\else\@ctrerr\fi}}
\makeatother

\author{Alberto Quaini\thanks{Corresponding author. Department of Econometrics,
Erasmus School of Economics, Erasmus University Rotterdam, and Tinbergen Institute. Email:
\href{mailto:Quaini@ese.eur.nl}{\texttt{Quaini@ese.eur.nl}}.}
\and
Fabio Trojani\thanks{Geneva Finance Research Institute, University of Geneva,
Swiss Finance Institute, and Collegio Carlo Alberto, University of Turin. Email:
\href{mailto:Fabio.Trojani@unige.ch}{\texttt{Fabio.Trojani@unige.ch}}.}}

\date{}
\maketitle
\thispagestyle{plain}

\begingroup
\renewcommand{\thefootnote}{}
\footnotetext{For useful comments, we thank the Co-Editors, three anonymous referees, Ming Yuan, Jianqing Fan, Elvezio Ronchetti, Chen Zhou, Olivier Scaillet, Whitney Newey, Patrick Gagliardini, Marco Avella Medina, Alberto Gonz\'alez-Sanz, Zoraida Fernandez Rico, Sofonias Alemu Korsaye, Federico Carlini, Nicola Gnecco, and Lorenzo Camponovo.}
\endgroup
\setcounter{footnote}{0}

\begin{abstract}
We develop a convex-analytic framework for constructing and analyzing a broad class of penalized estimators. Our method defines penalized estimators by applying proximal operators to benchmark estimators, rather than by penalizing the benchmark criterion itself. We establish general finite-sample properties, derive unified asymptotic distribution results, and provide tractable conditions for oracle selection and oracle asymptotic efficiency. We then apply the framework to linear regression with irregular designs, including singular and nearly-singular cases. In this setting, where the coefficient vector need not be point identified, we motivate the population Ridgeless estimand as the inferential target, construct a \(\sqrt{n}\)-consistent modified Ridgeless benchmark estimator, and combine it with adaptive proximal penalization to obtain, to the best of our knowledge, the first penalized estimator achieving the oracle property under both regular and irregular designs.
\end{abstract}

\noindent\textbf{Keywords:} proximal operator; convex penalization; oracle property; asymptotic theory; ridgeless estimator; irregular designs.

\noindent\textbf{JEL classification:} C13, C18, C51, C52.

\clearpage
\section{Introduction}

Penalization is a standard device in econometrics and statistics for stabilizing estimation, shrinking nuisance components, and enforcing low-dimensional structure such as sparsity; classical examples include Ridge, Lasso, SCAD, and related procedures in least-squares and likelihood settings \citep{hoerl1970ridge,tibshirani1996regression,fan2001variable,yuan2006model}. 
The usual approach is to append a penalty directly to the criterion that defines the estimator of interest, as in penalized least squares, penalized likelihood, or penalized GMM \citep{green1987penalized,frank1993statistical,caner2009lasso}. 
That approach is natural and often effective, but it also ties regularization to the specific global geometry of the criterion being optimized. 
As a consequence, existence, uniqueness, finite-sample stability, and asymptotic behavior typically have to be established for each penalized problem in a model- and criterion-specific way.

This paper develops a different approach.
Instead of regularizing by modifying the original criterion, we start from a benchmark estimator and then define a penalized estimator by applying a single proximal regularization step to it.
Throughout the paper, the parameter dimension \(p\) is fixed.
Given a benchmark estimator \(\hat{\bm\beta}_n^s\in\R^p\), a symmetric positive definite matrix \(\bm W_n\), a convex penalty function \(f_n\), and a tuning parameter \(\lambda_n>0\), we consider estimators of the form
\begin{equation}\label{eq:intro_prox}
\hat{\bm\beta}_n
:=
\argmin_{\bm\beta\in\R^p}
\left\{
\frac12\|\bm\beta-\hat{\bm\beta}_n^s\|_{\bm W_n}^2
+\lambda_n f_n(\bm\beta)
\right\},
\end{equation}
where \(\|\bm v\|_{\bm W_n}^2:=\bm v'\bm W_n \bm v\).
We call \(\hat{\bm\beta}_n\) a \emph{proximal estimator}.
In this construction, \(\hat{\bm\beta}_n^s\) carries the model-specific statistical content, \(\bm W_n\) determines the local geometry of the regularization step, and \(f_n\) determines the form of shrinkage, selection, or constraint enforcement.

The point of view in \eqref{eq:intro_prox} is not meant to replace penalize-the-criterion methods.
Rather, it provides a complementary and modular formulation of regularization.
The benchmark step and the regularization step are separated.
Once \(\hat{\bm\beta}_n^s\) is available, the proximal step is a strongly convex optimization problem under the sole requirement \(\bm W_n\succ\bm 0\), and is therefore globally well posed.
This separation turns out to be particularly useful.
It first produces desirable finite-sample regularity properties independently of the choice of benchmark estimator \(\hat{\bm\beta}_n^s\). Second, it
enables a unified asymptotic analysis for the regularized estimator, which seamlessly combines the limit law of the benchmark with the local geometry of the penalty.
In this sense, the proximal route isolates regularization from the model-specific geometry of the benchmark problem while preserving a transparent link between the two.

The connection with familiar penalized estimators is immediate in regular linear models.
When the benchmark is least squares and the proximal metric is chosen to be the sample design matrix, the proximal step reproduces the usual convex penalized least-squares estimator.
Thus, in standard settings, the proximal formulation does not create an unrelated class of procedures; it recovers familiar estimators from a different angle.
At the same time, the same construction remains meaningful well beyond least squares.
For instance, in nonlinear or multi-step problems such as two-step GMM, one may regularize an econometrically natural benchmark estimator without altering the original criterion, and hence without coupling the nonsmooth penalty to the full nonlinear objective.
This distinction is especially useful when one wants regularization to inherit the benchmark's local geometry while keeping the regularization step globally well-defined.

The paper develops a general convex-analysis framework for estimators of the form \eqref{eq:intro_prox}.
Our first contribution is to establish their finite-sample properties.
We show that proximal estimators are always uniquely defined under a positive definite metric, satisfy deterministic stability bounds, and admit equivalent primal and dual characterizations based on convex conjugacy.
For sublinear penalties, the dual characterization reduces to a projection-residual decomposition onto a closed convex set.
This geometry yields a transparent interpretation of shrinkage and selection and shows that familiar formulas, such as soft-thresholding for Lasso-type procedures under orthonormal designs, are special cases of a more general projection structure.

Our second contribution is asymptotic.
We derive a general limit theory for proximal estimators under a high-level assumption on the benchmark estimator and mild regularity conditions on the metric and the penalty sequence. All these asymptotic results are derived in the classical large-\(n\), fixed-\(p\) regime.
The resulting limit distributions are characterized in closed form by three objects only: the asymptotic distribution of the benchmark estimator, the local limit geometry of the penalty, and the limiting metric of the proximal step.
This framework nests classical asymptotic results for convex penalized least-squares estimators, including the benchmark results of \citet{knight2000asymptotics}, but it also applies in settings where the benchmark is not least squares and where penalize-the-criterion analyses require a dedicated local expansion of the penalized objective.
We further derive oracle results, in the sense of \citet{fan2001variable}, for proximal estimators.
The conditions are formulated directly in terms of the local behavior of the penalty along active and inactive coordinates and are therefore easy to verify for standard examples.
A useful implication is that the benchmark estimator and the penalty play conceptually distinct roles: the former supplies the local stochastic limit, while the latter determines whether shrinkage produces consistent selection and oracle-type asymptotic behavior.

These asymptotic and oracle results are especially transparent in dual form.
Through Moreau decomposition, the proximal step can be written as the benchmark estimator minus a dual correction determined by the conjugate penalty; for sublinear penalties, this correction is a projection onto a closed convex set.
This viewpoint isolates the mechanism through which selection occurs.
Exact zeros arise through dual feasibility: oracle selection fails when the dual correction cannot absorb the inactive benchmark fluctuation, and it holds when the local subdifferential geometry generates sufficient separation on inactive coordinates.
In this sense, the dual formulation turns support recovery and oracle efficiency into geometric questions about projections, dual sets, and local subgradients, rather than penalty-specific algebra.

Our third contribution concerns linear regression under \emph{irregular} designs, by which we mean that the population design matrix may be singular and the sample design may either inherit this singularity or approach it asymptotically.
This is the setting in which the proximal formulation is particularly useful.
Under a regular design, least squares and the corresponding penalized least-squares estimators are well behaved.
Under a singular or nearly-singular design, by contrast, the population normal equations need not identify a unique coefficient vector, and standard least-squares-type estimators may exhibit slower convergence rate, or otherwise nonstandard first-order behavior.
For Bridge-type penalized least-squares estimators under nearly-singular designs, slower rates and nonstandard limits are documented in \citet{knight2000asymptotics} and \citet{knight2008shrinkage}. Related irregularity also affects GMM and GEL procedures, as shown by \citet{caner2008nearly}, and more broadly motivates singularity-robust inference in moment condition models; see \citet{andrews2019identification}.

A central point, emphasized throughout our analysis, is that when the population design is singular, the coefficient vector of the linear model is generally not point identified. Accordingly, the paper does not claim regular inference for an unidentified coefficient vector. Instead, in the irregular-design problem we target the Ridgeless estimand, defined as the unique minimizer of the limiting population least-squares criterion over
the range of the limit design matrix. By eliminating the unidentified kernel component, this
target is the unique least-squares representative determined by the limiting
least-squares objects alone and, among linear selections from the identified set,
the unique one compatible with \(\sqrt n\)-regular inference without an auxiliary
normalization.

Building on this target, we construct a new modified Ridgeless benchmark estimator based on a rank-consistent estimate of the population design matrix. This estimator is \(\sqrt n\)-consistent and asymptotically normal for the Ridgeless target under regular, singular, and nearly-singular designs. Applying suitable proximal penalties to this benchmark then yields penalized estimators that remain well behaved across all three regimes and, under suitable adaptive penalties, can satisfy the oracle property of \citet{fan2001variable}.

The paper therefore has two intertwined messages. At a general level, it shows that a single proximal step applied to a benchmark estimator yields a tractable and flexible framework for penalized estimation and inference. At an applied theoretical level, it shows that this modular route is powerful enough to address a concrete problem in linear regression with irregular designs: one first stabilizes the benchmark estimator for a meaningful and identified least-squares target, and then introduces shrinkage through a proximal penalty while preserving a transparent first-order theory.

The remainder of the paper is organized as follows.
Section~\ref{sec:framework} develops the proximal framework, relates it to penalized least squares, and introduces the basic properties that drive the subsequent analysis.
Section~\ref{sec:asymptotics} derives the asymptotic distributions and oracle properties of proximal estimators.
Section~\ref{sec:irregular_designs} applies the general theory to linear regression under regular, singular, and nearly-singular designs, by introducing the modified Ridgeless benchmark and the associated oracle proximal estimators.
The convex-analysis background, supporting results, Monte Carlo analysis, and proofs are provided in the Supplementary Material (SM).

\section{Proximal estimation framework}\label{sec:framework}

Throughout the paper, \(n\) and \(p\) denote positive integers, with \(p\) fixed across the asymptotic analysis.
Real-valued vectors are written in boldface, for example \(\bm x,\bm y,\bm\beta\), and real-valued matrices in bold capitals, for example \(\bm X,\bm W,\bm Q\).
All random elements are defined on a complete probability space \((\Omega,\mathcal F,\Pr)\), and all stochastic convergences are understood with respect to \(\Pr\).
We use standard notation from convex analysis, summarized in SM Section~A.1.

Let \(\Gamma(\R^p)\) denote the class of proper, lower semicontinuous, convex functions from \(\R^p\) to \((-\infty,+\infty]\).\footnote{Allowing extended real values lets the framework cover both standard convex penalties and convex constraints through indicator functions.}
Throughout the paper, a proximal estimator is an estimator of the form \eqref{eq:intro_prox}, based on a benchmark estimator \(\hat{\bm\beta}_n^s\in\R^p\), a penalty \(f_n\in\Gamma(\R^p)\), a tuning parameter \(\lambda_n>0\), and a symmetric weight matrix \(\bm W_n\succ\bm 0\). Equivalently, 
\begin{equation}\label{eq:prox_def}
\hat{\bm\beta}_n
=
\prox_{\lambda_n f_n}^{\bm W_n}(\hat{\bm\beta}_n^s),
\end{equation}
where, for each \(\bm\theta\in\R^p\),
\begin{equation}\label{eq:prox_op_def}
\prox_{\lambda_n f_n}^{\bm W_n}(\bm\theta)
:=
\argmin_{\bm\beta\in\R^p}
\left\{
\frac12\|\bm\beta-\bm\theta\|_{\bm W_n}^2+\lambda_n f_n(\bm\beta)
\right\}
\end{equation}
denotes the proximal operator of \(\lambda_n f_n\) with respect to the inner product \(\langle\cdot,\cdot\rangle_{\bm W_n}\).\footnote{
    Proximal operators are also central in the literature on proximal algorithms in statistics
    and machine learning; see, for example, \citet{polson2015proximal}. That literature
    typically uses the proximal map as an iterative computational device for solving a
    pre-specified regularized optimization problem. Our use is different: the proximal map is
    applied once to a benchmark estimator and is part of the statistical definition of the
    estimator. The focus is consequently on well-posedness, first-order sampling theory, and
    oracle properties, rather than on algorithmic convergence.
}
This section relates the construction to penalized least squares, illustrates its use beyond linear models, and introduces basic properties that will play a central role in the sequel.

\subsection{Link to penalized least squares}\label{sec:link_plse}

We begin with linear regression, where the proximal construction can be compared directly with familiar penalized least-squares estimators.
Consider the model
\begin{equation}\label{linear_model_setting}
\bm Y=\bm X\bm\beta_0+\bm\varepsilon,
\end{equation}
where \(\bm X\in\R^{n\times p}\) is a possibly random design matrix and \(\bm\varepsilon\in\R^n\) is an error vector.
Given \(\lambda_n>0\) and \(f_n:\R^p\to(-\infty,+\infty]\), define the associated penalized least-squares estimator (PLSE) by
\begin{equation}\label{eq:plse_def}
\hat{\bm\beta}_n^{\mathrm{plse}}
\in
\argmin_{\bm\beta\in\R^p}
\left\{
\frac{1}{2n}\|\bm Y-\bm X\bm\beta\|_2^2+\lambda_n f_n(\bm\beta)
\right\}.
\end{equation}
This class includes many penalties from the literature: Ridge, Lasso, Adaptive Lasso, Group Lasso, Elastic Net, and indicator penalties corresponding to constrained least squares. We refer to these as our benchmark penalties and define them in SM Section~B.1.

When \(f_n\in\Gamma(\R^p)\), the criterion in \eqref{eq:plse_def} is convex, but uniqueness need not hold under singular designs or non-strictly convex penalties.
By contrast, the proximal estimator in \eqref{eq:prox_def} is always uniquely defined once \(\bm W_n\succ\bm 0\).
Under a regular design, the two constructions coincide exactly.

\begin{proposition}[PLSEs as proximal estimators under regular designs]\label{prop:prox_plse_regular}
In model~\eqref{linear_model_setting}, suppose
\(\bm Q_n:=\bm X'\bm X/n\succ\bm 0\).
Fix \(\lambda_n>0\) and \(f_n\in\Gamma(\R^p)\), and let
\begin{equation}\label{eq:ls_estimator_regular}
\hat{\bm\beta}_n^{\mathrm{ols}}:=\bm Q_n^{-1}\bm X'\bm Y/n.
\end{equation}
Then \(\hat{\bm\beta}_n^{\mathrm{plse}}\) in \eqref{eq:plse_def} is unique and satisfies
\[
\hat{\bm\beta}_n^{\mathrm{plse}}
=
\prox_{\lambda_n f_n}^{\bm Q_n}\big(\hat{\bm\beta}_n^{\mathrm{ols}}\big)
=
\argmin_{\bm\beta\in\R^p}
\left\{
\frac12\|\bm\beta-\hat{\bm\beta}_n^{\mathrm{ols}}\|_{\bm Q_n}^2+\lambda_n f_n(\bm\beta)
\right\}.
\]
\end{proposition}
Thus, under regular designs, every convex PLSE is a proximal regularization of least squares under the design metric \(\bm Q_n\).
In particular, the proximal formulation is not an alternative to familiar penalized least-squares procedures in such settings; it is an equivalent representation of them. When \(\bm Q_n\) is singular,
\(\|\cdot\|_{\bm Q_n}\) is only a seminorm rather than a norm. As a result, the proximal problem may not be well defined, due to the flatness of the quadratic term along \(\Kernel(\bm Q_n)\) in equation (\ref{eq:intro_prox}).
In these settings, we show that it is convenient
to additionally penalize
parametric directions in \(\Kernel(\bm Q_n)\).
Specifically, we consider the
Ridgeless estimator \citep{hastie2022surprises}:
\begin{equation}\label{eq:Ridgeless_est}
\hat{\bm\beta}_n^{\mathrm{rls}}
:=
\argmin_{\bm\beta\in\R^p}
\Big\{\|\bm\beta\|_2^2:\ \bm X'\bm X\bm\beta=\bm X'\bm Y\Big\}
=
\bm Q_n^+\bm X'\bm Y/n ,
\end{equation}
where $\bm Q_n^+$ is the  Moore--Penrose inverse of $\bm Q_n$, and define
\begin{equation}\label{eq:Qbar_def}
\overline{\bm Q}_n
:=
\bm Q_n+\big(\bm I_p-\bm Q_n\bm Q_n^+\big).
\end{equation}
Note that since \(\bm Q_n\bm Q_n^+\) is by construction the Euclidean orthogonal projector onto \(\Range(\bm Q_n)=\Range(\bm X')\),
\(\bm I_p-\bm Q_n\bm Q_n^+\) is the Euclidean orthogonal projector onto \(\Kernel(\bm Q_n)=\Kernel(\bm X)\).
Hence \(\overline{\bm Q}_n\) coincides with \(\bm Q_n\) on \(\Range(\bm Q_n)\), while adding an identity component on \(\Kernel(\bm Q_n)\),
thus penalizing parametric directions in
$\Kernel(\bm Q_n)$.

\begin{proposition}[Ridgeless PLSE representation of proximal estimators]\label{prop:prox_plse_Ridgeless}
In model~\eqref{linear_model_setting}, let \(\bm Q_n:=\bm X'\bm X/n\), fix \(\lambda_n>0\) and \(f_n\in\Gamma(\R^p)\), and define \(\hat{\bm\beta}_n^{\mathrm{rls}}\) and \(\overline{\bm Q}_n\) as in \eqref{eq:Ridgeless_est}--\eqref{eq:Qbar_def}.
Let
\begin{equation}\label{eq:extended_penalty_Ridgeless}
\overline f_n(\bm\beta)
:=
f_n(\bm\beta)
+\frac{1}{2\lambda_n}\bm\beta'
\big(\bm I_p-\bm Q_n\bm Q_n^+\big)\bm\beta.
\end{equation}
Then \(\overline{\bm Q}_n\succ\bm 0\), \(\overline f_n\in\Gamma(\R^p)\), and
\begin{equation}\label{eq:plse_prox_identity_Ridgeless}
\prox_{\lambda_n f_n}^{\overline{\bm Q}_n}\big(\hat{\bm\beta}_n^{\mathrm{rls}}\big)
=
\argmin_{\bm\beta\in\R^p}
\left\{
\frac{1}{2n}\|\bm Y-\bm X\bm\beta\|_2^2+\lambda_n \overline f_n(\bm\beta)
\right\}.
\end{equation}
\end{proposition}

Proposition~\ref{prop:prox_plse_Ridgeless} shows that singular designs can be handled by injecting curvature only along \(\Kernel(\bm X)\), i.e., along directions that do not change fitted values.
The quadratic adjustment in \eqref{eq:extended_penalty_Ridgeless} is convex and acts only on observationally irrelevant directions. It therefore removes indeterminacy along \(\Kernel(\bm X)\) while leaving the fitted-value component governed by the usual penalized least-squares tradeoff.
When \(\bm Q_n\succ\bm 0\), the projector \(\bm I_p-\bm Q_n\bm Q_n^+\) vanishes, \(\overline{\bm Q}_n=\bm Q_n\), \(\overline f_n=f_n\), and \eqref{eq:plse_prox_identity_Ridgeless} reduces to Proposition~\ref{prop:prox_plse_regular}.
A general PLSE representation based on an arbitrary positive semidefinite metric is given in SM Section~D.1.

\subsection{Beyond least squares}\label{sec:beyond_least_squares}

The proximal formulation is not restricted to least squares.
More generally, one may start from any econometrically natural benchmark estimator and regularize it through a separate proximal step.
This distinction is especially useful in nonlinear or multi-stage problems, where penalizing the original criterion directly couples the penalty to model-specific curvature and may complicate both computation and analysis.
The following example will serve as a running illustration.

\begin{example}[Penalized optimal two-step GMM in nonlinear moment models]\label{ex:penalized_2sgmm}
Let \(w_i\in\mathcal W\subseteq \R^{p'}\) denote the data, \(\bm\theta_0\in\R^d\) the structural parameter, and consider moment functions
\[
g_q(\cdot,\cdot):\mathcal W\times\R^d\to\R^q,
\qquad
g_k(\cdot,\cdot):\mathcal W\times\R^d\to\R^k,
\]
measurable in \(w\) and continuously differentiable in \(\bm\theta\).
Assume the first set is correctly specified, while the second may be misspecified:
\[
\E[g_q(w_i,\bm\theta_0)]=\bm 0,
\qquad
\E[g_k(w_i,\bm\theta_0)-\bm\gamma_0]=\bm 0,
\]
where \(\bm\gamma_0\in\R^k\) is a misspecification parameter.
Stack the parameters as
\[
\bm\beta_0:=(\bm\theta_0',\bm\gamma_0')'\in\R^p,
\qquad p:=d+k,
\]
and define
\[
 m(w,\bm\beta)
:=
\begin{pmatrix}
g_q(w,\bm\theta)\\
g_k(w,\bm\theta)-\bm\gamma
\end{pmatrix},
\qquad
\bm\beta=(\bm\theta,\bm\gamma),
\]
so that \(\E[m(w_i,\bm\beta_0)]=\bm 0\). Here, sparsity of \(\bm\gamma_0\) provides a convenient characterization of different levels of GMM misspecification and a natural starting point for moment selection.
Let \(\bar m_n(\bm\beta):=n^{-1}\sum_{i=1}^n m(w_i,\bm\beta)\).
Efficient estimation of \(\bm\beta_0\) is then obtained via the optimal two-step GMM estimator \(\hat{\bm\beta}_n^{\mathrm{2SGMM}}\), defined by
\[
\hat{\bm\beta}_n^{(1)}
\in
\argmin_{\bm\beta\in\R^p}\ \bar m_n(\bm\beta)'\bar m_n(\bm\beta),
\qquad
\widehat{\bm\Omega}_n
:=
\frac1n\sum_{i=1}^n
 m(w_i,\hat{\bm\beta}_n^{(1)}) m(w_i,\hat{\bm\beta}_n^{(1)})',
\]
and
\begin{equation}\label{eq:2sgmm}
\hat{\bm\beta}_n^{\mathrm{2SGMM}}
\in
\argmin_{\bm\beta\in\R^p}
\bar m_n(\bm\beta)'\widehat{\bm\Omega}_n^{-1}\bar m_n(\bm\beta).
\end{equation}
A standard approach in this setting, irrespective of moment selection, is to penalize the GMM criterion directly (\citet{caner2009lasso} and \citet{liao2013adaptive}).
In the present stacked formulation,
moment selection can be achieved with
an Adaptive Lasso penalty acting on the misspecification block.
Given a preliminary estimator \(\tilde{\bm\beta}_n=(\tilde{\bm\theta}_n,\tilde{\bm\gamma}_n)\), define:\footnote{With this definition, we cover situations where some components of the estimator $\tilde {\bm \gamma}_n$ may be exactly zero. Thus, whenever a preliminary coefficient \(\tilde\gamma_{nj}\) is zero, the corresponding coordinate of the resulting penalized estimator is hard-constrained to be exactly zero.} 
\begin{equation}\label{eq:adaptive_lasso_eta}
f_n(\bm\beta)
:=
\sum_{j:\tilde\gamma_{nj}\neq 0}\frac{|\gamma_j|}{|\tilde\gamma_{nj}|}
+
\sum_{j:\tilde\gamma_{nj}=0}\iota_{\{0\}}(\gamma_j),
\end{equation}
where \(\iota_C\) denotes the indicator function of \(C\subseteq\R^p\), defined by
\(\iota_C(\bm x)=0\) for \(x\in C\) and \(\iota_C(\bm x)=+\infty\) otherwise.
The corresponding penalize-the-criterion GMM estimator allowing for adaptive moment selection is
\begin{equation}\label{eq:pen_gmm_adal}
\tilde{\bm\beta}_n^{\mathrm{pen}}
\in
\argmin_{\bm\beta\in\R^p}
\Big\{
\bar m_n(\bm\beta)'\widehat{\bm\Omega}_n^{-1}\bar m_n(\bm\beta)
+\lambda_n f_n(\bm\beta)
\Big\}.
\end{equation}
Our proximal approach instead regularizes the benchmark estimator \(\hat{\bm\beta}_n^{\mathrm{2SGMM}}\) in the natural metric induced by the benchmark's local GMM geometry.
Define
\[
\widehat{\bm G}_n
:=
\frac{\partial \bar m_n(\bm\beta)}{\partial\bm\beta'}
\bigg|_{\bm\beta=\hat{\bm\beta}_n^{\mathrm{2SGMM}}},
\qquad
\bm W_n
:=
\widehat{\bm G}_n'\widehat{\bm\Omega}_n^{-1}\widehat{\bm G}_n.
\]
If needed, one may enforce positive definiteness by replacing \(\bm W_n\) with \(\bm W_n+\kappa_n\bm I_p\), with \(\kappa_n\downarrow0\).
The associated proximal GMM--Adaptive Lasso estimator is
\begin{equation}\label{eq:running_example_2sgmm_adal}
\hat{\bm\beta}_n
=
\prox_{\lambda_n f_n}^{\bm W_n}\big(\hat{\bm\beta}_n^{\mathrm{2SGMM}}\big)
=
\argmin_{\bm\beta\in\R^p}
\left\{
\frac12\|\bm\beta-\hat{\bm\beta}_n^{\mathrm{2SGMM}}\|_{\bm W_n}^2
+\lambda_n f_n(\bm\beta)
\right\}.
\end{equation}
Estimators \eqref{eq:pen_gmm_adal} and \eqref{eq:running_example_2sgmm_adal} use the same adaptive weights, the same tuning parameter, and induce sparsity in the same parameter block.
They differ only in how regularization enters the procedure. Indeed,
the former estimator modifies the nonlinear GMM criterion itself, whereas the latter leaves the benchmark criterion unchanged and applies a single proximal step with the efficient GMM metric tied to the benchmark's asymptotic covariance matrix.
This separation between benchmark construction and regularization is the central principle underlying our approach.
\end{example}

\subsection{Dual representation and projection geometry}

A useful feature of the proximal formulation is that it admits a dual characterization, and for sublinear penalties this characterization simplifies further to a projection formula. Let \(\Id:\R^p\to\R^p\) denote the identity map. Given \(g\in\Gamma(\R^p)\), define its convex conjugate under the inner product \(\langle\cdot,\cdot\rangle_{\bm W_n}\) by
\begin{equation}\label{eq:conjugate}
g^*(\bm u)
:=
\sup_{\bm\beta\in\R^p}
\big\{
\langle \bm u,\bm\beta\rangle_{\bm W_n}-g(\bm\beta)
\big\},
\qquad
\bm u\in\R^p.
\end{equation}
Recall that a function \(f:\R^p\to(-\infty,+\infty]\) is sublinear if it is positively homogeneous and subadditive; equivalently, \(f\) is the support function of a nonempty closed convex set.
\begin{proposition}[Dual representation and projection form]\label{prop:dual_repr}
The proximal operator~\eqref{eq:prox_op_def} satisfies
\begin{equation}\label{eq:moreau_W}
\prox_{\lambda_n f_n}^{\bm W_n}
=
\Id-\prox_{(\lambda_n f_n)^*}^{\bm W_n},
\end{equation}
where \((\lambda_n f_n)^*\) and \(f_n^*\) denote Fenchel conjugates with respect to the inner product \(\langle\cdot,\cdot\rangle_{\bm W_n}\). In particular,
\[
(\lambda_n f_n)^*(\bm\theta)
=
\lambda_n f_n^*(\bm\theta/\lambda_n),
\qquad
\bm\theta\in\R^p.
\]
Moreover, if \(f_n\) is sublinear, then \((\lambda_n f_n)^*=\iota_{C_n}\), where
\[
C_n
=
\bigcap_{\bm\beta\in\R^p}
\left\{
\bm\theta\in\R^p:\ 
\langle \bm\theta,\bm\beta\rangle_{\bm W_n}\le \lambda_n f_n(\bm\beta)
\right\},
\]
and therefore
\begin{equation}\label{eq:proj_formula_W}
\prox_{\lambda_n f_n}^{\bm W_n}
=
\Id-P_{C_n}^{\bm W_n},\qquad P_{C_n}^{\bm W_n}(\bm\theta)
:=
\argmin_{\tilde{\bm\theta}\in C_n}
\|\bm\theta-\tilde{\bm\theta}\|_{\bm W_n}.
\end{equation}
\end{proposition}
Proposition~\ref{prop:dual_repr} shows that proximal regularization can be written as the benchmark input minus a dual correction determined by the conjugate penalty. For sublinear penalties, the conjugate collapses to an indicator function, so the dual correction is simply the \(\bm W_n\)-metric projection onto \(C_n\). The proximal step therefore decomposes \(\bm\theta\) into a dual-feasible component \(P_{C_n}^{\bm W_n}(\bm\theta)\) and the residual \(\bm\theta-P_{C_n}^{\bm W_n}(\bm\theta)\), which is the regularized estimator. This representation is often more informative than the corresponding optimality conditions in the primal proximal estimation problem, because it identifies directly what the penalty subtracts from the benchmark estimator. It will also be useful below in the asymptotic and oracle analysis. SM Section~B.1 reports the conjugates of our benchmark penalties. For the sublinear cases, these conjugates are indicator functions, and the dual representation therefore reduces to the projection formula \eqref{eq:proj_formula_W}.
A particularly important case is provided by weighted \(\ell_1\) penalties.

\begin{corollary}[Lasso and Adaptive Lasso: projection sets]\label{cor:proj_formulas_lasso_adal}
Consider the Lasso and Adaptive Lasso penalties defined in SM Section~B.1.
\begin{enumerate}[label=(\roman*)]
\item For Lasso, the support set \(C_n\) is given by
\begin{equation}\label{eq:Cn_lasso}
C_n
=
\bigcap_{j=1}^p
\Bigl\{
\bm\theta\in\R^p:\ 
|( \bm W_n\bm\theta )_j|
\le
\lambda_n
\Bigr\}.
\end{equation}

\item For Adaptive Lasso, the support set \(C_n\) is given by
\begin{equation}\label{eq:Cn_adal}
C_n
=
\bigcap_{\{j:\tilde\beta_{nj}\neq 0\}}
\Bigl\{
\bm\theta\in\R^p:\ 
|( \bm W_n\bm\theta )_j|
\le
\lambda_n/|\tilde\beta_{nj}|
\Bigr\}.
\end{equation}
\end{enumerate}
\end{corollary}

The sets in \eqref{eq:Cn_lasso}--\eqref{eq:Cn_adal} are polyhedra.
Hence the Lasso and Adaptive Lasso proximal steps are projection residuals:
\[
\prox_{\lambda_n f_n}^{\bm W_n}(\bm\theta)
=
\bm\theta-P_{C_n}^{\bm W_n}(\bm\theta).
\]
When \(\bm W_n=\bm I_p\), the projection is coordinatewise: Lasso gives the usual soft-thresholding formula, while the Adaptive Lasso gives weighted soft-thresholding on coordinates with \(\tilde\beta_{nj}\neq0\) and sets coordinates with \(\tilde\beta_{nj}=0\) to zero under the extended-value convention above.
When \(\bm W_n\neq \bm I_p\), the same dual set induces a different shrinkage geometry because the projection is taken in the \(\bm W_n\)-metric.
SM Figure~B.1 illustrates these projection--residual geometries for Lasso and Adaptive Lasso in dimension \(p=2\).

\subsection{Finite--sample properties}\label{sec:finite_sample_properties}

The proximal construction \eqref{eq:prox_def} is well posed for a broad class of convex penalties and, importantly for econometric applications, enjoys deterministic regularity properties that hold pointwise in the data. Because these properties do not rely on probabilistic arguments, they apply uniformly in finite samples, irrespective of the stochastic structure of the benchmark estimator \(\hat{\bm\beta}_n^s\).
\begin{theorem}[Basic finite-sample properties of proximal estimators]\label{thm:properties}
Let \(\hat{\bm\beta}_n^s\) be a well-defined benchmark estimator.
Then the proximal estimator in \eqref{eq:prox_def}
exists and is unique, and the proximal operator in \eqref{eq:prox_op_def} satisfies:
\begin{enumerate}[label=(\roman*)]
\item\label{prop:properties:existence} Existence and uniqueness:
for every \(\bm\theta\in\R^p\), the problem
\[
\prox_{\lambda_n f_n}^{\bm W_n}(\bm\theta)
=
\argmin_{\bm\beta\in\R^p}
\left\{
\frac12\|\bm\beta-\bm\theta\|_{\bm W_n}^2+\lambda_n f_n(\bm\beta)
\right\}
\]
admits a unique solution.

\item\label{prop:properties:stability} Nonexpansiveness:
for any \(\bm\beta_1,\bm\beta_2\in\R^p\),
\begin{equation}\label{eq:stability_general}
\Big\|
\prox_{\lambda_n f_n}^{\bm W_n}(\bm\beta_1)
-
\prox_{\lambda_n f_n}^{\bm W_n}(\bm\beta_2)
\Big\|_{\bm W_n}
\le
\|\bm\beta_1-\bm\beta_2\|_{\bm W_n}.
\end{equation}

\item\label{prop:properties:differentiability} Almost-everywhere differentiability:
the map \(\prox_{\lambda_n f_n}^{\bm W_n}\) is differentiable Lebesgue--almost everywhere on \(\R^p\).
\end{enumerate}
\end{theorem}

Theorem~\ref{thm:properties} provides the deterministic foundation for treating proximal regularization as a disciplined second step.
Part~\ref{prop:properties:existence} guarantees that once a benchmark estimator is available, the regularization step is automatically well-defined and single-valued.
This is especially useful in econometric settings where the benchmark itself may be produced by a multi-stage or nonsmooth procedure: the proximal step cannot introduce multiplicity.
By contrast, uniqueness of the unpenalized benchmark need not carry over automatically to penalize-the-criterion formulations; SM Section~B.2 gives a simple GMM example in which the unpenalized criterion has a unique minimizer for every sample, while the corresponding \(\ell_1\)-penalized criterion has a non-singleton argmin with positive probability.

Part~\ref{prop:properties:stability} is the key finite-sample stability
property. It shows that proximal operators are nonexpansive in
the \(\bm W_n\)-norm: discrepancies in the benchmark estimator cannot be
amplified by the regularization step. Thus, if two benchmark inputs
\(\bm\beta_1\) and \(\bm\beta_2\) are close, the corresponding proximal
outputs are at least as close in the same metric. This property is useful
whenever the benchmark estimator is itself affected by preliminary
estimation, numerical approximation, nuisance estimation, or alternative
sample perturbations, because such perturbations are propagated through a
Lipschitz map with constant one.

Part~\ref{prop:properties:differentiability} gives a complementary local regularity property.
Although the penalty \(f_n\) may be nonsmooth, the proximal operator itself is still differentiable almost everywhere as a function of its input.

The dual representation developed in Proposition~\ref{prop:dual_repr} complements Theorem~\ref{thm:properties} by giving an explicit residual interpretation of the same regularization step.
Therefore, the deterministic properties in Theorem~\ref{thm:properties} apply directly to both the primal and dual expressions.

\begin{example}[Penalized optimal two-step GMM, continued]\label{ex:penalized_2sgmm_continued}
Consider again the proximal GMM--Adaptive Lasso estimator \eqref{eq:running_example_2sgmm_adal}, based on the benchmark \(\hat{\bm\beta}_n^{\mathrm{2SGMM}}\), the metric \(\bm W_n\), and the penalty \(f_n\) in \eqref{eq:adaptive_lasso_eta}, where only the misspecification block \(\bm\gamma\) is penalized. By applying Proposition~\ref{prop:dual_repr} to the block-separable penalty in \eqref{eq:adaptive_lasso_eta},
\begin{equation}\label{eq:2sgmm_adal_projection_merged}
\prox_{\lambda_n f_n}^{\bm W_n}\big(\hat{\bm\beta}_n^{\mathrm{2SGMM}}\big)
=
\hat{\bm\beta}_n^{\mathrm{2SGMM}}
-
P_{C_n}^{\bm W_n}\big(\hat{\bm\beta}_n^{\mathrm{2SGMM}}\big),
\end{equation}
where the support set is
\begin{align}\label{eq:Cn_adal_eta_merged}
C_n
=
\Big\{
\bm u\in\R^{d+k}\ :\
&\langle \bm e_\ell,\bm u\rangle_{\bm W_n}=0,\ \ell=1,\ldots,d,\ \text{and}\nonumber\\
&|\langle \bm e_{d+j},\bm u\rangle_{\bm W_n}|
\le
\lambda_n/|\tilde\gamma_{nj}|,
\ j:\tilde\gamma_{nj}\neq 0
\Big\}.
\end{align}
Thus the absence of a penalty on the structural block appears in the dual as the equality constraints
\(\langle \bm e_\ell,\bm u\rangle_{\bm W_n}=0\), \(\ell=1,\ldots,d\). The misspecification coordinates with \(\tilde\gamma_{nj}\neq0\) satisfy the displayed dual bounds, whereas coordinates with \(\tilde\gamma_{nj}=0\) are unconstrained in the dual because they are hard-constrained in the primal through \(\iota_{\{0\}}(\gamma_j)\).
If \(\bm W_n=\bm I_{d+k}\), the projection sets the structural dual correction to zero, clips the misspecification coordinates with \(\tilde\gamma_{nj}\neq0\), and leaves the dual coordinates with \(\tilde\gamma_{nj}=0\) unrestricted. Hence, the proximal estimator leaves \(\bm\theta\) unchanged, soft-thresholds the active preliminary \(\bm\gamma\)-coordinates, and sets the zero-preliminary \(\bm\gamma\)-coordinates to zero.
Under the efficient GMM metric \(\bm W_n\), by contrast, the correction subtracted from \(\hat{\bm\beta}_n^{\mathrm{2SGMM}}\) may mix \(\bm\theta\) and \(\bm\gamma\) through the off-diagonal blocks of \(\bm W_n\).
Moreover, Theorem~\ref{thm:properties} yields finite-sample implications that are especially useful in this multi-step setting. First, the proximal estimator is uniquely defined for every realization of the data, because the regularization step is a strongly convex problem once \(\bm W_n\succ\bm 0\). Second, the nonexpansiveness bound \eqref{eq:stability_general} implies that perturbations of the benchmark, for example those induced by replacing \(\hat{\bm\beta}_n^{\mathrm{2SGMM}}\) with an asymptotically equivalent or numerically stabilized version, cannot be amplified by the proximal step.
This separation between benchmark construction and regularization is a genuine advantage of the proximal formulation. For the penalize-the-criterion estimator \eqref{eq:pen_gmm_adal} in \citet{liao2013adaptive}, no analogous standalone regularization map is available because the penalty is built into the nonlinear GMM objective itself. The cost is that existence, uniqueness, and local analysis must be established for the penalty-augmented criterion as a whole: one must work with a case-specific local expansion of the penalized GMM objective,  by imposing conditions ensuring that the penalized minimizer is locally well behaved and that the penalty does not distort the first-order GMM geometry on the active coordinates. By contrast, in the proximal formulation the benchmark supplies the stochastic local approximation, while the regularization step is handled through generic properties of the proximal operator.
\end{example}
 \section{Asymptotic properties of proximal estimators}
\label{sec:asymptotics}

This section develops a fixed-\(p\) large-sample theory for proximal estimators
\[
\hat{\bm\beta}_n
=
\prox_{\lambda_n f_n}^{\bm W_n}(\hat{\bm\beta}_n^s),
\qquad n\to\infty,\qquad p\ \text{fixed}.
\]
The main advantage of the proximal formulation is its modularity.
The benchmark estimator \(\hat{\bm\beta}_n^s\), the metric \(\bm W_n\), and the penalty sequence \(\lambda_n f_n\) enter the asymptotic analysis through separate ingredients.
Probability limits follow from the stability of minimizers of strongly convex random functions.
Asymptotic distributions are then derived by recentering at \(\bm\beta_0\), rescaling at the benchmark rate \(r_n\), and identifying the corresponding first-order limit of the penalty. The weak convergence results in this section are pointwise; no uniform convergence over classes of data-generating processes satisfying the stated assumptions is asserted. Throughout, our asymptotic analysis is formulated using epigraph convergence, which is the appropriate notion of convergence for extended-real valued convex penalties, including the Adaptive Lasso and indicators of convex constraint sets.\footnote{
  SM Section~A.1 collects deterministic and stochastic definitions of epigraph convergence; see, e.g., \citet{rockafellarWets2009,salinetti1981convergence,salinetti1986convergence}.
}
This yields a unified fixed-\(p\) asymptotic theory covering deterministic and stochastic penalties, data-dependent metrics, and constraint-inducing penalties within the same framework.

\subsection{Probability limits}\label{subsec: preliminaries}

We begin with probability limits. Because a proximal estimator is the unique minimizer of a strongly convex criterion, consistency can be studied as a stability problem for minimizers of convex random functions. To this end, we impose high-level conditions on the benchmark estimator, the metric, and the penalty sequence.

\begin{assumption}\label{ass: beta}
\(\hat{\bm\beta}_n^s\to_{\Pr}\bm\beta_0\) for some \(\bm\beta_0\in\R^p\).
\end{assumption}

\begin{assumption}\label{ass: W}
\(\bm W_n\to_{\Pr}\bm W_0\) for some positive definite matrix \(\bm W_0\in\R^{p\times p}\).
\end{assumption}

\begin{assumption}\label{ass:f}
The following properties hold:
\begin{enumerate}[label=(\roman*)]
\item\label{epif ass}
\(f_n\to_{\Pr} f_0\) in epigraph for some proper function \(f_0:\R^p\to(-\infty,+\infty]\);
\item\label{domf ass}
\(\bm\beta_0\in \dom(f_0)\) and \(\bm\beta_0\in \dom(f_n)\) for all sufficiently large \(n\).
\end{enumerate}
\end{assumption}

These assumptions are deliberately modular. Assumptions~\ref{ass: beta} and~\ref{ass: W} concern only the benchmark estimator and the metric, and can therefore be verified model by model, independently of the subsequent regularization step. Assumption~\ref{ass:f}\ref{epif ass} requires epigraph convergence in probability of the penalty sequence. This condition is mild and accommodates our benchmark penalties; see SM Section~B.3. Assumption~\ref{ass:f}\ref{domf ass} ensures that the target \(\bm\beta_0\) remains feasible along the penalty sequence, including cases in which the limit penalty is extended-real valued and therefore encodes a convex constraint. 

The next proposition identifies the resulting probability limits of the proximal estimator under the two regimes \(\lambda_n\to\lambda_0>0\) and \(\lambda_n\to 0\).

\begin{proposition}[Probability limits of proximal estimators]
\label{prop:consistency}
Let Assumptions~\ref{ass: beta}, \ref{ass: W}, and~\ref{ass:f} hold, and let \(\lambda_n> 0\).
\begin{enumerate}[label=(\roman*)]

\item\label{prop_cons_i}
If \(\lambda_n\to\lambda_0>0\), then
\[
\prox_{\lambda_n f_n}^{\bm W_n}(\hat{\bm\beta}_n^s)
\to_{\Pr}
\prox_{\lambda_0 f_0}^{\bm W_0}(\bm\beta_0)=
\big(\Id-\prox_{(\lambda_0 f_0)^*}^{\bm W_0}\big)(\bm\beta_0),
\]
with \((\lambda_0 f_0)^*\) denoting the convex conjugate of \(\lambda_0 f_0\)
relative to the limit inner product \(\langle\cdot,\cdot\rangle_{\bm W_0}\).

\item\label{prop_cons_ii}
If \(\lambda_n\to 0\) and, in addition,
\begin{equation}\label{eq:condition:lambda_f}
\lambda_n f_n \to_{\Pr} \iota_{\dom(f_0)}
\qquad\text{in epigraph,}
\end{equation}
then
\[
\prox_{\lambda_n f_n}^{\bm W_n}(\hat{\bm\beta}_n^s)
\to_{\Pr}
\prox_{\iota_{\dom(f_0)}}^{\bm W_0}(\bm\beta_0)=
P_{\dom(f_0)}^{\bm W_0}(\bm\beta_0)
=
\bm\beta_0.
\]

\end{enumerate}
\end{proposition}

Proposition~\ref{prop:consistency} separates three conceptually distinct situations.
First, when \(\lambda_n\to 0\) and the scaled penalties vanish epigraphically, the proximal step is asymptotically inactive.
This covers the standard deterministic penalties that are finite everywhere, such as Ridge, Lasso, Elastic Net, and Group Lasso.
In that case \(\dom(f_0)=\R^p\), hence \(\iota_{\dom(f_0)}=0\), and \eqref{eq:condition:lambda_f} reduces to \(\lambda_n f_n\to_{\Pr}0\) in epigraph.
Consistency is then inherited directly from the benchmark:
\[
\prox_{\lambda_n f_n}^{\bm W_n}(\hat{\bm\beta}_n^s)-\hat{\bm\beta}_n^s
\to_{\Pr}\bm 0.
\]
Second, \(\lambda_n\to 0\) need not imply asymptotic inactivity.
For adaptive and constrained procedures, the scaled penalties may instead satisfy the constraint-type limit \(\lambda_n f_n\to_{\Pr}\iota_{\dom(f_0)}\).
Then the proximal step converges to the \(\bm W_0\)-projection onto \(\dom(f_0)\).
Because Assumption~\ref{ass:f}\ref{domf ass} ensures \(\bm\beta_0\in\dom(f_0)\), this projection equals \(\bm\beta_0\), so consistency is preserved even though the limiting effect of the penalty is nontrivial.
This regime is central for adaptive procedures, where vanishing tuning can coexist with asymptotically binding domain restrictions.

Third, part~\ref{prop_cons_i} records the nonvanishing regime \(\lambda_n\to\lambda_0>0\).
Here the estimator converges to a regularized population target,
\(\prox_{\lambda_0 f_0}^{\bm W_0}(\bm\beta_0)\),
rather than necessarily to \(\bm\beta_0\).
This regime is useful when one intentionally targets a regularized estimand, but it is not the default route to consistency for an underlying structural parameter unless \(\bm\beta_0\) is a fixed point of the limit proximal map.

\begin{example}[Probability limits for penalized optimal two-step GMM (Example~\ref{ex:penalized_2sgmm} cont'd)]\label{ex:consistency_2sgmm_adal}
Maintain the setup of Example~\ref{ex:penalized_2sgmm}, with stacked parameter
\(\bm\beta=(\bm\theta',\bm\gamma')'\) and benchmark \(\hat{\bm\beta}_n^{\mathrm{2SGMM}}\) defined by \eqref{eq:2sgmm}.
Under standard regularity conditions for nonlinear two-step GMM (see \citet{hansen1982large}),
\[
(\hat{\bm\beta}_n^{\mathrm{2SGMM}}-
\bm\beta_0)=O_{\Pr} (1/\sqrt{n}),
\qquad
\bm W_n=\widehat{\bm G}_n'\widehat{\bm\Omega}_n^{-1}\widehat{\bm G}_n \to_{\Pr}\bm W_0\succ \bm 0;
\]
hence, Assumptions~\ref{ass: beta}--\ref{ass: W} hold.
Let \(\mathcal A_{\bm \gamma}:=\{j\le k:\gamma_{0j}\neq 0\}\) and define the oracle subspace
\begin{eqnarray}
\mathcal S_0
:=
\big\{(\bm\theta,\bm\gamma)\in\R^{d+k}:\ \gamma_j=0\ \ \text{for every}\ j\notin \mathcal A_{\bm \gamma}\big\}.\label{eq: SD_0 def GMM example 1}
\end{eqnarray}
If \(\lambda_n\to0\), \(\sqrt n\lambda_n\to\infty\), and
\(\tilde{\bm\beta}_n-\bm\beta_0=O_{\Pr}(n^{-1/2})\), then the adaptive
weights in the penalty \eqref{eq:adaptive_lasso_eta} diverge on the truly
inactive coordinates, and the scaled penalties satisfy:\footnote{
For instance, $\tilde{\bm\beta}_n$ could be taken to be the first step GMM estimator $\hat{\bm\beta}^{(1)}_n$ based on the weighting matrix $\bm I_{m+k}$.}
\[
\lambda_n f_n\ \to_{\Pr}\ \iota_{\mathcal S_0}
\qquad\text{in epigraph;}
\]
see the verification in SM Section~B.3.
By construction
\(\bm\beta_0\in\mathcal S_0\), i.e.,
Assumption \ref{ass:f} also holds.
Proposition~\ref{prop:consistency}\ref{prop_cons_ii} then yields
\[
\hat{\bm\beta}_n
=
\prox_{\lambda_n f_n}^{\bm W_n}\big(\hat{\bm\beta}_n^{\mathrm{2SGMM}}\big)
\to_{\Pr}\bm\beta_0,
\qquad
\prox_{(\lambda_n f_n)^*}^{\bm W_n}\big(\hat{\bm\beta}_n^{\mathrm{2SGMM}}\big)\to_{\Pr}\bm 0.
\]
Thus, consistency is inherited from the benchmark estimator under the same high-level conditions used to establish consistency of optimal two-step GMM. The only additional ingredient is the epigraph convergence \(\lambda_nf_n\to_{\Pr}\iota_{\mathcal S_0}\), which, for the Adaptive Lasso penalty, follows from the asymptotic behavior of the adaptive weights on active and inactive misspecification coordinates.
\end{example}

\subsection{Asymptotic distribution}\label{sec: Asymptotic distribution}

We now obtain the asymptotic distribution of proximal estimators.
To this end, we need a weak limit law for the benchmark estimator under suitable centering and rescaling, together with a feasible first-order description of the penalty sequence in a shrinking neighborhood of \(\bm\beta_0\).

\begin{assumption}\label{ass: asy distr proximal}
\(r_n(\hat{\bm\beta}_n^s-\bm\beta_0)\to_d \bm\eta\)
for some random vector \(\bm\eta\) and some rate \(r_n\to\infty\).
\end{assumption}

Assumption~\ref{ass: asy distr proximal} postulates an asymptotic distribution for the benchmark estimator.
In particular, it implies
\[
r_n(\hat{\bm\beta}_n^s-\bm\beta_0)=O_{\Pr}(1),
\]
and hence strengthens Assumption~\ref{ass: beta} to the scale relevant for asymptotic distribution theory.
Once the benchmark fluctuation is fixed, the regularization step enters through the local behavior of the penalty around \(\bm\beta_0\).
Therefore, the relevant tuning quantity for characterizing the estimator's asymptotic distribution is the product \(\lambda_n r_n\), which determines whether the penalty contributes at the same order as the benchmark fluctuation or via
a hard feasibility constraint in the limit.
To describe these local contributions, we first introduce the directional derivative of the limit penalty \(f_0\) at \(\bm\beta_0\).
For \(\bm b\in\R^p\), define
\begin{equation}\label{eq:dirder_rho_intro}
\rho_{\bm\beta_0}(\bm b)
:=
\lim_{\alpha\downarrow 0}
\frac{f_0(\bm\beta_0+\alpha\bm b)-f_0(\bm\beta_0)}{\alpha}
\in(-\infty,+\infty].
\end{equation}
This object captures the first-order geometry of \(f_0\) at \(\bm\beta_0\).
To express the limit law in dual form, we also use the subdifferential of \(f_0\) at \(\bm\beta_0\), defined under the limit inner product \(\langle\cdot,\cdot\rangle_{\bm W_0}\),
\begin{equation}\label{eq:subgrad_W0_intro}
\partial f_0(\bm\beta_0)
:=
\bigcap_{\bm\beta\in\R^p}
\Big\{
\bm t\in\R^p:\ 
f_0(\bm\beta)-f_0(\bm\beta_0)-\langle \bm\beta-\bm\beta_0,\bm t\rangle_{\bm W_0}\ge 0
\Big\}.
\end{equation}
Under Assumption~\ref{ass:f}, \(\rho_{\bm\beta_0}\) is well-defined (\citet[Thm.~17.2]{bauschke2016convex}).
If, in addition, \(\partial f_0(\bm\beta_0)\neq\emptyset\),
the following tight duality relations hold (\citet[Prop.~17.17]{bauschke2016convex}):
\[
\rho_{\bm\beta_0}
=
(\iota_{\partial f_0(\bm\beta_0)})^*
=
\sigma_{\partial f_0(\bm\beta_0)},
\qquad
\sigma_C(\bm b):=\sup_{\bm t\in C}\langle \bm t,\bm b\rangle_{\bm W_0}.
\]
Note that when \(f_0=\iota_C\) for some nonempty closed convex set \(C\), then definition \eqref{eq:subgrad_W0_intro} simplifies to:
\[
\partial f_0(\bm\beta_0)=N_C(\bm\beta_0),\qquad N_C(\bm\beta_0)
:=
\Big\{
\bm t\in\R^p:\ \sup_{\bm\beta\in C}\langle\bm\beta-\bm\beta_0, \bm t\rangle_{\bm W_0}\le 0
\Big\}\ ,
\]
and
\[
\rho_{\bm\beta_0}
=
\sigma_{N_C(\bm\beta_0)}.
\]
SM Section~B.1 reports the directional derivatives and subgradients of our benchmark penalties.

Epi-convergence of \(f_n\) to \(f_0\) is enough for probability limits, but it does not by itself identify the local penalty contribution at the \(r_n\)-scale.
For asymptotic distribution theory, the relevant object is the rescaled local difference quotient
\begin{equation}
q_n(\bm b)
:=
r_n\left[
f_n\left(\bm\beta_0+\bm b/r_n\right)-f_n(\bm\beta_0)
\right],
\qquad
\bm b\in\R^p.\label{eq: q_n quotient}
\end{equation}
This sequence of functions records the first-order variation of the finite-sample penalties in shrinking neighborhoods of \(\bm\beta_0\).
The next theorem derives the asymptotic distribution for proximal estimators in the two regimes \(\lambda_n r_n\to\lambda_0>0\) and \(\lambda_n r_n\to\lambda_0=0\), with corresponding local epigraph convergence conditions on \(q_n\) toward the directional derivative of the relevant limit penalty.

\begin{theorem}[Weak limits of proximal estimators]
\label{prop:asymptotics:alt}
Let Assumptions~\ref{ass: W}, \ref{ass:f}, and
\ref{ass: asy distr proximal} hold, and let \(\lambda_n> 0\).
Then:
\begin{enumerate}[label=(\roman*)]
\item\label{prop_distr_i}
If \(\lambda_n r_n\to\lambda_0>0\) and
\begin{equation}\label{eq:local_regime}
q_n \to_{\Pr} \rho_{\bm\beta_0}
\qquad\text{in epigraph,}
\end{equation}
then
\begin{equation}\label{eq:asydistr_i}
r_n(\prox_{\lambda_n f_n}^{\bm W_n}(\hat{\bm\beta}_n^s)-\bm\beta_0)
\to_d
\prox_{\lambda_0\rho_{\bm\beta_0}}^{\bm W_0}(\bm\eta).
\end{equation}
If, in addition, \(\partial f_0(\bm\beta_0)\neq\emptyset\), then
\begin{equation}\label{eq:asym_moreau_i}
\prox_{\lambda_0\rho_{\bm\beta_0}}^{\bm W_0}(\bm\eta)
=
\big(\Id-P_{\lambda_0\partial f_0(\bm\beta_0)}^{\bm W_0}\big)(\bm\eta).
\end{equation}

\item\label{prop_distr_ii}
If \(\lambda_n r_n\to 0\) and
\begin{equation}\label{eq:local_scaled_zero_regime}
(\lambda_n r_n)q_n \to_{\Pr} \sigma_{N_{\dom(f_0)}(\bm\beta_0)}
\qquad\text{in epigraph,}
\end{equation}
then
\begin{equation}\label{eq:asydistr_ii}
r_n(\prox_{\lambda_n f_n}^{\bm W_n}(\hat{\bm\beta}_n^s)-\bm\beta_0)
\to_d
\prox_{\sigma_{N_{\dom(f_0)}(\bm\beta_0)}}^{\bm W_0}(\bm\eta),
\end{equation}
where
\begin{equation}\label{eq:asym_moreau_ii}
\prox_{\sigma_{N_{\dom(f_0)}(\bm\beta_0)}}^{\bm W_0}(\bm\eta)
=
\big(\Id-P_{N_{\dom(f_0)}(\bm\beta_0)}^{\bm W_0}\big)(\bm\eta).
\end{equation}
\end{enumerate}
\end{theorem}

Theorem~\ref{prop:asymptotics:alt} shows that the weak limits of proximal estimators are determined in a modular way by three ingredients: the benchmark limit \(\bm\eta\), the limit \(\bm W_0\)-metric, and the local behavior of the penalty around \(\bm\beta_0\). The difference between parts~\ref{prop_distr_i} and~\ref{prop_distr_ii} comes from the scaling of the penalty parameter relative to the benchmark rate, through the product \(\lambda_n r_n\).

In part~\ref{prop_distr_i}, \(\lambda_n r_n\to\lambda_0>0\), so the penalty remains active at the benchmark scale. Accordingly, the relevant local condition is that \(q_n\) converge in epigraph to \(\rho_{\bm\beta_0}\), the directional derivative of the limit penalty \(f_0\) at \(\bm\beta_0\). Thus the finite-sample penalties contribute at first order through the local slope of \(f_0\), and the resulting asymptotic distribution is \(\prox_{\lambda_0\rho_{\bm\beta_0}}^{\bm W_0}(\bm\eta)\). If, in addition, \(\partial f_0(\bm\beta_0)\neq\emptyset\), then \(\rho_{\bm\beta_0}=\sigma_{\partial f_0(\bm\beta_0)}\), so the same weak limit admits the dual characterization \eqref{eq:asym_moreau_i} as a projection--residual map.

In part~\ref{prop_distr_ii}, by contrast, \(\lambda_n r_n\to0\), so the penalty is asymptotically too small to contribute through the local slope of \(f_0\) itself. At that scale, only the domain restrictions encoded by the penalty can impact the limit. For this reason, the relevant condition is that the scaled quotient \((\lambda_n r_n)q_n\) converge in epigraph to \(\sigma_{N_{\dom(f_0)}(\bm\beta_0)}\), which is the directional derivative of the indicator function \(\iota_{\dom(f_0)}\) at \(\bm\beta_0\). Hence the first-order effect of regularization is no longer a local shrinkage effect, but a local feasibility effect. This gives rise to an asymptotic distribution given by \(\prox_{\sigma_{N_{\dom(f_0)}(\bm\beta_0)}}^{\bm W_0}(\bm\eta)\), or equivalently, by the residual of the  projection of \(\bm\eta\) onto the normal cone \(N_{\dom(f_0)}(\bm\beta_0)\).

The proof proceeds by recentering the proximal objective at \(\bm\beta_0\), introducing the local variable \(\bm\beta=\bm\beta_0+\bm b/r_n\), and analyzing the resulting random convex criterion in the local coordinate \(\bm b\). Under this localization, the benchmark contributes the local quadratic term centered at \(r_n(\hat{\bm\beta}_n^s-\bm\beta_0)\), while the penalty contributes through the local quotient \(q_n\). SM Section~B.3 verifies the corresponding local conditions in parts~\ref{prop_distr_i} and~\ref{prop_distr_ii} for our benchmark penalties.

\begin{remark}[General weak limit of proximal estimators]
The proof of Theorem~\ref{prop:asymptotics:alt} uses the specific
local penalty limits in the statement only to identify explicit forms of the
limit criterion.  The same argument yields a more general high-level result.
If \((\lambda_n r_n)q_n \to_{\Pr} h_0\) in epigraph, for some deterministic proper function \(h_0:\R^p\to(-\infty,+\infty]\),
then \(h_0\in\Gamma(\R^p)\) by the closedness of this class under
epigraph convergence \citep[Ch.~7]{rockafellarWets2009}. Consequently,
\[
r_n\left(
\prox_{\lambda_n f_n}^{\bm W_n}(\hat{\bm\beta}_n^s)-\bm\beta_0
\right)
\to_d
\prox_{h_0}^{\bm W_0}(\bm\eta)
=
\big(\Id-\prox_{h_0^*}^{\bm W_0}\big)(\bm\eta).
\]
Thus the scaled local penalty need not converge to
a directional derivative of a limiting penalty, nor to the support function of a
normal cone, for the proximal weak-limit representation to remain valid.
Theorem~\ref{prop:asymptotics:alt} is stated in this more structured form
because it isolates the local regimes arising for our benchmark penalties.
These regimes are organized
by the behavior of the local tuning sequence \(\lambda_n r_n\) relative to the
local difference quotient \(q_n\): case (i) corresponds to regularization that
remains active at the \(r_n\)-local scale, while case (ii) corresponds to locally vanishing
regularization that still leaves a tangent-cone or normal-cone restriction.
\end{remark}

\begin{remark}[Relation to benchmark results in the literature]
When specialized through the PLSE representation in Section~\ref{sec:link_plse}, Theorem~\ref{prop:asymptotics:alt} recovers several local weak limits already obtained in the literature for convex penalized least-squares estimators; SM Section~B.4 illustrates this reduction in the regular-design Lasso case by deriving the weak limit of \citet{knight2000asymptotics} together with its equivalent projection-residual representation.
\end{remark}
\begin{example}[Weak limit for penalized optimal two-step GMM; Example~\ref{ex:penalized_2sgmm} cont'd]\label{ex:asydistr_2sgmm_lasso}
In the setup of Example~\ref{ex:consistency_2sgmm_adal}, with benchmark estimator \(\hat{\bm\beta}_n^{\mathrm{2SGMM}}\) and metric
\[
\bm W_n
=
\widehat{\bm G}_n'\widehat{\bm\Omega}_n^{-1}\widehat{\bm G}_n
\to_{\Pr}\bm W_0\succ \bm 0, 
\]
the usual regularity conditions for optimal two-step GMM give (see \citet{hansen1982large}):
\[
\sqrt n\big(\hat{\bm\beta}_n^{\mathrm{2SGMM}}-\bm\beta_0\big)\to_d \bm\eta,
\qquad
\bm\eta\sim N(\bm 0,\bm W_0^{-1}),
\]
i.e., Assumption~\ref{ass: asy distr proximal} holds with \(r_n=\sqrt n\).
Consider again the proximal Adaptive Lasso penalty on the misspecification block under tuning such that
\(\lambda_n \sqrt{n}\to 0\) 
and \(n\lambda_n\to\infty\).
Then, the scaled local difference quotients (\ref{eq: q_n quotient}) produce the constraint-type epi-limit:
\[
\sqrt n\lambda_n q_n
\ \to_{\Pr}\
\sigma_{N_{\mathcal S_0}(\bm\beta_0)}
=
\iota_{\mathcal S_0}
\qquad\text{in epigraph,}
\]
with ${\cal S}_0$ defined in equation (\ref{eq: SD_0 def GMM example 1}).
Hence, Theorem~\ref{prop:asymptotics:alt}\ref{prop_distr_ii} yields:
\begin{equation}\label{eq:ex_2sgmm_adal_limit}
\sqrt n\big(\hat{\bm\beta}_n-\bm\beta_0\big)
\ \to_d\
\prox_{\sigma_{N_{\mathcal S_0}(\bm\beta_0)}}^{\bm W_0}(\bm\eta)
=
\big(\Id-P_{N_{\mathcal S_0}(\bm\beta_0)}^{\bm W_0}\big)(\bm\eta).
\end{equation}
Noting that \(\mathcal S_0\) is a linear subspace,
we further obtain
\(N_{\mathcal S_0}(\bm\beta_0)=\mathcal S_0^\perp\) (\citet[Prop.~6.22]{bauschke2016convex}).
As a consequence,
\begin{equation}\label{eq:ex_2sgmm_adal_projection}
\prox_{\sigma_{N_{\mathcal S_0}(\bm\beta_0)}}^{\bm W_0}(\bm\eta)
=\big(\Id-P_{\mathcal S_0^\perp}^{\bm W_0}\big)(\bm\eta) =
P_{\mathcal S_0}^{\bm W_0}(\bm\eta).
\end{equation}
Thus, the proximal Adaptive Lasso retains exactly the components of the benchmark Gaussian fluctuation lying in the oracle sparsity subspace ${\cal S}_0$, by removing
the orthogonal component to it.
This further highlights the simplification provided by the proximal approach in nonlinear GMM settings. Indeed, once the benchmark CLT and metric limit are established, the asymptotic distribution of the regularized estimator follows directly from the epi-limit of the penalty sequence, without requiring any expansion of a penalized nonlinear GMM criterion.
This contrasts, for example, with the penalize-the-criterion adaptive GMM framework  in \citet{liao2013adaptive}. There, the asymptotic law relies on a local quadratic expansion of the penalized GMM objective, requiring additional smoothness and derivative conditions to ensure that the penalty does not distort the first-order GMM geometry.
By contrast, the proximal approach dispenses with such smoothness requirements, as regularization is handled directly through the epi-limit
of \(\sqrt n\lambda_n f_n\) once the benchmark GMM asymptotic distribution is established. 
\end{example}

\subsection{Oracle property}\label{subsec:oracle}

We now turn to the classical oracle paradigm of \citet{fan2001variable} within our asymptotic framework.
To this end, let
\[
\mathcal A:=\{j\le p:\beta_{0j}\neq 0\},
\qquad
\mathcal A^c:=\{1,\dots,p\}\setminus\mathcal A,
\]
denote the active and inactive parameter subsets, respectively.
The infeasible oracle procedure knows \(\mathcal A\) in advance, sets \((\bm\beta_0)_{\mathcal A^c}=\bm 0\) and estimates only \((\bm\beta_0)_{\mathcal A}\).
Our main goal in this section is to identify primitive conditions under which a single proximal step
(i) selects \(\mathcal A\) consistently and
(ii) attains the same first-order weak limit on the active coordinates as the oracle estimator.
A key feature of the proximal framework is that these two requirements are governed by distinct objects.
Indeed, oracle selection is driven by the dual thresholding geometry induced by the limit penalty, through subgradient and normal cone properties.
Oracle distribution, by contrast, depends additionally on the limiting metric \(\bm W_0\), which determines how the proximal step asymptotically weights local directions around $\bm\beta_0$ in \(\R^p\).

To formulate the oracle distribution property, we next specialize the benchmark limit in Assumption~\ref{ass: asy distr proximal} to a (possibly degenerate) Gaussian law.

\begin{assumption}[Gaussian benchmark limit]\label{ass:gaussian initial estimator}
Given Assumption~\ref{ass: asy distr proximal}, suppose that
\[
\bm\eta=\bm M_0^{+}\bm Z,
\qquad
\bm Z\sim\mathcal N(\bm 0,\bm\Omega_0),
\]
for some \(p\times p\) matrices \(\bm M_0\) and \(\bm\Omega_0\), such that $\bm\Omega_0$ is positive semidefinite.
\end{assumption}

Assumption~\ref{ass:gaussian initial estimator} is satisfied under standard regularity conditions for, e.g., many M-- and GMM--type benchmark estimators.
It allows for singular Gaussian limits, since
\[
\Var(\bm\eta)=\bm M_0^{+}\bm\Omega_0\bm M_0^{+\prime}
\]
need not be full rank.
When \(\bm M_0\) and \(\bm\Omega_0\) are nonsingular, then \(\bm\eta\) is nondegenerate Gaussian and the assumption reduces to the usual asymptotic normal approximation.
To formulate the oracle property in this setting, we adopt in the next definition the notation \((\bm x)_{\mathcal A}\) and \((\bm M)_{\mathcal A}\) for subvectors and submatrices with coordinates indexed by \(\mathcal A\).
\begin{definition}[Oracle property]\label{oracle property}
Let
\[
\hat{\bm\beta}_n:=\prox_{\lambda_n f_n}^{\bm W_n}(\hat{\bm\beta}_n^s),
\qquad
\hat{\mathcal A}_n:=\{j\le p:\ (\hat{\bm\beta}_n)_j\neq 0\}.
\]
We say that \(\hat{\bm\beta}_n\) satisfies the oracle property if, as \(n\to\infty\):
\begin{enumerate}[label=(\alph*)]
\item\label{op1}
(\emph{Oracle selection}) \(\Pr(\hat{\mathcal A}_n=\mathcal A)\to 1\).

\item\label{op2}
(\emph{Oracle distribution}) Under Assumption~\ref{ass:gaussian initial estimator} and the constraint-type asymptotic regime~\ref{prop_distr_ii} of Theorem~\ref{prop:asymptotics:alt}, the active coordinates satisfy
\[
r_n\big(\hat{\bm\beta}_n-\bm\beta_0\big)_{\mathcal A}
\ \to_d\
[(\bm M_0)_{\mathcal A}]^{+}(\bm Z)_{\mathcal A},
\]
where the right-hand side is the target oracle Gaussian law associated with estimation on the active block alone.
\end{enumerate}
\end{definition}
\begin{remark}
It is useful to note that under
asymptotic regime~\ref{prop_distr_ii} in Theorem~\ref{prop:asymptotics:alt} the property \(N_{\dom(f_0)}(\bm\beta_0)=\mathcal S_0^\perp\)
directly implies that Oracle property~\ref{oracle property}\ref{op2} is equivalent to the condition
\[
\big(P_{\mathcal S_0}^{\bm W_0}(\bm\eta)\big)_{\mathcal A}
\stackrel{d}{=}
[(\bm M_0)_{\mathcal A}]^{+}(\bm Z)_{\mathcal A}.
\]
Thus, in our proximal framework the oracle distribution problem reduces to a projection question. Once the penalty enforces the correct sparsity subspace in the asymptotic distribution, the remaining issue is whether the \(\bm W_0\)--projection of $\bm\eta$ onto that subspace reproduces the efficient Gaussian law for the active block.
\end{remark}
We begin with Oracle selection, by providing
two results giving complementary criteria that are directly useful in applications.\footnote{For completeness, SM Section~B.5 provides an exact if-and-only-if characterization of Oracle property~\ref{oracle property}\ref{op1} based on the feasibility of the corresponding local optimality conditions. Although exact, this characterization is less transparent than the two criteria introduced below in this section, which are easier to verify.}
The first criterion is a necessary condition stated entirely in terms of the proximal limit law, and is therefore especially convenient for impossibility arguments.
The second gives a sufficient condition formulated entirely on the penalty side, facilitating verification.

By Theorem~\ref{prop:asymptotics:alt}, the
asymptotic distribution of the proximal estimator admits the representation
\[
r_n(\hat{\bm\beta}_n-\bm\beta_0)
\ \to_d\
(\Id-P_{B_0}^{\bm W_0})(\bm\eta),
\]
where
\[
B_0=
\left\{\begin{array}{cc}
 \lambda_0\partial f_0(\bm\beta_0)    &  \text{in asymptotic regime~\ref{prop_distr_i}}\\
N_{\dom(f_0)}(\bm\beta_0)     & 
\text{in asymptotic regime~\ref{prop_distr_ii}}
\end{array}  \right. .
\]
For support recovery, this decomposition is especially informative.
Oracle selection requires the inactive coordinates of \(\hat{\bm\beta}_n\) to be exactly zero with probability tending to one.
Hence, in the limit, the inactive block of the benchmark fluctuation \(\bm\eta\) must be cancelled entirely by the limiting dual correction \(P_{B_0}^{\bm W_0}(\bm\eta)\), yielding the next proposition.
\begin{proposition}[Necessary condition for Oracle property \ref{op1}]\label{prop:op1 necessary}
 Let Assumptions~\ref{ass: W}, \ref{ass:f}, and
\ref{ass: asy distr proximal} hold and
 asymptotic regimes
  \ref{prop_distr_i} or \ref{prop_distr_ii} in Theorem
  \ref{prop:asymptotics:alt} hold.
 If Oracle Property~\ref{oracle property}\ref{op1}  also holds then:
  \begin{equation}\label{op1 necessary}
    \Pr\left((\bm\eta)_{\A^c}=
    \left( P_{B_0}^{\bm W_0}(\bm\eta)\right)_{\A^c}\right)=1.
  \end{equation}
  \end{proposition}
Necessary condition~\eqref{op1 necessary}
shows that Oracle selection is possible only if the geometry of \(B_0\) is rich enough, in the inactive directions, to reproduce the inactive block of \(\bm\eta\) almost surely.
This condition is often easy to falsify.
Since \(P_{B_0}^{\bm W_0}(\bm\eta)\in B_0\) almost surely, if
$B_0$
is bounded in the inactive directions, condition~\eqref{op1 necessary} can hold only if \((\bm\eta)_{\mathcal A^c}\) is itself almost surely bounded.
This fails whenever the inactive block of \(\bm\eta\) has unbounded support, as in the standard case of a nondegenerate Gaussian benchmark limit.
\begin{example}
The Lasso illustrates a violation of condition~\eqref{op1 necessary} transparently. Indeed, whenever
\(\lambda_n r_n\to\lambda_0>0\) we obtain:
\begin{eqnarray*}
B_0&=&\lambda_0\partial \|\cdot\|_1(\bm\beta_0)
\\
&=&
\Big\{\bm t\in\R^p:
t_j=\lambda_0\sign(\beta_{0j}) \ \text{if }\beta_{0j}\neq 0,
\
t_j\in[-\lambda_0,\lambda_0] \ \text{if }\beta_{0j}=0
\Big\}.
\end{eqnarray*}
Thus, for every \(j\in\mathcal A^c\), the \(j\)th coordinate of any element of \(B_0\) is bounded in the interval \([-\lambda_0,\lambda_0]\):
\[
\big(P_{B_0}^{\bm W_0}(\bm\eta)\big)_j\in[-\lambda_0,\lambda_0]
\qquad\text{almost surely.}
\]
If \(\eta_j\) has unbounded support, then
\[
\Pr\left(\eta_j=\big(P_{B_0}^{\bm W_0}(\bm\eta)\big)_j\right)< 1 ,
\]
i.e., condition \eqref{op1 necessary} fails and the Lasso cannot achieve oracle selection under a finite local penalty level $\lambda_0 > 0$.
\end{example}

A positive oracle selection result requires a different approach from Proposition~\ref{prop:op1 necessary}. One must rule out, with probability tending to one, any locally optimal solution with nonzero inactive coordinates. This leads to a sufficient condition based on the Euclidean subdifferential of the penalty in an \(r_n^{-1}\)-neighborhood of \(\bm\beta_0\). For
any \(\bm\beta\in\R^p\), 
the Euclidean subdifferential is defined by:
\[
\partial^{e} f_n(\bm\beta)
:=
\Big\{
\bm t\in\R^p:\ 
f_n(\bm\beta')\ge f_n(\bm\beta)+\langle \bm t,\bm\beta'-\bm\beta\rangle
\ \text{ for all }\bm\beta'\in\R^p
\Big\},
\]
where \(\langle\cdot,\cdot\rangle\) denotes the standard Euclidean inner product. 
To quantify the strength of the penalty on inactive coordinates, fix
\(M>0\) and define:
\begin{equation*}\label{eq:Delta_oracle_selection}
\Delta_n(M)
:=
\inf\Big\{
r_n\|\bm v_{\mathcal A^c}\|_1:\ 
\|\bm\beta-\bm\beta_0\|_\infty\le M/r_n,\ 
\bm\beta_{\mathcal A^c}\neq \bm 0,\ 
\bm v\in \lambda_n\partial^{e} f_n(\bm\beta)
\Big\},
\end{equation*}
where the infimum is taken jointly over
\((\bm\beta,\bm v)\in\mathbb R^p\times\mathbb R^p\),
with the convention \(\inf\varnothing=+\infty\).
Thus, \(\Delta_n(M)\) is the minimal \(r_n\)-rescaled dual force on inactive coordinates
 compatible with locally optimal solutions that retain nonzero inactive components near
 \(\bm\beta_0\).
Since the benchmark estimator's fluctuations are of order $r_n^{-1}$, persistence of such solutions requires a dual force of the same order. The next proposition shows that this is excluded when $\Delta_n(M)\to\infty$, producing a sufficient condition for oracle selection.

\begin{proposition}[Sufficient condition for Oracle selection]\label{prop:op1 sufficient}
Suppose Assumption~\ref{ass: W} holds and, for some sequence \(r_n\to\infty\),
\begin{equation}\label{eq:oracle_localization}
\hat{\bm\beta}_n^s-\bm\beta_0 = O_{\Pr}(1/r_n),
\qquad
\hat{\bm\beta}_n-\bm\beta_0 = O_{\Pr}(1/r_n).
\end{equation}
If \(\Delta_n(M)\to_{\Pr}\infty\) for every fixed \(M>0\), then Oracle property~\ref{oracle property}\ref{op1} holds.
\end{proposition}

The localization requirement \eqref{eq:oracle_localization} is mild. Its benchmark component follows from Assumption~\ref{ass: asy distr proximal}, while the proximal component can be verified under the asymptotic regime of interest.
Given localization, Proposition~\ref{prop:op1 sufficient} reduces oracle selection to a condition on the penalty, ensuring that the local Euclidean subdifferential generates inactive-coordinate dual forces dominating benchmark estimator's fluctuations. This shifts the problem to a quantity that is directly tied to the penalty and therefore amenable to verification. In particular, for weighted \(\ell_1\)--type penalties the condition $\Delta (M)\to \infty$ takes a simple and easily checkable form.

\begin{corollary}[Weighted \(\ell_1\) penalties and Oracle selection]\label{cor:oracle_weighted_l1}
Consider the following adaptive penalty, based on
random weights \(w_{nj}\in[0,+\infty]\) for parameter components $\beta_j$:
\[
f_n(\bm\beta)
=
\sum_{j:w_{nj}<\infty} w_{nj}|\beta_j|
+
\sum_{j:w_{nj}=+\infty}\iota_{\{0\}}(\beta_j) ,
\]
and let Assumption~\ref{ass: W} and localization requirement \eqref{eq:oracle_localization} hold.
If in addition the following condition is satisfied:\footnote{By definition, we have $\min_{j\in {\cal A}^c } w_{nj}=+\infty$ when ${\cal A}^c = \emptyset$.}
\begin{equation}\label{eq:weighted_l1_oracle_condition}
\lambda_n r_n \min_{j\in\mathcal A^c} w_{nj}\ \to_{\Pr}\ +\infty\ ,
\end{equation}
then Oracle property~\ref{oracle property}\ref{op1} holds.
\end{corollary}
Corollary~\ref{cor:oracle_weighted_l1} reduces the sufficient condition of Proposition~\ref{prop:op1 sufficient} to a simple condition on the weights. For the Lasso, \(w_{nj}=1\), so condition \eqref{eq:weighted_l1_oracle_condition} becomes \(\lambda_n r_n \to \infty\), showing that a finite local penalty level is incompatible with support recovery in this setting.
For the Adaptive Lasso, the weights are
\[
w_{nj}=
\begin{cases}
|\tilde\beta_{nj}|^{-1}, & \tilde\beta_{nj}\neq 0,\\
+\infty, & \tilde\beta_{nj}=0.
\end{cases}
\]
Thus, condition \eqref{eq:weighted_l1_oracle_condition} becomes
\[
\lambda_n r_n \min_{j\in\mathcal A^c}|\tilde\beta_{nj}|^{-1}\to_{\Pr}\infty.
\]
A convenient sufficient condition for this to hold is:
\begin{equation}\label{eq:weighted_l1_oracle_condition_adal_sufficient}
(\tilde{\bm\beta}_n)_{\mathcal A^c}=O_{\Pr}(1/r_n)
\qquad\text{and}\qquad
\lambda_n r_n^2\to\infty.
\end{equation}
Under these conditions, the inactive weights diverge, so the penalty induces increasingly strong separation on \(\mathcal A^c\), yielding oracle selection.

We now turn to Oracle property~\ref{oracle property}\ref{op2}.
Once the local penalty limit enforces the oracle sparsity subspace, the remaining question is whether the active-block projection induced by \(\bm W_0\) reproduces the efficient Gaussian law.
The next proposition provides an exact answer

\begin{proposition}[Necessary and sufficient condition for Oracle distribution]\label{prop:op2}
Let Assumptions~\ref{ass: W}, \ref{ass:f}, and \ref{ass:gaussian initial estimator} hold, and define
\[
\mathcal S_0:=\{\bm \beta \in\R^p:\ \beta_j=0\  \text{ for all } j\in\mathcal A^c\}.
\]
Suppose that regime~\ref{prop_distr_ii} of Theorem~\ref{prop:asymptotics:alt} applies with \(N_{\dom(f_0)}(\bm\beta_0)=\mathcal S_0^\perp\).
Then, Oracle property~\ref{oracle property}\ref{op2} holds if and only if
\begin{equation}\label{eq:oracle_cov_identity}
(\bm W_0 \bm M_0^+ \bm \Omega_0 \bm M_0^{+\prime} \bm W_0 )_{\mathcal A}
=
(\bm W_0)_{\mathcal A}
[(\bm M_0)_{\mathcal A}]^{+}(\bm \Omega_0)_{\mathcal A}[(\bm M_0)_{\mathcal A}]^{+\prime}
(\bm W_0)_{\mathcal A}.
\end{equation}
\end{proposition}
Proposition~\ref{prop:op2} separates Oracle distribution from Oracle selection. It shows that under regime~\ref{prop_distr_ii} Oracle distribution depends on two distinct ingredients.
The first is geometric, requiring the local constraint induced by the penalty to coincide with the oracle sparsity subspace, i.e., \(N_{\dom(f_0)}(\bm\beta_0)=\mathcal S_0^\perp\).\footnote{Equivalently, the domain of the epigraphical penalty limit must impose the correct inactive-coordinate restrictions around \(\bm\beta_0\).}
The second is metric, ensuring that
conditional on that geometric identification the active-block projection induced by \(\bm W_0\) reproduces the efficient Gaussian law.

This second requirement is imposed through the covariance identity \eqref{eq:oracle_cov_identity}, which can be verified directly from the benchmark limit and the chosen metric. For regular designs, the identity follows automatically whenever
\[
\bm W_0={\bm \Omega}_0 =\bm M_0 .
\]
Under nearly singular and singular designs, however, the interaction between the population geometry and the metric requires additional analysis to ensure that condition \eqref{eq:oracle_cov_identity} remains satisfied. These issues are studied in detail in Section~\ref{sec:irregular_designs}.

\begin{example}[Oracle property for penalized optimal two-step GMM; Example~\ref{ex:penalized_2sgmm} cont'd]\label{ex:oracle_2sgmm_adal}
Recall
from Examples~\ref{ex:consistency_2sgmm_adal} and~\ref{ex:asydistr_2sgmm_lasso}
the notation \(\mathcal A_{\bm \gamma}:=\{j\le k:\gamma_{0j}\neq 0\}\) and:\footnote{For this block-penalized GMM example, the oracle coordinate subspace is indexed by
\[
\mathcal A^*
=
\{1,\ldots,d\}\cup\{d+j:\ j\in\mathcal A_\gamma\}.
\]
Hence, \(\mathcal S_0\) is the coordinate subspace associated with
\(\mathcal A^*\).}
\[
\mathcal S_0
:=
\big\{(\bm\theta,\bm\gamma)\in\R^{d+k}:\ \gamma_j=0\  \text{ for all } j\notin\mathcal A_{\bm \gamma}\big\}.
\]
To verify the oracle property in this setting, it suffices to check three conditions.
First, standard regularity conditions for optimal two-step GMM imply
\(\sqrt n\big(\hat{\bm\beta}_n^{\mathrm{2SGMM}}-\bm\beta_0\big)\to_d \bm\eta
\sim N(\bm 0,\bm W_0^{-1})\),
i.e., the benchmark estimator is \(\sqrt n\)-consistent and asymptotically efficient.
In the notation of Assumption~\ref{ass:gaussian initial estimator}, we therefore set
\(\bm M_0=\bm\Omega_0=\bm W_0\).
Second, 
choose the tuning sequence 
for the Adaptive Lasso penalty
so that
\[
\sqrt n\lambda_n\to 0,
\qquad
n\lambda_n\to+\infty,
\]
and assume that the preliminary estimator satisfies
\(\tilde{\bm\beta}_n-\bm\beta_0=O_{\Pr}(n^{-1/2})\).
Then the local adaptive penalty is asymptotically negligible on the truly
active misspecification coordinates and diverges on the truly inactive ones:
\[
\sqrt n\lambda_n|\tilde\gamma_{nj}|^{-1}\to_{\Pr}0,
\qquad j\in\mathcal A_\gamma,
\]
whereas
\[
\lambda_n\sqrt n|\tilde\gamma_{nj}|^{-1}
=
n\lambda_n\big(\sqrt n|\tilde\gamma_{nj}|\big)^{-1}
\to_{\Pr}+\infty,
\qquad j\notin\mathcal A_\gamma .
\]
Let \(\hat{\mathcal A}_{\gamma,n}
:=
\{j\le k:\hat\gamma_{nj}\neq 0\}\), where \(\hat{\bm\beta}_n=(\hat{\bm\theta}_n',\hat{\bm\gamma}_n')'\).
As a consequence, Corollary~\ref{cor:oracle_weighted_l1}, applied to the
misspecification block, yields oracle selection:
\[
\Pr\big(\hat{\mathcal A}_{\gamma,n}=\mathcal A_\gamma\big)\to 1.
\]
Equivalently,
\[
\Pr\big(\hat\gamma_{nj}=0\ \text{ for all }j\notin\mathcal A_\gamma\big)\to1,
\qquad
\Pr\big(\hat\gamma_{nj}\neq0\ \text{ for all }j\in\mathcal A_\gamma\big)\to1 .
\]
Finally, Oracle property~\ref{oracle property}\ref{op2} also follows.
Indeed, under the Adaptive Lasso penalty the epigraphical limit
domain is \(\dom(f_0)=\mathcal S_0\), which implies
\[
N_{\dom(f_0)}(\bm\beta_0)=N_{\mathcal S_0}(\bm\beta_0)=\mathcal S_0^\perp.
\]
Moreover, since \(\bm M_0=\bm\Omega_0=\bm W_0\), covariance identity \eqref{eq:oracle_cov_identity} trivially holds, with both sides of the identity reducing to \((\bm W_0)_{\mathcal A}\). Therefore,
Proposition~\ref{prop:op2} gives:
\[
\sqrt n\big(\hat{\bm\beta}_n-\bm\beta_0\big)_{\mathcal A}
\ \to_d\
[(\bm W_0)_{\mathcal A}]^{-1}(\bm Z)_{\mathcal A},
\qquad
\bm Z\sim N(\bm 0,\bm W_0),
\]
which is exactly the oracle Gaussian law obtained by restricting the misspecification block to
\(\gamma_j=0\) for all \(j\notin\mathcal A_{\bm\gamma}\).
Overall, this example confirms that oracle verification in nonlinear GMM settings is modular and directly applicable. Given the benchmark CLT and information metric, Oracle selection follows from the weight-growth condition in Corollary~\ref{cor:oracle_weighted_l1}, while Oracle distribution is characterized by the geometric condition
\(
N_{\dom(f_0)}(\bm\beta_0)=\mathcal S_0^\perp
\)
and the covariance identity in Proposition~\ref{prop:op2}. No penalty-specific expansion is required.
\end{example}
 \section{Proximal estimators for linear regression with irregular designs}\label{sec:irregular_designs}

This section studies proximal least-squares estimators for the linear regression model
\begin{equation}\label{linear model}
\bm Y=\bm X\bm\beta_0+\bm\varepsilon,
\end{equation}
under possibly singular or nearly-singular population designs.
Our objective is to construct the first estimators that remain well-defined and retain standard \(\sqrt n\)-Gaussian first-order behavior under both regular and irregular designs.

\subsection{Singular and nearly-singular designs}

The geometry of the problem is governed by
the sample design matrix
\[
\bm Q_n:=\bm X'\bm X/n.
\]
We allow \(\bm Q_n\) to converge to a possibly singular deterministic limit.
\begin{assumption}[Limit design]\label{ass:limit_design_irregular}
There exists a symmetric positive semidefinite matrix \(\bm Q_0\) such that
\(\bm Q_n\to_{\Pr}\bm Q_0\).
\end{assumption}

When \(\bm Q_0\) is singular, two asymptotically distinct situations may arise according to whether the finite-\(n\) population design, defined by
\[
\bm Q_{0n}:=\E[\bm Q_n]=\E[\bm X'\bm X/n],
\]
eventually has the same range as \(\bm Q_0\), or instead has strictly larger range for all large \(n\).

\begin{definition}[Singular and nearly-singular designs]\label{def:irregular designs}
Suppose Assumption~\ref{ass:limit_design_irregular} holds with singular limit design \(\bm Q_0\), and let \(\bm Q_{0n}\to\bm Q_0\).
\begin{enumerate}[label=(\roman*)]
\item The design is \emph{singular} if, for all sufficiently large \(n\),
\[
\Range(\bm Q_{0n})=\Range(\bm Q_0).
\]
\item The design is \emph{nearly-singular} if, for all sufficiently large \(n\),
\[
\Range(\bm Q_{0n})\supsetneq \Range(\bm Q_0).
\]
\end{enumerate}
\end{definition}

These regimes create a nested pair of difficulties.
First, if the limiting design is singular, the population normal equations do not identify a unique least-squares coefficient.
Second, if the design is nearly-singular, the finite-\(n\) population problem retains directions in its range that disappear in the limit. This further instability makes the Moore--Penrose inverse discontinuous along the sequence and implies that naive plug-in procedures need not converge to any suitable population target.
These are precisely the features behind the slower rates and non-Gaussian limits derived for classical penalized least-squares estimators under nearly-singular designs; see, for example, \citet{knight2000asymptotics,knight2008shrinkage}.

Our proximal estimation approach addresses these issues sequentially.
First, we resolve identification by fixing a natural population target, the Ridgeless estimand, defined as the unique solution of the population normal equations in \(\Range(\bm Q_0)\). Among all solutions, the Ridgeless estimand is completely determined by the population regression structure, without requiring external normalization.
Second, we address instability from near singularity by constructing a rank-consistent estimator of the limiting design matrix. This ensures convergence of the associated generalized inverse and enables consistent estimation of the Ridgeless estimand via a modified Ridgeless estimator, which can serve as the benchmark estimator in a subsequent proximal correction introducing shrinkage or adaptive regularization. 
We show that rank recovery stabilizes the benchmark estimator and isolates the effect of near singularity into an explicit drift term in the first-order asymptotics, without affecting the \(\sqrt n\) convergence rate or generating the non-Gaussian, slower-rate limiting distributions typically associated with existing nearly singular penalized least-squares estimators. The subsequent proximal step inherits this first-order asymptotic behavior and, under appropriate penalization schemes, can additionally achieve variable-selection consistency and oracle-type properties.

Throughout the section we use the notation collected in SM Section~A.2.
To characterize more precisely the instability induced by near singularity, we impose an explicit local parameterization of the finite-\(n\) population design around its limit, covering regular, singular, and nearly-singular designs.

\begin{assumption}[Local parameterization of \(\bm Q_{0n}\)]\label{ass: near singularity}
There exist a symmetric matrix \(\bm\Delta\in\R^{p\times p}\) and a sequence of positive scalars \(\tau_n\to\infty\) such that
\begin{equation}\label{eq: assumption nealy singular parametrization}
\bm Q_{0n}=
\bm Q_0+\tau_n^{-1}\bm\Delta.
\end{equation}
\end{assumption}

Let
\[
\bm P_0:=\bm I_p-\bm Q_0\bm Q_0^+,
\qquad
\bm P_0^\perp:=\bm Q_0\bm Q_0^+,
\]
denote the orthogonal projectors onto \(\Kernel(\bm Q_0)\) and \(\Range(\bm Q_0)\), respectively.
Under Assumption~\ref{ass: near singularity},
\[
\bm P_0\bm Q_{0n}\bm P_0
=
\tau_n^{-1}\bm P_0\bm\Delta\bm P_0.
\]
Thus \(\bm P_0\bm\Delta\bm P_0\neq \bm 0\) implies that the finite-\(n\) population design is nondegenerate on \(\Kernel(\bm Q_0)\), even though this component disappears asymptotically.
By contrast, if \(\bm P_0\bm\Delta\bm P_0=\bm 0\), no such nondegeneracy is created on \(\Kernel(\bm Q_0)\), i.e., no first-order near-singular perturbation acts within \(\Kernel(\bm Q_0)\).\footnote{When \(\bm Q_0\) is nonsingular, \(\bm P_0=\bm 0\), so no degeneracy is present in the limit.}

Nearly-singular designs have appeared in several forms in the previous literature; see, among others,
\citet{knight2000asymptotics,phillips2001regression,knight2008shrinkage,gabaix2011rank,aswani2011regression,phillips2016inference}.
The next example gives a simple triangular-array construction under which the finite-\(n\) population design converges to a singular limit and is nearly-singular when \(\bm P_0\bm\Delta\bm P_0\neq \bm 0\).

\begin{example}[Averaged noisy covariates and near singularity]
\label{ex:averaging_me_nearsing}
Consider a setting in which a covariate vector is observed only through
repeated noisy measurements, collected over time or across sources. A standard
approach is to average these measurements to reduce measurement error. This
strategy is common in empirical work: \citet{solon1992intergenerational} uses
multiyear averages of fathers' earnings to mitigate errors-in-variables bias,
while \citet{charles2003correlation} employ multiyear income averages as
proxies for lifetime resources and averages across survey waves to reduce
noise.
Let \(\bm x_i^\star\in\mathbb R^p\) be a latent regressor satisfying the
low-rank factor structure
\[
\bm x_i^\star=\bm B\bm f_i,
\qquad
\bm f_i\in\mathbb R^k,\quad
\bm B\in\mathbb R^{p\times k},\quad
k<p.
\]
Assume that
\[
\bm\Sigma_f:=\E[\bm f_i\bm f_i']\succ \bm0,
\qquad
\Rank(\bm B)=k.
\]
Then
\[
\bm Q_0
:=
\E[\bm x_i^\star\bm x_i^{\star\prime}]
=
\bm B\bm\Sigma_f\bm B'
\]
is positive semidefinite with
\[
\Rank(\bm Q_0)=k<p,
\qquad
\Range(\bm Q_0)=\Range(\bm B).
\]
Thus the limiting design is singular.
Suppose the outcome satisfies
\[
y_i=\bm x_i^{\star\prime}\bm\beta_0+u_i,
\qquad
\E[u_i\mid \bm x_i^\star]=0.
\]
For each \(n\), the econometrician observes \(m_n\) noisy replicates,
\[
\bm x_{i,t}=\bm x_i^\star+\bm\eta_{i,t},
\qquad
t=1,\ldots,m_n,
\]
where \(m_n\to\infty\) as $n\to \infty$. Assume that observations are independent across
\(i\), and that, conditional on \((\bm x_i^\star,u_i)\), the measurement
errors \(\bm\eta_{i,1},\ldots,\bm\eta_{i,m_n}\) are independent and identically
distributed with
\begin{equation}\label{eq:me_orthogonality}
\E[\bm\eta_{i,t}\mid \bm x_i^\star,u_i]=\bm0,
\qquad
\E[\bm\eta_{i,t}\bm\eta_{i,t}']=\bm\Delta,
\qquad
\E\|\bm\eta_{i,t}\|^2<\infty .
\end{equation}
These conditions imply
\[
\E[\bm\eta_{i,t}]=\bm0,
\qquad
\E[\bm\eta_{i,t}u_i]=\bm0,
\qquad
\E[\bm\eta_{i,t}\bm x_i^{\star\prime}]=\bm0,
\]
and, for \(s\neq t\), \(\E[\bm\eta_{i,t}\bm\eta_{i,s}']=\bm0\).
The averaged proxy for \(\bm x_i^\star\) is
\[
\bm x_i
=
\frac{1}{m_n}\sum_{t=1}^{m_n}\bm x_{i,t}
=
\bm x_i^\star+\bar{\bm\eta}_i,
\qquad
\bar{\bm\eta}_i
:=
\frac{1}{m_n}\sum_{t=1}^{m_n}\bm\eta_{i,t}.
\]
Therefore
\[
\E[\bar{\bm\eta}_i\bar{\bm\eta}_i']
=
\frac{1}{m_n}\bm\Delta,
\qquad
\E[\bar{\bm\eta}_i\bm x_i^{\star\prime}]=\bm0,
\]
and the finite-\(n\) population design of the observed regressors is
\[
\bm Q_{0n}
:=
\E[\bm x_i\bm x_i']
=
\bm Q_0+m_n^{-1}\bm\Delta.
\]
Suppose that the measurement-error covariance is
nondegenerate on the unidentified directions:
\[
\bm v'\bm\Delta\bm v>0
\quad
\text{for every nonzero }
\bm v\in\Kernel(\bm Q_0).
\]
Since \(\bm\Delta\succeq\bm0\), this condition implies:\footnote{Indeed, if \(\bm v\neq\bm0\) and
\(\bm P_0^\perp\bm v\neq\bm0\), then
\(\bm v'\bm Q_0\bm v>0\); while if
\(\bm P_0^\perp\bm v=\bm0\), then
\(\bm v\in\Kernel(\bm Q_0)\) and
\(\bm v'\bm\Delta\bm v>0\).}
\[
\bm Q_{0n}
=
\bm Q_0+m_n^{-1}\bm\Delta
\succ\bm0
\qquad
\text{for every } n.
\]
Consequently,
\[
\Range(\bm Q_{0n})=\mathbb R^p
\supsetneq
\Range(\bm Q_0) ,
\]
for every \(n\), while \(\bm Q_{0n}\to\bm Q_0\) as \( n\to\infty\).
Thus the design is nearly singular in the sense of
Definition~\ref{def:irregular designs}, with local parameterization
\[
\bm Q_{0n}
=
\bm Q_0+\tau_n^{-1}\bm\Delta,
\qquad
\tau_n=m_n.
\]
\end{example}

\subsection{The Ridgeless estimand}\label{subsubsec:Ridgeless_estimand}

\begin{figure}[!t]
\centering
\includegraphics[width=0.4\textwidth]{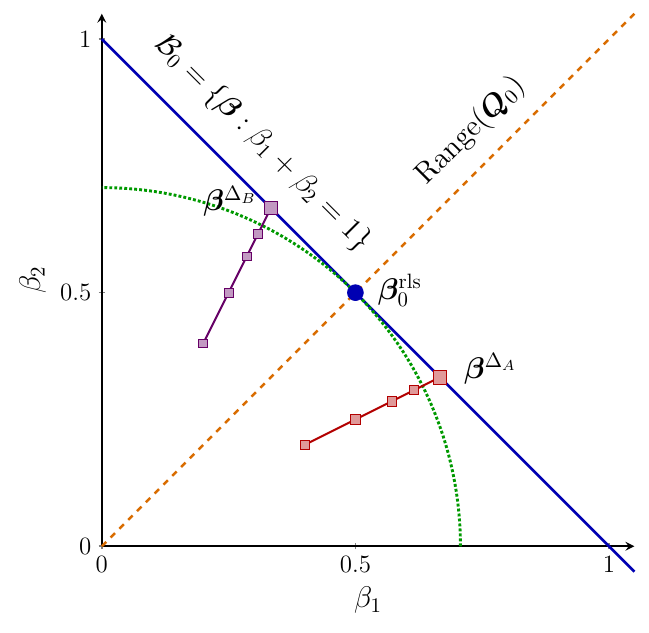}
\caption{Illustrative geometry of the Ridgeless estimand. Here \(\bm Q_0\) is a \(2\times2\) matrix of ones and \(\bm\delta_0=(1,1)'\), so \(\mathcal B_0=\{\bm\beta:\beta_1+\beta_2=1\}\) and \(\bm\beta_0^{\mathrm{rls}}=(1/2,1/2)'\). Under \(\bm Q_{0n}=\bm Q_0+\tau_n^{-1}\bm\Delta\), perturbations \(\bm\Delta_A=\diag(1,2)\) and \(\bm\Delta_B=\diag(2,1)\) yield limits \((2/3,1/3)'\) and \((1/3,2/3)'\), respectively, whereas the Ridgeless estimand depends only on \((\bm Q_0,\bm\delta_0)\).}
\label{fig:Ridgeless-path-dependence}
\end{figure}

Under irregular designs, the limiting least-squares relation need not identify
a unique regression parameter.
We first characterize the full set of population least-squares solutions and then isolate a representative that remains
point-identified and is suitable for first-order inference. To this end,
define the finite-\(n\) and limit population cross moments:
\begin{equation}\label{eq:delta0_def}
\bm\delta_{0n}:=\E[\bm X'\bm Y/n],
\qquad
\bm\delta_0:=\lim_{n\to\infty}\bm\delta_{0n},
\end{equation}
together with the limiting population least-squares risk:
\begin{equation}\label{eq:pop_ls_risk}
R_0(\bm\beta)
:=
\lim_{n\to\infty}\E\left[\|\bm Y-\bm X\bm\beta\|_2^2/n\right],
\qquad
\bm\beta\in\R^p.
\end{equation}
The following proposition characterizes the population least-squares solutions in this setting.
\begin{proposition}[Existence of population least-squares solutions]\label{prop:B0_nonempty}
Suppose Assumption~\ref{ass: near singularity} holds. If \(\bm\delta_0\) exists
and is finite and the limit \eqref{eq:pop_ls_risk} exists and is finite
for each fixed \(\bm\beta\in\R^p\), then \(\bm\delta_0\in\Range(\bm Q_0)\) and the population least-squares problem admits the set of solutions:
\begin{equation}\label{eq:B0_def}
\mathcal B_0
:=
\argmin_{\bm\beta\in\R^p}R_0(\bm\beta)
=
\{\bm\beta\in\R^p:\ \bm Q_0\bm\beta=\bm\delta_0\}\neq\emptyset.
\end{equation}
\end{proposition}
Under a regular design,
\(\bm Q_0\) is nonsingular
and Proposition \ref{prop:B0_nonempty} yields a unique limit population least squares solution:
\(\mathcal B_0=\{\bm Q_0^{-1}\bm\delta_0\}\). In contrast, when
\(\bm Q_0\) is singular then
\begin{equation}\label{eq:B0_param}
\mathcal B_0
=
\{\bm Q_0^+\bm\delta_0\}+\Kernel(\bm Q_0),
\end{equation}
so the limiting normal equations characterize only an affine solution set.
By construction, $\bm Q_0^+\bm\delta_0\in \Range(\bm Q_0)$ and
the indeterminacy in this set is exactly characterized by \(\Kernel(\bm Q_0)\), i.e., directions outside the limiting design span.
In particular, $\bm Q_0^+\bm\delta_0$ is the only element of ${\cal B}_0$ in $\Range(\bm Q_0)$.
Importantly, any element of $\mathcal B_0$ produces the 
same limiting fitted values relation:
\begin{equation*}\label{eq:identified_linear_predictor_Ridgeless}
\lim_{n\to\infty}
\E\left[\|\bm X(\bm\beta-\bm Q_0^{+}\bm\delta_0)\|_2^2/n\right]
=
(\bm\beta-\bm Q_0^{+}\bm\delta_0)'\bm Q_0(\bm\beta-\bm Q_0^{+}\bm\delta_0)
=
0,
\end{equation*}
for any \(\bm\beta\in\mathcal B_0\).
Thus, the limiting fit and prediction problem identifies the affine class \(\mathcal B_0\), but it does not distinguish among its elements. We therefore select a specific representative of \(\mathcal B_0\) as the target of estimation and inference in the irregular-design analysis.
\begin{definition}[Ridgeless estimand]\label{def:Ridgeless_estimand}
The Ridgeless estimand is the unique minimizer of the population
least-squares risk over \(\Range(\bm Q_0)\):
\begin{equation}\label{eq:Ridgeless_estimand_def}
\bm\beta_0^{\mathrm{rls}}
:=
\argmin_{\bm\beta\in\Range(\bm Q_0)}R_0(\bm\beta)
=
\bm Q_0^+\bm\delta_0.
\end{equation}
Equivalently,
it is 
the element of \(\mathcal B_0\) with minimum Euclidean norm:\footnote{Indeed, any
\(\bm\beta\in\mathcal B_0\) can be written as
\(\bm\beta=\bm Q_0^+\bm\delta_0+\bm v\) with
\(\bm v\in\Kernel(\bm Q_0)\).
Since \(\bm Q_0^+\bm\delta_0\in
\Range(\bm P_0^\perp)\)
and \(\bm v\in\Range(\bm P_0)\), this decomposition is orthogonal:
\[
\|\bm\beta\|_2^2
=
\|\bm Q_0^+\bm\delta_0\|_2^2+\|\bm v\|_2^2.
\]
Hence, the minimum norm over ${\cal B}_0$ is attained uniquely at \(\bm v=\bm 0\).}
\[
\bm\beta_0^{\mathrm{rls}}
=
\argmin_{\bm\beta\in\mathcal B_0}\|\bm\beta\|_2.
\]
\end{definition}
Hence, throughout the irregular-design analysis below, estimation and inference are directed at \(\bm\beta_0^{\mathrm{rls}}\), which may differ from the structural parameter \(\bm\beta_0\) when the latter is not point identified. The Ridgeless estimand is natural for the following reasons.
\begin{enumerate}[label=(\roman*)]
\item \textbf{Continuity with the regular case.}
If \(\bm Q_0\) is nonsingular, then \(\bm Q_0^+=\bm Q_0^{-1}\), so
\(\bm\beta_0^{\mathrm{rls}}\) coincides with the usual population least-squares coefficient.

\item \textbf{Intrinsic selection from the solution set and path invariance.}
The Ridgeless estimand is the unique element of \(\mathcal B_0\) that is fully determined by the limiting least-squares objects \((\bm Q_0,\bm\delta_0)\), without auxiliary normalization. Consequently, it is invariant to the particular nearly-singular sequence \((\bm Q_{0n}, \bm\delta_{0n})\) approaching \((\bm Q_0,\bm\delta_0)\). This path-invariant property is not shared, for example, by the finite-$n$ population least-squares estimand
\[
\bm\beta_{0n}^{\mathrm{ls}}:=\bm Q_{0n}^{+}\bm\delta_{0n}.
\]
Indeed, its limit can depend on the distortion matrix \(\bm\Delta\) under the local parameterization of Assumption~\ref{ass: near singularity},
even when the singular limit objects \((\bm Q_0,\bm\delta_0)\) are held fixed; see 
Figure~\ref{fig:Ridgeless-path-dependence} for an illustration.

\item \textbf{Compatibility with \(\sqrt n\)-regular inference.}
When \(\bm Q_0\) is singular, \citet[Thm.~2.1]{van1991differentiable} implies that a \(\sqrt n\)-regular estimand must be insensitive to
parameter perturbations along \(\Kernel(\bm Q_0)\). SM Section~B.6
shows that, among linear selection rules on \(\mathcal B_0\), the Ridgeless estimand is the unique one with this property.  It is therefore the unique least-squares representative compatible with regular inference, without requiring any auxiliary normalization.
\end{enumerate}

The preceding properties are complemented by a natural structural interpretation, after we impose the following asymptotic orthogonality condition.
\begin{assumption}[Asymptotic orthogonality of the cross-product term]\label{ass:orthogonality_beta0}
The cross moment \(\E[\bm X'\bm\varepsilon/n]\) is asymptotically negligible:
\(\E[\bm X'\bm\varepsilon/n]\to \bm 0\).
\end{assumption}
Together with Assumption~\ref{ass: near singularity}, this condition implies that
\(\bm\beta_0\) satisfies the limiting population least-squares normal equations:
\[
\bm\delta_0
=
\lim_{n\to\infty}\E[\bm X'\bm Y/n]
=
\lim_{n\to\infty}
\left\{
\bm Q_{0n}\bm\beta_0+\E[\bm X'\bm\varepsilon/n]
\right\}
=
\bm Q_0\bm\beta_0 .
\]
Hence \(\bm\beta_0\in\mathcal B_0\). Consequently,
\[
\bm\beta_0^{\mathrm{rls}}
=
\bm Q_0^+\bm\delta_0
=
\bm Q_0^+\bm Q_0\bm\beta_0
=
\bm P_0^\perp\bm\beta_0 .
\]
Thus the Ridgeless estimand is the component of the structural coefficient
identified by the limiting normal equations: it is the Euclidean projection of
\(\bm\beta_0\) onto \(\operatorname{Range}(\bm Q_0)\). Equivalently,
\[
\bm\beta_0
=
\bm\beta_0^{\mathrm{rls}}+\bm P_0\bm\beta_0,
\qquad
\bm P_0\bm\beta_0\in\Kernel(\bm Q_0),
\]
where the second component is left unrestricted by the limiting least-squares
relation.
Example~\ref{ex:averaging_me_nearsing} illustrates this interpretation: in that
setting, the Ridgeless estimand is the projection of the structural coefficient
\(\bm\beta_0\) onto \(\operatorname{Range}(\bm B)\), the span of the latent
factor loadings.

\begin{example}[Ridgeless estimand under vanishing measurement error]\label{ex:Ridgeless_me_projection}
Under the framework of Example~\ref{ex:averaging_me_nearsing},
the orthogonality conditions \eqref{eq:me_orthogonality} imply
\[
\bm\delta_{0n}
=
\E[(\bm x_i^\star+\bar{\bm\eta}_i)(\bm x_i^{\star\prime}\bm\beta_0+u_i)]
=
\bm Q_0\bm\beta_0
=:\bm\delta_0\in\Range(\bm Q_0).
\]
Hence the limiting population least-squares solutions satisfy
\[
\mathcal B_0
=
\{\bm\beta:\bm Q_0\bm\beta=\bm\delta_0\}
=
\{\bm\beta_0\}+\Kernel(\bm Q_0)
=
\{\bm\beta_0^{\mathrm{rls}}\}+\Kernel(\bm Q_0),
\]
where
\[
\bm\beta_0^{\mathrm{rls}}
=
\bm Q_0^+\bm\delta_0
=
\bm P_0^\perp\bm\beta_0.
\]
Thus \(\bm\beta_0^{\mathrm{rls}}\) is the projection of the structural
parameter \(\bm\beta_0\) onto
\(\Range(\bm Q_0)=\Range(\bm B)\). In the vanishing-measurement-error limit, the
least-squares relation therefore identifies only the component of \(\bm\beta_0\) that is spanned by the latent
factor loadings.
For each finite \(n\), let
\[
R_n(\bm\beta):=
\E\left[\|\bm Y-\bm X\bm\beta\|_2^2/n\right],
\qquad
\bm\beta_{0n}^{\mathrm{ls}}
:=
\argmin_{\bm\beta\in\R^p}R_n(\bm\beta).
\]
Since
\[
R_n(\bm\beta)
=
\E[\bm Y'\bm Y/n]-2\bm\delta_{0n}'\bm\beta+\bm\beta'\bm Q_{0n}\bm\beta ,
\]
where \(\bm\delta_{0n}=\bm\delta_0\) and
\(\bm Q_{0n}=\bm Q_0+m_n^{-1}\bm\Delta\succ \bm 0\), we further have:
\[
\bm\beta_{0n}^{\mathrm{ls}}
=
\bm Q_{0n}^{-1}\bm\delta_0,\qquad
R_n(\bm\beta)
=
R_0(\bm\beta)+\frac{1}{m_n}\bm\beta'\bm\Delta\bm\beta + c_n,
\]
where \(c_n\) is a constant not depending on \(\bm\beta\).
Thus, \(\bm\beta_{0n}^{\mathrm{ls}}\) minimizes the limiting risk
\(R_0(\bm\beta)\) contaminated by a vanishing quadratic 
term. As a result, its limit depends on
the perturbation direction \(\bm\Delta\), not only
on \((\bm Q_0,\bm\delta_0)\). This path dependence shows that the finite-\(n\) least-squares coefficient does not induce a path-invariant population target for regular asymptotic analysis. The Ridgeless estimand removes this dependence by selecting the
unique representative in \(\mathcal B_0\) that lies in \(\Range(\bm Q_0)\).
\end{example}

\subsection{Rank-consistent design estimation}\label{subsec:stability_mp_rank}

When \(\bm Q_0\) is nonsingular, inversion is continuous: if
\(\bm Q_n\to_{\Pr}\bm Q_0\), then
\(\bm Q_n^{-1}\to_{\Pr}\bm Q_0^{-1}\).
When \(\bm Q_0\) is singular, the Moore--Penrose inverse need not enjoy the same continuity:
\(\bm Q_n\to_{\Pr}\bm Q_0\) does not in general imply
\(\bm Q_n^+\to_{\Pr}\bm Q_0^+\)
\citep{Madanetal1984, andrews1987asymptotic}.
The next proposition shows that, for the design matrices considered here, stochastic continuity of the Moore--Penrose inverse is equivalent to asymptotic rank stability.

\begin{proposition}[Stochastic continuity of the Moore--Penrose inverse]\label{prop:random_MP_continuity}
Suppose Assumption~\ref{ass:limit_design_irregular} holds. Then
\[
\bm Q_n^+\to_{\Pr}\bm Q_0^+
\qquad\Longleftrightarrow\qquad
\Pr\big(\Rank(\bm Q_n)=\Rank(\bm Q_0)\big)\to 1.
\]
In particular:
\begin{enumerate}[label=(\roman*)]
\item If the design is singular, then \(\bm Q_n^+\to_{\Pr}\bm Q_0^+\).

\item If the design is nearly-singular and
\(\Pr\big(\Rank(\bm Q_n)=\Rank(\bm Q_{0n})\big)\to 1\),
then \(\bm Q_n^+\not\to_{\Pr}\bm Q_0^+\).
\end{enumerate}
\end{proposition}

Thus the key issue is not singularity itself, but whether the sample rank stabilizes at the limiting rank. Under singular designs, the finite-\(n\) population rank already coincides with the limiting rank, so the sample pseudoinverse consistently estimates \(\bm Q_0^+\). Under near singularity, by contrast, \(\bm Q_{0n}\) contains additional range directions that vanish asymptotically. Hence, if the sample rank tracks this larger finite-\(n\) rank, then it cannot stabilize at \(\Rank(\bm Q_0)\), and the Moore--Penrose plug-in is not consistent for \(\bm Q_0^+\).

This motivates replacing direct pseudoinversion with a rank-recovery step, retaining the \(r_0:=\Rank(\bm Q_0)\) directions with nonvanishing eigenvalues and discarding those that vanish asymptotically.\footnote{This reduced-rank truncation is related in spirit to \citet{lutkepohl1997modified}.}
This idea leads naturally to hard-thresholding the empirical spectrum of \(\bm Q_n\), which yields a rank-consistent estimator of \(\bm Q_0\).
Since \(\bm Q_n\) is symmetric positive semidefinite, it admits the spectral decomposition
\begin{equation}\label{eq:spectral_Qn}
\bm Q_n=\bm E_n\diag(\bm\sigma_n)\bm E_n',
\end{equation}
where \(\bm E_n\in\R^{p\times p}\) is an orthogonal matrix whose columns are the eigenvectors of the sample design matrix \(\bm Q_n\), and
\[
\bm\sigma_n=(\sigma_{n1},\dots,\sigma_{np})',
\qquad
\sigma_{n1}\ge\cdots\ge \sigma_{np}\ge 0,
\]
is the vector of ordered eigenvalues.
Given a threshold \(\nu_n>0\), define
\begin{equation}\label{eq:check_Qn_def}
\check{\bm Q}_n
:=
\bm E_n\diag(\check{\bm\sigma}_n)\bm E_n',
\end{equation}
where
\[
\check{\bm\sigma}_n
:=
\big(\sigma_{n1}\mathbf 1\{\sigma_{n1}>\nu_n\},\dots,\sigma_{np}\mathbf 1\{\sigma_{np}>\nu_n\}\big)'.
\]
Thus \(\check{\bm Q}_n\) is obtained from \(\bm Q_n\) by setting to zero all empirical eigenvalues not exceeding \(\nu_n\), while leaving the remaining eigenvalues unchanged.
The purpose of this truncation is to retain the \(r_0\) directions corresponding to the nonzero eigenvalues of \(\bm Q_0\) and discard those associated with eigenvalues that vanish asymptotically.

To make this separation operational, we complement Assumption~\ref{ass:limit_design_irregular} with a condition controlling the deviation of the sample design from its finite-\(n\) population counterpart at the usual \(\sqrt n\) scale.

\begin{assumption}\label{ass:bounded_Qn_op}
Let \(\vech(\bm A)\) denote the vector of lower-triangular entries of a symmetric matrix \(\bm A\), including the diagonal. Assume
\[
\sqrt{n}\vech(\bm Q_n-\bm Q_{0n})=O_{\Pr}(1).
\]
\end{assumption}

Assumption~\ref{ass:bounded_Qn_op} effectively sets the stochastic scale for thresholding. The threshold must dominate the sampling noise in \(\bm Q_n\), yet remain small enough to preserve the nondegenerate eigenstructure of \(\bm Q_0\).

\begin{proposition}[Consistency and rank recovery for \(\check{\bm Q}_n\)]\label{pro:cons_rank_Qn_check}
Let Assumptions~\ref{ass: near singularity} and~\ref{ass:bounded_Qn_op} hold.
Define \(\check{\bm Q}_n\) by \eqref{eq:check_Qn_def}, and let \(\nu_n>0\) satisfy
\begin{equation}\label{eq:nu_rate_conditions}
\nu_n\to 0,
\qquad
\sqrt n\nu_n\to\infty,
\qquad
\tau_n\nu_n\to\infty.
\end{equation}
Then
\[
\check{\bm Q}_n\to_{\Pr}\bm Q_0,
\qquad
\Pr\big(\Rank(\check{\bm Q}_n)=\Rank(\bm Q_0)\big)\to 1,
\qquad
\check{\bm Q}_n^+\to_{\Pr}\bm Q_0^+.
\]
\end{proposition}

The rate conditions \eqref{eq:nu_rate_conditions} have a transparent interpretation.
The requirement \(\nu_n\to 0\) ensures that no strictly positive eigenvalue of \(\bm Q_0\) is asymptotically discarded.
The condition \(\sqrt n\nu_n\to\infty\) makes the threshold dominate the stochastic fluctuation of the sample spectrum, while \(\tau_n\nu_n\to\infty\) makes it dominate the asymptotically vanishing eigenvalues generated by near singularity.
Hence \(\check{\bm Q}_n\) recovers the limiting rank and yields a pseudoinverse that consistently estimates \(\bm Q_0^+\), even when the direct Moore--Penrose plug-in \(\bm Q_n^+\) fails to do so.

\subsection{Modified Ridgeless estimator}\label{sec: Parameter of interest and modified Ridgeless estimator}

After recovering the limiting rank through the hard-thresholded design estimator \(\check{\bm Q}_n\), we can define a sample analogue of the Ridgeless estimand that remains well-defined under both singular and nearly-singular designs.
The construction replaces the Moore--Penrose plug-in based on \(\bm Q_n\), which need not converge under near singularity, with the rank-consistent design \(\check{\bm Q}_n\), thereby restricting the estimator to the stable range of the design and removing components along the estimated null space.
This leads to a modified Ridgeless estimator targeting \(\bm\beta_0^{\mathrm{rls}}\).

\begin{definition}[Modified Ridgeless estimator]\label{def:modified_Ridgeless}
The \emph{modified Ridgeless estimator} is defined as:\footnote{
   A proof of the equality in \eqref{eq:modified_Ridgeless_def} is provided in the proof of the following proposition, which reports an alternative equivalent expression of the modified Ridgeless estimator.
   \begin{proposition}\label{prop:modified_Ridgeless_repr}
Estimator \(\check{\bm\beta}_n^{\mathrm{rls}}\) in Definition~\ref{def:modified_Ridgeless} is equivalently given by
\begin{align*}
\check{\bm\beta}_n^{\mathrm{rls}}
=
\argmin_{\bm\beta\in\R^p}
\left\{
\|\bm\beta\|_2
:\ 
\check{\bm Q}_n\bm\beta
=
\check{\bm Q}_n\check{\bm Q}_n^+\bm X'\bm Y/n
\right\}.
\end{align*}
\end{proposition}
}
\begin{equation}\label{eq:modified_Ridgeless_def}
\check{\bm\beta}_n^{\mathrm{rls}}
:=\argmin_{\bm\beta\in\Range(\check{\bm Q}_n)}
\|\bm Y-\bm X\bm\beta\|_2^2/n=
\check{\bm Q}_n^+\bm X'\bm Y/n.
\end{equation}
\end{definition}

The next proposition establishes the consistency of the modified Ridgeless estimator with respect to the Ridgeless estimand, under a suitable law of large number.

\begin{proposition}[Consistency of the modified Ridgeless estimator]
\label{prop:consistency_modified_Ridgeless}
Suppose the assumptions of Proposition~\ref{pro:cons_rank_Qn_check} hold, with
\(\nu_n\) satisfying \eqref{eq:nu_rate_conditions}. Suppose further that
Assumption~\ref{ass:orthogonality_beta0} holds and that \(\bm X'\bm\varepsilon/n
\to_{\Pr}\bm 0\).
Then
\[
\check{\bm\beta}_n^{\mathrm{rls}}
\to_{\Pr}
\bm\beta_0^{\mathrm{rls}} .
\]
\end{proposition}

To derive the asymptotic distribution of \(\check{\bm\beta}_n^{\mathrm{rls}}\), we impose two additional conditions. The first is a joint weak limit for the centered sample design matrix and the cross-product term \(\bm X'\bm\varepsilon\). The second requires the kernel--range cross block of the centered sample design matrix to be negligible at the \(\sqrt n\) scale.

\begin{assumption}[Joint weak limit for key terms]\label{ass: joint_clt_Q_score}
There exist a symmetric random matrix \(\bm\Theta\in\R^{p\times p}\) and a random vector \(\bm Z\in\R^p\) such that
\[
\left(
\sqrt n\vech(\bm Q_n-\bm Q_{0n}),
\ \bm X'\bm\varepsilon/\sqrt n
\right)
\to_d
\left(
\vech(\bm\Theta),
\ \bm Z
\right),
\]
where \(\big(\vech(\bm\Theta)',\bm Z'\big)'\) is 
multivariate normal.
\end{assumption}

\begin{assumption}[\(\sqrt n\)-negligibility of the kernel--range cross block]\label{ass: vanishing stochastic term in modified Ridgeless asymptotics}
\[
\big\|\bm P_0(\bm Q_n-\bm Q_{0n})\bm P_0^\perp\big\|_{\op}=o_{\Pr}(n^{-1/2}).
\]
\end{assumption}

Assumption~\ref{ass: vanishing stochastic term in modified Ridgeless asymptotics} implies that the cross-block term \(\bm P_0(\bm Q_n-\bm Q_{0n})\bm P_0^\perp\) is negligible at the \(\sqrt{n}\)-scale. This condition will be used to simplify the asymptotic distribution of \(\check{\bm\beta}_n^{\mathrm{rls}}\) in the nearly-singular setting.
As the next example illustrates, it holds in triangular-array settings where these cross second moments arise from an averaged noise component.

\begin{example}[Verification in Example~\ref{ex:averaging_me_nearsing}]
\label{rem:verify_crossblock_me}
In the framework of Example~\ref{ex:averaging_me_nearsing},
\(\bm x_i^\star=\bm B\bm f_i\in\Range(\bm B)=\Range(\bm Q_0)\) almost surely. Therefore, \(\bm P_0\bm x_i^\star=\bm 0\) almost surely and
\[
\bm P_0(\bm Q_n-\bm Q_{0n})\bm P_0^\perp
=
\frac1n\sum_{i=1}^n
\bm P_0\bar{\bm\eta}_i\bm x_i^{\star\prime}\bm P_0^\perp
+
\frac1n\sum_{i=1}^n
\bm P_0\Bigl(\bar{\bm\eta}_i\bar{\bm\eta}_i'
-m_n^{-1}\bm\Delta
\Bigr)\bm P_0^\perp.
\]
Under the preceding conditional mean-zero and conditional independence
conditions, and assuming finite fourth moments, the first term is an average
of mean-zero matrices with entrywise variance of order \(m_n^{-1}\):
\begin{eqnarray}
 \frac1n\sum_{i=1}^n
\bm P_0\bar{\bm\eta}_i\bm x_i^{\star\prime}\bm P_0^\perp
= O_{\Pr} ((n m_n)^{-1/2})   .
\end{eqnarray}
The second term is a sample mean of mean zero matrices with variance of order $m_n^{-2}$. Therefore,
\begin{eqnarray}
 \frac1n\sum_{i=1}^n
\bm P_0\Bigl(\bar{\bm\eta}_i\bar{\bm\eta}_i'
-m_n^{-1}\bm\Delta
\Bigr)\bm P_0^\perp = O_{\Pr} (n^{-1/2}m_n^{-1})   .
\end{eqnarray}
We conclude that, under the cross-sectional independence and conditional
i.i.d.\ measurement error assumptions imposed in
Example~\ref{ex:averaging_me_nearsing}, it follows:
\[
\big\|\bm P_0(\bm Q_n-\bm Q_{0n})\bm P_0^\perp\big\|_{\op}
=
O_{\Pr}\bigl((nm_n)^{-1/2}\bigr).
\]
Specifically,
Assumption~\ref{ass: vanishing stochastic term in modified Ridgeless asymptotics}
holds whenever \(m_n\to\infty\).
\end{example}

The next theorem states the asymptotic distribution of the modified Ridgeless estimator.

\begin{theorem}\label{thm:asy_modified_Ridgeless}
Suppose that Assumptions~\ref{ass:limit_design_irregular},
\ref{ass: near singularity},
\ref{ass:orthogonality_beta0}, and
\ref{ass: joint_clt_Q_score} hold, and that the conditions of
Proposition~\ref{pro:cons_rank_Qn_check} are satisfied. Then the following statements hold:
\begin{enumerate}[label=(\roman*)]
\item If the design is singular,
\[
\sqrt n\big(\check{\bm\beta}_n^{\mathrm{rls}}-\bm\beta_0^{\mathrm{rls}}\big)
\to_d
\bm Q_0^+\bm Z.
\]

\item If the design is nearly-singular and \(\sqrt n/\tau_n\to c\in[0,\infty)\),
\[
\sqrt n\big(\check{\bm\beta}_n^{\mathrm{rls}}-\bm\beta_0^{\mathrm{rls}}\big)
\to_d
\bm P_0(\bm\Theta+c\bm\Delta)\bm Q_0^+\bm\beta_0^{\mathrm{rls}}
+
\bm Q_0^+(\bm\Theta+c\bm\Delta)\bm P_0\bm\beta_0
+
\bm Q_0^+\bm Z.
\]

\item If
Assumption~\ref{ass: vanishing stochastic term in modified Ridgeless asymptotics}
holds, then the limit in part~(ii) simplifies to:
\[
\sqrt n\big(\check{\bm\beta}_n^{\mathrm{rls}}-\bm\beta_0^{\mathrm{rls}}\big)
\to_d
c(\bm P_0\bm\Delta\bm Q_0^+\bm\beta_0^{\mathrm{rls}}
+
\bm Q_0^+\bm\Delta\bm P_0\bm\beta_0 )
+
\bm Q_0^+\bm Z.
\]
\end{enumerate}
\end{theorem}

Theorem~\ref{thm:asy_modified_Ridgeless} shows that \(\check{\bm\beta}_n^{\mathrm{rls}}\) admits a \(\sqrt n\)-weak limit approximation under both singular and nearly-singular designs.
Under singularity, the limit law is given by the Gaussian term \(\bm Q_0^+\bm Z\), which is entirely concentrated on \(\Range(\bm Q_0)\)
and is centered at $\bm 0$ because typically $\mathbb E [{\bm X' \bm\epsilon}/\sqrt{n}]
= o (1)$ under a singular design.
Thus, when the population design is singular, the modified Ridgeless estimator retains a regular first-order approximation fully supported on the identified subspace. This asymptotic behavior is fully analogous to the standard nonsingular design case, except that it takes place on the lower-dimensional identified subspace \(\Range(\bm Q_0)\).

Under near singularity, this limit law
is more involved and not fully supported on \(\Range(\bm Q_0)\). In general, it is perturbed by two additional first-order components, one in the kernel and one in the range of
$\bm Q_0$.
The term
\[
\bm P_0(\bm\Theta+c\bm\Delta)\bm Q_0^+\bm\beta_0^{\mathrm{rls}}
\]
belongs to \(\Kernel(\bm Q_0)\). It arises from the interaction of the nearly-singular perturbation and the identified Ridgeless estimand \(\bm\beta_0^{\mathrm{rls}}\in\Range(\bm Q_0)\). Conversely, the term
\[
\bm Q_0^+(\bm\Theta+c\bm\Delta)\bm P_0\bm\beta_0
\]
belongs to \(\Range(\bm Q_0)\) and is induced by the kernel component
\(\bm P_0\bm\beta_0\) of the structural parameter \(\bm\beta_0\in\mathcal B_0\).
Furthermore, under a nearly singular design the limit quantity $\bm Q_0^+ \bm Z$ may involve a Gaussian term $\bm Z$ with non-zero mean,
because the design may imply $\mathbb E [{\bm X' \bm\epsilon}/{\sqrt{n}}]
\not \to \bm 0$.

Under Assumption~\ref{ass: vanishing stochastic term in modified Ridgeless asymptotics}, the stochastic cross-block contributions generated by the sample fluctuation \(\bm Q_n-\bm Q_{0n}\) vanish asymptotically.
In this case, the first-order effect of near singularity comes from the deterministic local perturbation $\bm \Delta = \tau_n (\bm Q_{0n}-\bm Q_0)$.
Accordingly,
the limit law reduces to the sum of the deterministic drift
\[c (\bm P_0\bm\Delta\bm Q_0^+\bm\beta_0^{\mathrm{rls}}+ \bm Q_0^+\bm\Delta\bm P_0\bm\beta_0) ,\]
which captures only the impact of the nearly-singular population design, and the Gaussian component \(\bm Q_0^+\bm Z\).
Hence the identified directions continue to exhibit a Gaussian first-order approximation, which is complemented by an additional deterministic displacement. This contrasts sharply with classical nearly-singular PLSE asymptotics (e.g., \citet{knight2000asymptotics, knight2008shrinkage}), where the convergence rate is slower and the asymptotic behavior is governed entirely by the kernel directions; a detailed comparison is found in
SM Section~B.7.

\begin{example}[Asymptotic distribution of modified Ridgeless estimator in Example~\ref{ex:averaging_me_nearsing}]\label{remark: Shift of the -Gaussian component under near singularity}
Under near singularity, the local perturbation of the population design may induce a nonzero first-order mean in the limit vector \(\bm Z\). In that case, the term
\(
\bm Q_0^+\bm Z
\)
inherits a non-zero mean component in \(\Range(\bm Q_0)\).
This is illustrated by the benchmark triangular-array model of Example~\ref{ex:averaging_me_nearsing}, where
\[
\mathbb E[\bm X_i\epsilon_i]
=
-\mathbb E[\bm\eta_i\bm\eta_i']\bm\beta_0
=
-\frac{\bm\Delta}{\tau_n}\bm\beta_0,
\qquad
\tau_n:=m_n.
\]
Hence, whenever \(\sqrt n/\tau_n\to c>0\),
\[
\frac{1}{\sqrt n}\bm X'\bm\epsilon
=
\frac{1}{\sqrt n}\sum_{i=1}^n \bm X_i\epsilon_i
\to_d
\bm Z,
\]
where \(\bm Z\) is Gaussian with mean \(\mathbb E[\bm Z]
=
-c\bm\Delta\bm\beta_0\).
Accordingly,
\(
\bm Q_0^+\bm Z
\)
has mean
\(
-c\bm Q_0^+\bm\Delta\bm\beta_0
\), so this Gaussian component itself carries a first-order shift generated by the nearly-singular sequence. 
In this case, the asymptotic distribution
in Theorem \ref{thm:asy_modified_Ridgeless} 
simplifies as
\begin{eqnarray*}
\sqrt n\big(\check{\bm\beta}_n^{\mathrm{rls}}-\bm\beta_0^{\mathrm{rls}}\big)
&\to_d&
c (\bm P_0 \bm\Delta\bm Q_0^+\bm\beta_0^{\mathrm{rls}}
-  \bm Q_0^+ \bm\Delta \bm P_0^\perp \bm\beta_0 
)
+
\bm Q_0^+\tilde {\bm Z}\\
&=&
c (\bm P_0 \bm\Delta\bm Q_0^+
-  \bm Q_0^+ \bm\Delta 
)\bm\beta_0^{\mathrm{rls}}
+
\bm Q_0^+\tilde {\bm Z}\ ,
\end{eqnarray*}
with zero-mean Gaussian variable $\tilde {\bm Z}=\bm Z-\E[\bm Z]$,
since Assumption \ref{ass: vanishing stochastic term in modified Ridgeless asymptotics} holds in this setting (see again Example \ref{rem:verify_crossblock_me}).
The limit distribution therefore consists of a zero-mean Gaussian
stochastic component in \(\Range(\bm Q_0)\) and a deterministic
drift term generated by the limiting effect of the near singularity.
This drift vanishes whenever the singularity disappears sufficiently
fast, namely \(\sqrt n=o(\tau_n)\). Otherwise, the near-singular
perturbation generally induces a persistent first-order shift through the
interaction between the singular and regular subspaces of \(\bm Q_0\).
In particular, since \(\bm\Delta\) is positive definite on
\(\Kernel (\bm Q_0)\), the drift disappears only when the perturbation is
orthogonal to \(\Range(\bm Q_0)\), that is,
\(
\bm \Delta \bm P_0^\perp=\bm 0
\).
In this case, the perturbation acts only along kernel directions and
therefore does not generate any first-order asymptotic shift in the
regular component of the estimator.
\end{example}

\subsection{Oracle proximal estimation under irregular designs}\label{sec: Oracle proximal estimation for irregular designs}

The modified Ridgeless estimator
\(\check{\bm\beta}_n^{\mathrm{rls}}\) provides a benchmark estimator
for the Ridgeless estimand \(\bm\beta_0^{\mathrm{rls}}\) under singular
and nearly-singular designs. We now study a further proximal
regularization of this benchmark. We first record the resulting
asymptotic distribution inherited from the general proximal asymptotic
theory. We then discuss how to obtain the oracle property in both regular and irregular designs, under an additional active-block compatibility condition
when the limiting design is irregular.

To this end, let \(f_n\in\Gamma(\mathbb R^p)\), \(\lambda_n>0\), and define
\[
\overline{\check{\bm Q}}_n
:=
\check{\bm Q}_n+\bm I_p-\check{\bm Q}_n\check{\bm Q}_n^+
=
\check{\bm Q}_n+\check{\bm P}_n .
\]
Since \(\check{\bm P}_n\) is the orthogonal projector onto
\(\Kernel(\check{\bm Q}_n)\), the matrix
\(\overline{\check{\bm Q}}_n\) is positive definite. The corresponding
proximal estimator is
\begin{equation}\label{PLSERR2}
\check{\bm\beta}_n^{+}
:=
\prox_{\lambda_n f_n}^{\overline{\check{\bm Q}}_n}
\big(\check{\bm\beta}_n^{\mathrm{rls}}\big)
=
\argmin_{\bm\beta\in\mathbb R^p}
\left\{
\frac12
\|\check{\bm\beta}_n^{\mathrm{rls}}-\bm\beta\|_{\overline{\check{\bm Q}}_n}^2
+
\lambda_n f_n(\bm\beta)
\right\}.
\end{equation}
By Proposition~\ref{pro:cons_rank_Qn_check},
\[
\overline{\check{\bm Q}}_n
\to_{\Pr}
\overline{\bm Q}_0
:=
\bm Q_0+\bm P_0,
\]
and \(\overline{\bm Q}_0\) is positive definite. Therefore the general
results of Section~\ref{sec:asymptotics} apply with
\[
\hat{\bm\beta}_n^s=\check{\bm\beta}_n^{\mathrm{rls}},
\qquad
\bm W_n=\overline{\check{\bm Q}}_n,
\qquad
\bm\beta_0=\bm\beta_0^{\mathrm{rls}},
\qquad
r_n=\sqrt n,
\]
provided the benchmark modified Ridgeless limit law holds; here, the equality \(\bm\beta_0=\bm\beta_0^{\mathrm{rls}}\) means that the generic target \(\bm\beta_0\) in Section~\ref{sec:asymptotics} is the Ridgeless estimand in the present irregular-design application. This gives
the following corollary.

\begin{corollary}[Proximal weak limit under irregular designs]\label{cor:asymptotics_irregular}
Let Assumption~\ref{ass:f} hold, and assume the conditions of Theorem~\ref{thm:asy_modified_Ridgeless} under which
\[
\sqrt n\big(\check{\bm\beta}_n^{\mathrm{rls}}-\bm\beta_0^{\mathrm{rls}}\big)
\to_d
\bm Q_0^+\bm Z,
\]
where \(\bm Z\) is as in Assumption~\ref{ass: joint_clt_Q_score}.\footnote{
    For singular designs, this conclusion is exactly Theorem~\ref{thm:asy_modified_Ridgeless}(i).
    For nearly-singular designs, Theorem~\ref{thm:asy_modified_Ridgeless}(iii)
    gives the displayed limit whenever the additional first-order drift terms vanish.
    This holds, for instance, if \(\sqrt n=o(\tau_n)\). It also holds if
    \(\bm Q_0\bm\Delta=\bm0\), since then, by symmetry,
    \(\bm\Delta\bm Q_0^+=\bm0\) and \(\bm Q_0^+\bm\Delta=\bm0\).
    In the example of Remark~\ref{remark: Shift of the -Gaussian component under near singularity},
    \(\sqrt n=o(\tau_n)\) also implies \(\mathbb E[\bm Z]=\bm0\).
}
Then:
\begin{enumerate}[label=(\roman*)]
\item\label{prop_distr_i_irregular}
If \(\lambda_n\sqrt n\to\lambda_0>0\) and \(q_n \to_{\Pr} \rho_{\bm\beta_0^{\mathrm{rls}}}\) in epigraph,
then
\begin{equation}\label{eq:asydistr_i_irregular}
\sqrt n\big(\check{\bm\beta}_n^{+}-\bm\beta_0^{\mathrm{rls}}\big)
\to_d
\prox_{\lambda_0\rho_{\bm\beta_0^{\mathrm{rls}}}}^{\overline{\bm Q}_0}(\bm Q_0^+\bm Z).
\end{equation}
If, in addition, \(\partial f_0(\bm\beta_0^{\mathrm{rls}})\neq\emptyset\), then
\begin{equation}\label{eq:asym_moreau_i_irregular}
\prox_{\lambda_0\rho_{\bm\beta_0^{\mathrm{rls}}}}^{\overline{\bm Q}_0}(\bm Q_0^+\bm Z)
=
\big(\Id-P_{\lambda_0\partial f_0(\bm\beta_0^{\mathrm{rls}})}^{\overline{\bm Q}_0}\big)(\bm Q_0^+\bm Z).
\end{equation}

\item\label{prop_distr_ii_irregular}
If \(\lambda_n\sqrt n\to 0\) and \((\lambda_n\sqrt n)q_n \to_{\Pr} \sigma_{N_{\dom(f_0)}(\bm\beta_0^{\mathrm{rls}})}\) in epigraph,
then
\begin{equation}\label{eq:asydistr_ii_irregular}
\sqrt n\big(\check{\bm\beta}_n^{+}-\bm\beta_0^{\mathrm{rls}}\big)
\to_d
\prox_{\sigma_{N_{\dom(f_0)}(\bm\beta_0^{\mathrm{rls}})}}^{\overline{\bm Q}_0}(\bm Q_0^+\bm Z),
\end{equation}
where
\begin{equation}\label{eq:asym_moreau_ii_irregular}
\prox_{\sigma_{N_{\dom(f_0)}(\bm\beta_0^{\mathrm{rls}})}}^{\overline{\bm Q}_0}(\bm Q_0^+\bm Z)
=
\big(\Id-P_{N_{\dom(f_0)}(\bm\beta_0^{\mathrm{rls}})}^{\overline{\bm Q}_0}\big)(\bm Q_0^+\bm Z).
\end{equation}
\end{enumerate}
\end{corollary}

We next turn to oracle estimation. Let
\[
\mathcal A
:=
\{j\in\{1,\ldots,p\}:\beta_{0j}^{\mathrm{rls}}\neq0\}.
\]
Under an irregular design, the active block \((\bm Q_0)_{\mathcal A}\)
need not be nonsingular. Hence the natural active-block target is not
defined through an ordinary inverse, but through the \emph{reduced Ridgeless
target}
\begin{eqnarray}
[(\bm Q_0)_{\mathcal A}]^+(\bm\delta_0)_{\mathcal A}.\label{eq: def reduced ridgeless target}
\end{eqnarray}
For oracle efficiency relative to the full Ridgeless estimand, this
reduced target must coincide with the active coordinates of
\(\bm\beta_0^{\mathrm{rls}}\). We therefore impose the following
compatibility condition.

\begin{assumption}\label{ass: ridgless oracle compatibility condition}
The active coordinates of the Ridgeless estimand satisfy
\[
(\bm\beta_0^{\mathrm{rls}})_{\mathcal A}
\in
\Range\big((\bm Q_0)_{\mathcal A}\big).
\]
\end{assumption}

From our earlier findings, oracle behavior may follow in the vanishing-penalty regime of
Corollary~\ref{cor:asymptotics_irregular}\ref{prop_distr_ii_irregular}
using a suitable adaptive choice of penalty weights. To this end,
we first introduce the Adaptive Lasso penalty
\begin{equation}\label{eq:penalty_Ridgeless_adaptive}
f_n(\bm\beta)
:=
\sum_{j:\check\beta_{n,j}^{\mathrm{rls}}\neq 0}
|\check\beta_{n,j}^{\mathrm{rls}}|^{-1}|\beta_j|
+
\sum_{j:\check\beta_{n,j}^{\mathrm{rls}}=0}
\iota_{\{0\}}(\beta_j),
\end{equation}
and its corresponding limit penalty:
\begin{equation}\label{eq:limit_penalty_Ridgeless_adaptive}
f_0(\bm\beta)
:=
\sum_{j:\beta_{0j}^{\mathrm{rls}}\neq 0}|\beta_{0j}^{\mathrm{rls}}|^{-1}|\beta_j|
+
\sum_{j:\beta_{0j}^{\mathrm{rls}}=0}\iota_{\{0\}}(\beta_j).
\end{equation}
The next corollary characterizes the conditions under which this adaptive proximal construction recovers the oracle law relative to the Ridgeless target
using a proximal estimator of the form \eqref{PLSERR2}.

\begin{corollary}[Oracle proximal estimation under irregular designs]
\label{cor:oracle_prox_est}
Assume the conditions of Theorem~\ref{thm:asy_modified_Ridgeless} under which
\[
\sqrt n\big(\check{\bm\beta}_n^{\mathrm{rls}}-\bm\beta_0^{\mathrm{rls}}\big)
\to_d
\bm Q_0^+\bm Z,
\]
where \(\bm Z\) is as in Assumption~\ref{ass: joint_clt_Q_score} and such that
\[
\bm Z\sim\mathcal N(\bm 0,\sigma^2\bm Q_0)
\]
for some \(\sigma>0\). Let \(\check{\bm\beta}_n^+\) be the proximal estimator
\eqref{PLSERR2} based on the adaptive penalty
\eqref{eq:penalty_Ridgeless_adaptive}, and define
\[
\check{\mathcal A}_n^+
:=
\{j\in\{1,\ldots,p\}:(\check{\bm\beta}_n^+)_j\neq 0\},
\qquad
\mathcal A
:=
\{j\in\{1,\ldots,p\}:(\bm\beta_0^{\mathrm{rls}})_j\neq 0\}.
\]
If \(\lambda_n\sqrt n\to 0\) and \(\lambda_n n\to\infty\),
then
\begin{equation}\label{eq:oracle_limit_oracle_prox_consistency}
\Pr(\check{\mathcal A}_n^+=\mathcal A)\to 1.
\end{equation}
Moreover,
\begin{equation}\label{eq:oracle_limit_oracle_prox}
\sqrt n(\check{\bm\beta}_n^+-\bm\beta_0^{\mathrm{rls}})_{\mathcal A}
\to_d
[(\bm Q_0)_{\mathcal A}]^+(\bm Z)_{\mathcal A}
\end{equation}
holds if and only if the covariance identity \eqref{eq:oracle_cov_identity}
of Proposition~\ref{prop:op2} holds with
\(\bm M_0=\bm Q_0\) and \(\bm W_0=\overline{\bm Q}_0\).
If, in addition, Assumption~\ref{ass: ridgless oracle compatibility condition}
holds, then the limiting distribution in
\eqref{eq:oracle_limit_oracle_prox} is the efficient oracle distribution for
the active coordinates of the full Ridgeless estimand.
\end{corollary}
Corollary~\ref{cor:oracle_prox_est} requires additional conditions beyond the benchmark \(\sqrt n\)-limit. The restriction
\[
\bm Z\sim\mathcal N(\bm0,\sigma^2\bm Q_0)
\]
has two roles. First, it places the benchmark term
\(\bm Q_0^+\bm Z\) within the centered Gaussian framework of
Proposition~\ref{prop:op2}. Second, it aligns the covariance of the score
limit with the limit design. The remaining covariance identity
\eqref{eq:oracle_cov_identity}, specialized to
\(\bm M_0=\bm Q_0\) and \(\bm W_0=\overline{\bm Q}_0\), is then the exact
condition ensuring that the \(\overline{\bm Q}_0\)-projection limit in
Corollary~\ref{cor:asymptotics_irregular}\ref{prop_distr_ii_irregular}
coincides with the oracle Gaussian law on the active block.
This last requirement is substantive: although the proximal step with limiting metric \(\overline{\bm Q}_0\) achieves consistent support recovery under the stated tuning conditions, it need not preserve the efficient covariance structure on the selected coordinates.

\begin{remark}\label{remark: remark on challenge for oracle under near singularity 1}
A sufficient structural condition for oracle efficiency in Corollary~\ref{cor:oracle_prox_est} is that the unidentified directions are supported only on inactive coordinates, namely
\begin{equation}\label{eq:kernel_inactive_condition}
\Kernel(\bm Q_0)\subseteq
\{\bm\beta\in\mathbb R^p:\bm\beta_{\mathcal A}=\bm0\}.
\end{equation}
Mathematically, condition \eqref{eq:kernel_inactive_condition} requires the
kernel of the limit design matrix to be confined to inactive coordinates. Hence, no singular direction of the population design involves components relevant for the reduced Ridgeless target. Equivalently, the active coordinate subspace is contained in the identified space \(\Range(\bm Q_0)\), so that any linear combination supported on \(\mathcal A\) remains identifiable at first order.
Intuitively, the lack of identification generated by the nearly-singular design is confined to regressors that are asymptotically irrelevant for the full Ridgeless estimand. The regressors carrying signal for the oracle parameter are therefore unaffected by the degeneracy, while singularity arises only through redundant or weakly relevant inactive variables.
Under condition \eqref{eq:kernel_inactive_condition},
\[
(\bm P_0)_{\mathcal A}=\bm0,
\qquad
(\overline{\bm Q}_0)_{\mathcal A}=(\bm Q_0)_{\mathcal A}.
\]
Thus the kernel correction vanishes on the active block, and the canonical proximal metric coincides there with the oracle design geometry. Consequently, the additional regularization needed to stabilize singular directions acts only on inactive coordinates and does not distort the first-order law of the active components.
\end{remark}
Remark~\ref{remark: remark on challenge for oracle under near singularity 1} highlights that, unless \((\bm Q_0)_{\mathcal A}\) is positive definite, constructing a single proximal estimator that jointly achieves oracle selection and oracle efficiency is challenging. The reason is that one cannot, in general, build a limiting weighting matrix that is simultaneously positive definite and exactly aligned with the geometry induced by \(\bm Q_0\) on the active block.

When \((\bm Q_0)_{\mathcal A}\succ\bm 0\), this difficulty can be overcome through a two-step procedure based on the support-adapted limiting matrix \(\bm W_0\), defined by
\[
\bm W_0 = \bm Q_0 + \bm D ({\cal A}) ,
\]
where $\bm D ({\cal A})$ is a diagonal matrix such that
\[ (\bm D ({\cal A}))_{\cal A} =\bm 0 \text{ and } 
(\bm D ({\cal A}))_{{\cal A}^c}=\bm I_{|{\cal A}^c|}.\]
By construction, matrix $\bm W_0$ is positive definite whenever \((\bm Q_0)_{\mathcal A}\succ\bm 0\). It preserves the oracle geometry on the active block while regularizing inactive directions, and can therefore deliver an efficient oracle limit law once the support has been consistently recovered.

A
consistent estimator of \(\bm W_0\) is obtained by a plug-in construction \(\check{\bm W}_n\), replacing \(\bm Q_0\) and the unknown active set \(\mathcal A\) with their consistent estimators \({\bm Q}_n\) and \(\check{\mathcal A}_n^+\), namely:
\[
\check{\bm W}_n
=
{\bm Q}_n
+\bm D {(\check{\mathcal A}_n^+)}.
\]
Building on this consistent weighting matrix estimator, the next corollary establishes an efficient oracle proximal estimator for the Ridgeless target under a nearly singular full design and a positive definite active block \((\bm Q_0)_{\mathcal A}\). In particular,
the active-block compatibility Assumption \ref{ass: ridgless oracle compatibility condition}  is then automatically satisfied.

\begin{corollary}[Efficient two-stage oracle proximal estimation]
\label{cor:oracle_prox_est_double}
Assume the conditions of Theorem~\ref{thm:asy_modified_Ridgeless} under which
\[
\sqrt n\big(\check{\bm\beta}_n^{\mathrm{rls}}-\bm\beta_0^{\mathrm{rls}}\big)
\to_d
\bm Q_0^+\bm Z,
\qquad
\bm Z\sim\mathcal N(\bm 0,\sigma^2\bm Q_0),
\]
for some \(\sigma>0\). Let \(\check{\bm\beta}_n^+\) be the first-stage adaptive proximal estimator and suppose that
\[
\Pr(\check{\mathcal A}_n^+=\mathcal A)\to1,
\qquad
\mathcal A:=\{j:(\bm\beta_0^{\mathrm{rls}})_j\neq0\}.
\]
Assume moreover that \((\bm Q_0)_{\mathcal A}\succ0\). Let \(\check{\bm W}_n\) be the support-adapted weighting matrix defined above, and define the second-stage proximal estimator
\[
\check{\check{\bm\beta}}_n^{+}
:=
\prox_{\lambda_n f_n}^{\check{\bm W}_n}\big(\check{\bm\beta}_n^{\mathrm{rls}}\big)
=
\argmin_{\bm\beta\in\mathbb R^p}
\left\{
\frac12\|\check{\bm\beta}_n^{\mathrm{rls}}-\bm\beta\|_{\check{\bm W}_n}^{2}
+
\lambda_n f_n(\bm\beta)
\right\},
\]
where \(f_n\) is the adaptive penalty in \eqref{eq:penalty_Ridgeless_adaptive}. Define
\[
\check{\check{\mathcal A}}_n^+
:=
\{j\in\{1,\ldots,p\}:(\check{\check{\bm\beta}}_n^{+})_j\neq 0\}.
\]
If \(\lambda_n\sqrt n\to0\) and \(\lambda_n n\to\infty\),
then
\[
\Pr\left(\check{\check{\mathcal A}}_n^+=\mathcal A\right)\to1,\qquad
\sqrt n\big(\check{\check{\bm\beta}}_n^{+}-\bm\beta_0^{\mathrm{rls}}\big)_{\mathcal A}
\to_d
[(\bm Q_0)_{\mathcal A}]^{+}(\bm Z)_{\mathcal A},
\]
where $[(\bm Q_0)_{\mathcal A}]^{+}=[(\bm Q_0)_{\mathcal A}]^{-1}$.
Hence, the estimator achieves oracle support recovery and the oracle active-block limit distribution.
\end{corollary}
Despite the oracle result in Corollary~\ref{cor:oracle_prox_est_double}, the conclusion remains sensitive to the requirement that \((\bm Q_0)_{\mathcal A}\) be positive definite. In nearly singular designs, the active block may itself be singular, which complicates the construction of a proximal metric that is both positive definite and preserves the oracle geometry on the selected coordinates.

To avoid this restriction, we finally consider a more robust two-step procedure. In the first stage, the adaptive proximal estimator is used solely for consistent selection of the active coordinates. In the second stage, conditional on the selected support, a modified Ridgeless estimator is refitted on the selected block.
By construction, such a refitted estimator aims to estimate the reduced Ridgeless target 
(\ref{eq: def reduced ridgeless target}), which under
compatibility Assumption~\ref{ass: ridgless oracle compatibility condition} coincides with the reduced Ridgeless parameter.
In particular, this requires a consistent estimator of the active-block Moore--Penrose inverse $[(\bm Q_0)_{\cal A}]^{+}$.
However, even when $\check{\bm Q}_n$ is a consistent and rank-consistent estimator of $\bm Q_0$, 
the corresponding submatrix $(\check{\bm Q}_n)_{\cal A}$ need not be rank-consistent for $(\bm Q_0)_{\cal A}$, as rank is not preserved under restriction. 

To address this issue, we introduce a spectrally thresholded estimator of the active block. With this purpose in mind, consider for any candidate index set $\cal S$ the spectral decomposition
\begin{eqnarray}
(\bm Q_{n})_{\cal S} = \bm E_n ({\cal S}) \diag(\bm \sigma_n ({\cal S}))
\bm E_n ({\cal S})'\ ,
\end{eqnarray}
and the associated probability limit
\begin{eqnarray}
(\bm Q_{0})_{\cal S} = \bm E_0 ({\cal S}) \diag (\bm \sigma_0 ({\cal S}))
\bm E_0 ({\cal S})'.
\end{eqnarray}
We can then define the thresholded estimator
\begin{eqnarray}
 \check {\bm Q}_n({\cal S}):=
 \bm E_n ({\cal S}) \diag(\check {\bm \sigma}_n ({\cal S}))
\bm E_n ({\cal S})'  \ ,
\end{eqnarray}
where $\check {\bm \sigma}_n ({\cal S}):=(
{ \sigma}_{n1} ({\cal S}){\bf 1}\{{ \sigma}_{n1} ({\cal S})\ge \nu_n\},\ldots, { \sigma}_{n|{\cal S}|} ({\cal S}){\bf 1} \{{ \sigma}_{n |{\cal S}|} ({\cal S})\ge \nu_n\}))'$.
Under the same assumptions as in Proposition \ref{pro:cons_rank_Qn_check} it follows that:
\begin{eqnarray}
 \check {\bm Q}_n({\cal A})\to_{\Pr} (\bm Q_0)_{\cal A}  
 \text{ and } \Pr \big(\Rank ( \check {\bm Q}_n({\cal A}))= \Rank((\bm Q_0)_{\cal A})\big)\to 1 .
\end{eqnarray}
Furthermore, under the conditions of Corollary \ref{cor:oracle_prox_est}:
\begin{eqnarray*}
 \Pr(\check{\mathcal A}_n^+=\mathcal A)\to 1.   
\end{eqnarray*}
Together, this yields the following corollary.

\begin{corollary}[Two-step oracle estimation by post-selection modified Ridgeless refit]

\label{cor:oracle_post_Ridgeless}

Assume the conditions of
Theorem~\ref{thm:asy_modified_Ridgeless}, and let
\(\check{\bm\beta}_n^+\) be the adaptive proximal estimator based on
\eqref{eq:penalty_Ridgeless_adaptive}. Suppose that
\[
\Pr(\check{\mathcal A}_n^+=\mathcal A)\to1.
\]
Define the post-selection estimator
\(\tilde{\bm\beta}_n\) by
\[
(\tilde{\bm\beta}_n)_{\check{\mathcal A}_n^+}
=
\big[\check{\bm Q}_n(\check{\mathcal A}_n^+)\big]^+
\left(\bm X'\bm Y/n\right)_{\check{\mathcal A}_n^+},
\qquad
(\tilde{\bm\beta}_n)_{(\check{\mathcal A}_n^+)^c}
=
\bm0.
\]
Assume additionally
Assumption~\ref{ass: ridgless oracle compatibility condition}. Then
\[
\Pr\left(
\{j:(\tilde{\bm\beta}_n)_j\neq0\}
=
\mathcal A
\right)\to1.
\]
Moreover, in the nearly-singular regime
\(\sqrt n/\tau_n\to c\in[0,\infty)\),
\[
\sqrt n\big(
(\tilde{\bm\beta}_n)_{\mathcal A}
-
(\bm\beta_0^{\mathrm{rls}})_{\mathcal A}
\big)
\]
admits the same first-order asymptotic expansion as in
Theorem~\ref{thm:asy_modified_Ridgeless} for the modified Ridgeless estimator on the active block.
In particular, if $\sqrt{n}=o(\tau_n)$ and
Assumption~\ref{ass: vanishing stochastic term in modified Ridgeless asymptotics}
holds for the active subblock, namely
\[
\Big\|
\bm P_0({\cal A})
(\bm Q_n-\bm Q_{0n})_{\cal A}
\bm P_0^\perp({\cal A})
\Big\|_{\op}
=
o_{\Pr}(n^{-1/2}),
\]
where $\bm P_0 ({\cal A}) := \bm I_{|{\cal A}|}- ({\bm Q_0})_{\cal A}[({\bm Q_0})_{\cal A}]^+$,
then
\[
\sqrt n\big(
\tilde{\bm\beta}_n-\bm\beta_0^{\mathrm{rls}}
\big)_{\mathcal A}
\to_d
[(\bm Q_0)_{\mathcal A}]^+
(\bm Z)_{\mathcal A},
\qquad
\bm Z  \sim\mathcal N(\bm 0,\sigma^2 \bm Q_0).
\]
Hence, the estimator achieves oracle support recovery and the oracle active-block limit distribution without requiring
\((\bm Q_0)_{\mathcal A}\succ0\).
\end{corollary}

In Corollary~\ref{cor:oracle_post_Ridgeless}, Assumption~\ref{ass: ridgless oracle compatibility condition} implies \((\bm\beta_0^{\mathrm{rls}})_{\mathcal A}
=
[(\bm Q_0)_{\mathcal A}]^+(\bm\delta_0)_{\mathcal A}\),
so that the active coordinates of the full Ridgeless estimand coincide with the reduced target.
When this condition fails, the post-selection estimator remains well defined and continues to estimate the reduced quantity $[(\bm Q_0)_{\mathcal A}]^+(\bm\delta_0)_{\mathcal A}$. It should then be interpreted as a refit within the geometry of the selected submodel, rather than as recovery of the active coordinates of the full Ridgeless estimand. The discrepancy arises because removing inactive variables alters the underlying singular structure.

% \section{Concluding remarks}\label{sec: Conclusions}
%
% By construction, the proximal framework developed in this paper extends naturally to several open problems. These include, for example, the analysis of Least Absolute Deviations proximal estimators under irregular designs and of modified Ridgeless-type instrumental variable estimators under weak identification. Another promising direction is the study of modified Ridgeless-type estimators in high-dimensional settings with nearly singular designs, as well as the role of the induced regularization under alternative structures of near singularity. Finally, the smooth dependence of proximal estimators on their initial estimator makes the framework well suited for constructing robust procedures with bounded influence functions, by combining it with appropriately chosen bounded-influence initial estimators.

\clearpage
\section*{Supplementary Material}
\setcounter{section}{0}
\renewcommand{\thesection}{\Alph{section}}
\renewcommand{\thesubsection}{\thesection.\arabic{subsection}}
\numberwithin{equation}{section}
\numberwithin{figure}{section}
\numberwithin{table}{section}
\numberwithin{theorem}{subsection}

\noindent This supplementary material is organized into five sections. Section~\ref{sec:notation} introduces the notation and conventions used throughout the article. Section~\ref{sec:appendix_theory} collects supporting results, including benchmark penalty calculations and additional asymptotic arguments. Section~\ref{sec:simulations} presents the Monte Carlo design, results, and figures. Section~\ref{sec:proofs} contains proofs of the results stated in the main article, and Section~\ref{sec:proofs_appendix} contains proofs of the supporting results established here.

\noindent\textbf{Data and code availability.}
The Monte Carlo study uses simulated data and requires no external dataset.
The R package, replication code, generated results, and figures are available at
\url{https://github.com/a91quaini/replicateProxEst}.

\section{Notation}\label{sec:notation}

\subsection{Notation for Proximal Estimators}

We collect here the notation and convex-analysis notions used
in the main text.
Standard references include \citet{SM:bauschke2016convex}, \citet{SM:rockafellarWets2009}, \citet{SM:attouch1984variational}, \citet{SM:hiriart2004fundamentals}.

\noindent \textbf{Linear algebra.}
For a matrix \(\bm M\), \(\bm M'\) denotes its transpose and \(\bm M^+\) its Moore--Penrose inverse.
We write \(\Range(\bm M)\), \(\Kernel(\bm M)\), and \(\Rank(\bm M)\) for the range, kernel, and rank of \(\bm M\), respectively.
For symmetric matrices \(\bm A,\bm B\in\R^{p\times p}\), \(\bm A\succeq \bm B\) means that \(\bm A-\bm B\) is positive semidefinite, and \(\bm A\succ \bm B\) means that \(\bm A-\bm B\) is positive definite.
The standard Euclidean inner product on \(\R^p\) is
\[
\langle \bm x,\bm y\rangle := \bm x'\bm y,
\qquad \bm x,\bm y\in\R^p,
\]
with induced Euclidean norm
\[
\|\bm x\|_2 := \sqrt{\langle \bm x,\bm x\rangle}.
\]
If \(\bm W\succ\bm 0\), we define the \(\bm W\)-weighted inner product and norm by
\[
\langle \bm x,\bm y\rangle_{\bm W} := \bm x'\bm W\bm y,
\qquad
\|\bm x\|_{\bm W} := \sqrt{\langle \bm x,\bm x\rangle_{\bm W}}.
\]
If \(\mathcal S\subseteq\R^p\) is a linear subspace, we write
\[
\mathcal S^\perp
:=
\{\bm v\in\R^p:\ \langle \bm v,\bm s\rangle_{\bm W}=0\ \text{for all }\bm s\in\mathcal S\}.
\]
When the Euclidean inner product is intended instead, we state this explicitly.
For a matrix \(\bm M\in\R^{m\times p}\), we write
\[
\|\bm M\|
:=
\sup_{\bm x\neq \bm 0}\frac{\|\bm M\bm x\|_2}{\|\bm x\|_2}
=
\sup_{\|\bm x\|_2\le 1}\|\bm M\bm x\|_2
\]
for its spectral norm, that is, the operator norm induced by the Euclidean norm.
Whenever \(\|\cdot\|\) denotes a norm on \(\R^p\), its dual norm is
\[
\|\bm y\|_*:=\sup\{\langle \bm y,\bm x\rangle:\ \|\bm x\|\le 1\},
\qquad \bm y\in\R^p.
\]
For \(r\in[1,\infty)\), the \(\ell_r\)-norm is
\[
\|\bm x\|_r:=\left(\sum_{i=1}^p |x_i|^r\right)^{1/r},
\]
and \(\|\bm x\|_\infty:=\max_{1\le i\le p}|x_i|\).
Its dual norm is \(\|\cdot\|_q\), where \(1/r+1/q=1\).
We write
\[
\mathcal B_r:=\{\bm x\in\R^p:\ \|\bm x\|_r\le 1\}
\]
for the closed unit ball of \(\|\cdot\|_r\); in particular,
\(\mathcal B_\infty=\{\bm x\in\R^p:\ \|\bm x\|_\infty\le 1\}\).

\noindent \textbf{Functions, domains, and convex conjugates.}
For a function \(f:\R^p\to(-\infty,+\infty]\), its domain is
\[
\dom(f):=\{\bm x\in\R^p:\ f(\bm x)<+\infty\}.
\]
The function \(f\) is proper if \(\dom(f)\neq\varnothing\) and \(f(\bm x)>-\infty\) for all \(\bm x\in\R^p\).
It is convex if
\[
f(\lambda \bm x+(1-\lambda)\bm y)\le \lambda f(\bm x)+(1-\lambda)f(\bm y)
\]
for all \(\bm x,\bm y\in\R^p\) and all \(\lambda\in[0,1]\).
It is lower semicontinuous if all its sublevel sets are closed.
We write
\[
\Gamma(\R^p)
\]
for the class of proper, lower semicontinuous, convex functions from \(\R^p\) to \((-\infty,+\infty]\).
A convex function \(f\) is sublinear if it is positively homogeneous and subadditive, that is,
\[
f(\alpha\bm x)=\alpha f(\bm x)\quad \text{for all }\alpha\ge 0,
\qquad
f(\bm x+\bm y)\le f(\bm x)+f(\bm y)
\]
for all \(\bm x,\bm y\in\R^p\).
Given \(f\in\Gamma(\R^p)\) and an inner product \(\langle\cdot,\cdot\rangle_{\bm W}\), the convex conjugate of \(f\) under \(\langle\cdot,\cdot\rangle_{\bm W}\) is
\[
f^{*,\bm W}(\bm y)
:=
\sup_{\bm x\in\R^p}
\big\{
\langle \bm y,\bm x\rangle_{\bm W}-f(\bm x)
\big\},
\qquad \bm y\in\R^p.
\]
Whenever the underlying inner product is clear from the context, we simply write \(f^*\).
For a set \(C\subseteq\R^p\), the indicator function \(\iota_C:\R^p\to\{0,+\infty\}\) is defined by
\[
\iota_C(\bm x)=0 \ \text{if }\bm x\in C,
\qquad
\iota_C(\bm x)=+\infty \ \text{otherwise}.
\]
The identity operator on \(\R^p\) is denoted by \(\Id\).
Given the working inner product \(\langle\cdot,\cdot\rangle_{\bm W}\), the support function of a nonempty set \(C\subseteq\R^p\) is
\[
\sigma_C(\bm x)
:=
\sup_{\bm y\in C}\langle \bm y,\bm x\rangle_{\bm W},
\qquad \bm x\in\R^p.
\]
When the Euclidean inner product is intended instead, we state this explicitly.

\noindent \textbf{Projections, proximal operators, subgradients, and normal cones.}
Given a nonempty closed convex set \(C\subseteq\R^p\) and \(\bm W\succ\bm 0\), the weighted projection of \(\bm x\in\R^p\) onto \(C\) is
\[
P_C^{\bm W}(\bm x)
:=
\argmin_{\bm y\in C}\|\bm x-\bm y\|_{\bm W}.
\]
The minimizer exists and is unique.
Given \(f\in\Gamma(\R^p)\) and \(\bm W\succ\bm 0\), the proximal operator of \(f\) with respect to \(\bm W\) is
\[
\prox_f^{\bm W}(\bm x)
:=
\argmin_{\bm y\in\R^p}
\left\{
f(\bm y)+\frac12\|\bm x-\bm y\|_{\bm W}^2
\right\},
\qquad
\bm x\in\R^p.
\]
The associated Moreau envelope is
\[
e_f^{\bm W}(\bm x)
:=
\inf_{\bm y\in\R^p}
\left\{
f(\bm y)+\frac12\|\bm x-\bm y\|_{\bm W}^2
\right\}.
\]
For \(f\in\Gamma(\R^p)\), \(\bm x\in\dom(f)\), and the working inner product \(\langle\cdot,\cdot\rangle_{\bm W}\), the subdifferential of \(f\) at \(\bm x\) is
\[
\partial f(\bm x)
:=
\Big\{
\bm v\in\R^p:\ 
f(\bm y)\ge f(\bm x)+\langle \bm v,\bm y-\bm x\rangle_{\bm W}
\ \text{for all }\bm y\in\R^p
\Big\}.
\]
Likewise, for a nonempty closed convex set \(C\subseteq\R^p\), the normal cone to \(C\) at \(\bm x\in C\) is
\[
N_C(\bm x)
:=
\Big\{
\bm v\in\R^p:\ 
\langle \bm v,\bm y-\bm x\rangle_{\bm W}\le 0
\ \text{for all }\bm y\in C
\Big\},
\]
and \(N_C(\bm x):=\varnothing\) if \(\bm x\notin C\).
When the Euclidean inner product is required explicitly, we write \(\partial^e f(\bm x)\).
For \(f:\R^p\to(-\infty,+\infty]\) and \(\bm x\in\dom(f)\), the one-sided directional derivative of \(f\) at \(\bm x\) in direction \(\bm y\in\R^p\) is
\[
\rho_{\bm x}(\bm y)
:=
\lim_{\alpha\downarrow 0}
\frac{f(\bm x+\alpha\bm y)-f(\bm x)}{\alpha},
\]
whenever the limit exists; for \(f\in\Gamma(\R^p)\), this limit exists in \((-\infty,+\infty]\).
For a nonempty set \(C\subseteq\R^p\), the Euclidean distance from \(\bm x\in\R^p\) to \(C\) is
\[
d_C(\bm x):=\inf_{\bm y\in C}\|\bm x-\bm y\|_2.
\]
When \(\bm W\succ\bm 0\), the weighted distance is
\[
d_C^{\bm W}(\bm x):=\inf_{\bm y\in C}\|\bm x-\bm y\|_{\bm W}.
\]
For a cone \(K\subseteq\R^p\), its polar cone is
\[
K^\circ
:=
\{\bm v\in\R^p:\ \langle \bm v,\bm w\rangle_{\bm W}\le 0
\ \text{for all }\bm w\in K\}.
\]
Again, this notation is understood relative to the inner product \(\langle\cdot,\cdot\rangle_{\bm W}\).

\noindent \textbf{Epigraph convergence and stochastic convergence.}
For extended-real-valued random functions \(f_n,f\in\Gamma(\R^p)\), we write
\(f_n\to_{\Pr} f\) in epigraph and \(f_n\to_d f\) in epigraph for convergence in
probability and in distribution, respectively, with respect to any metric that
metrizes epigraph convergence on \(\Gamma(\R^p)\), such as the epi-distance
\(d_l\); see \citet[Thm.~7.58]{SM:rockafellarWets2009} and
\citet{SM:salinetti1981convergence,SM:salinetti1986convergence,SM:pflug1995asymptotic,SM:vogel2003continuous,SM:lachout2006epi}.
Throughout, \(\to_{\Pr}\) and \(\to_d\) denote ordinary convergence in probability and in distribution, respectively, unless we explicitly specify that the convergence is in epigraph.
Likewise, \(X\stackrel{d}{=} Y\) means that random variables \(X\) and \(Y\) are equal in distribution.

\subsection{Notation for linear regression with irregular designs }\label{subsec:notation_irregular_lr}

Here, we collect the linear-algebraic notation used throughout the analysis of linear regression problems with irregular designs in the main text.
We continue to write \(\Range(\cdot)\), \(\Kernel(\cdot)\), \(\Rank(\cdot)\), and \((\cdot)^+\) for the range, kernel, rank, and Moore--Penrose inverse of a matrix, respectively.

\paragraph{Symmetric positive semidefinite matrices.}
Let \(\bm M\in\R^{p\times p}\) be symmetric and positive semidefinite. Here and throughout this part of the notation, \({}^\perp\), orthogonality, and orthogonal projections are understood with respect to the Euclidean inner product \(\langle\cdot,\cdot\rangle\).
Then
\[
\Range(\bm M)=\Kernel(\bm M)^\perp,
\]
and every \(\bm x\in\R^p\) can be written uniquely as
\[
\bm x=\bm x_{\Range}+\bm x_{\Kernel},
\qquad
\bm x_{\Range}\in\Range(\bm M),\quad
\bm x_{\Kernel}\in\Kernel(\bm M).
\]
Moreover,
\[
\bm M\bm M^+=\bm M^+\bm M,\qquad \bm I_p-\bm M\bm M^+
\]
are the orthogonal projectors onto \(\Range(\bm M)\) and \(\Kernel(\bm M)\), respectively.
Thus every \(\bm x\in\R^p\) admits the orthogonal decomposition
\[
\bm x=\bm M\bm M^+\bm x+(\bm I_p-\bm M\bm M^+)\bm x,
\]
where the first term lies in \(\Range(\bm M)\) and the second in \(\Kernel(\bm M)\).

\paragraph{Spectral decompositions.}
If \(\bm M\in\R^{p\times p}\) is symmetric, we write its spectral decomposition as
\[
\bm M=\bm E\,\diag(\lambda_1,\ldots,\lambda_p)\,\bm E',
\qquad
\lambda_1\ge\cdots\ge\lambda_p,
\]
where \(\bm E\) is orthogonal and \(\lambda_1,\ldots,\lambda_p\in\R\) are the eigenvalues of \(\bm M\).
If \(\bm M\succeq \bm 0\), then \(\lambda_j\ge 0\) for every \(j\), and its Moore--Penrose inverse is
\[
\bm M^+
=
\bm E\,\diag(\lambda_1^+,\ldots,\lambda_p^+)\,\bm E',
\qquad
\lambda_j^+:=
\begin{cases}
\lambda_j^{-1}, & \lambda_j>0,\\
0, & \lambda_j=0.
\end{cases}
\]
In particular, \(\Rank(\bm M)\) is the number of strictly positive eigenvalues of \(\bm M\).

\paragraph{Half-vectorization.}
For a symmetric matrix \(\bm M\in\R^{p\times p}\), we write
\[
\vech(\bm M)
\]
for the vector obtained by stacking the lower-triangular entries of \(\bm M\), including the diagonal.

\paragraph{Linear spans and affine spaces.}
For vectors \(\bm v_1,\ldots,\bm v_k\in\R^p\), we write
\[
\Span\{\bm v_1,\ldots,\bm v_k\}
\]
for their linear span.
Given a vector \(\bm x_0\in\R^p\) and a linear subspace \(\mathcal S\subseteq\R^p\), we write
\[
\bm x_0+\mathcal S:=\{\bm x_0+\bm s:\bm s\in\mathcal S\}
\]
for the corresponding affine subspace.

\paragraph{Linear systems and minimum-norm solutions.}
Let \(\bm M\in\R^{m\times p}\) and \(\bm b\in\R^m\).
The linear system
\[
\bm M\bm x=\bm b
\]
is solvable if and only if \(\bm b\in\Range(\bm M)\).
When solvable, its solution set is
\[
\{\bm x\in\R^p:\bm M\bm x=\bm b\}
=
\bm M^+\bm b+\Kernel(\bm M).
\]
The unique minimum-Euclidean-norm solution is
\[
\bm x^\dagger:=\bm M^+\bm b.
\]
Equivalently, \(\bm x^\dagger\) is the unique solution satisfying
\[
\bm x^\dagger\in\Range(\bm M')
\qquad\text{and}\qquad
\bm M\bm x^\dagger=\bm b.
\]

\paragraph{Least-squares problems.}
For arbitrary \(\bm M\in\R^{m\times p}\) and \(\bm b\in\R^m\), the least-squares problem
\[
\min_{\bm x\in\R^p}\|\bm b-\bm M\bm x\|_2^2
\]
has solution set
\[
\argmin_{\bm x\in\R^p}\|\bm b-\bm M\bm x\|_2^2
=
\{\bm x\in\R^p:\bm M'\bm M\bm x=\bm M'\bm b\}
=
\bm M^+\bm b+\Kernel(\bm M).
\]
Its unique minimum-Euclidean-norm solution is again
\[
\bm x^\dagger=\bm M^+\bm b.
\]
The associated fitted value is the orthogonal projection of \(\bm b\) onto \(\Range(\bm M)\):
\[
\bm M\bm x^\dagger
=
\bm M\bm M^+\bm b
=
P_{\Range(\bm M)}(\bm b).
\]
If \(\bm M\in\R^{p\times p}\) is symmetric positive semidefinite and \(\bm b\in\Range(\bm M)\), then
\[
\argmin_{\bm x\in\R^p}\left\{\frac12\,\bm x'\bm M\bm x-\bm b'\bm x\right\}
=
\{\bm x\in\R^p:\bm M\bm x=\bm b\}.
\]

\section{Supporting results}\label{sec:appendix_theory}

\subsection{Benchmark penalties}\label{app:sec:benchmark_penalties}

For the Adaptive Lasso, let \(\mathcal A_n:=\{j:\tilde\beta_{nj}\neq 0\}\), \(\mathcal A_0:=\{j:\beta_{0j}\neq 0\}\), and
\[
\mathcal S_0:=\{\bm\beta\in\R^p:\ \beta_j=0 \text{ for all } j\notin\mathcal A_0\}.
\]
For the Group Lasso, write \(\bm\beta=(\bm\beta_1',\ldots,\bm\beta_K')'\) for the given group partition.
For the Elastic Net, fix \(w\in(0,1)\).

Tables~\ref{tab:penalties}--\ref{tab:subgradient} collect, for several commonly used penalties \(f_n\in\Gamma(\R^p)\), the formulas used repeatedly in the paper.
Tables~\ref{tab:conjugates_euclidean} and \ref{tab:conjugates} report convex conjugates under, respectively, the Euclidean inner product and the \(\langle\cdot,\cdot\rangle_{\bm W_n}\)-inner product.
Directional derivatives in Table~\ref{tab:directional} are intrinsic.
In Tables~\ref{tab:directional} and \ref{tab:subgradient}, the objects \(\sigma_C\), \(N_C\), and \(\partial f_0\) are understood under the limit inner product \(\langle\cdot,\cdot\rangle_{\bm W_0}\), in line with the asymptotic theory.
Whenever Euclidean subdifferentials are needed, we write \(\partial^e\).

\begin{table}[ht]
\centering
\begin{tabular}{c l p{0.68\textwidth}}
\hline
ID & Penalty & \(f_n(\bm\beta)\) \\
\hline
(r)  & Ridge &
\(\dfrac12\|\bm\beta\|_2^2\) \\[0.4em]

(l)  & Lasso &
\(\|\bm\beta\|_1\) \\[0.4em]

(al) & Adaptive Lasso &
\(
\displaystyle
\sum_{j:\tilde\beta_{nj}\neq 0}\frac{|\beta_j|}{|\tilde\beta_{nj}|}
+
\sum_{j:\tilde\beta_{nj}=0}\iota_{\{0\}}(\beta_j)
\) \\[0.7em]

(gl) & Group Lasso &
\(\displaystyle \sum_{k=1}^K \|\bm\beta_k\|_2\) \\[0.4em]

(en) & Elastic Net &
\(\displaystyle w\|\bm\beta\|_1+\frac{1-w}{2}\|\bm\beta\|_2^2\) \\[0.4em]

(cs) & Constraint set &
\(\iota_C(\bm\beta)\) \\
\hline
\end{tabular}
\caption{\bf Benchmark penalties.}
\label{tab:penalties}
\end{table}

\begin{table}[ht]
\centering
\begin{tabular}{c l p{0.45\textwidth}}
\hline
ID & \(\dom(f_0)\) & Limit penalty \(f_0(\bm\beta)\) \\
\hline
(r) &
\(\R^p\) &
\(\dfrac12\|\bm\beta\|_2^2\) \\[0.8em]

(l) &
\(\R^p\) &
\(\|\bm\beta\|_1\) \\[0.8em]

(al) &
\(\mathcal S_0\) &
\(
\displaystyle
\sum_{j\in\mathcal A_0}\frac{|\beta_j|}{|\beta_{0j}|}
+
\iota_{\mathcal S_0}(\bm\beta)
\) \\[1.0em]

(gl) &
\(\R^p\) &
\(\displaystyle \sum_{k=1}^K \|\bm\beta_k\|_2\) \\[0.8em]

(en) &
\(\R^p\) &
\(\displaystyle w\|\bm\beta\|_1+\frac{1-w}{2}\|\bm\beta\|_2^2\) \\[0.8em]

(cs) &
\(C\) &
\(\iota_C(\bm\beta)\) \\
\hline
\end{tabular}
\caption{\bf Epigraphical limit penalties \(f_0\) and their domains for the benchmark penalties.}
\label{tab:limit_penalties}
\end{table}

\begin{table}[ht]
\centering
\begin{tabular}{c p{0.9\textwidth}}
\hline
ID & \(\lambda_n f_n^*(\bm\theta/\lambda_n)\), conjugation under the Euclidean inner product \\
\hline
(r) &
\(\dfrac{1}{2\lambda_n}\|\bm\theta\|_2^2\) \\[1.0em]

(l) &
\(
\iota_{B_n}(\bm\theta),
\qquad
B_n:=\bigl\{\bm\theta\in\R^p:\ |\theta_j|\le \lambda_n,\ j=1,\dots,p\bigr\}
\) \\[1.0em]

(al) &
\(
\iota_{\widetilde B_n}(\bm\theta),
\qquad
\widetilde B_n
:=
\Bigl\{
\bm\theta\in\R^p:
|\theta_j|\le \lambda_n/|\tilde\beta_{nj}|
\ \forall j\in\mathcal A_n
\Bigr\}
\) \\[1.0em]

(gl) &
\(
\iota_{G_n}(\bm\theta),
\qquad
G_n:=\bigl\{\bm\theta\in\R^p:\ \|\bm\theta_k\|_2\le \lambda_n,\ k=1,\dots,K\bigr\}
\) \\[0.8em]

(en) &
\(
\displaystyle
\frac{1}{2\lambda_n(1-w)}
\sum_{j=1}^p\bigl(|\theta_j|-\lambda_n w\bigr)_+^2
\) \\[1.0em]

(cs) &
\(\sigma_C(\bm\theta)\) \\
\hline
\end{tabular}
\caption{\bf Convex conjugates under the Euclidean inner product.}
\label{tab:conjugates_euclidean}
\end{table}

\begin{table}[ht]
\centering
\begin{tabular}{c p{0.9\textwidth}}
\hline
ID & \(\lambda_n f_n^*(\bm\theta/\lambda_n)\), conjugation under the \(\langle\cdot,\cdot\rangle_{\bm W_n}\) inner product \\
\hline
(r) &
\(\dfrac{1}{2\lambda_n}\|\bm W_n\bm\theta\|_2^2\)
\\[1.0em]

(l) &
\(
\iota_{B_n}(\bm\theta),
\qquad
B_n
:=
\Bigl\{
\bm\theta\in\R^p:\ \|\bm W_n\bm\theta\|_\infty\le \lambda_n
\Bigr\}
\)
\\[1.0em]

(al) &
\(
\iota_{\widetilde B_n}(\bm\theta),
\qquad
\widetilde B_n
:=
\Bigl\{
\bm\theta\in\R^p:\ 
\|\tilde{\bm\beta}_n\circ(\bm W_n\bm\theta)\|_\infty\le \lambda_n
\Bigr\}
\)
\\[1.0em]

(gl) &
\(
\iota_{G_n}(\bm\theta),
\qquad
G_n
:=
\Bigl\{
\bm\theta\in\R^p:\ 
\|(\bm W_n\bm\theta)_k\|_2\le \lambda_n,\ k=1,\dots,K
\Bigr\}
\)
\\[1.0em]

(en) &
\(
\displaystyle
\frac{1}{2\lambda_n(1-w)}
\sum_{j=1}^p\bigl(|(\bm W_n\bm\theta)_j|-\lambda_n w\bigr)_+^2
\)
\\[1.0em]

(cs) &
\(\sigma_C(\bm\theta)\)
\\
\hline
\end{tabular}
\caption{\bf Convex conjugates under the \(\bm W_n\)-inner product.}
\label{tab:conjugates}
\end{table}

\begin{table}[ht]
\centering
\begin{tabular}{c p{0.9\textwidth}}
\hline
ID & {\bf Directional derivative} \(\rho_{\bm\beta_0}(\bm b)\) \\
\hline
(r) &
\(\bm\beta_0'\bm b\) \\[0.6em]

(l) &
\(
\displaystyle
\sum_{j\in\mathcal A_0}\sign(\beta_{0j})b_j
+
\sum_{j\notin\mathcal A_0}|b_j|
\) \\[0.8em]

(al) &
\(
\displaystyle
\sum_{j\in\mathcal A_0}\frac{\sign(\beta_{0j})}{|\beta_{0j}|}b_j
+
\iota_{\mathcal S_0}(\bm b)
\) \\[0.8em]

(gl) &
\(
\displaystyle
\sum_{k:\bm\beta_{0k}\neq \bm 0}
\left(\frac{\bm\beta_{0k}}{\|\bm\beta_{0k}\|_2}\right)'\bm b_k
+
\sum_{k:\bm\beta_{0k}= \bm 0}\|\bm b_k\|_2
\) \\[0.9em]

(en) &
\(
\displaystyle
w\left[
\sum_{j\in\mathcal A_0}\sign(\beta_{0j})b_j
+
\sum_{j\notin\mathcal A_0}|b_j|
\right]
+(1-w)\bm\beta_0'\bm b
\) \\[0.9em]

(cs) &
\(
\iota_{T_C(\bm\beta_0)}(\bm b)
=
\sigma_{N_C(\bm\beta_0)}(\bm b)
\) \\
\hline
\end{tabular}
\caption{\bf Directional derivatives \(\rho_{\bm\beta_0}(\bm b)\).}
\label{tab:directional}
\end{table}

\begin{table}[ht]
\centering
\begin{tabular}{c p{0.9\textwidth}}
\hline
ID & {\bf Subgradient} \(\partial f_0(\bm\beta_0)\) under \(\langle\cdot,\cdot\rangle_{\bm W_0}\) \\
\hline
(r) &
\(\{\bm W_0^{-1}\bm\beta_0\}\) \\[1.0em]

(l) &
\(
\Bigl\{
\bm t\in\R^p:\ 
(\bm W_0\bm t)_j=\sign(\beta_{0j})\ \forall j\in\mathcal A_0,
\quad
(\bm W_0\bm t)_j\in[-1,1]\ \forall j\notin\mathcal A_0
\Bigr\}
\) \\[1.0em]

(al) &
\(
\Bigl\{
\bm t\in\R^p:\ 
(\bm W_0\bm t)_j=\beta_{0j}^{-1}\ \forall j\in\mathcal A_0
\Bigr\}
\) \\[1.0em]

(gl) &
\(
\Bigl\{
\bm t\in\R^p:\ 
(\bm W_0\bm t)_k=\dfrac{\bm\beta_{0k}}{\|\bm\beta_{0k}\|_2}
\ \forall k \text{ s.t. } \bm\beta_{0k}\neq \bm 0,
\quad
\|(\bm W_0\bm t)_k\|_2\le 1
\ \forall k \text{ s.t. } \bm\beta_{0k}= \bm 0
\Bigr\}
\) \\[1.0em]

(en) &
\(
\Bigl\{
\bm t\in\R^p:\ 
(\bm W_0\bm t)_j=w\sign(\beta_{0j})+(1-w)\beta_{0j}
\ \forall j\in\mathcal A_0,
\quad
(\bm W_0\bm t)_j\in[-w,w]
\ \forall j\notin\mathcal A_0
\Bigr\}
\) \\[1.0em]

(cs) &
\(N_C(\bm\beta_0)\) \\
\hline
\end{tabular}
\caption{\bf Subgradients \(\partial f_0(\bm\beta_0)\) under the \(\bm W_0\)-inner product.}
\label{tab:subgradient}
\end{table}

\begin{figure}[htbp]
\centering
\includegraphics[width=0.35\linewidth]{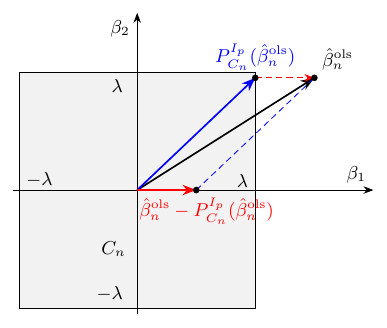}\includegraphics[width=0.35\linewidth]{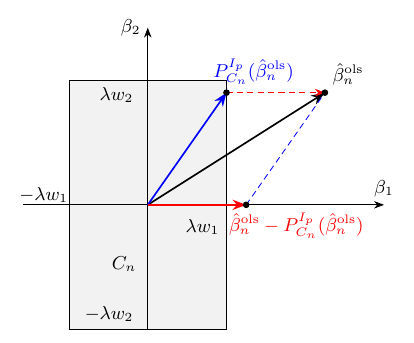}\\
\includegraphics[width=0.35\linewidth]{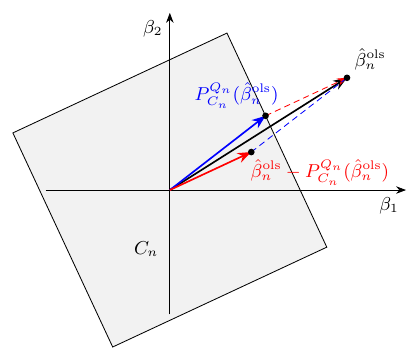}\includegraphics[width=0.35\linewidth]{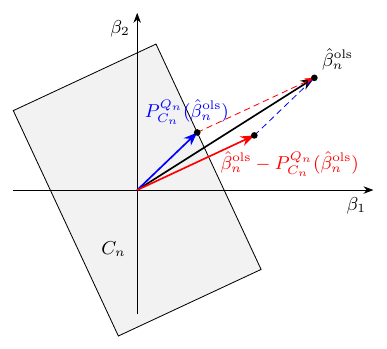}
\caption{Projection--residual representation of proximal regularization for polyhedral sets \(C_n\) (here \(p=2\)).
Columns correspond to Lasso (left) and Adaptive Lasso (right), with adaptive weights \(w_1=1/|\hat\beta_{n1}^{\mathrm{ols}}|, w_2=1/|\hat\beta_{n2}^{\mathrm{ols}}|\). Rows correspond to the Euclidean metric \(\bm W_n=\bm I_p\) (top) and the rotated geometry \(\bm W_n=\bm Q_n\neq\bm I_p\) (bottom).
Each panel shows the benchmark input \(\hat{\bm\beta}_n^{\mathrm{ols}}\) (black), the projection \(P_{C_n}^{\bm W_n}(\hat{\bm\beta}_n^{\mathrm{ols}})\) (blue), and the residual \(\prox_{\lambda_n f_n}^{\bm W_n}(\hat{\bm\beta}_n^{\mathrm{ols}})=\hat{\bm\beta}_n^{\mathrm{ols}}-P_{C_n}^{\bm W_n}(\hat{\bm\beta}_n^{\mathrm{ols}})\) (red).}
\label{fig:sm_prox_dual_lasso_adal}
\end{figure}

Figure~\ref{fig:sm_prox_dual_lasso_adal} illustrates the projection--residual geometry associated with the Lasso and Adaptive Lasso conjugates in dimension \(p=2\).

\FloatBarrier
\subsection{Non-unique solution for penalized GMM}\label{app:nonunique_penGMM}

This section presents a simple GMM example where the sample GMM criterion admits a unique minimizer almost surely, but the associated \(\ell_1\)-penalized criterion has a non-singleton set of minimizers with strictly positive probability. The point is purely finite-sample and concerns well-posedness: unlike the proximal construction, direct penalization of the GMM criterion may require additional conditions to ensure uniqueness of the solution.

Let \((w_i)_{i=1}^n\) be a sequence of i.i.d. random variables such that:
\[
\Pr(w_i=1)=\Pr(w_i=-1)=1/2.
\]
For a a scalar parameter \(\beta\in\R\), consider the parameterization
\begin{equation}\label{eq:h_def_app}
h(\beta):=
\begin{cases}
\sqrt{1-\beta}, & \beta\le 1,\\[0.3em]
-(\beta-1), & \beta\ge 1.
\end{cases}
\end{equation}
and
\begin{equation}\label{eq:g_def_app}
g(w_i,\beta):=h(\beta)+w_i.
\end{equation}
This induces the moment condition
\[
\E[g(w_i,\beta_0)]=0\ ,
\]
which admits the unique solution \(\beta_0=1\).
Let
\[
\bar g_n(\beta):=\frac1n\sum_{i=1}^n g(w_i,\beta)
= h (\beta) +\bar w_n
,
\]
where $
\bar w_n:=\frac1n\sum_{i=1}^n w_i$, and
consider the one-step GMM objective function:
\begin{equation}\label{eq:Qn_def_app}
Q_n(\beta):=\big(\bar g_n(\beta)\big)^2
=
\big(h(\beta)+\bar w_n\big)^2.
\end{equation}
Finally, for a fixed tuning parameter \(\lambda>0\), define the \(\ell_1\)-penalized GMM criterion
\begin{equation}\label{eq:pen_obj_app}
F_{n,\lambda}(\beta):=Q_n(\beta)+\lambda|\beta|.
\end{equation}

\begin{proposition}[Unique unpenalized benchmark and non-unique penalized criterion]\label{prop:nonunique_penGMM_app}
The following statements hold.

\begin{enumerate}[label=(\roman*)]
\item
For every realization of \((w_i)_{i=1}^n\), the unpenalized sample GMM criterion \(Q_n\) 
in equation (\ref{eq:Qn_def_app}) has a unique global minimizer \(\hat\beta_n\), because \(h\) is a continuous strictly decreasing bijection from \(\R\) onto \(\R\). Equivalently,
\[
h(\hat\beta_n)=-\bar w_n.
\]

\item
Suppose \(n\) is even and \(\lambda=1\) in equation \eqref{eq:pen_obj_app}. Then, on the event \(\{\bar u_n=0\}\),
\[
\argmin_{\beta\in\R} F_{n,1}(\beta)=[0,1].
\]
Consequently, for event \(E:=\big\{\argmin_{\beta\in\R} F_{n,1}(\beta)
\text{ is not a singleton}\big\}\),
we have
\[
\Pr(E)\ge \Pr(\bar u_n=0)=\binom{n}{n/2}2^{-n}>0.
\]
\end{enumerate}
\end{proposition}

Part~(i) of Proposition \ref{prop:nonunique_penGMM_app} shows that the benchmark GMM problem is well posed in every sample, because the unpenalized criterion always has a unique minimizer.
Part~(ii) shows that this property need not survive direct penalization of the criterion.
On the event \(\{\bar u_n=0\}\), the \(\ell_1\) penalty exactly offsets the slope of the unpenalized criterion on the whole interval \([0,1]\), thereby creating a flat region of global minimizers.
Hence, for penalize-the-criterion GMM estimators, uniqueness is a separate property that must be verified for the penalized objective itself and does not follow from uniqueness of the unpenalized benchmark.
By contrast, in the proximal framework the regularization step is strongly convex whenever the metric is positive definite, so existence and uniqueness are ensured by construction.

\subsection{Benchmark epigraph convergence results}\label{app:sec:epigraph_benchmark}

For the benchmark penalties reported in Table~\ref{tab:penalties},
we verify in this section the epigraph-convergence conditions invoked in
Proposition~\ref{prop:consistency} and
Theorem~\ref{prop:asymptotics:alt}.
Specifically, we verify the following four limit conditions:

\begin{enumerate}[label=(\roman*)]
    \item Assumption~\ref{ass:f}\ref{epif ass}:
    \(f_n\to_{\Pr} f_0\) in epigraph,
    for some proper function \(f_0:\R^p\to(-\infty,+\infty]\).

    \item Scaled epi-limit condition in
    Proposition~\ref{prop:consistency}\ref{prop_cons_ii}:
    \(\lambda_n f_n \to_{\Pr} \iota_{\dom(f_0)}\)
    in epigraph.

    \item Local epi-limit condition in
    Theorem~\ref{prop:asymptotics:alt}\ref{prop_distr_i}:
    \(q_n \to_{\Pr} \rho_{\bm\beta_0}\)
    in epigraph, where
    \[
    q_n(\bm b)
    :=
    r_n\left[
    f_n\left(\bm\beta_0+\bm b/r_n\right)-f_n(\bm\beta_0)
    \right],
    \qquad
    \bm b\in\R^p,
    \]
    and \(\rho_{\bm\beta_0}\) is the directional derivative of \(f_0\) at
    \(\bm\beta_0\), defined in~\eqref{eq:dirder_rho_intro}.

    \item Scaled local epi-limit condition in
    Theorem~\ref{prop:asymptotics:alt}\ref{prop_distr_ii}:
    \[
    (\lambda_n r_n)q_n
    \to_{\Pr}
    \sigma_{N_{\dom(f_0)}(\bm\beta_0)}
    \qquad\text{in epigraph.}
    \]
\end{enumerate}
Note that under Assumption~\ref{ass:f}\ref{domf ass}, the local difference quotient \(q_n\) in items (iii)--(iv) is well-defined for all sufficiently large \(n\).
Accordingly, we impose Assumption~\ref{ass:f}\ref{domf ass} throughout the sequel, although in some cases it is trivially satisfied.

Only the Adaptive Lasso is genuinely sample-dependent. All other benchmark penalties are fixed convex functions, so their epi-limits and local epi-limits follow from generic convex-analytic facts.
For the Adaptive Lasso penalty, we make use of the notation
\[
\mathcal A_0:=\{j:\beta_{0j}\neq 0\},
\qquad
\mathcal S_0:=\{\bm \beta\in\R^p:\ \beta_j=0 \text{ for all } j\notin\mathcal A_0\},
\]
and
\[
f_0^{\mathrm{al}}(\bm\beta)
:=
\sum_{j\in\mathcal A_0}\frac{|\beta_j|}{|\beta_{0j}|}
+\iota_{\mathcal S_0}(\bm\beta).
\]
Therefore, since \(\bm\beta_0\in\mathcal S_0\),
\[
\dom(f_0^{\mathrm{al}})=\mathcal S_0,
\qquad
\sigma_{N_{\dom(f_0^{\mathrm{al}})}(\bm\beta_0)}
=
\sigma_{N_{\mathcal S_0}(\bm\beta_0)}
=
\iota_{\mathcal S_0}.
\]
We first verify Assumption~\ref{ass:f}\ref{epif ass} for the benchmark penalties in Table~\ref{tab:penalties}, by characterizing explicitly the corresponding epi-limits \(f_0\).

\begin{proposition}[Epi-limits for benchmark penalties]
\label{prop:epi_table1}
Consider the penalties \(f_n\) listed in Table~\ref{tab:penalties}. Then
\(f_n\to_{\Pr} f_0\) in epigraph, where
\(f_0:\R^p\to(-\infty,+\infty]\)
is a proper lower semicontinuous convex function given in explicit form as follows:
\begin{enumerate}[label=(\roman*)]
    \item\label{prop:epi_table1:ridge} Ridge: \(f_0(\bm\beta)=\frac12\|\bm\beta\|_2^2\),
    \item\label{prop:epi_table1:lasso} Lasso: \(f_0(\bm\beta)=\|\bm\beta\|_1\),
    \item\label{prop:epi_table1:glasso} Group Lasso: \(f_0(\bm\beta)=\sum_{k=1}^K\|\bm\beta_k\|_2\),
    \item\label{prop:epi_table1:elnet} Elastic Net: \(f_0(\bm\beta)=w\|\bm\beta\|_1+\frac{1-w}{2}\|\bm\beta\|_2^2\), \quad \(w\in(0,1)\),
    \item\label{prop:epi_table1:constr} Constrained LS: \(f_0(\bm\beta)=\iota_C(\bm\beta)\),
    \item\label{prop:epi_table1:alasso} Adaptive Lasso:
    \(f_0 (\bm \beta)=f_0^{\mathrm{al}}(\bm\beta)\), whenever \(\tilde{\bm\beta}_n\to_{\Pr}\bm\beta_0\).
\end{enumerate}
\end{proposition}
Second, we verify the condition \(\lambda_n f_n\to_{\Pr}\iota_{\dom(f_0)}\) in epigraph 
for
\(\lambda_n\to 0\) (see Proposition~\ref{prop:consistency}\ref{prop_cons_ii}).

\begin{proposition}[Scaled epi-limits for vanishing tuning]
\label{prop:epi_scaled_table1}
Let \(\lambda_n>0\), \(\lambda_n\to 0\)
and \(f_0\) denote the epi-limit in Proposition~\ref{prop:epi_table1}.
Then \(\lambda_n f_n\to_{\Pr}\iota_{\dom(f_0)}\) in epigraph in the following cases:
\begin{enumerate}[label=(\roman*)]
    \item \label{prop:epi_scaled_table1:i} Ridge, Lasso, Group Lasso, Elastic Net.

    \item \label{prop:epi_scaled_table1:ii} Constrained least squares.

    \item \label{prop:epi_scaled_table1:iii} Adaptive Lasso: If \(\tilde{\bm\beta}_n\to_{\Pr}\bm\beta_0\) and
    $
    |\tilde\beta_{nj}|/\lambda_n\to_{\Pr}0$ for every
    $j\notin\mathcal A_0$.
\end{enumerate}
\end{proposition}

Third, we verify the local epi-limit condition in
Theorem~\ref{prop:asymptotics:alt}\ref{prop_distr_i},
\(q_n \to_{\Pr} \rho_{\bm\beta_0}\) in epigraph.

\begin{proposition}[Local epi-limits for benchmark penalties]
\label{prop:local_epi_table1}
Let \(r_n\to\infty\)
and
\(f_0\) denote the epi-limit
in Proposition~\ref{prop:epi_table1}.
Then, 
\(q_n \to_{\Pr} \rho_{\bm\beta_0}\) in epigraph
in the following benchmark cases:
\begin{enumerate}[label=(\roman*)]
    \item\label{prop:local_epi_table1:ridge}
    Ridge:
    \[
    \rho_{\bm\beta_0}(\bm b)=\bm\beta_0'\bm b.
    \]

    \item\label{prop:local_epi_table1:lasso}
    Lasso:
    \[
    \rho_{\bm\beta_0}(\bm b)
    =
    \sum_{j:\beta_{0j}\neq 0}\sign(\beta_{0j})b_j
    +
    \sum_{j:\beta_{0j}=0}|b_j|.
    \]

    \item\label{prop:local_epi_table1:glasso}
    Group Lasso: if \(\bm\beta=(\bm\beta_1',\ldots,\bm\beta_K')'\) is partitioned according to the groups in Table~\ref{tab:penalties}, then
    \[
    \rho_{\bm\beta_0}(\bm b)
    =
    \sum_{k:\bm\beta_{0k}\neq \bm 0}
    \left(\frac{\bm\beta_{0k}}{\|\bm\beta_{0k}\|_2}\right)'\bm b_k
    +
    \sum_{k:\bm\beta_{0k}=\bm 0}\|\bm b_k\|_2.
    \]

    \item\label{prop:local_epi_table1:elnet}
    Elastic Net:
    \[
    \rho_{\bm\beta_0}(\bm b)
    =
    w\left(
    \sum_{j:\beta_{0j}\neq 0}\sign(\beta_{0j})b_j
    +
    \sum_{j:\beta_{0j}=0}|b_j|
    \right)
    +(1-w)\bm\beta_0'\bm b,
    \qquad w\in(0,1).
    \]

    \item\label{prop:local_epi_table1:constr}
    Constrained least squares: if \(\bm\beta_0\in C\), then
    \[
    \rho_{\bm\beta_0}(\bm b)
    =
    \iota_{T_C(\bm\beta_0)}(\bm b)
    =
    \sigma_{N_C(\bm\beta_0)}(\bm b),
    \]
    where \(T_C(\bm\beta_0)\) and \(N_C(\bm\beta_0)\) denote, respectively,
    the tangent cone and the normal cone of \(C\) at \(\bm\beta_0\).

    \item\label{prop:local_epi_table1:alasso}
    Adaptive Lasso: if \(\tilde{\bm\beta}_n\to_{\Pr}\bm\beta_0\), then
    \[
    \rho_{\bm\beta_0}(\bm b)
    =
    \sum_{j\in\mathcal A_0}
    \frac{\sign(\beta_{0j})}{|\beta_{0j}|}b_j
    +
    \iota_{\mathcal S_0}(\bm b).
    \]
\end{enumerate}
\end{proposition}

Finally, we address the scaled local epi-limit condition
\((\lambda_n r_n) q_n \to_{\Pr} \sigma_{N_{\dom(f_0)}(\bm\beta_0)}\)
in epigraph
(see Theorem~\ref{prop:asymptotics:alt}\ref{prop_distr_ii}).

\begin{proposition}[Scaled epi-limits for benchmark penalties]
\label{prop:local_scaled_epi_table1}
Let \(\lambda_n>0\) and \(r_n\to\infty\) be such that
$\lambda_n r_n\to 0$.
Let further \(f_0\) denote the epi-limit 
in Proposition~\ref{prop:epi_table1}.
Then, 
\[
\lambda_n r_n
q_n
\to_{\Pr}
\sigma_{N_{\dom(f_0)}(\bm\beta_0)}
\qquad\text{in epigraph,}
\]
in the following benchmark cases:
\begin{enumerate}[label=(\roman*)]
    \item\label{prop:local_scaled_epi_table1:finite}
    Ridge, Lasso, Group Lasso, and Elastic Net:
\[
    \dom(f_0)=\R^p,
    \qquad
    \sigma_{N_{\dom(f_0)}(\bm\beta_0)}
    =
    \sigma_{\{\bm 0\}}
    =0.
    \]

    \item\label{prop:local_scaled_epi_table1:constr}
    Constrained least squares: if \(\bm\beta_0\in C\), then
    \[
    \dom(f_0)=C,
    \qquad
    \sigma_{N_{\dom(f_0)}(\bm\beta_0)}
    =
    \sigma_{N_C(\bm\beta_0)}
    =
    \iota_{T_C(\bm\beta_0)}.
    \]

    \item\label{prop:local_scaled_epi_table1:alasso}
    Adaptive Lasso: if \(\tilde{\bm\beta}_n\to_{\Pr}\bm\beta_0\) and
    \[
    \frac{|\tilde\beta_{nj}|}{\lambda_n r_n}\to_{\Pr}0
    \qquad\text{for every }j\notin\mathcal A_0,
    \]
    where
    \[
    \dom(f_0^{\mathrm{al}})=\mathcal S_0,
\qquad
    \sigma_{N_{\dom(f_0^{\mathrm{al}})}(\bm\beta_0)}
=
    \iota_{\mathcal S_0}.
    \]
\end{enumerate}
\end{proposition}

\subsection{Standard Lasso weak limit results}
\label{app:KF_lasso_dual}

Theorem~\ref{prop:asymptotics:alt} is formulated for proximal estimators with penalties in
\(\Gamma(\R^p)\). Through the PLSE embedding of Section~\ref{sec:link_plse}, it covers convex penalized least-squares estimators as a special case.
This section shows how the general weak-limit theory in Theorem~\ref{prop:asymptotics:alt} produces the standard limit results of \citet{SM:knight2000asymptotics}
for this special case.
Their analysis is focused on least-squares estimation with Bridge penalties
\[
f_n(\bm\beta)=f(\bm\beta):=\sum_{j=1}^p |\beta_j|^\gamma,
\qquad \gamma>0,
\]
covering both convex \((\gamma\ge 1)\) and nonconvex \((0<\gamma<1)\) specifications. Our proximal framework is instead based on generic convex (possibly sample-dependent and extended-real-valued)
penalties for estimation problems beyond least squares.

For brevity, we focus in the sequel on the regular-design Lasso setting
of \citet{SM:knight2000asymptotics}. To this end,
consider the linear model~\eqref{linear_model_setting} under a regular design, i.e., using
assumptions~(3)--(4) in \citet{SM:knight2000asymptotics}:
\[
\bm Q_n:=\bm X'\bm X/n \to_{\Pr} \bm Q_0,
\qquad
\bm Q_0\succ \bm 0,
\]
and
\[
\sqrt n\big(\hat{\bm\beta}_n^{\mathrm{ls}}-\bm\beta_0\big)
\to_d
\bm\eta,
\qquad
\bm\eta\sim N(\bm 0,\sigma^2\bm Q_0^{-1}),
\]
where
\[
\hat{\bm\beta}_n^{\mathrm{ls}}
=
\bm Q_n^{-1}\bm X'\bm Y/n.
\]
For the Lasso penalty \(f_n(\bm\beta)=\|\bm\beta\|_1\), they further assume
\[
\lambda_n\sqrt n\to \lambda_0 > 0,
\]
and show that
\[
\sqrt n\big(\hat{\bm\beta}_n^{\mathrm{lasso}}-\bm\beta_0\big)
\to_d
\argmin_{\bm u\in\R^p}
\left\{
-2\bm u'\bm W+\bm u'\bm Q_0\bm u
+\lambda_0\rho_{\bm\beta_0}(\bm u)
\right\},
\]
where
\[
\bm W\sim N(\bm 0,\sigma^2\bm Q_0)
\]
and
\[
\rho_{\bm\beta_0}(\bm u)
=
\sum_{j:\beta_{0j}\neq 0}\sign(\beta_{0j})u_j
+
\sum_{j:\beta_{0j}=0}|u_j|.
\]
This result is an immediate special case of our more general proximal limit
theorem. Indeed, the regular-design PLSE embedding of
Proposition~\ref{prop:prox_plse_regular}
directly yields:
\[
\prox_{\lambda_n\|\cdot\|_1}^{\bm Q_n}
\big(\hat{\bm\beta}_n^{\mathrm{ls}}\big)
=
\argmin_{\bm\beta\in\R^p}
\left\{
\frac{1}{2n}\|\bm Y-\bm X\bm\beta\|_2^2+\lambda_n\|\bm\beta\|_1
\right\}=\hat{\bm\beta}_n^{\mathrm{lasso}}.
\]
Furthermore, given \(f_0(\bm\beta)=\|\bm\beta\|_1\) the local penalty in
Theorem~\ref{prop:asymptotics:alt}\ref{prop_distr_i} is precisely the
directional derivative
\[
\rho_{\bm\beta_0}(\bm u)
=
\sum_{j:\beta_{0j}\neq 0}\sign(\beta_{0j})u_j
+
\sum_{j:\beta_{0j}=0}|u_j|.
\]
Therefore, Theorem~\ref{prop:asymptotics:alt}\ref{prop_distr_i} gives
\[
\sqrt n\Big(
\hat{\bm\beta}_n^{\mathrm{lasso}}
-\bm\beta_0
\Big)
\to_d
\prox_{\lambda_0\rho_{\bm\beta_0}}^{\bm Q_0}(\bm\eta).
\]
With the reparameterization
\(\bm W:=\bm Q_0\bm\eta\) we further have
\(\bm W\sim N(\bm 0,\sigma^2\bm Q_0)\)
and
\begin{eqnarray*}
\prox_{\lambda_0\rho_{\bm\beta_0}}^{\bm Q_0}(\bm\eta)
&=&
\argmin_{\bm u\in\R^p}
\left\{
\frac12\bm u'\bm Q_0\bm u-\bm u'\bm W
+\lambda_0\rho_{\bm\beta_0}(\bm u)
\right\}\\
&=&
 \argmin_{\bm u\in\R^p}
\left\{
-2\bm u'\bm W+\bm u'\bm Q_0\bm u
+2\lambda_0\rho_{\bm\beta_0}(\bm u)
\right\}.
\end{eqnarray*}
This is precisely the Lasso weak limit in
\citet[Thm.~2, case \(\gamma=1\)]{SM:knight2000asymptotics}.
Our framework further yields an equivalent dual representation of this limit.
Indeed, since
\[
\rho_{\bm\beta_0}
=
\sigma_{\partial\|\cdot\|_1(\bm\beta_0)},
\]
Moreau decomposition implies:
\[
\prox_{\lambda_0\rho_{\bm\beta_0}}^{\bm Q_0}(\bm\eta)
=
\big(
\Id-
P_{\lambda_0\partial\|\cdot\|_1(\bm\beta_0)}^{\bm Q_0}
\big)(\bm\eta).
\]
Accordingly, the Lasso weak limit is equivalently the residual 
of a projection of the limit OLS fluctuation $\bm\eta$ onto the scaled Lasso subgradient set under the
\(\bm Q_0\)-metric.
Whenever \(\bm Q_0=\bm I_p\), this limit reduces to the familiar coordinatewise
soft-thresholding representation of the Lasso weak limit.

Finally, note that when \(\lambda_n\sqrt n\to 0\) the Lasso penalty is asymptotically negligible at the \(\sqrt n\)-scale. Accordingly,
Theorem~\ref{prop:asymptotics:alt}\ref{prop_distr_ii} yields the OLS limit law for the Lasso estimator in this case (recall that \(\dom(\|\cdot\|_1)=\R^p\) and \(N_{\R^p}(\bm\beta_0)=\{\bm 0\}\)):
\[
\prox_{\sigma_{N_{\dom(\|\cdot\|_1)}(\bm\beta_0)}}^{\bm W_0}(\bm\eta)
=
\bigl(\Id-P_{N_{\dom(\|\cdot\|_1)}(\bm\beta_0)}^{\bm W_0}\bigr)(\bm\eta)
=
\bigl(\Id-P_{N_{\R^p}(\bm\beta_0)}^{\bm W_0}\bigr)(\bm\eta)
=
\bm\eta.
\]

\subsection{Characterizing conditions for  oracle selection}\label{app:oracle_selection}

This section provides an exact characterization of Oracle property~\ref{oracle property}\ref{op1}. This builds on the fact that support recovery holds if and only if, with probability tending to one, the proximal problem admits no locally optimal solution within an \(O(1/r_n)\)-neighborhood of \(\bm\beta_0\) that retains a nonzero inactive coordinate.
Thus, once \(\hat{\bm\beta}_n\) is \(r_n\)-consistent, Oracle selection reduces to a local exclusion condition for such solutions. This characterization underlies the more primitive sufficient and necessary conditions developed below.

\begin{proposition}[Exact characterization of oracle selection]\label{prop:op1_iff}
Let
\[
\hat{\bm\beta}_n:=\prox_{\lambda_n f_n}^{\bm W_n}(\hat{\bm\beta}_n^s), .\]
Suppose further that \(\bm W_n\succ \bm 0\) for every \(n\), and that
as $r_n\to\infty$:
\[
\hat{\bm\beta}_n-\bm\beta_0 = O_{\Pr}(1/r_n),
\]
For \(M>0\), define the set:
\[
\mathcal K_n(M)
:=
\Big\{
\bm\beta\in\R^p:
\|\bm\beta-\bm\beta_0\|_\infty\le \frac{M}{r_n},\ 
\bm\beta_{\mathcal A^c}\neq \bm 0,\ 
\bm W_n(\hat{\bm\beta}_n^s-\bm\beta)\in \lambda_n\partial^e f_n(\bm\beta)
\Big\}.
\]
Then, Oracle property~\ref{oracle property}\ref{op1} holds if and only if
\[
\lim_{M\to\infty}\limsup_{n\to\infty}
\Pr\big(\mathcal K_n(M)\neq\emptyset\big)=0.
\]
\end{proposition}

\subsection{Ridgeless estimand: identification and regularity}
\label{app:Ridgeless_differentiability}

We formalize in this section the regularity properties of the Ridgeless
estimand.
We work on a fixed-law Gaussian submodel with singular second-moment
matrix \(\bm Q_0\), because any target that is \(\sqrt n\)-regular on a larger
model must also be \(\sqrt n\)-regular on each regular parametric submodel;
see \citet[Thm.~2.1]{SM:van1991differentiable}.
Within this submodel, singularity of \(\bm Q_0\) implies that coefficient
perturbations in \(\Kernel(\bm Q_0)\) do not change the law of the observables
and therefore cannot be identified from the data.
The proposition below records the resulting null-space restriction for
a regular parameter target and then shows that the Ridgeless representative is the
unique linear selection into
\(\Range(\bm Q_0)\) from the affine set \(\bm\beta+\Kernel(\bm Q_0)\).
Let
\[
\bm P_0:=\bm I_p-\bm Q_0\bm Q_0^+,
\qquad
\bm P_0^\perp:=\bm Q_0\bm Q_0^+,
\]
be the orthogonal projectors onto \(\Kernel(\bm Q_0)\) and
\(\Range(\bm Q_0)\), respectively.
Suppose that \(\bm x\) has a fixed distribution satisfying
\(\E[\bm x\bm x']=\bm Q_0\), and that, conditional on \(\bm x\),
\begin{equation}\label{eq:appendix_rls_gaussian_submodel}
y=\bm x'\bm\beta+\varepsilon,
\qquad
\varepsilon\sim N(0,\sigma_0^2),
\qquad
\varepsilon\perp \bm x,
\end{equation}
with \(\sigma_0^2>0\).
Let \(P_{\bm\beta}\) denote the induced law.
For every \(\bm v\in\Kernel(\bm Q_0)\), we obtain
\[
\E[(\bm x'\bm v)^2]=\bm v'\bm Q_0\bm v=0,
\]
i.e., \(\bm x'\bm v=0\) almost surely. Hence
\[
y=\bm x'(\bm\beta+\bm v)+\varepsilon
=
\bm x'\bm\beta+\varepsilon,
\]
almost surely,
and
$
P_{\bm\beta+\bm v}=P_{\bm\beta}
$.
That is, \(\bm\beta\) and \(\bm\beta+\bm v\) parametrize the same joint law of the
observable pair \((y,\bm x)\).
Hence, the model is indexed only up to the quotient
\(\R^p/\Kernel(\bm Q_0)\).
In particular, any coefficient target \(K(P_{\bm\beta})=\psi(\bm\beta)\) that
is well-defined on this quotient must satisfy
for any $\bm v\in\Kernel(\bm Q_0)$:
\[
\psi(\bm\beta+\bm v)=\psi(\bm\beta)
\qquad
.
\]
Whenever \(\psi\) is differentiable at \(\bm\beta_0\), then this necessarily implies
\begin{equation}\label{eq:appendix_identification_nullspace_condition}
\Kernel(\bm Q_0)\subseteq \Kernel(\dot\psi_0),
\end{equation}
where \(\dot\psi_0\) denotes the Jacobian matrix of \(\psi\) at
\(\bm\beta_0\).

\begin{proposition}[Regular estimand under singular information]
\label{prop:Ridgeless_regular_targets}
Consider the Gaussian submodel \eqref{eq:appendix_rls_gaussian_submodel},
and let \(\psi:\R^p\to\R^m\) be differentiable at \(\bm\beta_0\), with
derivative matrix \(\dot\psi_0\in\R^{m\times p}\); that is, as
\(\bm h\to \bm 0\),
\[
\psi(\bm\beta_0+\bm h)
=
\psi(\bm\beta_0)+\dot\psi_0\bm h+o(\|\bm h\|_2).
\]
Then the following statements hold.

\begin{enumerate}[label=(\roman*)]
\item
If the coefficient target \(K(P_{\bm\beta})=\psi(\bm\beta)\) is
\(\sqrt n\)-regularly estimable at \(P_{\bm\beta_0}\) in the sense of
\citet[Thm.~2.1]{SM:van1991differentiable}, then
\begin{equation}\label{eq:appendix_rls_nullspace_condition}
\Kernel(\bm Q_0)\subseteq \Kernel(\dot\psi_0).
\end{equation}
In particular, when \(\bm Q_0\) is singular, the full coefficient vector
\(\psi(\bm\beta)=\bm\beta\) cannot be a regular target.

\item
For every \(\bm\beta\in\R^p\),
\begin{equation}\label{eq:appendix_rls_singleton}
\bigl(\bm\beta+\Kernel(\bm Q_0)\bigr)\cap\Range(\bm Q_0)
=
\{\bm P_0^\perp\bm\beta\}.
\end{equation}
Accordingly, define
\begin{equation}\label{eq:appendix_rls_map}
K^{\mathrm{rls}}(P_{\bm\beta})
:=
\bm P_0^\perp\bm\beta .
\end{equation}
This functional is well-defined on the quotient model. Moreover, for every
path
\[
\bm\beta_t=\bm\beta_0+t\bm h+o(t),
\qquad
\bm h\in\R^p,
\]
we have, as \(t\to 0\),
\begin{equation}\label{eq:appendix_rls_derivative}
K^{\mathrm{rls}}(P_{\bm\beta_t})
=
K^{\mathrm{rls}}(P_{\bm\beta_0})
+
t\bm P_0^\perp\bm h
+
o(t).
\end{equation}
Equivalently,
\[
\left.\frac{d}{dt}K^{\mathrm{rls}}(P_{\bm\beta_t})\right|_{t=0}
=
\bm P_0^\perp\bm h.
\]

\item
Let \(L:\R^p\to\R^p\) be linear.
If for every \(\bm\beta\in\R^p\),
\begin{equation}\label{eq:appendix_rls_selection_conditions}
L\bm\beta\in
\bigl(\bm\beta+\Kernel(\bm Q_0)\bigr)\cap\Range(\bm Q_0),
\end{equation}
then
\begin{equation}\label{eq:appendix_rls_unique_projection}
L=\bm P_0^\perp.
\end{equation}
Consequently, among linear coefficient selections that choose one
representative from each observational-equivalence class
\(\bm\beta+\Kernel(\bm Q_0)\) and require the selected coefficient to lie in
\(\Range(\bm Q_0)\), the unique candidate is the Ridgeless rule
\(\bm\beta\mapsto \bm P_0^\perp\bm\beta\).
This rule applied to any \(\tilde{\bm\beta}_0\in\mathcal B_0\) yields
\(\bm P_0^\perp\tilde{\bm\beta}_0=\bm\beta_0^{\mathrm{rls}}\).
\end{enumerate}
\end{proposition}
Proposition~\ref{prop:Ridgeless_regular_targets} isolates the two main issues exploited in the main text.
First, regular coefficient targets must be locally insensitive to perturbations
in \(\Kernel(\bm Q_0)\).
Second, among linear rules that select one coefficient from each
observational-equivalence class, the Ridgeless rule is the unique one that
keeps the selected coefficient in \(\Range(\bm Q_0)\).
This is why \(\bm\beta_0^{\mathrm{rls}}\) is the natural least-squares target
for regular inference under singularity.

This corresponds to reparameterizing the model on the identified component of the regressors, i.e., the projection onto \(\Range(\bm Q_0)\).
Let \(\bm U_0\in\R^{p\times r_0}\) have orthonormal columns spanning
\(\Range(\bm Q_0)\)
and \(\bm z:=\bm U_0'\bm x\). Then \(\bm x=\bm U_0\bm z\) almost surely with
and
\(\E[\bm z\bm z']=\bm U_0'\bm Q_0\bm U_0\succ\bm 0\).
Using the reparameterization \(\bm\gamma:=\bm U_0'\bm\beta\), the model becomes
\[
y=\bm z'\bm\gamma+\varepsilon,
\qquad
K^{\mathrm{rls}}(P_{\bm\beta})=\bm U_0\bm\gamma.
\]
Standard least-squares theory therefore applies on the identified
span and the null-space component is not an inferential target, because it is
neither identified by the law of \((y,\bm x)\) nor compatible with regular
estimation.

\subsection{Relation to nearly-singular penalized least-squares asymptotics}\label{app:comparison_classical_nearlysingular}

In this section, we relate the limit results of Section~\ref{sec:irregular_designs} -- in particular Theorem~\ref{thm:asy_modified_Ridgeless} and Corollary~\ref{cor:oracle_prox_est} -- to the literature on nearly-singular penalized least-squares regression.
The closest benchmark is the theory of \citet[Sect.~5, Thm.~6]{SM:knight2000asymptotics} and \citet[Thm.~1]{SM:knight2008shrinkage} for penalized least-squares estimators under deterministic nearly-singular designs. While this literature is seminal for Bridge- and Lasso-type estimators in such settings, its first-order conclusions differ markedly from ours.
To facilitate comparison, we rewrite their setup in our notation. Consider the Bridge estimator
\[
\hat{\bm\beta}_n^{\mathrm{br}}
\in
\argmin_{\bm\beta\in\R^p}
\Big\{
\|\bm Y-\bm X\bm\beta\|_2^2
+
\lambda_n\sum_{j=1}^p |\beta_j|^\gamma
\Big\},
\qquad
\gamma>0,
\]
under a deterministic triangular-array design, so that
\[
\bm Q_n=\bm Q_{0n}=\bm X'\bm X/n.
\]
Their nearly-singular framework is defined by
\[
\bm Q_{0n}\succ \bm 0,\qquad
\bm Q_{0n}\to \bm Q_0,
\qquad
\tau_n(\bm Q_{0n}-\bm Q_0)\to \bm\Delta,
\]
where \(\tau_n\to\infty\) and \(\bm\Delta\) is positive definite on \(\Kernel(\bm Q_0)\). This corresponds to our Assumption~\ref{ass: near singularity} in the special case of a deterministic design without sampling fluctuations.

Given the normalization \(b_n:=(n/\tau_n)^{1/2}\), the main result of \citet{SM:knight2000asymptotics, SM:knight2008shrinkage} is that a nondegenerate limit arises at the \(b_n\)-scale, rather than at the standard \(\sqrt n\)-scale. Specifically,
\begin{equation}
b_n\big(\hat{\bm\beta}_n^{\mathrm{br}}-\bm\beta_0\big)
\to_d
\argmin_{\bm u\in\Kernel(\bm Q_0)} V_\gamma(\bm u),
\label{eq: classical PLSE limit law}
\end{equation}
under the tuning regimes
\[
\lambda_n/b_n\to\lambda_0
\quad\text{if }\gamma\ge 1,
\qquad
\lambda_n/b_n^\gamma\to\lambda_0
\quad\text{if }0<\gamma<1,
\]
where
\[
V_\gamma(\bm u)
:=
\bm u'\bm\Delta \bm u
-
2\bm u'\bm W
+
\lambda_0 \rho_0(\bm u),
\]
with \(\bm W\sim N(\bm 0,\sigma^2\bm\Delta)\), and
\[
\rho_0(\bm u)
:=
\begin{cases}
\sum_{j=1}^p u_j\sign(\beta_{0j})|\beta_{0j}|^{\gamma-1},
& \gamma>1,\\[1ex]
\sum_{j=1}^p
\Big[
u_j\sign(\beta_{0j})\mathbf 1\{\beta_{0j}\neq 0\}
+
|u_j|\mathbf 1\{\beta_{0j}=0\}
\Big],
& \gamma=1,\\[1ex]
\sum_{j=1}^p |u_j|^\gamma\mathbf 1\{\beta_{0j}=0\},
& 0<\gamma<1.
\end{cases}
\]
Three features of this theory are particularly relevant in a comparison with our proximal estimation approach:

\begin{enumerate}[label=(\roman*)]

\item \emph{Nonstandard convergence rate.}  
Since \(b_n=(n/\tau_n)^{1/2}=o(\sqrt n)\), the convergence rate is slower than \(\sqrt n\). Under the local regime \(\sqrt n/\tau_n\to c\in(0,\infty)\), one has \(b_n\asymp n^{1/4}\), so classical PLSEs achieve rate \(n^{-1/4}\). In contrast, our estimator retains the standard \(n^{-1/2}\) rate.

\item \emph{Support of the limit law.}  
The limit \eqref{eq: classical PLSE limit law} is supported on \(\Kernel(\bm Q_0)\), and thus reflects only directions that are not identified in the limit. The identified component in \(\Range(\bm Q_0)\) is negligible at this scale. Since our object of interest,
\[
\bm\beta_0^{\mathrm{rls}}=\bm Q_0^+\bm\delta_0,
\]
lies in \(\Range(\bm Q_0)\), this limit does not provide a first-order description of the coordinates relevant for inference.

\item \emph{Gaussian limits and identification.}  
Gaussian limits may arise under weak penalization, but remain supported on \(\Kernel(\bm Q_0)\). As such, they do not yield a Gaussian approximation for the identified component.

\end{enumerate}
Our proximal estimation approach relies on a different mechanism. We separate rank recovery from regularization by first constructing the modified Ridgeless estimator
\[
\check{\bm\beta}_n^{\mathrm{rls}}
=
\check{\bm Q}_n^+\bm X'\bm Y/n,
\]
which removes non-identified directions prior to penalization.

\begin{enumerate}[label=(\roman*)]

\item \emph{Standard \(\sqrt n\)-limit.}  
The modified Ridgeless estimator admits a nondegenerate \(\sqrt n\)-limit, with near singularity entering through a drift term. The Gaussian component lies in \(\Range(\bm Q_0)\), so the identified component retains standard asymptotic behavior.

\item \emph{Oracle property.}  
Proximal estimators built on this benchmark preserve the \(\sqrt n\)-limit and can enforce sparsity in a second step. Under suitable tuning, they achieve the oracle property relative to the Ridgeless target. This contrasts with classical PLSE theory, where oracle estimation is not directly achievable, as no regular Oracle estimator for the relevant target admits a standard first-order approximation.

\end{enumerate}
Summarizing, we conclude that
classical nearly-singular PLSE theory operates at a slower scale and focuses on non-identified directions, whereas our approach retains a standard \(\sqrt n\)-approximation for the identified component.
These differences in limiting behavior arise, in part, from the distinct centering targets. In the classical PLSE framework, the analysis is centered at a fixed prelimit coefficient \(\bm\beta_0\), which satisfies \(\bm Q_0 \bm\beta_0 = \bm\delta_0\),
and hence belongs to the solution set \(\mathcal B_0\). Since this set is not a singleton, the limit behavior depends on the particular representative \(\bm\beta_0\).
By contrast, our analysis is centered on the Ridgeless estimand \(\bm\beta_0^{\mathrm{rls}}=\bm Q_0^+\bm\delta_0\). This distinction in centering contributes directly to the different asymptotic conclusions.

\clearpage

\section{Monte Carlo analysis}\label{sec:simulations}

\subsection{Design}\label{subsec:mc_design}

We study by Monte Carlo simulation the finite-sample behavior of our proximal estimators under irregular designs.
Starting from the regular Gaussian design in \citet[Ex.~1]{SM:tibshirani1996regression}, we construct a data-generating scheme that jointly encompasses a regular design, a singular design, and a nearly-singular triangular-array design in the sense of Definition~\ref{def:irregular designs}.

We simulate $n$ 
i.i.d.\ realizations of observations from linear model~\eqref{linear model}:
\[
y_i=\bm X_i'\bm\beta_0+\varepsilon_i,
\]
where \(\varepsilon_i\) is independent of \(\bm X_i\), \(\bm X_i\sim\mathcal N(\bm 0,\bm Q_{0n})\),
\(\varepsilon_i\sim\mathcal N(0,\sigma_0^2)\),
$\sigma_0^2=2$, $p=8$ and
\begin{equation}\label{eq: beta_star_numbers_sim}
\bm\beta_0=(3,\ 1.5,\ 0,\ 0,\ 2,\ 0,\ 0,\ 0)'.
\end{equation}
To incorporate regular, singular and nearly-singular designs,
the population second moment matrix \(\bm Q_{0n}:=\E[\bm X_i\bm X_i']\) is specified as follows.

\begin{figure}[H]
\centering
\includegraphics[scale=0.4]{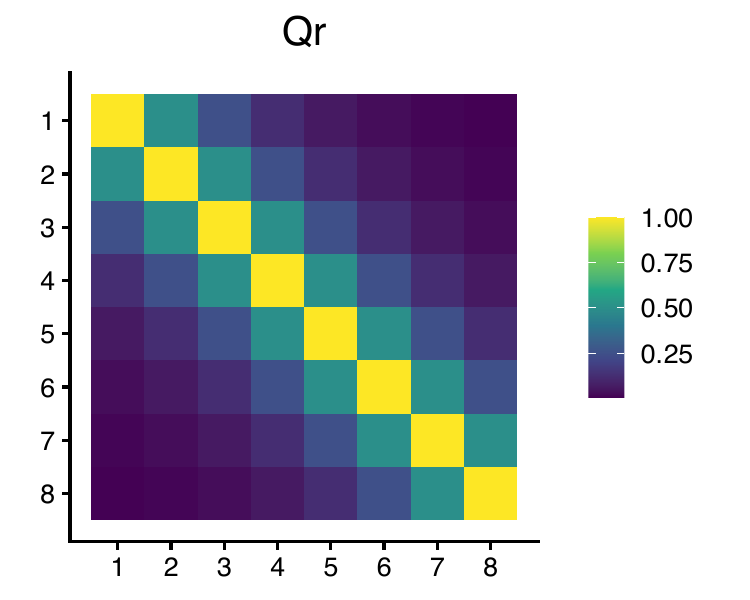}\includegraphics[scale=0.4]{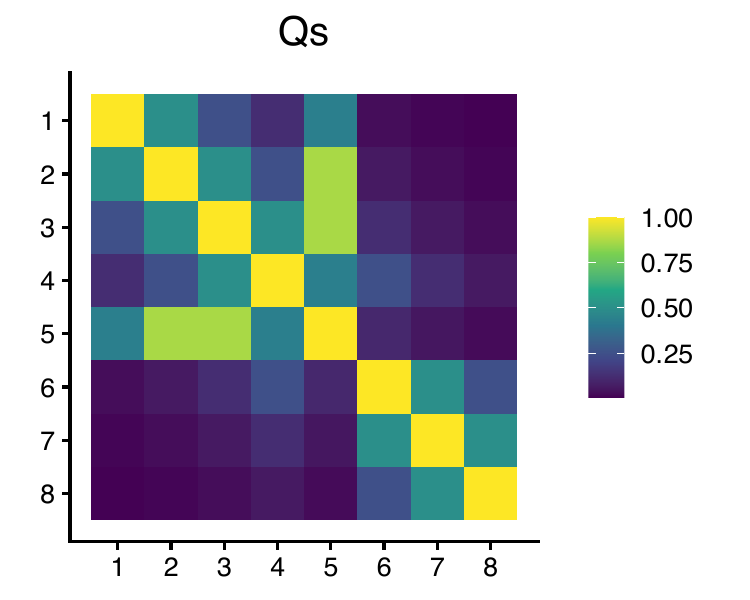}
\caption{{Heatmap of population design matrices $\bm Q_r$ and $\bm Q_0$ in our Monte Carlo simulation settings}.}\label{fig:QrQ0}
\end{figure}

For the regular design, we set \(\bm Q_{0n}=\bm Q_r\) using a positive definite matrix \(\bm Q_r\) defined as in \citet{SM:tibshirani1996regression} with the Toeplitz entries
\[
(\bm Q_r)_{jk}=0.5^{|j-k|},\qquad j,k\in\{1,\ldots,8\}.
\]
For the singular design, we set \(\bm Q_{0n}=\bm Q_0\), where
\(\bm Q_0\) is a rank-deficient matrix built from \(\bm Q_r\).
Specifically, \(\bm Q_0\) coincides with \(\bm Q_r\) except in the fifth row
and fifth column. The fifth regressor is constructed, at the population
level, as the normalized linear combination of the second and third
coordinates under the regular design. Let
\[
C^2
:=
(\bm Q_r)_{22}+(\bm Q_r)_{33}+2(\bm Q_r)_{23},
\qquad
\widetilde X_i
:=
\frac{X_{i2}+X_{i3}}{C},
\]
where the covariances are computed under the regular design
\(\bm X_i\sim N(\bm 0,\bm Q_r)\). Thus \(\Var(\widetilde X_i)=1\).
We define
\[
(\bm Q_0)_{kj}
=
\begin{cases}
(\bm Q_r)_{kj}, & k,j\neq 5,\\
1, & k=j=5,\\
\{(\bm Q_r)_{2j}+(\bm Q_r)_{3j}\}/C, & k=5,\ j\neq5,\\
\{(\bm Q_r)_{k2}+(\bm Q_r)_{k3}\}/C, & k\neq5,\ j=5.
\end{cases}
\]
With this construction, \(\bm Q_0\) is positive semidefinite with
\(\Rank(\bm Q_0)=p-1\).
For comparison, the different correlation structures induced by
\(\bm Q_r\) and \(\bm Q_0\) are illustrated in Figure~\ref{fig:QrQ0}.

For the nearly-singular design, we interpolate between the singular limit
\(\bm Q_0\) and the regular design \(\bm Q_r\). Specifically, we set
\[
\bm Q_{0n}
=
\bm Q_0+\frac{1}{\sqrt n}\bm\Delta,
\qquad
\bm\Delta:=\bm Q_r-\bm Q_0 .
\]
Equivalently,
\[
\bm Q_{0n}
=
\Big(1-\frac{1}{\sqrt n}\Big)\bm Q_0
+
\frac{1}{\sqrt n}\bm Q_r .
\]
Therefore, \(\bm Q_{0n}\) is positive definite for every \(n\), because
\(\bm Q_r\succ0\) and \(\bm Q_0\succeq0\), while
\(\bm Q_{0n}\to\bm Q_0\). The near-singularity vanishes at rate
\(\tau_n^{-1}=n^{-1/2}\).

Given these definitions,
the population cross-moment under the nearly-singular designs is given by:
\[
\bm\delta_{0n}=\bm Q_{0n}\bm\beta_0 \to 
\bm\delta_0=\bm Q_0\bm\beta_0
,
\]
where $\bm\delta_0$ is the population cross-moment under the singular design.
Hence, the limit population least-squares solutions form the affine class
\[
\mathcal B_0=\{\bm\beta\in\R^p:\bm Q_0\bm\beta=\bm\delta_0\}
=\bm\beta_0+\Kernel(\bm Q_0).
\]
Accordingly,
the identifiable regular parameter of interest in this class is the Ridgeless estimand
\[
\bm\beta_0^{\mathrm{rls}}:=\bm Q_0^+\bm\delta_0
=\bm Q_0\bm Q_0^+\bm\beta_0
=\bm P_0^\perp\bm\beta_0,
\qquad
\bm P_0^\perp:=\bm Q_0\bm Q_0^+ .
\]
Under the regular design, 
\(\bm\beta_0^{\mathrm{rls}}=\bm\beta_0\) by construction,
while under an irregular design:
\begin{equation}\label{eq: beta0_rls_numbers_sim}
\bm\beta_0^{\mathrm{rls}}
=\bm P_0^\perp\bm\beta_0
\approx (3,\ 1.893,\ 0.393,\ 0,\ 1.32,\ 0,\ 0,\ 0)'.
\end{equation}
In each design, the oracle active set used to compute selection probabilities is
\[
\mathcal A
:=
\{j\le p:(\bm\beta_0^{\mathrm{rls}})_j\neq0\}.
\]
Thus, under the regular design \(\mathcal A=\{1,2,5\}\), whereas under
the singular and nearly-singular designs
\(\mathcal A=\{1,2,3,5\}\) for the numerical design in
\eqref{eq: beta0_rls_numbers_sim}.

To estimate parameter $\bm\beta_0^{\mathrm{rls}}$,
we build the modified Ridgeless estimator~\eqref{eq:modified_Ridgeless_def} by spectrally hard-thresholding the sample design matrix \(\bm Q_n\), thereby
obtaining the thresholded matrix \(\check{\bm Q}_n\)
in equation \eqref{eq:check_Qn_def}.
We present most of our simulation results using the tuning parameter choice $\nu_n=n^{-3/8}$ for eigenvalue thresholding,
but the Monte Carlo evidence is essentially unchanged for any
parametrization \(\nu_n=n^{-\alpha}\) with \(\alpha\in[3/8,1/2)\);\footnote{The admissible tuning parameter range implied by
Proposition~\ref{pro:cons_rank_Qn_check} for our simulation setting
with \(\tau_n=\sqrt n\)
is \(0<\alpha<1/2\).}
see
Figure~\ref{fig:QSERidgeless} below.

\subsection{Estimators and results}\label{subsec:mc_results}

We report Monte Carlo results for four estimators. As benchmarks, we consider the standard Ridgeless estimator (RL) and our modified Ridgeless estimator (MRL),
\[
\hat{\bm\beta}_n^{\mathrm{rls}}:=\bm Q_n^+(\bm X'\bm Y/n),\qquad
\check{\bm\beta}_n^{\mathrm{rls}}:=\check{\bm Q}_n^+(\bm X'\bm Y/n) .\]
We also consider proximal estimators obtained by applying the Adaptive Lasso penalty to these benchmarks. Specifically, we apply the adaptive
\(\ell_1\) proximal step with tuning parameter \(\lambda_n\) to \(\hat{\bm\beta}_n^{\mathrm{rls}}\) and \(\check{\bm\beta}_n^{\mathrm{rls}}\), yielding the estimators RLAL \(\hat{\bm\beta}_n^+\) and MRLAL \(\check{\bm\beta}_n^+\). For each procedure, the adaptive weights are computed from its corresponding benchmark, and the proximal criterion is centered at that same benchmark. Thus, RLAL is centered at \(\hat{\bm\beta}_n^{\mathrm{rls}}\) with metric \(\overline{\bm Q}_n\), whereas MRLAL is centered at \(\check{\bm\beta}_n^{\mathrm{rls}}\) with metric \(\overline{\check{\bm Q}}_n\).

\begin{figure}[H]
\centering
\includegraphics[scale=0.35]{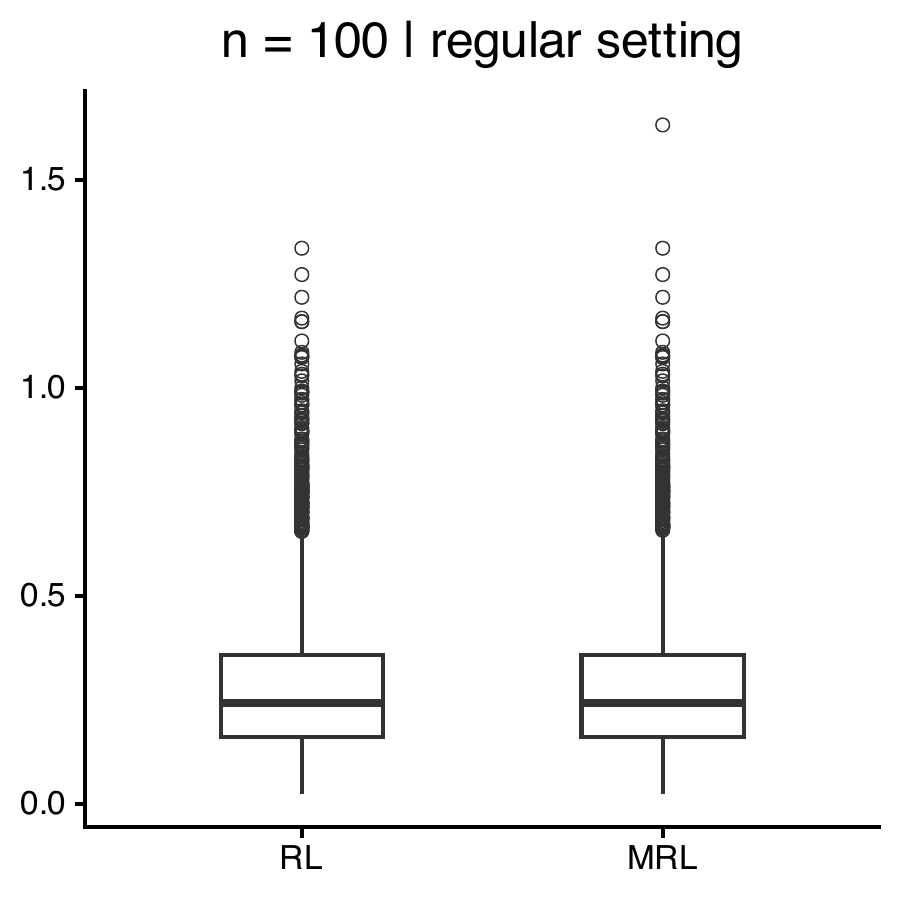}\includegraphics[scale=0.35]{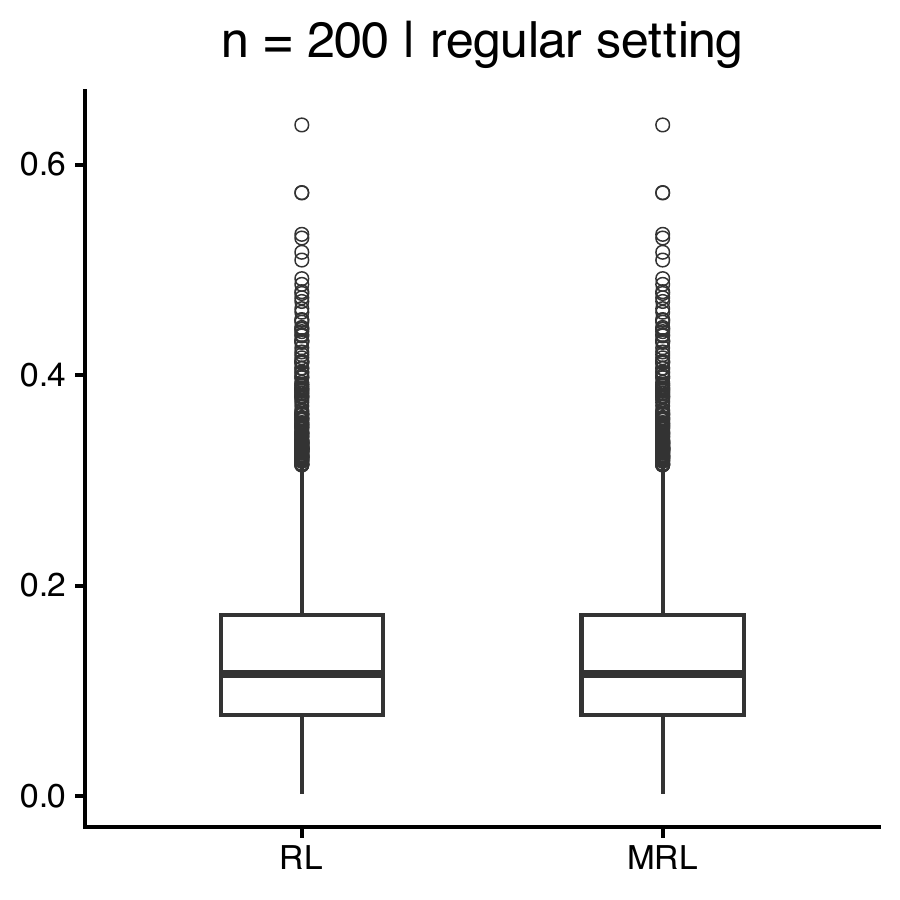}\\
\includegraphics[scale=0.35]{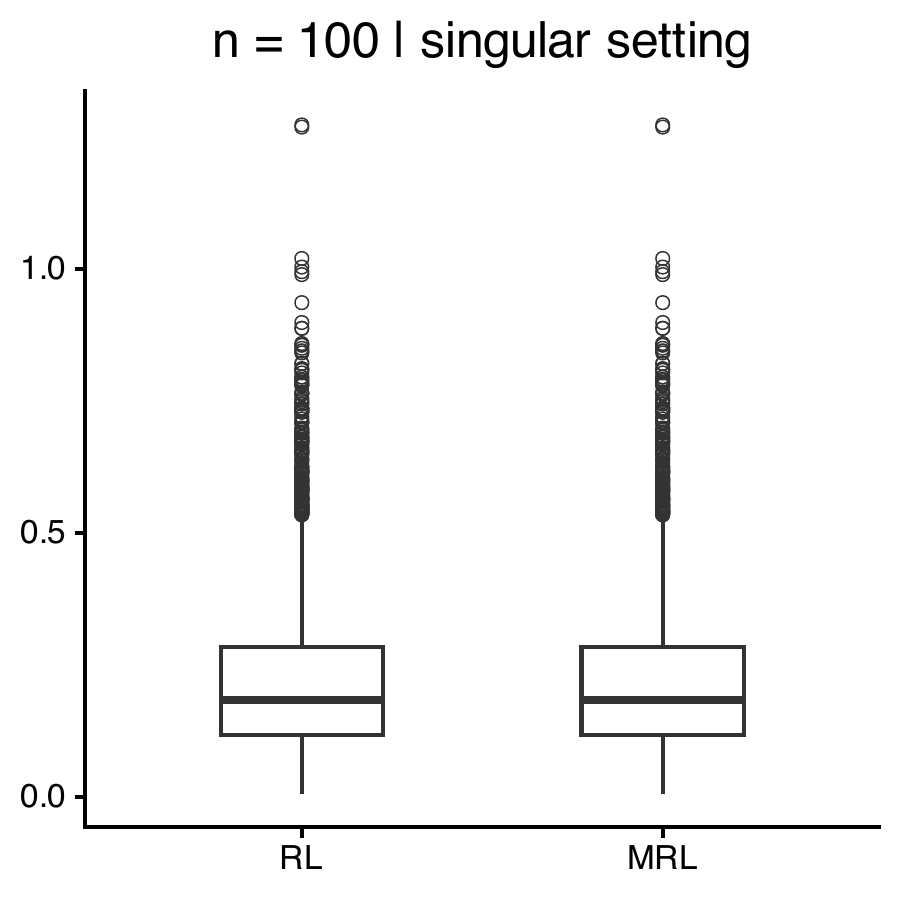}\includegraphics[scale=0.35]{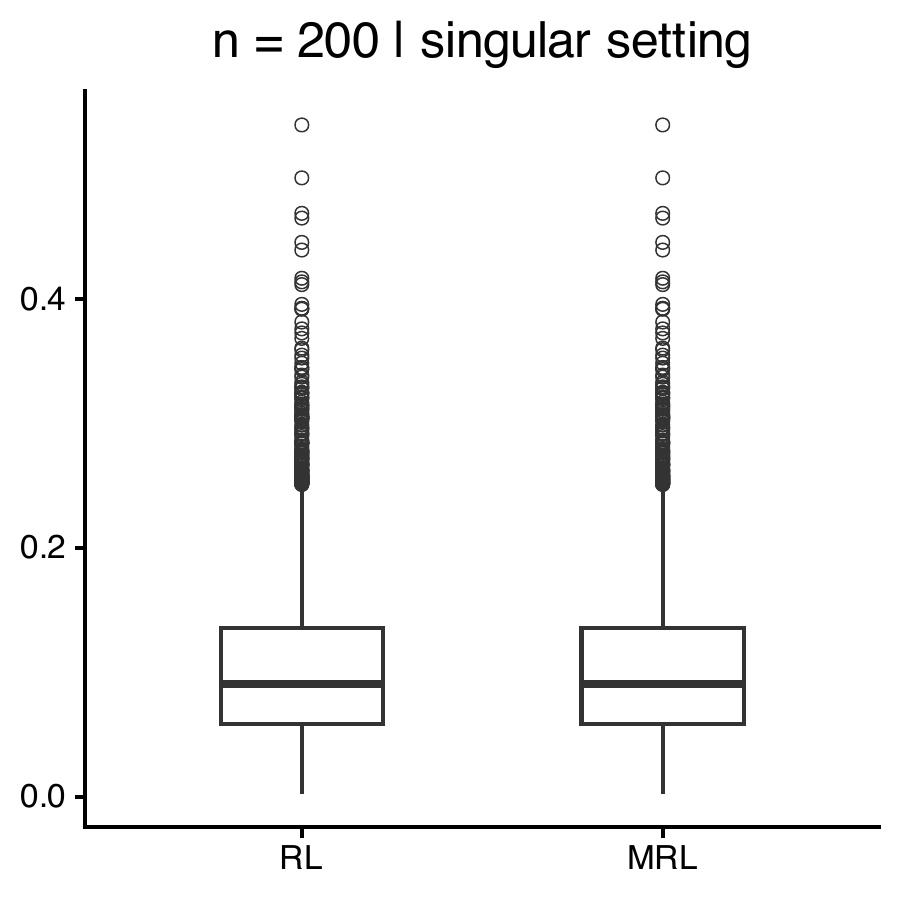}\\
\includegraphics[scale=0.35]{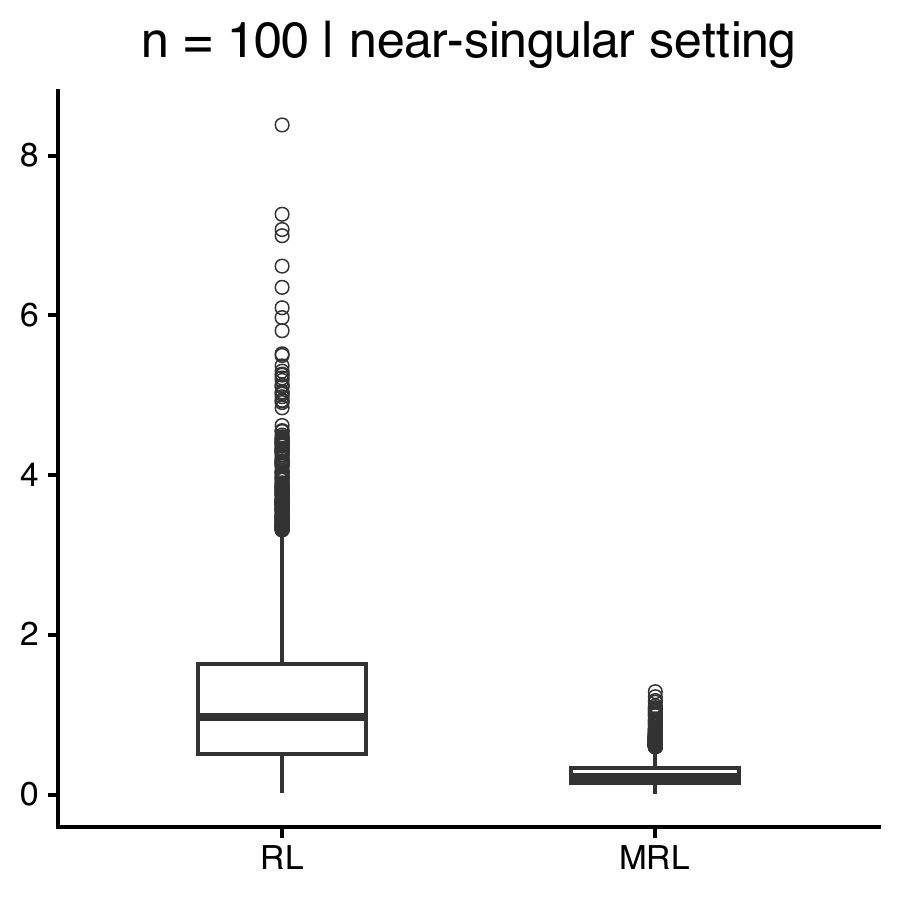}\includegraphics[scale=0.35]{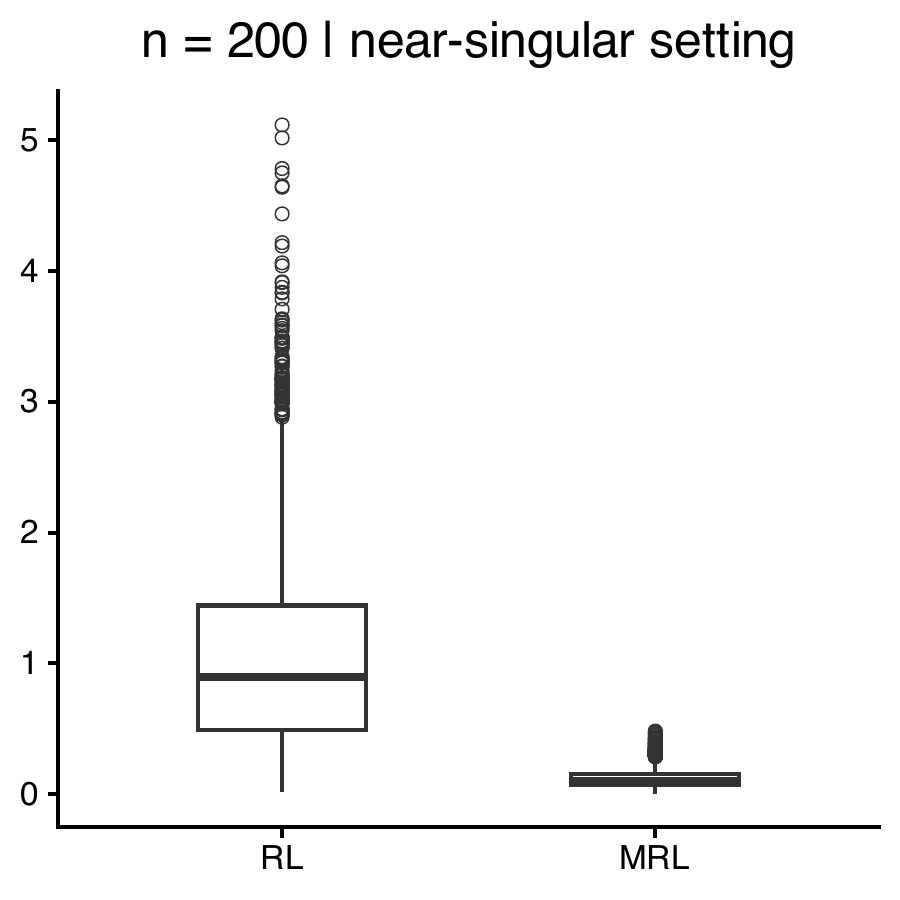}\\
\caption{Monte Carlo boxplots of sample squared errors $\|\hat{\bm\beta}_n-\bm\beta_0^{\mathrm{rls}}\|_2^2$
for Ridgeless \(\hat{\bm\beta}_n=\hat{\bm\beta}_n^{rls}\) (RL) and
modified Ridgeless \(\hat{\bm\beta}_n=\check{\bm\beta}_n^{rls}\) (MRL), using a tuning parameter $\nu_n= n^{-3/8}$
and sample sizes $n=100, 200$,
under the regular, singular, and nearly-singular designs, respectively.}\label{fig:SERidgeless}
\end{figure}

\begin{figure}[H]
\centering
\includegraphics[scale=0.35]{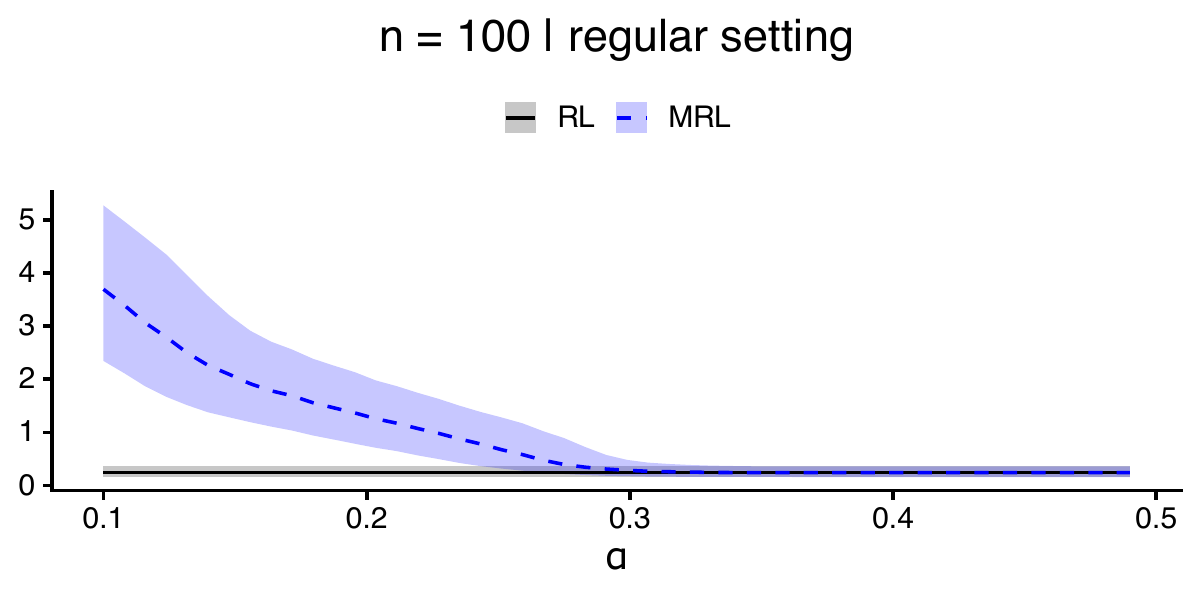}\includegraphics[scale=0.35]{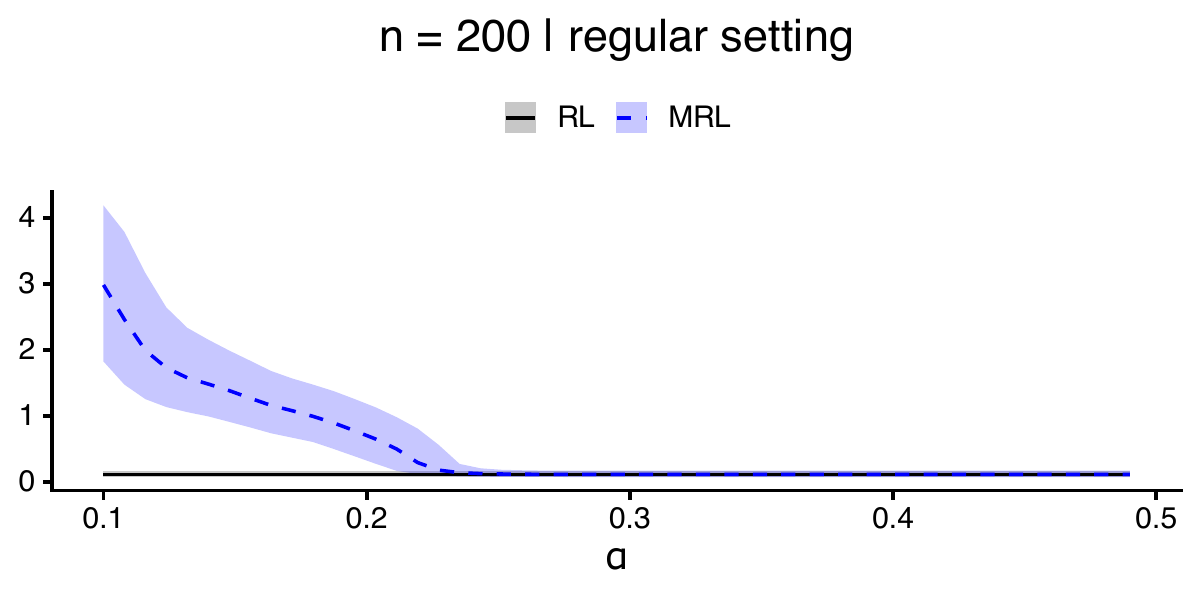}\\
\includegraphics[scale=0.35]{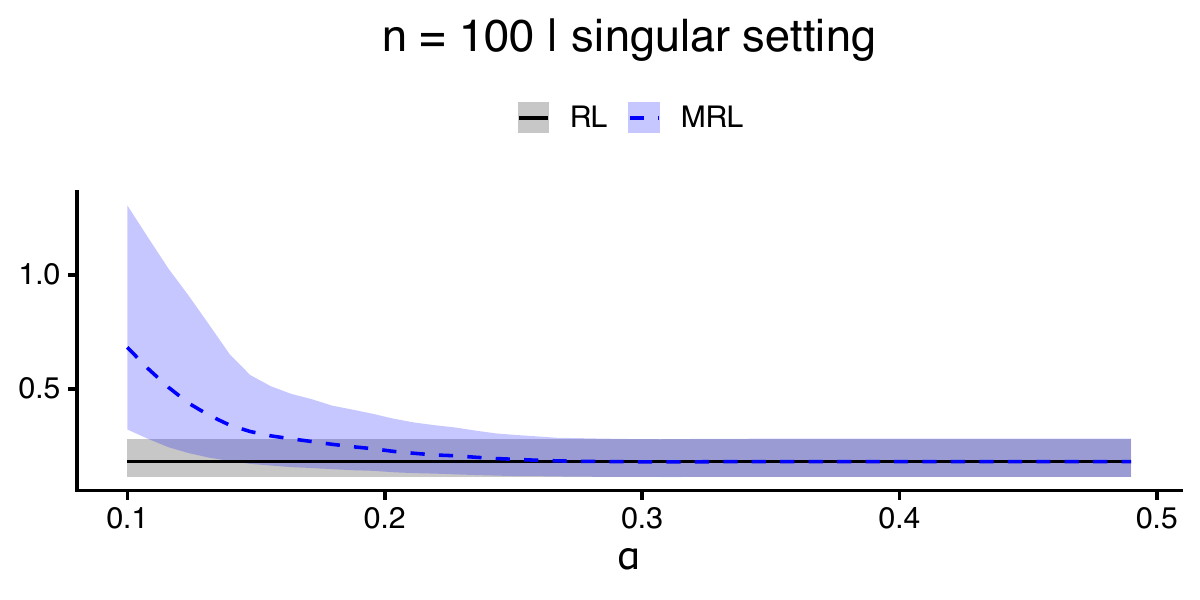}\includegraphics[scale=0.35]{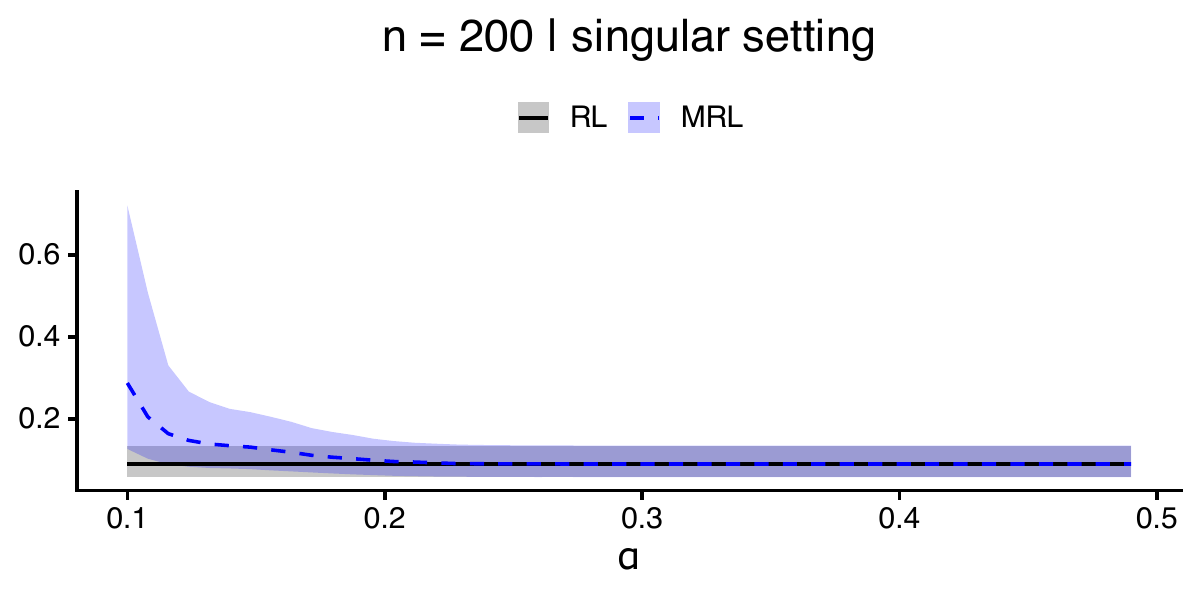}\\
\includegraphics[scale=0.35]{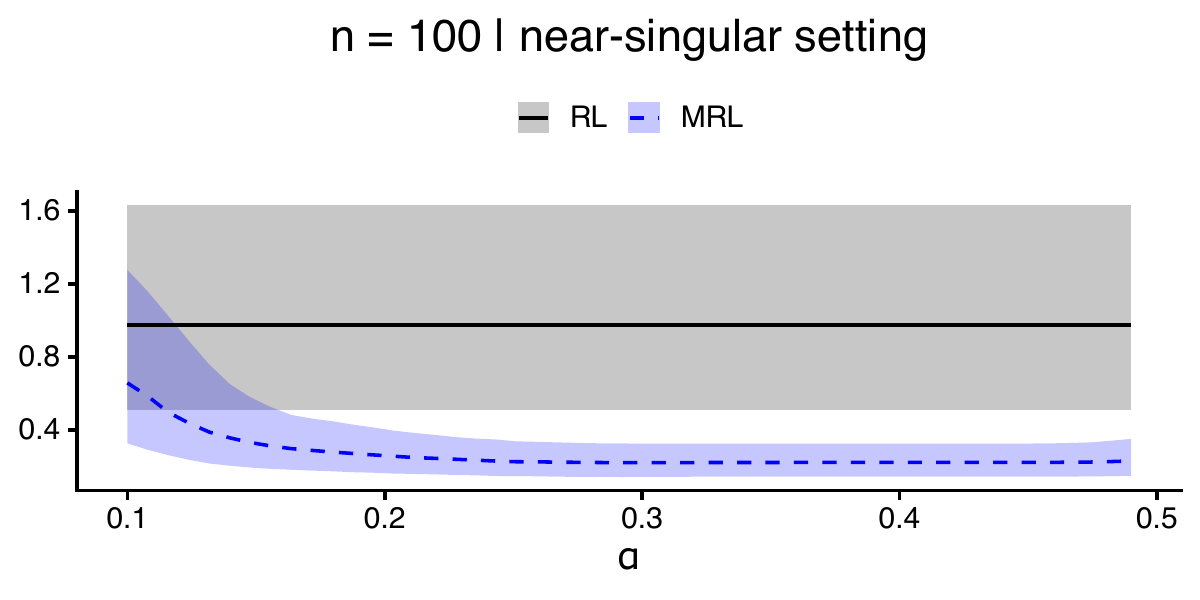}\includegraphics[scale=0.35]{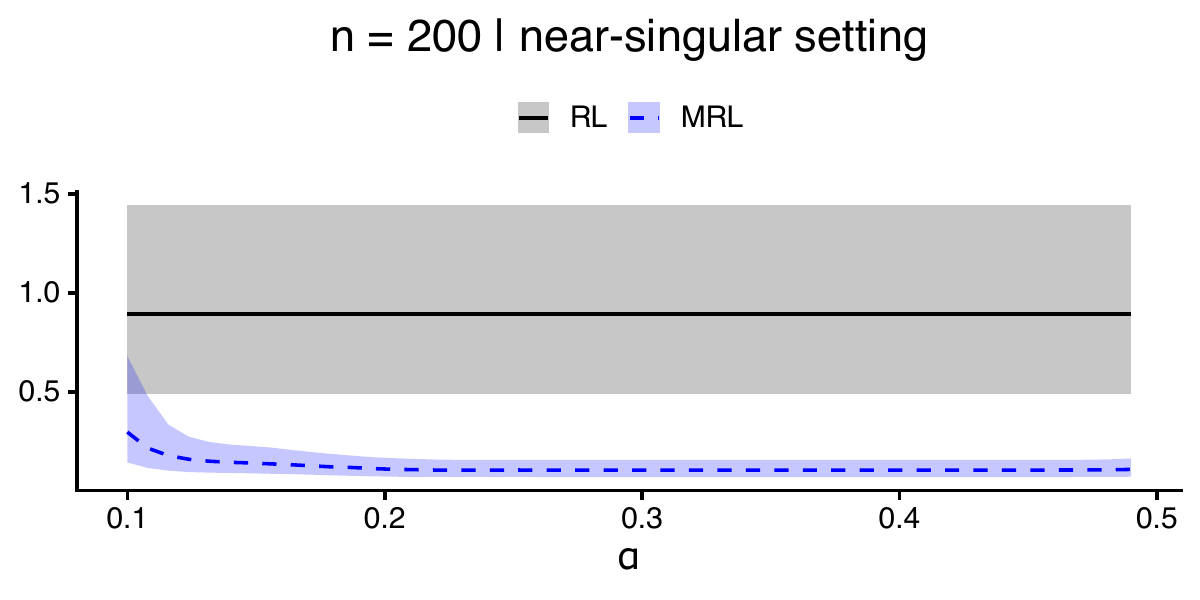}\\
\caption{Monte Carlo quartiles of sample squared errors $\|\hat{\bm\beta}_n-\bm\beta_0^{\mathrm{rls}}\|_2^2$
for Ridgeless \(\hat{\bm\beta}_n=\hat{\bm\beta}_n^{rls}\) (RL, black solid lines and areas) and
modified Ridgeless \(\hat{\bm\beta}_n=\check{\bm\beta}_n^{rls}\) (MRL, blue dashed lines and areas), using a tuning parameter $\nu_n=n^{-\alpha}$ [$\alpha\in[0.1,0.5)$]
and sample sizes $n=100, 200$,
under the regular, singular and nearly-singular design, respectively.}\label{fig:QSERidgeless}
\end{figure}

\begin{figure}[H]
\centering
\includegraphics[scale=0.35]{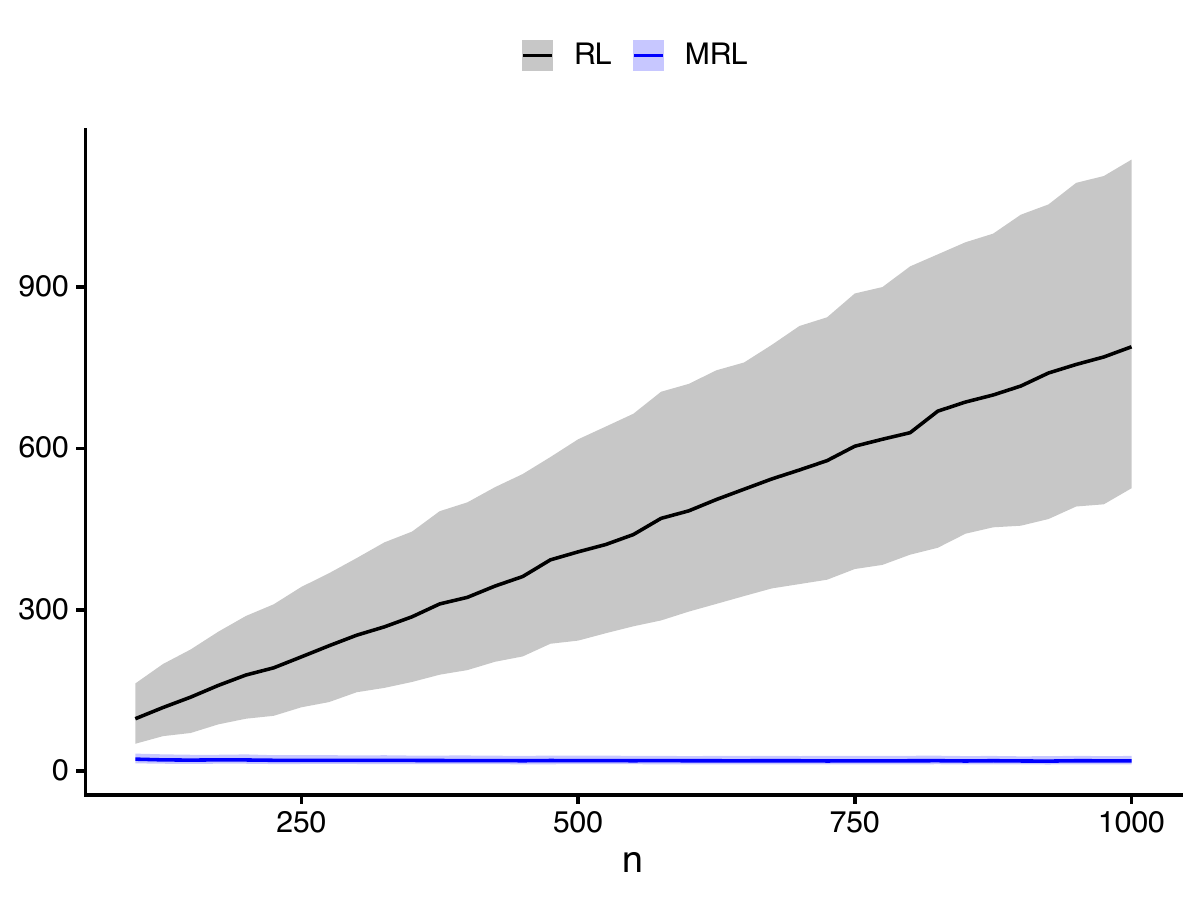}
\caption{Monte Carlo quartiles of normalized sample squared errors $n \|\hat{\bm\beta}_n-\bm\beta_0^{\mathrm{rls}}\|_2^2$
for Ridgeless \(\hat{\bm\beta}_n=\hat{\bm\beta}_n^{rls}\) (RL, black line and area) and modified Ridgeless \(\hat{\bm\beta}_n=\check{\bm\beta}_n^{rls}\) (MRL, blue line and area)
initial estimators, using a tuning parameter $\nu_n= n^{-3/8}$ and sample sizes $n\in [100, 1000]$,
under the nearly-singular design.}
\label{fig:QSELinearEstimatorsNearSingular}
\end{figure}

\begin{figure}[H]
\centering
\includegraphics[scale=0.35]{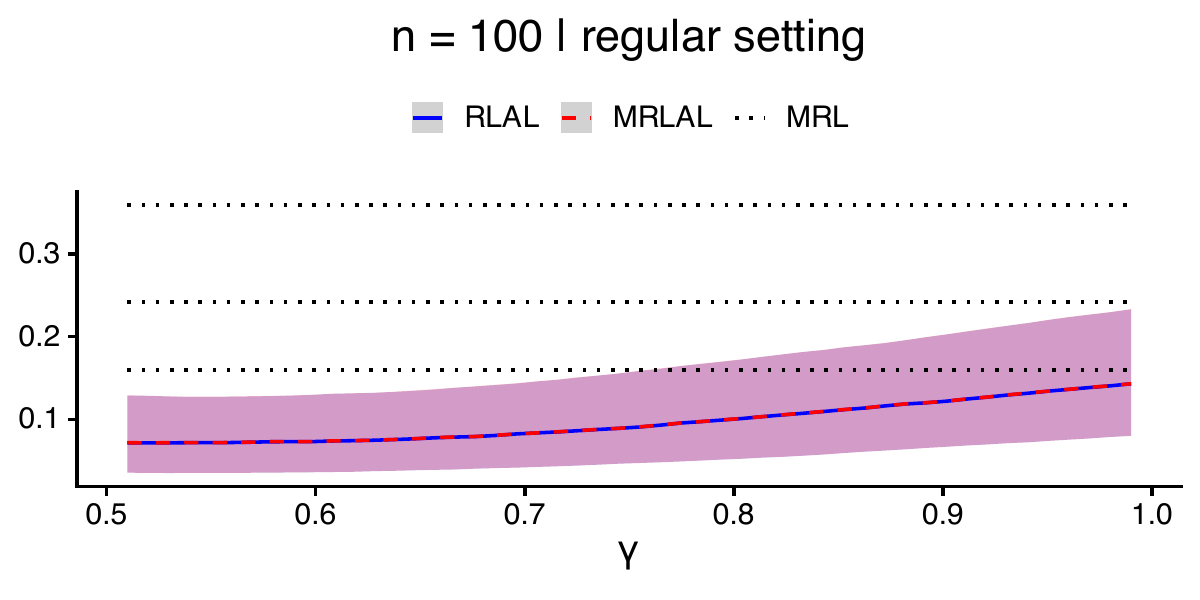}\includegraphics[scale=0.35]{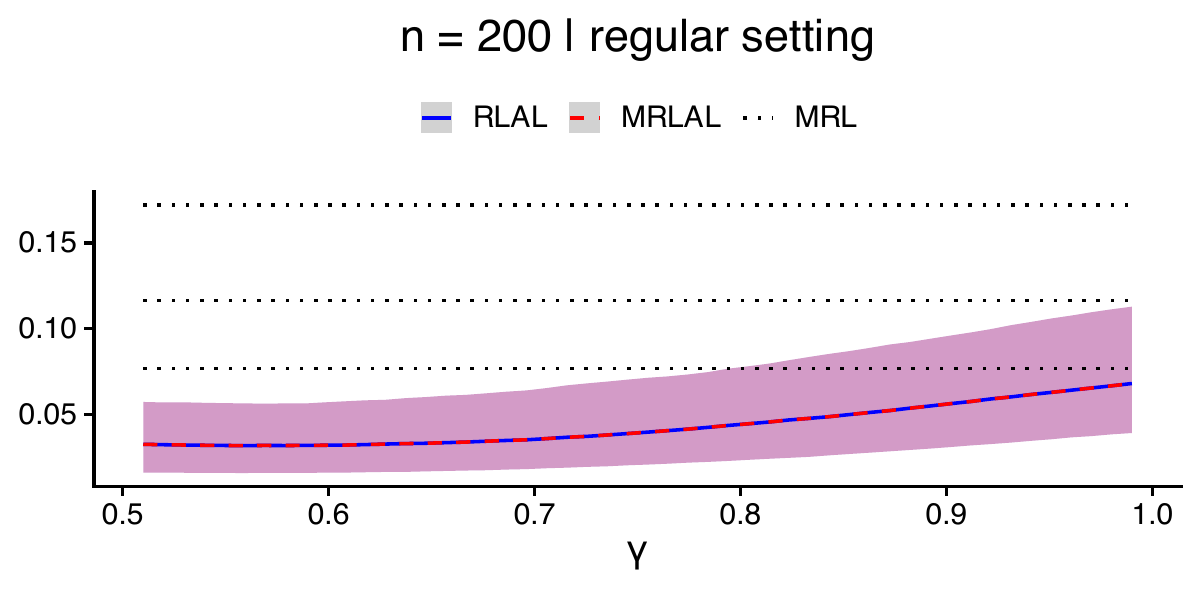}\\
\includegraphics[scale=0.35]{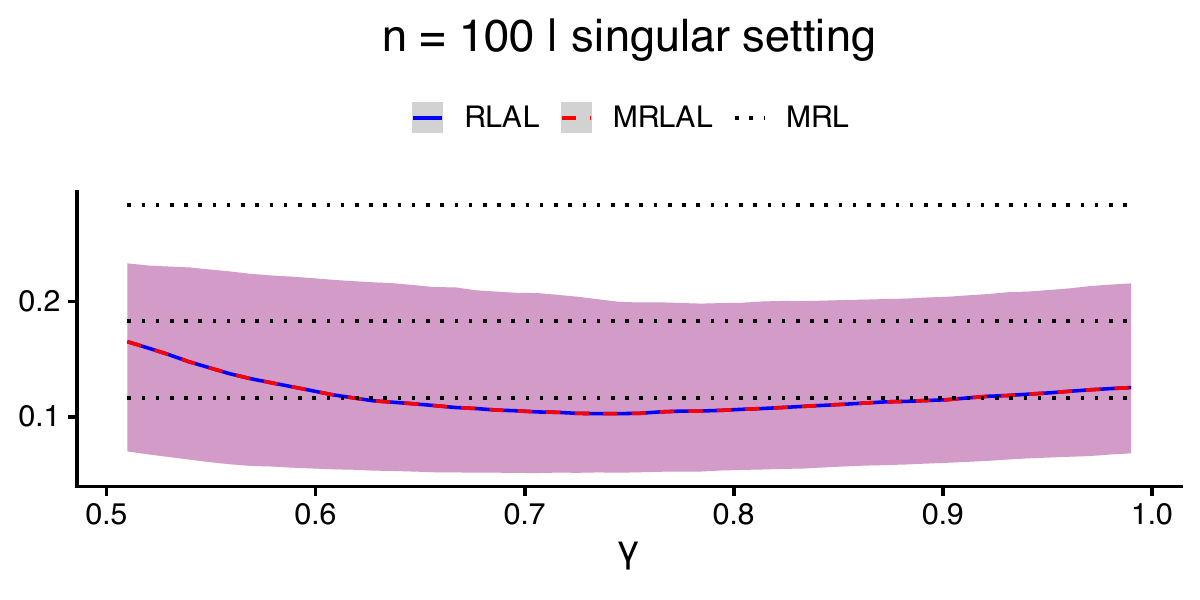}\includegraphics[scale=0.35]{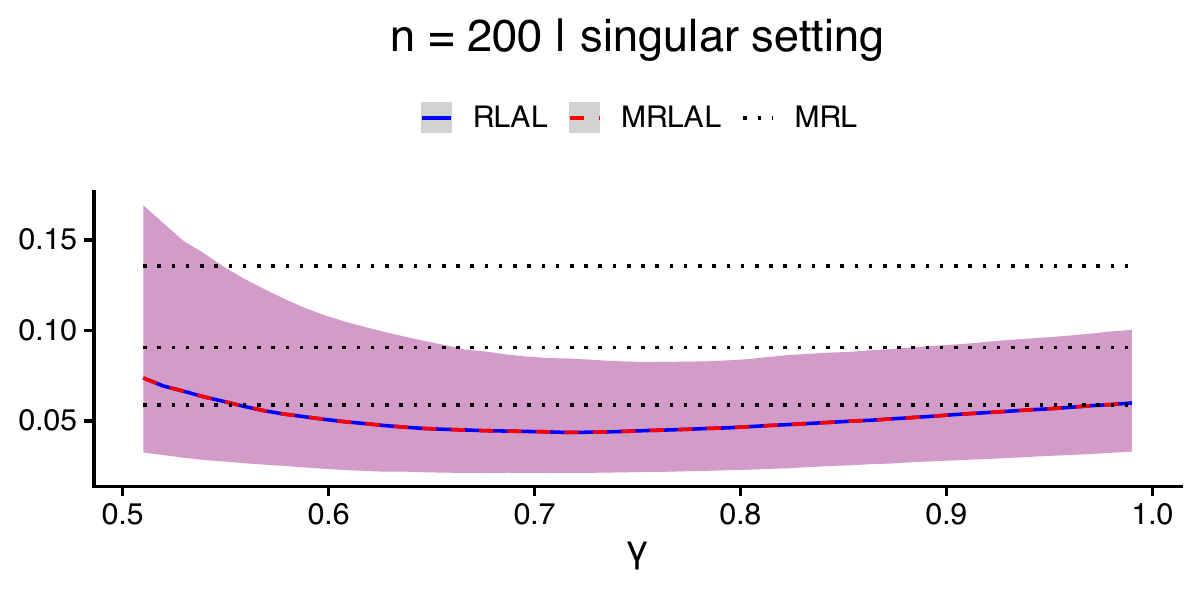}\\
\includegraphics[scale=0.35]{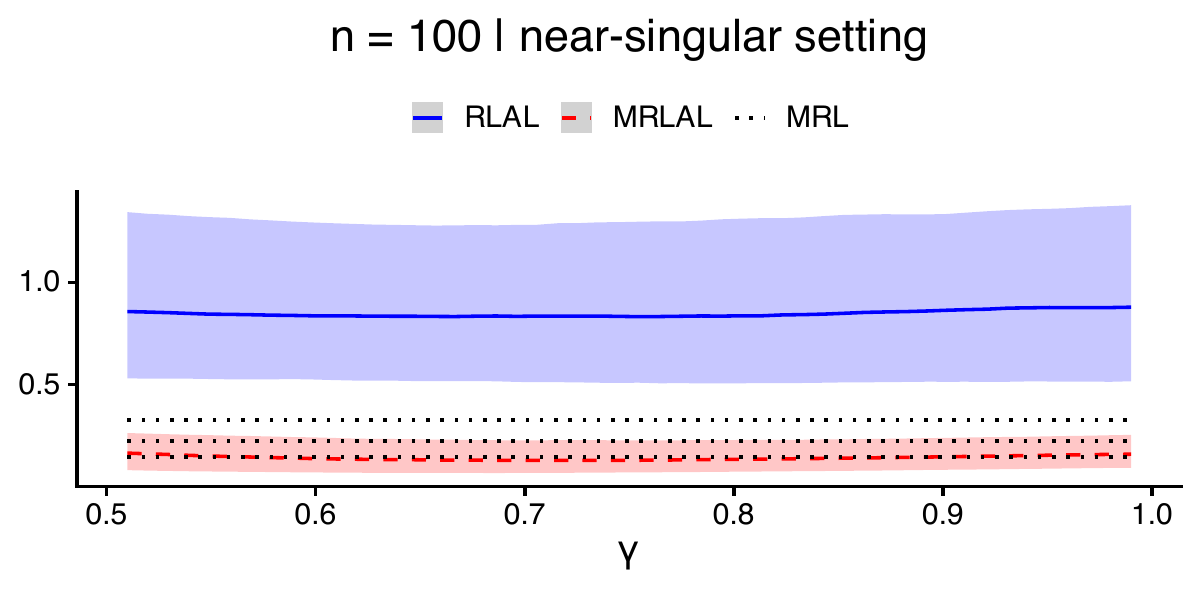}\includegraphics[scale=0.35]{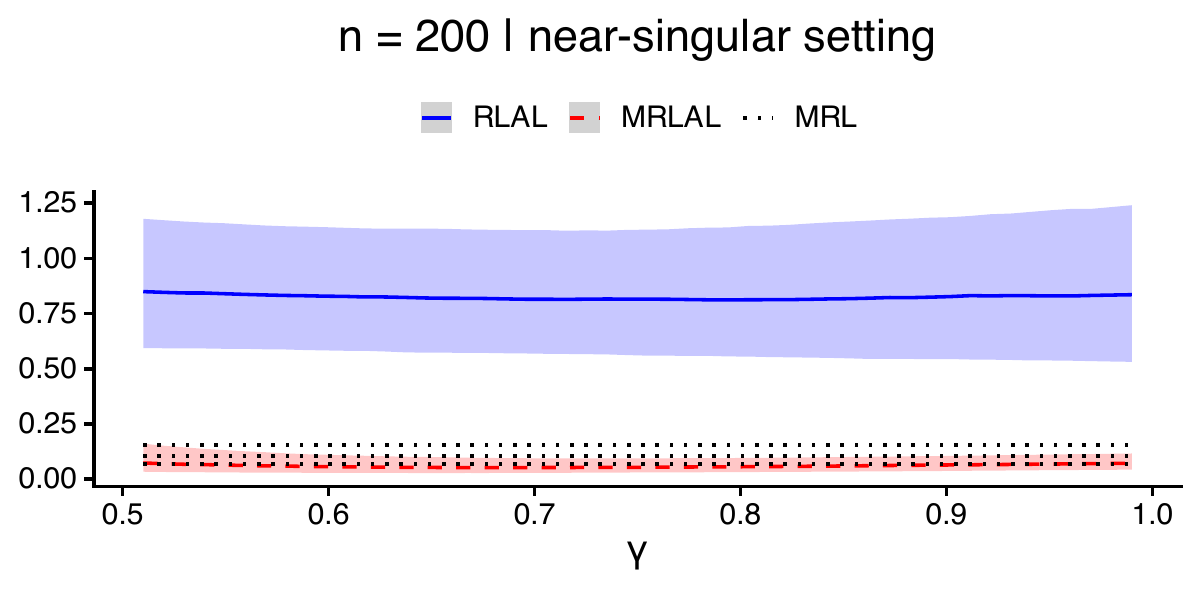}\\
\caption{Monte Carlo quartiles of sample squared errors $\|\hat{\bm\beta}_n-\bm\beta_0^{\mathrm{rls}}\|_2^2$
for modified Ridgeless \(\hat{\bm\beta}_n=\check{\bm\beta}_n^{\mathrm{rls}}\) (MRL, black dotted lines) using tuning parameter $\nu_n=n^{-3/8}$,
Ridgeless Adaptive Lasso \(\hat{\bm\beta}_n=\hat{\bm\beta}_n^{+}\) (RLAL, blue solid lines and areas) and modified Ridgeless
Adaptive Lasso \(\hat{\bm\beta}_n=\check{\bm\beta}_n^{+}\) (MRLAL, red dashed lines and areas) proximal estimators, using tuning parameters $\lambda_n =n^{-\gamma}$ [$\gamma\in(0.5,1)$] and sample sizes $n=100, 200$, under
the regular, singular and nearly-singular design, respectively.}\label{fig:QSEPLSEs}
\end{figure}

\begin{figure}[H]
\centering
\includegraphics[scale=0.35]{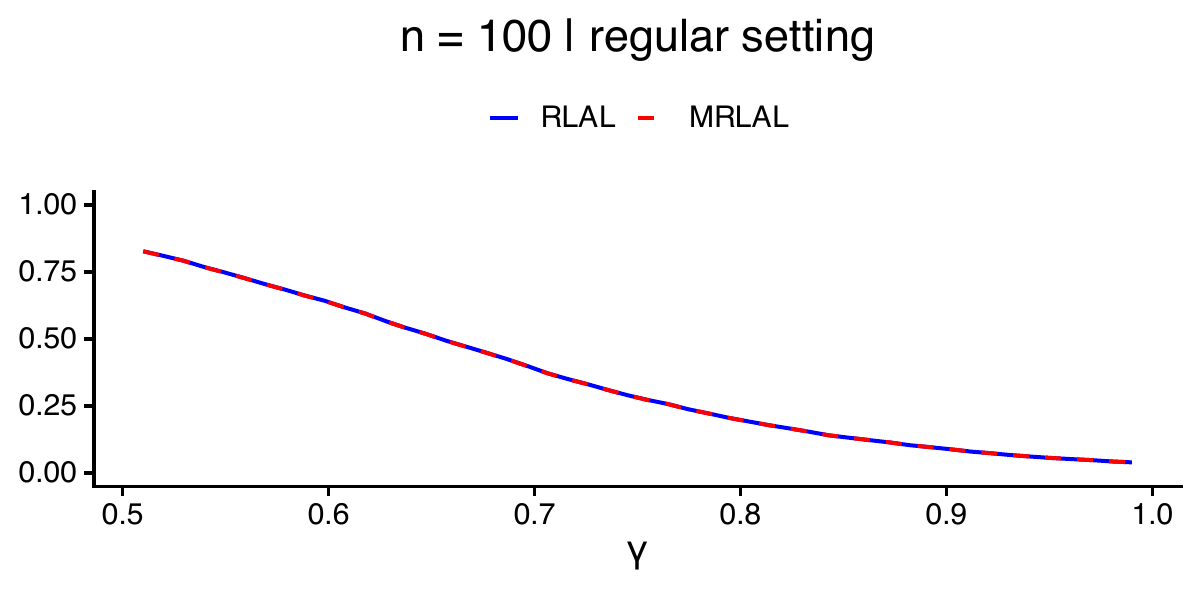}\includegraphics[scale=0.35]{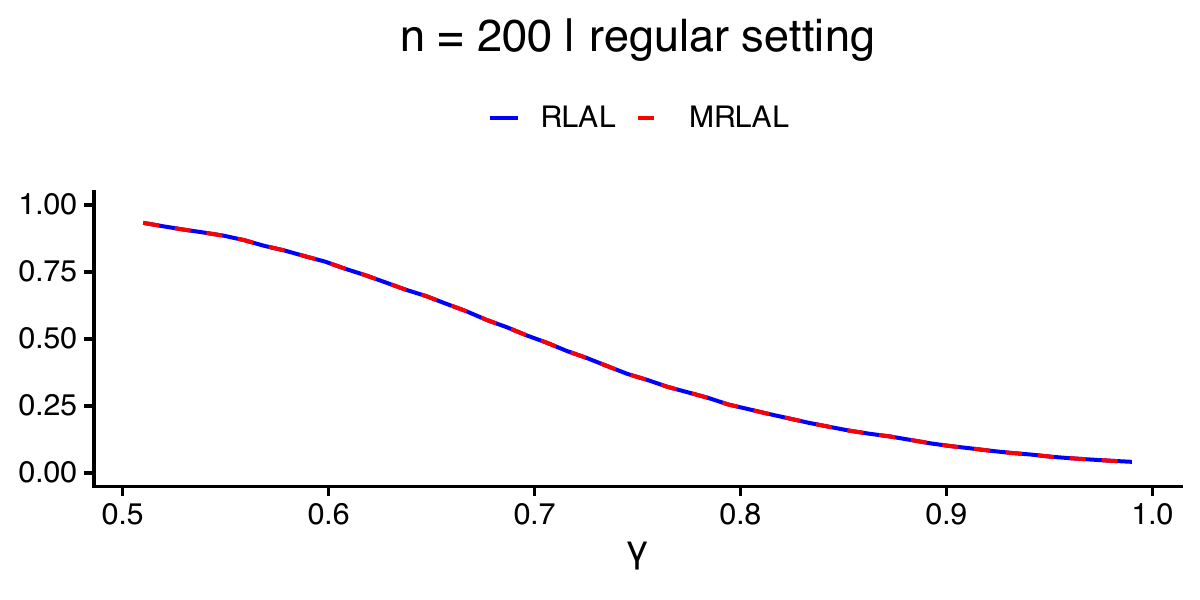}\\
\includegraphics[scale=0.35]{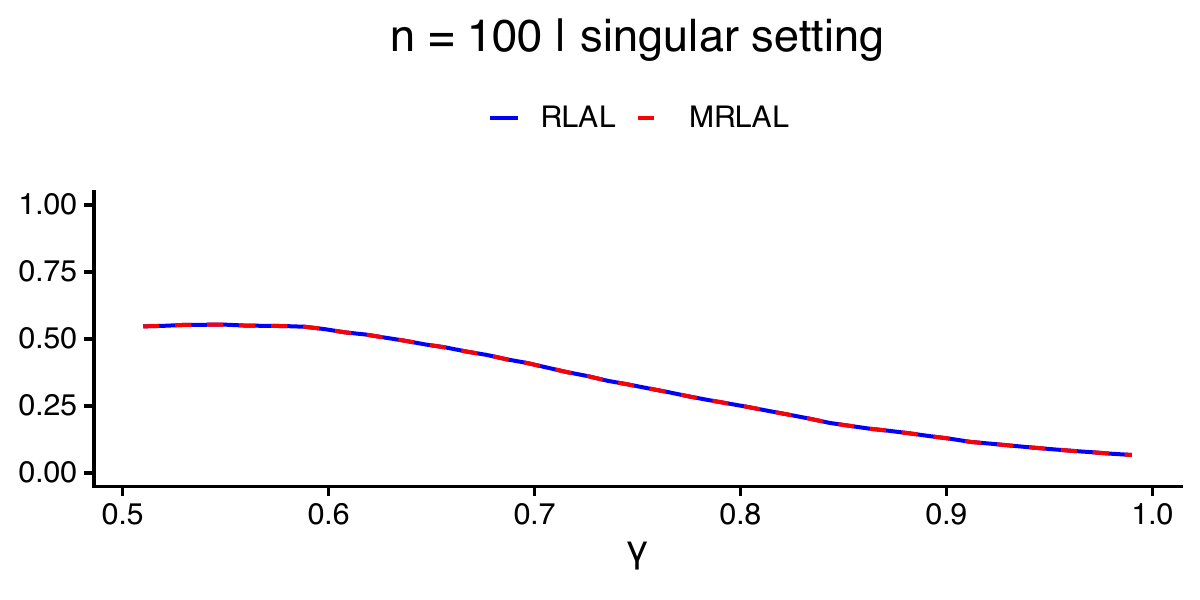}\includegraphics[scale=0.35]{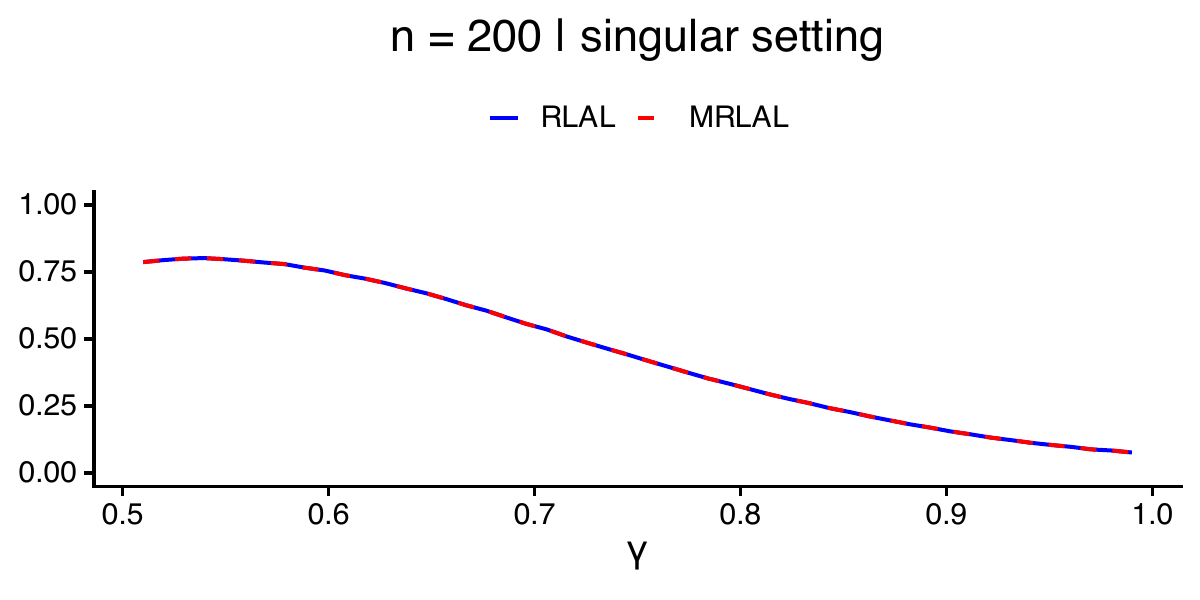}\\
\includegraphics[scale=0.35]{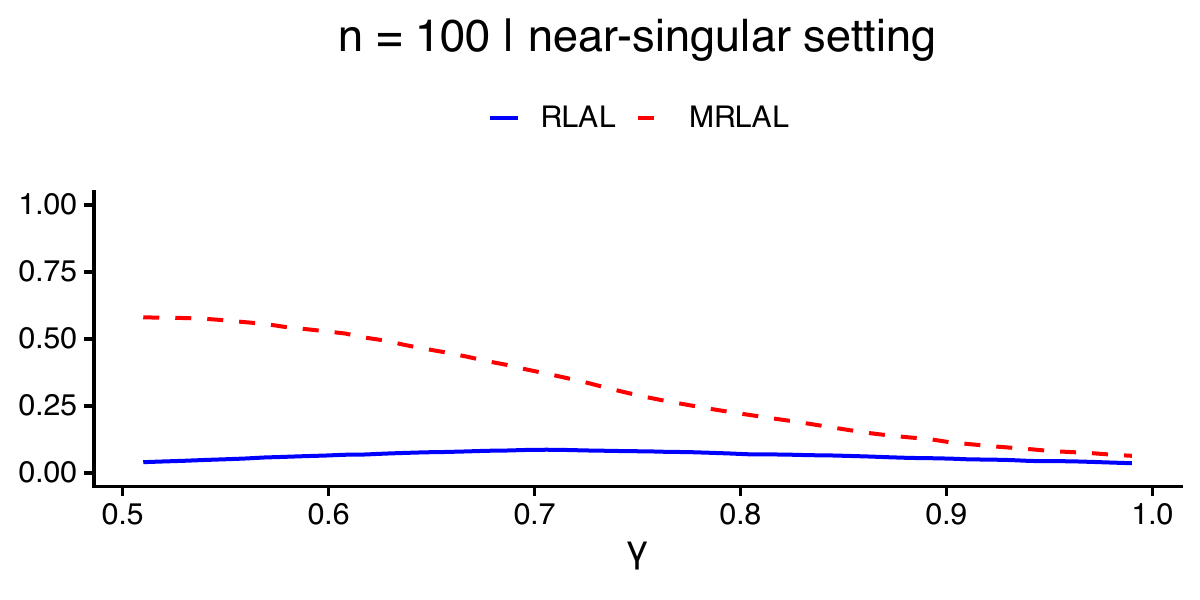}\includegraphics[scale=0.35]{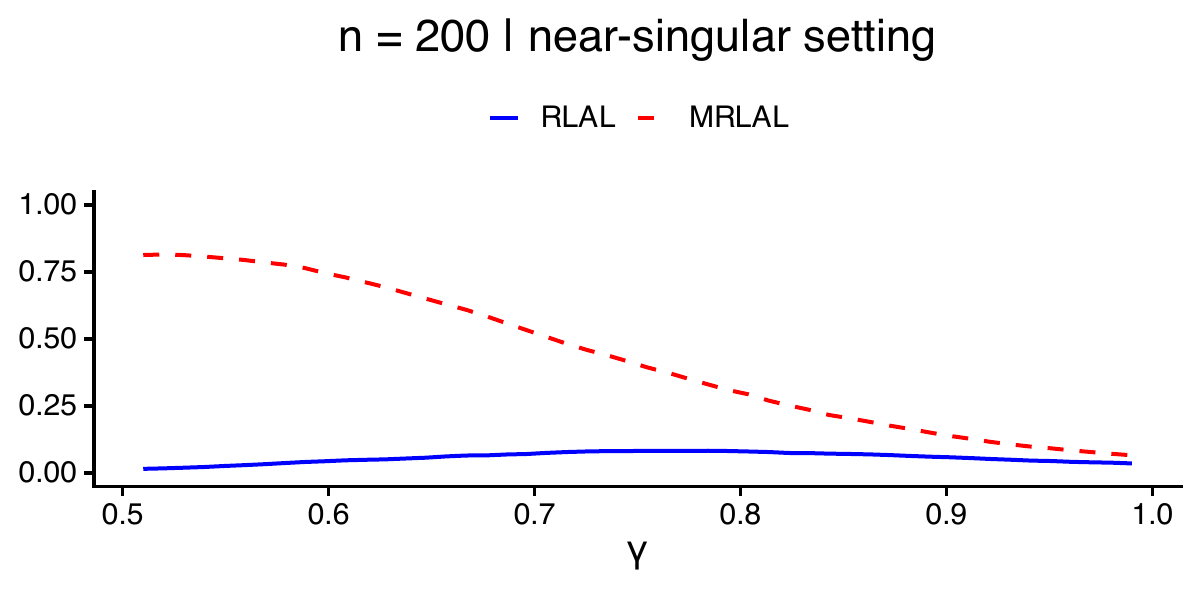}\\
\caption{Monte Carlo detection probabilities
$\Pr(\hat{\mathcal A}_n = \mathcal A)$
for Ridgeless Adaptive Lasso \(\hat{\bm\beta}_n^{+}\) (RLAL, blue solid line) and modified Ridgeless
Adaptive Lasso \(\check{\bm\beta}_n^{+}\) (MRLAL, red dashed line) proximal estimators, using tuning parameters $\lambda_n= n^{-\gamma}$
[$\gamma\in(0.5,1)$] and sample sizes $n=100, 200$,
under the regular, singular and nearly-singular design, respectively.}\label{fig:VSPLSEs}
\end{figure}

\begin{figure}[H]
\centering
\includegraphics[scale=0.35]{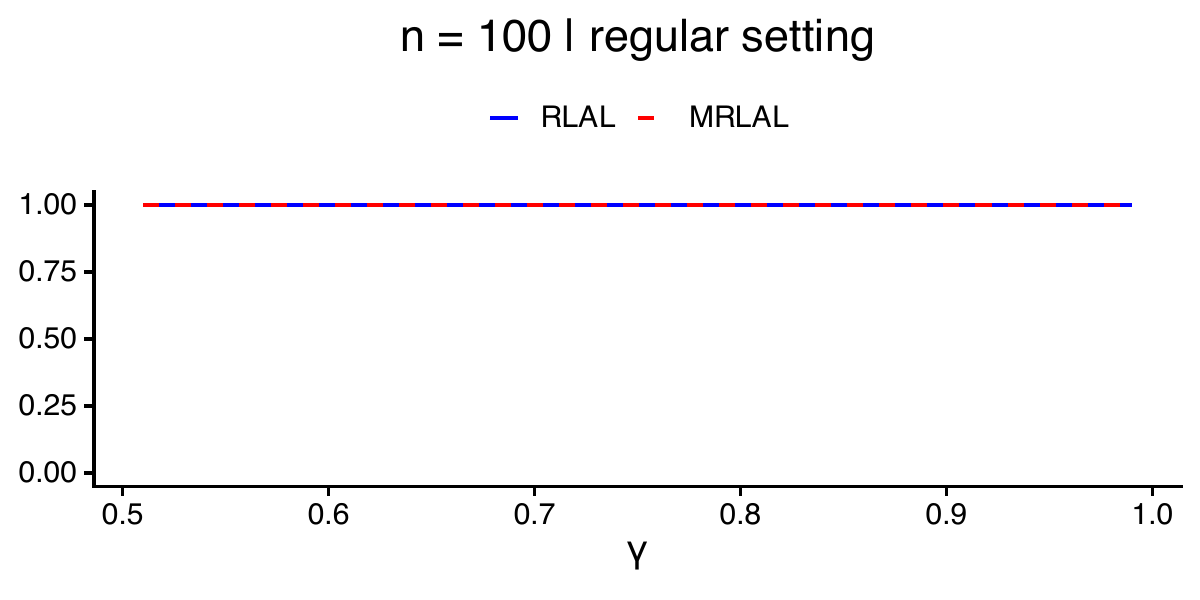}\includegraphics[scale=0.35]{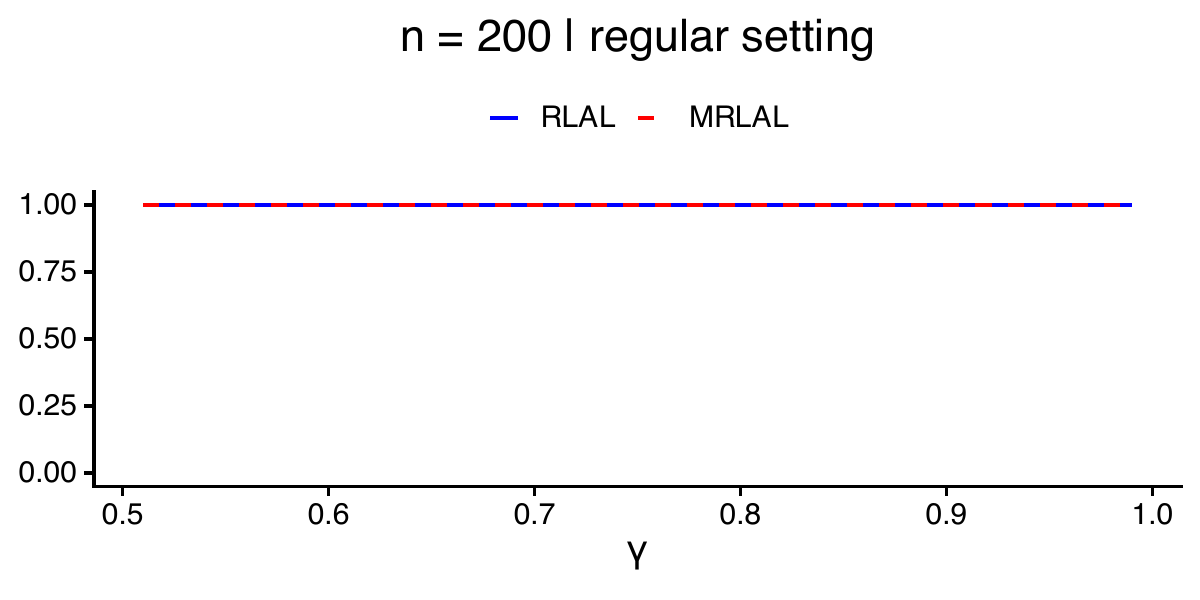}\\
\includegraphics[scale=0.35]{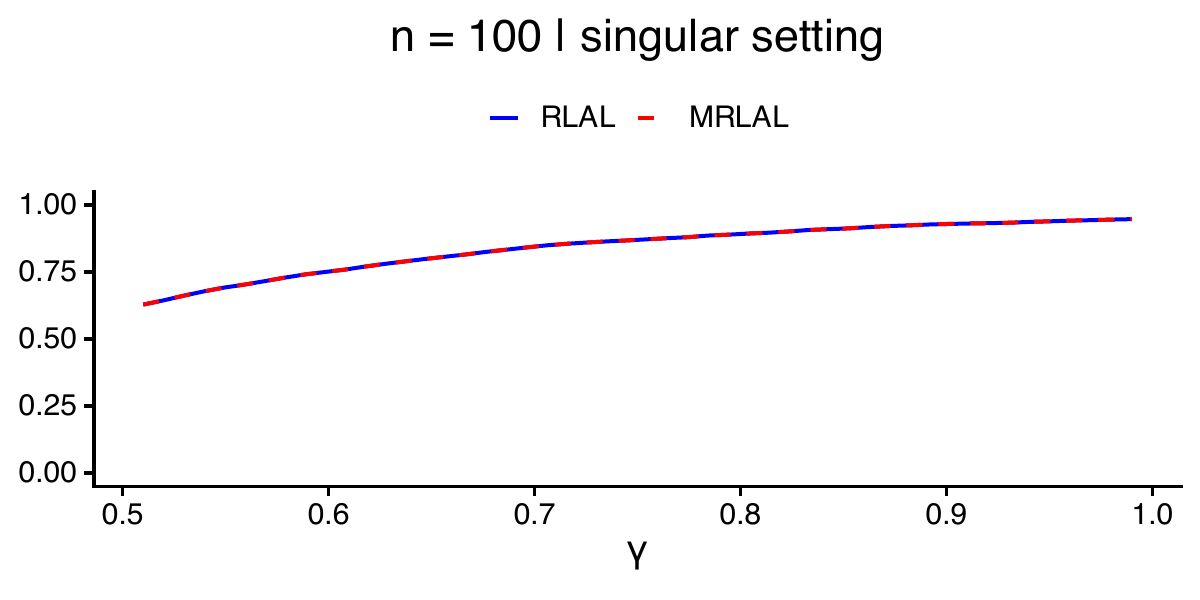}\includegraphics[scale=0.35]{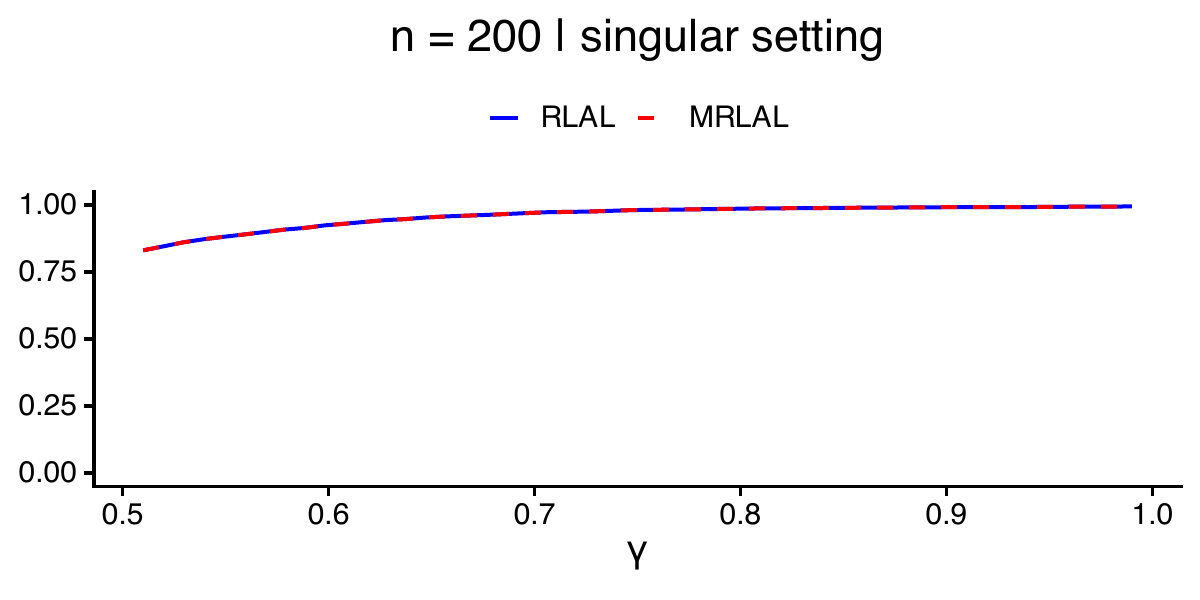}\\
\includegraphics[scale=0.35]{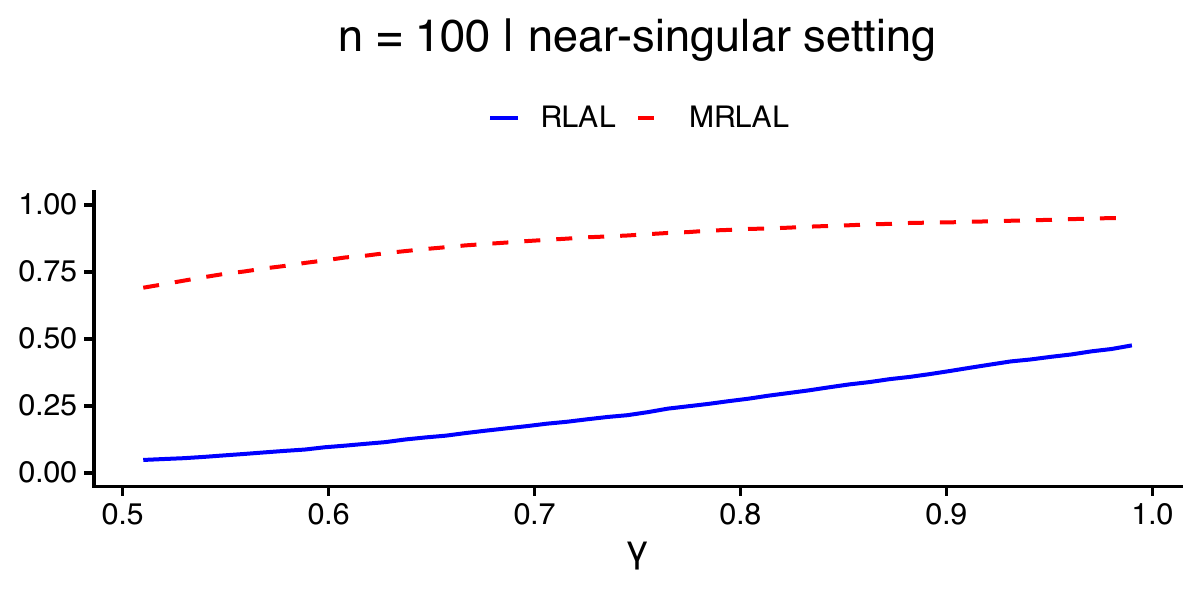}\includegraphics[scale=0.35]{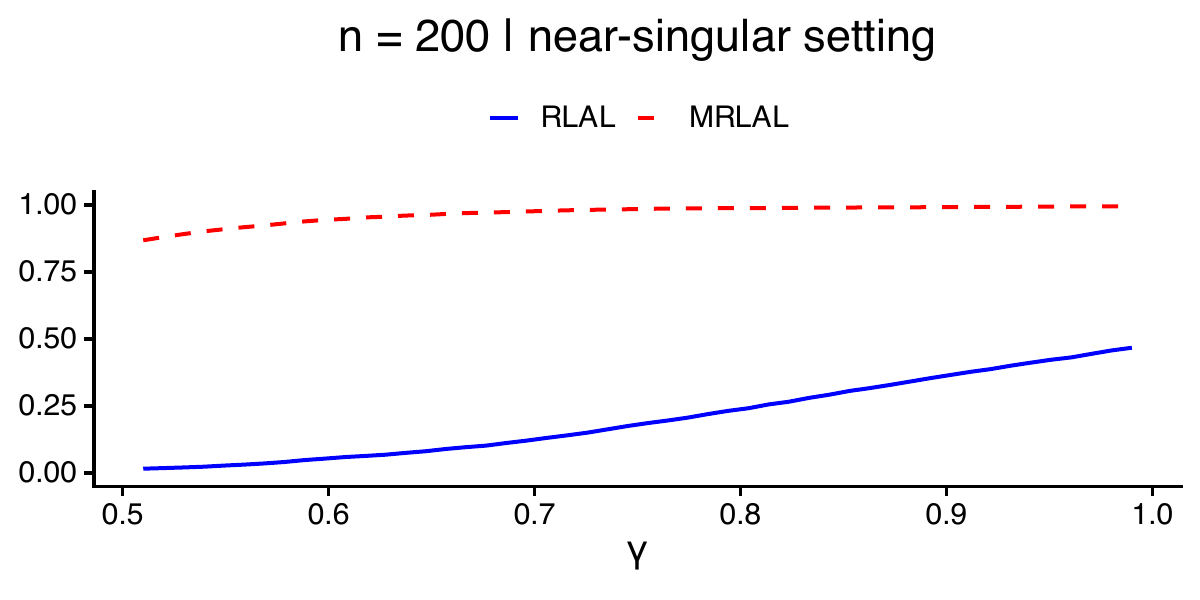}\\
\caption{Monte Carlo inclusion probabilities
$\Pr(\hat{\mathcal A}_n\supset\mathcal A)$ for Ridgeless Adaptive Lasso \(\hat{\bm\beta}_n^{+}\) (RLAL, blue solid line) and modified Ridgeless
Adaptive Lasso \(\check{\bm\beta}_n^{+}\) (MRLAL, red dashed line) proximal estimators, using tuning parameters $\lambda_n= n^{-\gamma}$
[$\gamma\in(0.5,1)$] and sample sizes $n=100, 200$,
under the regular, singular and nearly-singular design, respectively.}\label{fig:VSinclusionPLSEs}
\end{figure}

All results are based on 5,000 Monte Carlo replications. The threshold sensitivity analysis uses 50 equally spaced values of \(\alpha\) from \(0.10\) to \(0.49\), while the adaptive estimators use 50 equally spaced values of \(\gamma\) from \(0.51\) to \(0.99\). A coordinate is recorded as selected when its absolute value exceeds \(10^{-8}\).
Figures~\ref{fig:SERidgeless}--\ref{fig:VSinclusionPLSEs} below report the full set of results. We summarize here the main patterns.

Figure~\ref{fig:SERidgeless} reports boxplots of the squared estimation error for RL and MRL at sample sizes \(n=100, 200\).
Under regular and singular designs, the two estimators exhibit nearly identical distributions in finite samples.
Under the nearly-singular design, by contrast, the squared error distribution of RL is substantially more dispersed and shifted upward.
In this design, \(\bm Q_{0n}\) is positive definite for every \(n\), so RL continues to target the prelimit coefficient \(\bm\beta_0\), rather than the limiting Ridgeless estimand \(\bm\beta_0^{\mathrm{rls}}\). Because these parameters differ in our numerical design, RL is not consistent for \(\bm\beta_0^{\mathrm{rls}}\). By contrast, MRL restores the standard \(\sqrt n\)-rate for the Ridgeless estimand through a rank recovery step that removes directions whose eigenvalues vanish asymptotically before pseudoinversion.

Figure~\ref{fig:QSERidgeless} examines sensitivity to the eigenvalue threshold over the stated range of \(\alpha\). Small values of \(\alpha\) produce relatively large thresholds and can discard stable directions at the sample sizes considered. As \(\alpha\) approaches the baseline choice \(3/8\), the MRL curves become essentially flat across all three designs, while its advantage under near singularity remains evident. This supports the baseline choice and illustrates that asymptotic admissibility need not imply identical performance in finite samples when the threshold is relatively large.

This difference is even more apparent in Figure~\ref{fig:QSELinearEstimatorsNearSingular}, which reports quartiles of the normalized squared error
\[
n\|\hat{\bm\beta}_n-\bm\beta_0^{\mathrm{rls}}\|_2^2,
\]
where \(\hat{\bm\beta}_n=\hat{\bm\beta}_n^{\mathrm{rls}}\) or \(\hat{\bm\beta}_n=\check{\bm\beta}_n^{\mathrm{rls}}\),
under the nearly-singular design as \(n\) increases.
For RL, the normalized squared error remains highly dispersed and does not stabilize, whereas the corresponding distribution for MRL stabilizes rapidly. This behavior is consistent with Theorem~\ref{thm:asy_modified_Ridgeless}, which establishes a standard \(\sqrt n\)-limit for the modified Ridgeless estimator under singular and nearly-singular designs.

The properties of the initial estimators propagate to the associated proximal estimators.
Figure~\ref{fig:QSEPLSEs} reports Monte Carlo quartiles of squared estimation errors for RLAL and MRLAL and also displays the MRL benchmark, over a grid of tuning parameters \(\lambda_n=n^{-\gamma}\) such that \(\gamma\in(1/2,1)\). This is the admissible range in Corollary~\ref{cor:oracle_prox_est} for the conditions
\(\lambda_n\sqrt n\to 0\) and \(\lambda_n n\to\infty\) to hold.

Under regular and singular designs, RLAL and MRLAL perform similarly, as do their benchmark estimators. Under the nearly-singular design, by contrast, RLAL inherits the instability of RL, while MRLAL remains well behaved and systematically improves upon MRL across essentially all tuning choices. Thus, under near singularity, the proximal step alone does not repair the instability of RL; the favorable behavior of MRLAL relies on first stabilizing the benchmark through rank recovery.

Figure~\ref{fig:VSPLSEs} reports the Monte Carlo variable-selection probabilities
\[
\Pr(\hat{\mathcal A}_n=\mathcal A) ,
\]
for RLAL and MRLAL, as functions of \(\gamma\in(1/2,1)\).
Under regular and singular designs, the two procedures yield similar selection probabilities. These probabilities increase when the sample size rises from \(100\) to \(200\) and generally decline as \(\gamma\) approaches one, in line with the theory.
Under the nearly-singular design, MRLAL is substantially more reliable. Its selection probability is higher at every reported value of \(\gamma\) and increases when the sample size rises from \(100\) to \(200\), whereas RLAL performs markedly worse and does not improve uniformly with the sample size.

Figure~\ref{fig:VSinclusionPLSEs} reports the corresponding inclusion probabilities
\[
\Pr(\hat{\mathcal A}_n\supset \mathcal A), 
\]
confirming the same pattern. Under regular and singular designs the two procedures behave similarly, while under near singularity the inclusion probability of MRLAL is higher at every reported value of \(\gamma\).

\section{Proofs}\label{sec:proofs}

Throughout, assertions involving random variables are understood to hold almost surely, unless stated otherwise.

\subsection{Proofs of the results in Section~\ref{sec:framework}}

To derive Propositions~\ref{prop:prox_plse_regular} and \ref{prop:prox_plse_Ridgeless}, we first establish the following more general representation result.

\begin{proposition}[General PLSE representation of proximal estimators]\label{prop:prox_plse_general}
In model~\eqref{linear_model_setting}, let
\(\bm Q_n:=\bm X'\bm X/n\),
and let $\bm A_n$ be a symmetric positive semidefinite $p\times p$ matrix with Moore--Penrose inverse $\bm A_n^{+}$.
Define the linear benchmark estimator
\[
\hat{\bm\beta}_n^{s}:=\bm A_n^{+}\bm X'\bm Y/n,
\]
and the positive definite extension of \(\bm A_n\),
\[
\overline{\bm A}_n
:=\bm A_n+\big(\bm I_p-\bm A_n\bm A_n^{+}\big).
\]
Fix $\lambda_n>0$ and $f_n\in\Gamma(\R^p)$, and define the adjusted penalty
\begin{equation}\label{eq:extended_penalty_general}
\overline f_n(\bm\beta)
:=f_n(\bm\beta)+\frac{1}{2\lambda_n}\bm\beta'(\overline{\bm A}_n-\bm Q_n)\bm\beta.
\end{equation}
If
\begin{equation}\label{eq:kernel_condition_general}
\Kernel(\bm A_n)\subseteq \Kernel(\bm Q_n)
\qquad\text{and}\qquad
\overline{\bm A}_n-\bm Q_n\succeq \bm 0,
\end{equation}
then $\overline{\bm A}_n\succ\bm 0$, $\overline f_n\in\Gamma(\R^p)$, and
\begin{equation}\label{eq:plse_prox_identity_general}
\prox^{\overline{\bm A}_n}_{\lambda_n f_n}\big(\hat{\bm\beta}_n^{s}\big)
=
\argmin_{\bm\beta\in\R^p}
\left\{
\frac{1}{2n}\|\bm Y-\bm X\bm\beta\|_2^2+\lambda_n \overline f_n(\bm\beta)
\right\}.
\end{equation}
In particular, the minimizer on the right-hand side is unique.
\end{proposition}

\begin{proof}[Proof of Proposition~\ref{prop:prox_plse_general}]
Fix $\lambda_n>0$ and $f_n\in\Gamma(\R^p)$.
Consider first an arbitrary symmetric matrix $\bm W_n\succ\bm 0$.
Dropping terms independent of $\bm\beta$, we have
\[
\prox^{\bm W_n}_{\lambda_n f_n}\big(\hat{\bm\beta}_n^{s}\big)
=
\argmin_{\bm\beta\in\R^p}
\left\{
\frac12\bm\beta'\bm W_n\bm\beta
-\bm\beta'\bm W_n\hat{\bm\beta}_n^{s}
+\lambda_n f_n(\bm\beta)
\right\}.
\]
If $\bm W_n-\bm Q_n\succeq\bm 0$, define
\[
\widetilde f_n^{\bm W}(\bm\beta)
:=
f_n(\bm\beta)
+\frac{1}{2\lambda_n}\bm\beta'(\bm W_n-\bm Q_n)\bm\beta.
\]
Then $\widetilde f_n^{\bm W}\in\Gamma(\R^p)$, and
\begin{equation}\label{eq:prox_expanded_check_revised}
\prox^{\bm W_n}_{\lambda_n f_n}\big(\hat{\bm\beta}_n^{s}\big)
=
\argmin_{\bm\beta\in\R^p}
\left\{
\frac12\bm\beta'\bm Q_n\bm\beta
-\bm\beta'\bm W_n\hat{\bm\beta}_n^{s}
+\lambda_n \widetilde f_n^{\bm W}(\bm\beta)
\right\}.
\end{equation}
On the other hand, for any $\widetilde f\in\Gamma(\R^p)$,
\begin{equation}\label{eq:plse_expanded_check_revised}
\argmin_{\bm\beta\in\R^p}
\left\{
\frac{1}{2n}\|\bm Y-\bm X\bm\beta\|_2^2+\lambda_n \widetilde f(\bm\beta)
\right\}
=
\argmin_{\bm\beta\in\R^p}
\left\{
\frac12\bm\beta'\bm Q_n\bm\beta
-\bm\beta'\bm X'\bm Y/n
+\lambda_n \widetilde f(\bm\beta)
\right\},
\end{equation}
since the two objectives differ only by the constant $\|\bm Y\|_2^2/(2n)$.
Hence the minimizers in \eqref{eq:prox_expanded_check_revised} and
\eqref{eq:plse_expanded_check_revised} coincide whenever
\[
\bm W_n-\bm Q_n\succeq\bm 0
\qquad\text{and}\qquad
\bm W_n\hat{\bm\beta}_n^{s}=\bm X'\bm Y/n.
\]

Now set
\[
\bm W_n=\overline{\bm A}_n
:=\bm A_n+\big(\bm I_p-\bm A_n\bm A_n^+\big),
\qquad
\hat{\bm\beta}_n^{s}=\bm A_n^+\bm X'\bm Y/n.
\]
Because $\bm A_n$ is symmetric positive semidefinite, $\bm A_n\bm A_n^+$ is the orthogonal projector onto $\Range(\bm A_n)$ and $\bm I_p-\bm A_n\bm A_n^+$ is the orthogonal projector onto $\Kernel(\bm A_n)$.
Writing any $\bm v\in\R^p$ as
$\bm v=\bm v_R+\bm v_K$ with
$\bm v_R\in\Range(\bm A_n)$ and $\bm v_K\in\Kernel(\bm A_n)$, we obtain
\[
\bm v'\overline{\bm A}_n\bm v
=
\bm v_R'\bm A_n\bm v_R+\|\bm v_K\|_2^2.
\]
Since $\bm A_n$ is positive definite on $\Range(\bm A_n)$, the right-hand side is strictly positive for every $\bm v\neq\bm 0$. Thus $\overline{\bm A}_n\succ\bm 0$.
Next,
\[
\overline{\bm A}_n\bm A_n^+
=
(\bm A_n+\bm I_p-\bm A_n\bm A_n^+)\bm A_n^+
=
\bm A_n\bm A_n^+,
\]
because $\Range(\bm A_n^+)=\Range(\bm A_n)$ implies
$(\bm I_p-\bm A_n\bm A_n^+)\bm A_n^+=\bm 0$.
Therefore
\[
\overline{\bm A}_n\hat{\bm\beta}_n^s
=
\bm A_n\bm A_n^+\bm X'\bm Y/n.
\]
Now $\bm X'\bm Y/n\in\Range(\bm X')=\Range(\bm Q_n)$.
Since $\bm A_n$ and $\bm Q_n$ are symmetric,
$\Kernel(\bm A_n)\subseteq\Kernel(\bm Q_n)$ is equivalent to
$\Range(\bm Q_n)\subseteq\Range(\bm A_n)$.
Hence $\bm X'\bm Y/n\in\Range(\bm A_n)$, and so
\[
\overline{\bm A}_n\hat{\bm\beta}_n^s
=
\bm X'\bm Y/n.
\]
Finally, the condition $\overline{\bm A}_n-\bm Q_n\succeq\bm 0$ implies that
\[
\overline f_n(\bm\beta)
=
f_n(\bm\beta)
+\frac{1}{2\lambda_n}\bm\beta'(\overline{\bm A}_n-\bm Q_n)\bm\beta
\]
belongs to $\Gamma(\R^p)$.
Applying the preceding argument with $\bm W_n=\overline{\bm A}_n$ and
$\widetilde f_n^{\bm W}=\overline f_n$ yields
\[
\prox^{\overline{\bm A}_n}_{\lambda_n f_n}\big(\hat{\bm\beta}_n^{s}\big)
=
\argmin_{\bm\beta\in\R^p}
\left\{
\frac{1}{2n}\|\bm Y-\bm X\bm\beta\|_2^2+\lambda_n \overline f_n(\bm\beta)
\right\}.
\]
The minimizer on the right-hand side is unique because
$\overline{\bm A}_n\succ\bm 0$, so the left-hand side is a well-defined proximal estimator.
\end{proof}

\begin{proof}[Proof of Proposition~\ref{prop:prox_plse_regular}]
Apply Proposition~\ref{prop:prox_plse_general} with \(\bm A_n=\bm Q_n\).
Since \(\bm Q_n\succ\bm 0\), we have \(\bm Q_n^+=\bm Q_n^{-1}\) and
\[
\hat{\bm\beta}_n^s
=
\bm A_n^+\bm X'\bm Y/n
=
\bm Q_n^{-1}\bm X'\bm Y/n
=
\hat{\bm\beta}_n^{\mathrm{ls}}.
\]
Moreover,
\[
\overline{\bm A}_n
=
\bm Q_n+\big(\bm I_p-\bm Q_n\bm Q_n^+\big)
=
\bm Q_n,
\]
because \(\bm Q_n\bm Q_n^+=\bm I_p\), and
\[
\overline f_n(\bm\beta)
=
f_n(\bm\beta)+\frac{1}{2\lambda_n}\bm\beta'(\overline{\bm A}_n-\bm Q_n)\bm\beta
=
f_n(\bm\beta).
\]
The conditions in \eqref{eq:kernel_condition_general} are immediate: \(\Kernel(\bm Q_n)=\{\bm 0\}\), so
\[
\Kernel(\bm A_n)\subseteq\Kernel(\bm Q_n),
\qquad
\overline{\bm A}_n-\bm Q_n=\bm 0\succeq\bm 0.
\]
Hence Proposition~\ref{prop:prox_plse_general} yields
\[
\prox_{\lambda_n f_n}^{\bm Q_n}\big(\hat{\bm\beta}_n^{\mathrm{ls}}\big)
=
\argmin_{\bm\beta\in\R^p}
\left\{
\frac{1}{2n}\|\bm Y-\bm X\bm\beta\|_2^2+\lambda_n f_n(\bm\beta)
\right\}.
\]
In particular, the minimizer is unique.
\end{proof}

\begin{proof}[Proof of Proposition~\ref{prop:prox_plse_Ridgeless}]
Apply Proposition~\ref{prop:prox_plse_general} with \(\bm A_n=\bm Q_n\).
Then
\[
\hat{\bm\beta}_n^{s}
=\bm Q_n^{+}\bm X'\bm Y/n
=\hat{\bm\beta}_n^{\mathrm{rls}},
\qquad
\overline{\bm A}_n
=
\bm Q_n+(\bm I_p-\bm Q_n\bm Q_n^{+})
=
\overline{\bm Q}_n.
\]
Moreover,
\[
\Kernel(\bm A_n)=\Kernel(\bm Q_n)\subseteq\Kernel(\bm Q_n),
\qquad
\overline{\bm A}_n-\bm Q_n
=
\bm I_p-\bm Q_n\bm Q_n^{+}\succeq\bm 0,
\]
since \(\bm Q_n\bm Q_n^+\) is the orthogonal projector onto \(\Range(\bm Q_n)\).
Hence Proposition~\ref{prop:prox_plse_general} yields \(\overline{\bm Q}_n\succ\bm 0\),
\(\overline f_n\in\Gamma(\R^p)\), and
\[
\prox_{\lambda_n f_n}^{\overline{\bm Q}_n}\big(\hat{\bm\beta}_n^{\mathrm{rls}}\big)
=
\argmin_{\bm\beta\in\R^p}
\left\{
\frac{1}{2n}\|\bm Y-\bm X\bm\beta\|_2^2+\lambda_n\overline f_n(\bm\beta)
\right\},
\]
where
\[
\overline f_n(\bm\beta)
=
f_n(\bm\beta)
+\frac{1}{2\lambda_n}\bm\beta'(\bm I_p-\bm Q_n\bm Q_n^+)\bm\beta.
\]
This is exactly \eqref{eq:extended_penalty_Ridgeless}--\eqref{eq:plse_prox_identity_Ridgeless}.
\end{proof}

\begin{proof}[Proof of Proposition~\ref{prop:dual_repr}]
Since \(f_n\in\Gamma(\R^p)\) and \(\lambda_n>0\), the function \(\lambda_n f_n\) is proper, lower semicontinuous, and convex.
Therefore Moreau's decomposition theorem applies; see, e.g., \citet[Thm.~14.3]{SM:bauschke2016convex}. It yields, for every \(\bm\theta\in\R^p\),
\[
\prox_{\lambda_n f_n}^{\bm W_n}(\bm\theta)
+
\prox_{(\lambda_n f_n)^*}^{\bm W_n}(\bm\theta)
=
\bm\theta,
\]
where \((\lambda_n f_n)^*\) is the convex conjugate of \(\lambda_n f_n\) with respect to scalar product \(\langle\cdot,\cdot\rangle_{\bm W_n}\). Rearranging gives
\[
\prox_{\lambda_n f_n}^{\bm W_n}(\bm\theta)
=
\bm\theta-\prox_{(\lambda_n f_n)^*}^{\bm W_n}(\bm\theta).
\]
Since this holds for every \(\bm\theta\in\R^p\), equation \eqref{eq:moreau_W} follows.
Finally, with conjugacy taken with respect to \(\langle\cdot,\cdot\rangle_{\bm W_n}\),
\[
(\lambda_n f_n)^*(\bm\theta)
=
\lambda_n f_n^*(\bm\theta/\lambda_n),
\qquad \bm\theta\in\R^p.
\]
When \(f_n\) is sublinear, because \(f_n\in\Gamma(\R^p)\), the function \(\lambda_n f_n\) is proper, lower semicontinuous, convex, and sublinear. Hence, by \citet[Thm.~3.1.1]{SM:hiriart2004fundamentals}, there exists a unique nonempty closed convex set \(C_n\subseteq\R^p\) such that
\[
\lambda_n f_n(\bm\beta)
=
\sigma_{C_n}(\bm\beta)
:=
\sup_{\bm\theta\in C_n}\langle \bm\theta,\bm\beta\rangle_{\bm W_n},
\qquad \bm\beta\in\R^p,
\]
and
\[
C_n
=
\Big\{
\bm\theta\in\R^p:\ 
\langle \bm\theta,\bm\beta\rangle_{\bm W_n}\le \lambda_n f_n(\bm\beta)
\ \text{for all }\bm\beta\in\R^p
\Big\}.
\]
Moreover, by \citet[Ex.~13.3(i)]{SM:bauschke2016convex}, the \(\bm W_n\)-Fenchel conjugate of a support function is the indicator of its support set. Therefore
\[
(\lambda_n f_n)^*
=
(\sigma_{C_n})^*
=
\iota_{C_n}.
\]
Applying Proposition~\ref{prop:dual_repr} gives
\[
\prox_{\lambda_n f_n}^{\bm W_n}
=
\Id-\prox_{\iota_{C_n}}^{\bm W_n}.
\]
Finally, for each \(\bm\theta\in\R^p\),
\[
\prox_{\iota_{C_n}}^{\bm W_n}(\bm\theta)
=
\argmin_{\tilde{\bm\theta}\in C_n}
\frac12\|\tilde{\bm\theta}-\bm\theta\|_{\bm W_n}^2
=
P_{C_n}^{\bm W_n}(\bm\theta),
\]
where the minimizer is unique because \(\bm W_n\succ\bm 0\) and \(C_n\) is nonempty, closed, and convex. Substituting into the previous display yields \eqref{eq:proj_formula_W}.
\end{proof}

\begin{proof}[Proof of Corollary~\ref{cor:proj_formulas_lasso_adal}]
Both penalties in (i) and (ii) are sublinear, so the dual representation~\ref{eq:proj_formula_W} applies. Hence
\begin{align*}
C_n
&=
\Bigl\{
\bm\theta\in\R^p:\ 
\langle \bm\theta,\bm\beta\rangle_{\bm W_n}\le \lambda_n f_n(\bm\beta)
\ \text{for all }\bm\beta\in\R^p
\Bigr\}\\
&=
\Bigl\{
\bm\theta\in\R^p:\ 
\langle \bm W_n\bm\theta,\bm\beta\rangle\le \lambda_n f_n(\bm\beta)
\ \text{for all }\bm\beta\in\R^p
\Bigr\}.
\end{align*}

\begin{enumerate}[label=(\roman*)]

\item
If \(f_n(\bm\beta)=\|\bm\beta\|_1\), then
\[
\langle \bm W_n\bm\theta,\bm\beta\rangle\le \lambda_n\|\bm\beta\|_1
\quad\text{for all }\bm\beta\in\R^p
\]
holds if and only if
\(\|\bm W_n\bm\theta\|_\infty\le \lambda_n\),
that is,
\(|( \bm W_n\bm\theta )_j|\le \lambda_n\)
for every \(j=1,\dots,p\).
Therefore
\[
C_n
=
\bigcap_{j=1}^p
\Bigl\{
\bm\theta\in\R^p:\ 
|( \bm W_n\bm\theta )_j|
\le
\lambda_n
\Bigr\}.
\]

\item
For the Adaptive Lasso, use the extended-value formulation
\[
f_n(\bm\beta)
=
\sum_{j:\tilde\beta_{nj}\neq 0}\frac{|\beta_j|}{|\tilde\beta_{nj}|}
+
\sum_{j:\tilde\beta_{nj}=0}\iota_{\{0\}}(\beta_j).
\]
Let
\[
A_n:=\{j:\tilde\beta_{nj}\neq 0\},
\qquad
A_n^c:=\{1,\ldots,p\}\setminus A_n=\{j:\tilde\beta_{nj}=0\}.
\]
Then the defining inequality for \(C_n\) becomes
\[
\langle \bm W_n\bm\theta,\bm\beta\rangle
\le
\lambda_n
\left(
\sum_{j\in A_n}\frac{|\beta_j|}{|\tilde\beta_{nj}|}
+
\sum_{j\in A_n^c}\iota_{\{0\}}(\beta_j)
\right)
\qquad
\text{for all }\bm\beta\in\R^p.
\]
If \(\beta_j\neq 0\) for some \(j\in A_n^c\), then the right-hand side equals \(+\infty\), so the inequality is automatic. Hence it is enough to consider vectors \(\bm\beta\) such that \(\beta_j=0\) for all \(j\in A_n^c\). For such \(\bm\beta\),
\[
\langle \bm W_n\bm\theta,\bm\beta\rangle
=
\sum_{j\in A_n}(\bm W_n\bm\theta)_j\beta_j,
\]
and the condition becomes
\[
\sum_{j\in A_n}(\bm W_n\bm\theta)_j\beta_j
\le
\lambda_n\sum_{j\in A_n}\frac{|\beta_j|}{|\tilde\beta_{nj}|}
\qquad
\text{for all }\bm\beta \text{ with }\bm\beta_{A_n^c}=\bm 0.
\]
This holds if and only if
\(|(\bm W_n\bm\theta)_j|
\le
\lambda_n/|\tilde\beta_{nj}|\)
for every \(j\in A_n\),
that is,
\[
C_n
=
\bigcap_{\{j:\tilde\beta_{nj}\neq 0\}}
\Bigl\{
\bm\theta\in\R^p:\ 
|( \bm W_n\bm\theta )_j|
\le
\lambda_n/|\tilde\beta_{nj}|
\Bigr\}.
\]
This concludes the proof.
\end{enumerate}
\end{proof}

\begin{proof}[Proof of Theorem~\ref{thm:properties}]
By Definition~\eqref{eq:intro_prox} of proximal estimators,  \(\bm W_n\succ\bm 0\). It follows that the bilinear form
\(\langle \bm u,\bm v\rangle_{\bm W_n}:=\bm u'\bm W_n\bm v\)
defines an inner product on \(\R^p\), with associated norm \(\|\cdot\|_{\bm W_n}\). Hence \(\R^p\), equipped with this inner product, is a Hilbert space.

\begin{enumerate}[label=(\roman*)]
\item \emph{Existence and uniqueness.}
Fix \(\bm\theta\in\R^p\) and define
\[
\phi_{\bm\theta}(\bm\beta)
:=
\frac12\|\bm\beta-\bm\theta\|_{\bm W_n}^2+\lambda_n f_n(\bm\beta).
\]
Since \(f_n\in\Gamma(\R^p)\), the function \(\phi_{\bm\theta}\) is proper and lower semicontinuous. Because \(\bm W_n\succ\bm 0\), the quadratic term is continuous, coercive, and \(1\)-strongly convex with respect to \(\|\cdot\|_{\bm W_n}\). Hence \(\phi_{\bm\theta}\) is coercive and strongly convex, so it admits a unique minimizer. By the definition of the proximal mapping under the inner product induced by \(\bm W_n\)
this minimizer is
\(\prox^{\bm W_n}_{\lambda_n f_n}(\bm\theta)\).
Evaluating at \(\bm\theta=\hat{\bm\beta}_n^s\) gives existence and uniqueness of the proximal estimator.

\item \emph{Stability.}
Since \(\lambda_n f_n\in\Gamma(\R^p)\), the proximal mapping
\(\prox^{\bm W_n}_{\lambda_n f_n}\) is firmly nonexpansive, hence nonexpansive, on the Hilbert space \((\R^p,\langle\cdot,\cdot\rangle_{\bm W_n})\); see \citet[Prop.~12.28]{SM:bauschke2016convex}. Therefore, for any \(\bm\beta_1,\bm\beta_2\in\R^p\),
\[
\Big\|
\prox^{\bm W_n}_{\lambda_n f_n}(\bm\beta_1)
-
\prox^{\bm W_n}_{\lambda_n f_n}(\bm\beta_2)
\Big\|_{\bm W_n}
\le
\|\bm\beta_1-\bm\beta_2\|_{\bm W_n},
\]
which is property \eqref{eq:stability_general}.

\item \emph{Almost-everywhere differentiability.}
By part~(ii), \(\prox^{\bm W_n}_{\lambda_n f_n}\) is \(1\)-Lipschitz under \(\|\cdot\|_{\bm W_n}\). Since \(\bm W_n\succ\bm 0\) and the space is finite-dimensional, \(\|\cdot\|_{\bm W_n}\) and the Euclidean norm are equivalent. Hence \(\prox^{\bm W_n}_{\lambda_n f_n}\) is Lipschitz continuous with respect to the Euclidean norm, hence strictly continuous.  Therefore, by Rademacher's theorem, \(\prox_{\lambda_n f_n}^{\bm W_n}\) is Fr\'echet differentiable Lebesgue--almost everywhere on \(\R^p\); see \citet[Thm.~9.60]{SM:rockafellarWets2009}.
\end{enumerate}
\end{proof}

\subsection{Proofs of the results in Section~\ref{sec:asymptotics}}

The next lemma isolates the epi-convergence in probability and minimizer-convergence in probability arguments used in the proof of Proposition~\ref{prop:consistency}.

\begin{lemma}[Stochastic epi-convergence of sums and exact minimizers]
\label{lem:epi_sum_argmin_prob}
Let \(g_n,h_n\) be random functions on \(\R^p\), and let \(g_0,h_0\) be deterministic functions on \(\R^p\). Define
\[
G_n:=g_n+h_n,
\qquad
G_0:=g_0+h_0.
\]
Assume that
\begin{enumerate}[label=(\roman*)]
\item \(g_n\to_{\Pr} g_0\) uniformly on every compact subset of \(\R^p\), and each \(g_n\) and \(g_0\) is finite-valued and continuous;
\item \(h_n\to_{\Pr} h_0\) in epigraph, and each \(h_n\) and \(h_0\) is proper, lower semicontinuous, and convex;
\item for each \(n\), \(G_n\) is proper, lower semicontinuous, and convex, and \(\hat{\bm x}_n\) is a measurable exact minimizer of \(G_n\);
\item \(G_0\) is proper, lower semicontinuous, convex, coercive, and has a unique minimizer \(\bm x_0\).
\end{enumerate}
Then
\[
G_n\to_{\Pr} G_0
\qquad\text{in epigraph,}\qquad \hat{\bm x}_n\to_{\Pr}\bm x_0.
\]
\end{lemma}

\begin{proof}[Proof of Lemma \ref{lem:epi_sum_argmin_prob}]
Because uniform convergence on compact subsets of \(\R^p\) (local uniform convergence) and epigraph convergence can both be characterized through suitable distances on complete separable metric spaces, with the latter characterized by epi-distance \citep[Thm.~7.58]{SM:rockafellarWets2009}), we can use the usual subsequence criterion for convergence in probability. Let \(\{n_k\}\) index an arbitrary subsequence. Since \(g_n\to_{\Pr} g_0\) locally uniformly on \(\R^p\), there exists a further subsequence, again indexed by \(\{n_k\}\) for simplicity, such that
\[
g_{n_k}\to g_0
\qquad\text{almost surely locally uniformly.}
\]
From this subsequence, since \(h_n\to_{\Pr} h_0\) in epigraph, we may extract a further subsequence, not relabeled, such that
\[
h_{n_k}\to h_0
\qquad\text{almost surely in epigraph.}
\]
Now, fix a generic outcome \(\omega\) in the intersection of the two almost sure events above. Then
\[
G_{n_k}(\cdot,\omega)
=
g_{n_k}(\cdot,\omega)+h_{n_k}(\cdot,\omega)
\to
g_0+h_0
=
G_0
\qquad\text{in epigraph,}
\]
by \citet[Thm.~2.15]{SM:attouch1984variational}. Since the original subsequence was arbitrary, the subsequence characterization of convergence in probability yields
\[
G_n\to_{\Pr} G_0
\qquad\text{in epigraph.}
\]
On the same almost sure event, \(G_0\) is coercive, hence level-bounded. In addition, each \(G_{n_k}(\cdot,\omega)\) is convex, so every sublevel set is connected. In particular, \(\{G_{n_k}(\cdot,\omega)\}\) is level-bounded by \citet[Ex.~7.32(c)]{SM:rockafellarWets2009}. Because \(\hat{\bm x}_{n_k}(\omega)\in\argmin G_{n_k}(\cdot,\omega)\) and \(G_0\) has the unique minimizer \(\bm x_0\), the deterministic epi-argmin theorem \citep[Thm.~7.33]{SM:rockafellarWets2009} gives
\[
\hat{\bm x}_{n_k}(\omega)\to \bm x_0.
\]
Thus every subsequence of \(\{\hat{\bm x}_n\}\) admits a further subsequence converging almost surely, hence in probability, to \(\bm x_0\). Therefore,
\[
\hat{\bm x}_n\to_{\Pr}\bm x_0.
\]
This concludes the proof of the lemma.
\end{proof}

\begin{proof}[Proof of Proposition \ref{prop:consistency}]
Let
\[
g_n(\bm\beta):=\frac{1}{2}\norm{\hat{\bm\beta}_n^s-\bm\beta}_{\bm W_n}^2,
\qquad
g_0(\bm\beta):=\frac{1}{2}\norm{\bm\beta_0-\bm\beta}_{\bm W_0}^2,
\]
and define
\[
\mathcal P_n(\bm\beta):=g_n(\bm\beta)+\lambda_n f_n(\bm\beta),
\qquad
\mathcal P_0(\bm\beta):=g_0(\bm\beta)+h_0(\bm\beta),
\]
where
\[
h_0:=
\begin{cases}
\lambda_0 f_0, & \text{under regime \ref{prop_cons_i}},\\[3pt]
\iota_{\dom(f_0)}, & \text{under regime \ref{prop_cons_ii}}.
\end{cases}
\]
Assumptions~\ref{ass: beta} and~\ref{ass: W} imply \(g_n\to_{\Pr} g_0\) locally uniformly on \(\R^p\). Under regime~\ref{prop_cons_i}, Assumption~\ref{ass:f}\ref{epif ass} and \(\lambda_n\to\lambda_0>0\) yield
\[
\lambda_n f_n\to_{\Pr}\lambda_0 f_0=h_0
\qquad\text{in epigraph.}
\]
Conversely, under regime~\ref{prop_cons_ii} the epi-limit is given by \eqref{eq:condition:lambda_f}. Thus, in both regimes,
\[
\lambda_n f_n\to_{\Pr} h_0
\qquad\text{in epigraph.}
\]
By Assumption~\ref{ass:f}\ref{domf ass}, after deleting finitely many initial terms we may assume \(\bm\beta_0\in\dom(f_n)\) for every \(n\), which does not affect convergence in probability. Then \(\mathcal P_n\) is proper for every \(n\). Since \(g_n\) and \(g_0\) are finite-valued, continuous, strictly convex, and coercive, while \(\lambda_n f_n\) and \(h_0\) are proper lower semicontinuous convex functions, it follows that \(\mathcal P_n\) and \(\mathcal P_0\) are proper, lower semicontinuous, and convex, and that \(\mathcal P_0\) is strictly convex and coercive. Hence, it has the unique minimizer
\[
\argmin_{\bm\beta\in\R^p}\mathcal P_0(\bm\beta)
=
\prox_{h_0}^{\bm W_0}(\bm\beta_0).
\]
Moreover, by definition of the proximal map,
\[
\hat{\bm\beta}_n
=
\prox_{\lambda_n f_n}^{\bm W_n}(\hat{\bm\beta}_n^s)
\in
\argmin_{\bm\beta\in\R^p}\mathcal P_n(\bm\beta).
\]
Hence, Lemma~\ref{lem:epi_sum_argmin_prob} yields
\[
\hat{\bm\beta}_n
\to_{\Pr}
\prox_{h_0}^{\bm W_0}(\bm\beta_0).
\]
Under regime~\ref{prop_cons_i}, \(h_0=\lambda_0 f_0\), so
\[
\prox_{\lambda_n f_n}^{\bm W_n}(\hat{\bm\beta}_n^s)
\to_{\Pr}
\prox_{\lambda_0 f_0}^{\bm W_0}(\bm\beta_0).
\]
Moreover, Moreau's decomposition under the limit inner product
\(\langle\cdot,\cdot\rangle_{\bm W_0}\) yields
\[
\prox_{\lambda_0 f_0}^{\bm W_0}(\bm\beta_0)
=
\bm\beta_0-\prox_{(\lambda_0 f_0)^*}^{\bm W_0}(\bm\beta_0),
\]
where \((\lambda_0 f_0)^*\) is the convex conjugate of \(\lambda_0 f_0\)
relative to \(\langle\cdot,\cdot\rangle_{\bm W_0}\)
\citep[Thm.~14.3]{SM:bauschke2016convex}.
Under regime~\ref{prop_cons_ii}, \(h_0=\iota_{\dom(f_0)}\) and \(\bm\beta_0\in\dom(f_0)\). Hence
\[
\prox_{\lambda_n f_n}^{\bm W_n}(\hat{\bm\beta}_n^s)
\to_{\Pr}
\prox_{\iota_{\dom(f_0)}}^{\bm W_0}(\bm\beta_0)
=
P_{\dom(f_0)}^{\bm W_0}(\bm\beta_0)
=
\bm\beta_0.
\]
This concludes the proof of Proposition \ref{prop:consistency}.
\end{proof}

In the next lemma, we collect the epi-convergence in distribution and minimizer-convergence in distribution arguments used in the proof of Theorem~\ref{prop:asymptotics:alt}.

\begin{lemma}[Weak epi-convergence and minimizer convergence for random convex criteria]
\label{lem:weak_epi_stability_distr}
Let \(C(\R^p)\) denote the space of real-valued continuous functions on \(\R^p\),
equipped with the topology of uniform convergence on compact subsets, and let
\(d_l\) be the epi-distance on \(\Gamma(\R^p)\); see \citet[Thm.~7.58]{SM:rockafellarWets2009}.
For each \(n\in\mathbb N\), let \(g_n\) be a random finite-valued continuous convex function on \(\R^p\),
and let \(h_n\) be a random element of \(\Gamma(\R^p)\).
Let \(g_0\) be a random finite-valued continuous convex function on \(\R^p\),
and let \(h_0\) be a deterministic element of \(\Gamma(\R^p)\).
Define
\[
G_n:=g_n+h_n,
\qquad
G_0:=g_0+h_0.
\]
Assume that:
\begin{enumerate}[label=(\roman*)]
\item \(g_n\to_d g_0\) in \(C(\R^p)\);
\item \(h_n\to_{\Pr} h_0\) in epigraph;
\item for every \(n\), \(\hat{\bm x}_n\) is a measurable exact minimizer of \(G_n\);
\item \(G_0\) is almost surely coercive and has an almost surely unique measurable minimizer \(\hat{\bm x}_0\).
\end{enumerate}
Then \(G_n\to_d G_0\) in epigraph, and \(\hat{\bm x}_n\to_d \hat{\bm x}_0\).
\end{lemma}

\begin{proof}
Let \(d_C\) be a metric on \(C(\R^p)\) that generates the topology of uniform
convergence on compact subsets. By \citet[Thm.~7.58]{SM:rockafellarWets2009},
the epi-distance \(d_l\) metrizes epigraph convergence on \(\Gamma(\R^p)\), and
\((\Gamma(\R^p),d_l)\) is a Polish space. Since \(h_0\) is deterministic,
assumption (ii) implies \(h_n\to_d h_0\) in \((\Gamma(\R^p),d_l)\). Hence,
together with assumption (i), Slutsky's theorem yields
\[
(g_n,h_n)\to_d (g_0,h_0)
\]
in \(C(\R^p)\times \Gamma(\R^p)\).
Next, by \citet[Thm.~2.15]{SM:attouch1984variational}, the addition map
\[
(g,h)\mapsto g+h
\]
is continuous from \(C(\R^p)\times\Gamma(\R^p)\), equipped respectively with the topology of uniform convergence on compact subsets and the topology of epigraph convergence,
into \(\Gamma(\R^p)\) endowed with the topology of epigraph convergence. Therefore, by the
continuous mapping theorem,
\[
G_n=g_n+h_n \to_d g_0+h_0=G_0
\qquad\text{in epigraph.}
\]
Finally, to prove convergence of the minimizers, let \(\{n_k\}\) index an arbitrary subsequence.
Since \((\Gamma(\R^p),d_l)\) is a Polish space and \(G_{n_k}\to_d G_0\) in epigraph,
Skorokhod representation theorem \citep[Thm.~6.7]{SM:billingsley2013convergence} yields a further subsequence,
not relabeled for simplicity, and random criteria \(\tilde G_k,\tilde G_0\) defined on a common probability
space such that \(\tilde G_k\stackrel{d}{=}G_{n_k}\) for every \(k\), \(\tilde G_0\stackrel{d}{=}G_0\), and \(\tilde G_k \to \tilde G_0\) almost surely in epigraph.
Using a regular conditional distribution of \(\hat{\bm x}_{n_k}\) given \(G_{n_k}\), we may construct measurable exact minimizers \(\tilde{\bm x}_k\) of \(\tilde G_k\) such that, for every \(k\),
\((\tilde G_k,\tilde{\bm x}_k)\stackrel{d}{=}(G_{n_k},\hat{\bm x}_{n_k})\).
Let \(\tilde{\bm x}_0\) denote the unique minimizer of \(\tilde G_0\). Since \(\tilde G_0\) is almost
surely coercive, \(\{\tilde G_k\}\) is eventually level-bounded almost surely by
\citet[Ex.~7.32(c)]{SM:rockafellarWets2009}. Therefore
\citet[Thm.~7.33]{SM:rockafellarWets2009} implies \(\tilde{\bm x}_k \to \tilde{\bm x}_0\) almost surely.
Hence \(\hat{\bm x}_{n_k}\to_d \hat{\bm x}_0\) along the subsequence. Since the subsequence
\(\{n_k\}\) was arbitrary, it follows that
\(\hat{\bm x}_n\to_d \hat{\bm x}_0\).
\end{proof}

\begin{proof}[Proof of Theorem \ref{prop:asymptotics:alt}]
Let
\[
\bm\xi_n:=r_n(\hat{\bm\beta}_n^s-\bm\beta_0),
\qquad
\hat{\bm b}_n:=r_n\big(\prox_{\lambda_n f_n}^{\bm W_n}(\hat{\bm\beta}_n^s)-\bm\beta_0\big).
\]
By Assumption~\ref{ass:f}\ref{domf ass}, after deleting finitely many initial terms and relabeling if necessary,
we may assume that \(\bm\beta_0\in\dom(f_n)\) for every \(n\).

For \(\bm b\in\R^p\), define
\[
g_n(\bm b):=\frac12\norm{\bm\xi_n-\bm b}_{\bm W_n}^2,
\qquad
q_n(\bm b):=
r_n\left[
f_n\left(\bm\beta_0+\bm b/r_n\right)-f_n(\bm\beta_0)
\right],
\]
and
\[
h_n(\bm b):=(\lambda_n r_n)q_n(\bm b),
\qquad
\mathcal P_n(\bm b):=g_n(\bm b)+h_n(\bm b).
\]
Similarly, define
\[
g_0(\bm b):=\frac12\norm{\bm\eta-\bm b}_{\bm W_0}^2,
\qquad
\mathcal P_0(\bm b):=g_0(\bm b)+h_0(\bm b),
\]
where
\[
h_0(\bm b):=
\begin{cases}
\lambda_0\rho_{\bm\beta_0}(\bm b), & \text{under regime \ref{prop_distr_i}},\\[3pt]
\sigma_{N_{\dom(f_0)}(\bm\beta_0)}(\bm b), & \text{under regime \ref{prop_distr_ii}}.
\end{cases}
\]
Because
\[
\bm b\mapsto f_n\left(\bm\beta_0+\bm b/r_n\right)-f_n(\bm\beta_0)
\]
is an affine reparametrization of \(f_n\) minus a finite constant,
\(q_n\in\Gamma(\R^p)\), hence \(h_n\in\Gamma(\R^p)\).
Under regime~\ref{prop_distr_i}, assumption~\eqref{eq:local_regime} and \(q_n\in\Gamma(\R^p)\) imply
\(\rho_{\bm\beta_0}\in\Gamma(\R^p)\), so \(h_0=\lambda_0\rho_{\bm\beta_0}\in\Gamma(\R^p)\).
Under regime~\ref{prop_distr_ii}, Assumption~\ref{ass:f}\ref{domf ass} gives
\(\bm\beta_0\in\dom(f_0)\), so \(N_{\dom(f_0)}(\bm\beta_0)\) is a nonempty closed convex cone.
Therefore
\[
h_0=\sigma_{N_{\dom(f_0)}(\bm\beta_0)}\in\Gamma(\R^p)
\]
by \citet[Ex.~13.3(i)]{SM:bauschke2016convex}.
Since \(\bm W_n\) is positive definite by construction of the proximal estimator and
\(\bm W_0\) is positive definite by Assumption~\ref{ass: W},
each \(g_n\) is finite-valued, continuous, strongly convex, and coercive, and the same holds for \(g_0\) almost surely.
Hence, for each \(n\), \(\mathcal P_n\) is in \(\Gamma(\R^p)\) and is almost surely strongly convex and coercive, while \(\mathcal P_0\) is in \(\Gamma(\R^p)\) and is strongly convex and coercive.
Therefore, for each \(n\), \(\mathcal P_n\) admits an almost surely unique minimizer, and \(\mathcal P_0\) admits a unique minimizer \citep[Prop.~11.14]{SM:bauschke2016convex}.

We next characterize these minimizers. Let
\[
\bm\beta=\bm\beta_0+\bm b/r_n.
\]
Then
\begin{align*}
\frac12\norm{\hat{\bm\beta}_n^s-\bm\beta}_{\bm W_n}^2+\lambda_n f_n(\bm\beta)
&=
\frac12\norm{\hat{\bm\beta}_n^s-\bm\beta_0-\bm b/r_n}_{\bm W_n}^2
+\lambda_n f_n\left(\bm\beta_0+\bm b/r_n\right)\\
&=
\frac{1}{2r_n^2}\norm{\bm\xi_n-\bm b}_{\bm W_n}^2
+\lambda_n f_n\left(\bm\beta_0+\bm b/r_n\right)\\
&=
\frac{1}{r_n^2}\mathcal P_n(\bm b)+\lambda_n f_n(\bm\beta_0).
\end{align*}
Hence
\[
\hat{\bm b}_n
=
\argmin_{\bm b\in\R^p}\mathcal P_n(\bm b).
\]
In particular, \(\hat{\bm b}_n\) is a measurable exact minimizer of \(\mathcal P_n\).
Likewise, almost surely,
\[
\argmin_{\bm b\in\R^p}\mathcal P_0(\bm b)
=
\prox_{h_0}^{\bm W_0}(\bm\eta).
\]
Since the proximal map \(\prox_{h_0}^{\bm W_0}\) is continuous, \(\prox_{h_0}^{\bm W_0}(\bm\eta)\) is measurable.
We next verify the assumptions of Lemma~\ref{lem:weak_epi_stability_distr}.
By Assumptions~\ref{ass: W} and \ref{ass: asy distr proximal},
Slutsky's theorem gives
\[
(\bm\xi_n,\bm W_n)\to_d(\bm\eta,\bm W_0).
\]
Since the map
\[
(\bm u,\bm W)\mapsto\left[\bm b\mapsto \frac12\norm{\bm u-\bm b}_{\bm W}^2\right]
\]
is continuous from \(\R^p\times\mathbb S_{++}^n\), where \(\mathbb S_{++}^p\) denotes the set of symmetric positive definite \(p\times p\) matrices, into \(C(\R^p)\), equipped with the topology of uniform convergence on compact subsets,
the continuous mapping theorem yields
\[
g_n\to_d g_0
\qquad\text{in } C(\R^p).
\]
For the penalty term, under regime~\ref{prop_distr_i},
assumption~\eqref{eq:local_regime} together with \(\lambda_n r_n\to\lambda_0>0\) yields
\[
h_n=(\lambda_n r_n)q_n\to_{\Pr}\lambda_0\rho_{\bm\beta_0}=h_0
\qquad\text{in epigraph.}
\]
Under regime~\ref{prop_distr_ii}, assumption~\eqref{eq:local_scaled_zero_regime} gives
\[
h_n\to_{\Pr}h_0
\qquad\text{in epigraph.}
\]
Therefore all assumptions of Lemma~\ref{lem:weak_epi_stability_distr} are satisfied, and thus
\[
\mathcal P_n\to_d \mathcal P_0
\qquad\text{in epigraph,}
\qquad
\hat{\bm b}_n\to_d \prox_{h_0}^{\bm W_0}(\bm\eta).
\]
This proves \eqref{eq:asydistr_i} under regime~\ref{prop_distr_i}
and \eqref{eq:asydistr_ii} under regime~\ref{prop_distr_ii}.
It remains to derive the projection forms. Here and below, \(^{*}\) denotes Fenchel conjugacy relative to the limit inner product \(\langle\cdot,\cdot\rangle_{\bm W_0}\).
Under regime~\ref{prop_distr_i}, if \(\partial f_0(\bm\beta_0)\neq\emptyset\), then
\citet[Prop.~17.17]{SM:bauschke2016convex} gives
\[
\rho_{\bm\beta_0}
=
\sigma_{\partial f_0(\bm\beta_0)},
\qquad
\lambda_0\rho_{\bm\beta_0}
=
\sigma_{\lambda_0\partial f_0(\bm\beta_0)}.
\]
Hence
\[
(\lambda_0\rho_{\bm\beta_0})^*
=
\iota_{\lambda_0\partial f_0(\bm\beta_0)}
\]
by \citet[Ex.~13.3(i)]{SM:bauschke2016convex}.
Moreau's decomposition under \(\langle\cdot,\cdot\rangle_{\bm W_0}\) therefore yields
\[
\prox_{\lambda_0\rho_{\bm\beta_0}}^{\bm W_0}(\bm\eta)
=
\big(\Id-P_{\lambda_0\partial f_0(\bm\beta_0)}^{\bm W_0}\big)(\bm\eta),
\]
which is \eqref{eq:asym_moreau_i}.
Under regime~\ref{prop_distr_ii},
\[
h_0=\sigma_{N_{\dom(f_0)}(\bm\beta_0)},
\qquad
h_0^*=\iota_{N_{\dom(f_0)}(\bm\beta_0)},
\]
again by \citet[Ex.~13.3(i)]{SM:bauschke2016convex}.
Thus Moreau's decomposition under \(\langle\cdot,\cdot\rangle_{\bm W_0}\) yields
\[
\prox_{\sigma_{N_{\dom(f_0)}(\bm\beta_0)}}^{\bm W_0}(\bm\eta)
=
\big(\Id-P_{N_{\dom(f_0)}(\bm\beta_0)}^{\bm W_0}\big)(\bm\eta),
\]
which is \eqref{eq:asym_moreau_ii}. This concludes the proof.
\end{proof}

\begin{proof}[Proof of Proposition \ref{prop:op1 necessary}]
  Condition $\A=\hat{\A}_n$ holds if and only if $\A^c=\hat{\A}_n^c$, which implies 
  $$\hat{\beta}_{nj}^s=
  \left(\prox_{(\lambda_nf_n)^*}^{\bm W_n}(\hat{\bm\beta}_n^s)\right)_j$$
  for all $j\in\A^c$.
Therefore, if Oracle property~\ref{oracle property}\ref{op1} holds, we obtain
  \begin{align*}
    1=\limsup_{n\to\infty}\Pr\left(\hat{\A}_n=\A\right)&=
    \limsup_{n\to\infty}\Pr\left(\hat{\A}_n^c=\A^c\right)\\
    &\le\limsup_{n\to\infty}\Pr\left(r_n(\hat{\bm\beta}_n^s)_{\A^c}=
    r_n\left(\prox_{(\lambda_nf_n)^*}^{\bm W_n}(\hat{\bm\beta}_n^s)\right)_{\A^c}\right)\\
    &\le\Pr\left((\bm\eta)_{\A^c}=\left(P_{B_0}^{\bm W_0}(\bm\eta)\right)_{\A^c}\right)\ ,
\end{align*}
  where the last inequality follows from Theorem \ref{prop:asymptotics:alt} and Portmanteau's theorem. 
  This concludes the proof.
\end{proof}

\begin{proof}[Proof of Proposition~\ref{prop:op1 sufficient}]
For any \(M>0\) we define the localization event
\[
L_n(M):=
\left\{
\|\hat{\bm\beta}_n-\bm\beta_0\|_\infty\le \frac{M}{r_n}
\right\}.
\]
By condition \eqref{eq:oracle_localization}, for every \(\varepsilon>0\) there exists
\(M\) sufficiently large such that
\begin{eqnarray}
\liminf_{n\to\infty}\Pr\big(L_n(M)\big)\ge 1-\varepsilon.\label{eq: large L_M probability}
\end{eqnarray}
We next show that under the stated assumptions the event
\[
\{(\hat{\bm\beta}_n)_{\mathcal A^c}\neq \bm 0\}
\]
is asymptotically negligible on the localization event $L_n(M)$. To this end, recall that since \(\hat{\bm\beta}_n\) minimizes
the objective function
\[
\ell_n(\bm\beta)+\lambda_n f_n(\bm\beta),
\qquad
\ell_n(\bm\beta):=\frac12\|\bm\beta-\hat{\bm\beta}_n^s\|_{\bm W_n}^2,
\]
there exists
$
\bm v_n\in \lambda_n\partial^e f_n(\hat{\bm\beta}_n)
$
such that
\[
\nabla \ell_n(\hat{\bm\beta}_n)+\bm v_n=\bm 0.
\]
Specifically, since
$
\nabla \ell_n(\hat{\bm\beta}_n)
=
\bm W_n(\hat{\bm\beta}_n-\hat{\bm\beta}_n^s),
$
we have
\[
\bm v_n
=
-\bm W_n(\hat{\bm\beta}_n-\hat{\bm\beta}_n^s).
\]
Now suppose that event \(L_n(M)\) occurs and that
$
(\hat{\bm\beta}_n)_{\mathcal A^c}\neq \bm 0.
$
Then, the definition of \(\Delta_n(M)\) applies at \(\bm\beta=\hat{\bm\beta}_n\), and the following inclusion holds
\begin{eqnarray}
L_n(M)\cap   \{(\hat{\bm\beta}_n)_{\mathcal A^c}\neq \bm 0\} \subseteq 
\{ r_n\|(\bm v_n)_{\mathcal A^c}\|_1\ge \Delta_n(M)\}.
\label{eq:oracle_lower_bound_vn}
\end{eqnarray}
On the other hand, given that
\[
\hat{\bm\beta}_n-\hat{\bm\beta}_n^s
=
(\hat{\bm\beta}_n-\bm\beta_0)-(\hat{\bm\beta}_n^s-\bm\beta_0),
\]
condition \eqref{eq:oracle_localization} implies
\[
\hat{\bm\beta}_n-\hat{\bm\beta}_n^s=O_{\Pr}(1/r_n).
\]
Furthermore, since \(\bm W_n=O_{\Pr}(1)\) by Assumption~\ref{ass: W}, we also obtain
\begin{eqnarray}
r_n\|(\bm v_n)_{\mathcal A^c}\|_1
=
r_n\big\|(\bm W_n(\hat{\bm\beta}_n-\hat{\bm\beta}_n^s))_{\mathcal A^c}\big\|_1
=
O_{\Pr}(1).
\label{eq:oracle_upper_bound_vn}
\end{eqnarray}
Recalling that
under the main hypothesis of the proposition one has for every $M>0$:
\[
\Delta_n(M)\to_{\Pr}\infty,
\]
we can combine
relations \eqref{eq:oracle_lower_bound_vn} and
\eqref{eq:oracle_upper_bound_vn}
to get for any fixed $M$ as $n\to\infty$:
\[
\Pr\left(
L_n(M)\cap\{(\hat{\bm\beta}_n)_{\mathcal A^c}\neq \bm 0\}
\right)\le \Pr ( r_n\|(\bm v_n)_{\mathcal A^c}\|_1\ge \Delta_n(M))\to 0 .
\]
Therefore,
\begin{align*}
\limsup_{n\to \infty} \Pr\big((\hat{\bm\beta}_n)_{\mathcal A^c}\neq \bm 0\big)
\le&
\limsup_{n\to \infty} \Pr\big(L_n(M)^c\big)
,
\end{align*}
where we can use property \eqref{eq: large L_M probability} and choose \(M\) large enough so that
\[\limsup_{n\to \infty} \Pr(L_n(M)^c)<\epsilon.\]
Since $\epsilon$ was arbitrary, this finally yields
\[
\Pr\big((\hat{\bm\beta}_n)_{\mathcal A^c}\neq \bm 0\big)\to 0,
\]
i.e., Oracle property~\ref{oracle property}\ref{op1}. This concludes the proof.
\end{proof}

\begin{proof}[Proof of Corollary~\ref{cor:oracle_weighted_l1}]
If \(\mathcal A^c=\varnothing\), then Oracle property~\ref{oracle property}\ref{op1} trivially holds. Hence suppose \(\mathcal A^c\neq\varnothing\).
Fix \(M>0\), and recall \(\Delta_n(M)\) from \eqref{eq:Delta_oracle_selection}.
Let \(\bm\beta\in\R^p\) satisfy
\[
\|\bm\beta-\bm\beta_0\|_\infty\le M/r_n,
\qquad
\bm\beta_{\mathcal A^c}\neq \bm 0,
\]
and let \(\bm v\in\lambda_n\partial^e f_n(\bm\beta)\).
For the above extended-valued weighted \(\ell_1\) penalty, the Euclidean subdifferential is coordinatewise:
\[
u_j\in
\begin{cases}
\{w_{nj}\sign(\beta_j)\}, & w_{nj}<\infty,\ \beta_j\neq 0,\\[3pt]
[-w_{nj},w_{nj}], & w_{nj}<\infty,\ \beta_j=0,\\[3pt]
\varnothing, & w_{nj}=+\infty,\ \beta_j\neq 0,\\[3pt]
\R, & w_{nj}=+\infty,\ \beta_j=0,
\end{cases}
\qquad j=1,\dots,p.
\]
Equivalently, since \(\bm v\in\lambda_n\partial^e f_n(\bm\beta)\),
\[
v_j\in
\begin{cases}
\{\lambda_n w_{nj}\sign(\beta_j)\}, & w_{nj}<\infty,\ \beta_j\neq 0,\\[3pt]
[-\lambda_n w_{nj},\lambda_n w_{nj}], & w_{nj}<\infty,\ \beta_j=0,\\[3pt]
\varnothing, & w_{nj}=+\infty,\ \beta_j\neq 0,\\[3pt]
\R, & w_{nj}=+\infty,\ \beta_j=0.
\end{cases}
\]
Because \(\bm\beta_{\mathcal A^c}\neq \bm 0\), there exists \(j^\star\in\mathcal A^c\) such that \(\beta_{j^\star}\neq 0\).
Since \(\bm v\in\lambda_n\partial^e f_n(\bm\beta)\), the subdifferential must be nonempty at \(j^\star\), and therefore necessarily \(w_{n j^\star}<\infty\). For this coordinate,
\[
|v_{j^\star}|=\lambda_n w_{n j^\star}
\ge
\lambda_n \min_{j\in\mathcal A^c} w_{nj}.
\]
Hence
\[
\|\bm v_{\mathcal A^c}\|_1
\ge
|v_{j^\star}|
\ge
\lambda_n \min_{j\in\mathcal A^c} w_{nj},
\]
and therefore
\[
r_n\|\bm v_{\mathcal A^c}\|_1
\ge
\lambda_n r_n \min_{j\in\mathcal A^c} w_{nj}.
\]
Since this lower bound holds for every admissible pair \((\bm\beta,\bm v)\) in the definition of \(\Delta_n(M)\), and \(\inf\varnothing=+\infty\), it follows that
\[
\Delta_n(M)
\ge
\lambda_n r_n \min_{j\in\mathcal A^c} w_{nj}.
\]
By \eqref{eq:weighted_l1_oracle_condition}, the right-hand side converges in probability to \(+\infty\). Hence, for every fixed $M>0$:
\[
\Delta_n(M)\to_{\Pr}\infty
.
\]
Proposition~\ref{prop:op1 sufficient} therefore yields oracle property~\ref{oracle property}\ref{op1}.
\end{proof}

\begin{proof}[Proof of Proposition~\ref{prop:op2}]
Since regime~\ref{prop_distr_ii} of Theorem~\ref{prop:asymptotics:alt} applies, Theorem~\ref{prop:asymptotics:alt}\ref{prop_distr_ii} yields
\[
r_n(\hat{\bm\beta}_n-\bm\beta_0)
\ \to_d\
\prox_{\sigma_{N_{\dom(f_0)}(\bm\beta_0)}}^{\bm W_0}(\bm\eta),
\qquad
\bm\eta=\bm M_0^+\bm Z.
\]
By assumption,
\(N_{\dom(f_0)}(\bm\beta_0)=\mathcal S_0^\perp\).
Therefore
\[
r_n(\hat{\bm\beta}_n-\bm\beta_0)
\ \to_d\
\prox_{\sigma_{\mathcal S_0^\perp}}^{\bm W_0}(\bm\eta).
\]
Since \(\mathcal S_0^\perp\) is a linear subspace, its support function is the indicator of its polar cone. As
\((\mathcal S_0^\perp)^\circ=\mathcal S_0\), it follows that
\(\sigma_{\mathcal S_0^\perp}=\iota_{\mathcal S_0}\).
Hence
\[
\prox_{\sigma_{\mathcal S_0^\perp}}^{\bm W_0}(\bm\eta)
=
\prox_{\iota_{\mathcal S_0}}^{\bm W_0}(\bm\eta)
=
P_{\mathcal S_0}^{\bm W_0}(\bm\eta).
\]
It follows that
\[
r_n(\hat{\bm\beta}_n-\bm\beta_0)
\ \to_d\
P_{\mathcal S_0}^{\bm W_0}(\bm\eta).
\]
We next compute the active block of this weighted projection. Let \(\bm x\in\R^p\), and write
\(\bm c:=P_{\mathcal S_0}^{\bm W_0}(\bm x)\).
Then \(\bm c_{\mathcal A^c}=\bm 0\), and \(\bm c\) minimizes
\[
\frac12\|\bm u-\bm x\|_{\bm W_0}^2
\qquad\text{over }\bm u\in\mathcal S_0.
\]
The first-order condition with respect to the active block is
\[
\big(\bm W_0(\bm c-\bm x)\big)_{\mathcal A}=\bm 0.
\]
Since \(\bm c_{\mathcal A^c}=\bm 0\), this becomes
\[
(\bm W_0)_{\mathcal A}\bm c_{\mathcal A}
=
(\bm W_0\bm x)_{\mathcal A}.
\]
Because \((\bm W_0)_{\mathcal A}\) is a principal submatrix of the positive definite matrix \(\bm W_0\), it is itself positive definite and hence invertible. Therefore
\begin{equation}\label{eq:projection_active_block_oracle2}
\big(P_{\mathcal S_0}^{\bm W_0}(\bm x)\big)_{\mathcal A}
=
[(\bm W_0)_{\mathcal A}]^{-1}(\bm W_0\bm x)_{\mathcal A},
\qquad
\big(P_{\mathcal S_0}^{\bm W_0}(\bm x)\big)_{\mathcal A^c}
=
\bm 0.
\end{equation}
Applying \eqref{eq:projection_active_block_oracle2} with \(\bm x=\bm\eta=\bm M_0^+\bm Z\), we obtain
\[
\big(P_{\mathcal S_0}^{\bm W_0}(\bm\eta)\big)_{\mathcal A}
=
[(\bm W_0)_{\mathcal A}]^{-1}
(\bm W_0\bm M_0^+\bm Z)_{\mathcal A}.
\]
Hence \(\big(P_{\mathcal S_0}^{\bm W_0}(\bm\eta)\big)_{\mathcal A}\) is a centered Gaussian vector with covariance matrix
\begin{equation}\label{eq:cov_projected_active_oracle2}
[(\bm W_0)_{\mathcal A}]^{-1}
(\bm W_0\bm M_0^+\bm\Omega_0\bm M_0^{+\prime}\bm W_0)_{\mathcal A}
[(\bm W_0)_{\mathcal A}]^{-1}.
\end{equation}
On the other hand, by assumption,
\[
[(\bm M_0)_{\mathcal A}]^{+}(\bm Z)_{\mathcal A}
\]
is the efficient Gaussian limit for the active block. Since \((\bm Z)_{\mathcal A}\sim N(\bm 0,(\bm\Omega_0)_{\mathcal A})\), this vector is centered Gaussian with covariance matrix
\begin{equation}\label{eq:cov_oracle_active_oracle2}
[(\bm M_0)_{\mathcal A}]^{+}(\bm\Omega_0)_{\mathcal A}[(\bm M_0)_{\mathcal A}]^{+\prime}.
\end{equation}
By Definition~\ref{oracle property}\ref{op2}, Oracle property~\ref{oracle property}\ref{op2} holds if and only if
\[
\big(P_{\mathcal S_0}^{\bm W_0}(\bm\eta)\big)_{\mathcal A}
\ \stackrel{d}{=}\
[(\bm M_0)_{\mathcal A}]^{+}(\bm Z)_{\mathcal A}.
\]
Since both random vectors are centered Gaussian, this is equivalent to equality of the covariance matrices in
\eqref{eq:cov_projected_active_oracle2}
and
\eqref{eq:cov_oracle_active_oracle2}. That is, Oracle property~\ref{oracle property}\ref{op2} holds if and only if
\[
[(\bm W_0)_{\mathcal A}]^{-1}
(\bm W_0\bm M_0^+\bm\Omega_0\bm M_0^{+\prime}\bm W_0)_{\mathcal A}
[(\bm W_0)_{\mathcal A}]^{-1}
=
[(\bm M_0)_{\mathcal A}]^{+}(\bm\Omega_0)_{\mathcal A}[(\bm M_0)_{\mathcal A}]^{+\prime}.
\]
Multiplying on the left and on the right by \((\bm W_0)_{\mathcal A}\) yields the stated matrix identity. This proves the claim.
\end{proof}

\subsection{Proofs of the results in Section~\ref{sec:irregular_designs}}

\begin{proof}[Proof of Proposition~\ref{prop:B0_nonempty}]
Let
\[
a_n:=\E[\|\bm Y\|_2^2/n],
\qquad
\bm\delta_{0n}:=\E[\bm X'\bm Y/n],
\qquad
\bm Q_{0n}:=\E[\bm X'\bm X/n].
\]
By assumption, \(a_n \to a_0:=R_0(\bm 0)\in\R\), \(\bm\delta_{0n}\to \bm\delta_0\) and \(\bm Q_{0n}\to \bm Q_0\).
We first show that \(\bm\delta_0\in\Range(\bm Q_0)\).
Let \(\bm v\in\Kernel(\bm Q_0)\) be arbitrary, and for \(t\in\R\) define
\[
r_n(t):=\E\left[\|\bm Y-\bm X(t\bm v)\|_2^2/n\right].
\]
Expanding the square gives
\[
r_n(t)
=
a_n
-2t\bm v'\bm\delta_{0n}
+t^2\bm v'\bm Q_{0n}\bm v.
\]
Since \(\bm v\in\Kernel(\bm Q_0)\), we have \(\bm v'\bm Q_{0n}\bm v \to \bm v'\bm Q_0\bm v=0\).
Hence, for every fixed \(t\in\R\),
\[
R_0(t\bm v)
=
\lim_{n\to\infty} r_n(t)
=
a_0-2t\bm v'\bm\delta_0.
\]
Because each \(r_n(t)\ge 0\), also \(R_0(t\bm v)\ge 0\) for every \(t\in\R\). The affine function
\[
t\mapsto a_0-2t\bm v'\bm\delta_0
\]
can be nonnegative for all \(t\in\R\) only if \(\bm v'\bm\delta_0=0\). Since \(\bm v\in\Kernel(\bm Q_0)\) was arbitrary, it follows that \(\bm\delta_0\perp\Kernel(\bm Q_0)\).
As \(\bm Q_0\) is symmetric, \(\Kernel(\bm Q_0)^\perp=\Range(\bm Q_0)\),
and therefore \(\bm\delta_0\in\Range(\bm Q_0)\).
Now let \(\bar{\bm\beta}\) satisfy \(\bm Q_0\bar{\bm\beta}=\bm\delta_0\);
for instance, one may take \(\bar{\bm\beta}=\bm Q_0^+\bm\delta_0\).
For any \(\bm\beta\in\R^p\),
\[
R_0(\bm\beta)
=
\lim_{n\to\infty}
\E[\|\bm Y-\bm X\bm\beta\|_2^2/n]
=
a_0-2\bm\beta'\bm\delta_0+\bm\beta'\bm Q_0\bm\beta.
\]
Therefore,
\begin{align*}
R_0(\bm\beta)-R_0(\bar{\bm\beta})
&=
-2(\bm\beta-\bar{\bm\beta})'\bm\delta_0
+\bm\beta'\bm Q_0\bm\beta-(\bar{\bm\beta})'\bm Q_0\bar{\bm\beta} \\
&=
-2(\bm\beta-\bar{\bm\beta})'\bm Q_0\bar{\bm\beta}
+\bm\beta'\bm Q_0\bm\beta-(\bar{\bm\beta})'\bm Q_0\bar{\bm\beta} \\
&=
(\bm\beta-\bar{\bm\beta})'\bm Q_0(\bm\beta-\bar{\bm\beta})
\ge0,
\end{align*}
since \(\bm Q_0\) is positive semidefinite.
Thus \(\bar{\bm\beta}\) is a minimizer and \(\mathcal B_0\neq\emptyset\). Moreover,
\[
R_0(\bm\beta)=R_0(\bar{\bm\beta})
\quad\Longleftrightarrow\quad
(\bm\beta-\bar{\bm\beta})'\bm Q_0(\bm\beta-\bar{\bm\beta})=0
\quad\Longleftrightarrow\quad
\bm Q_0(\bm\beta-\bar{\bm\beta})=\bm 0,
\]
where the last equivalence uses symmetry and positive semidefiniteness of \(\bm Q_0\).
Hence
\[
\mathcal B_0
=
\{\bm\beta\in\R^p\ :\ \bm Q_0\bm\beta=\bm\delta_0\}.
\]
This proves the claim.
\end{proof}

\begin{proof}[Proof of Proposition~\ref{prop:random_MP_continuity}]
The equivalence
\[
\bm Q_n^+\to_{\Pr}\bm Q_0^+
\qquad\Longleftrightarrow\qquad
\Pr\big(\Rank(\bm Q_n)=\Rank(\bm Q_0)\big)\to 1
\]
follows directly from \citet[Thm.~2(b)]{SM:andrews1987asymptotic}, applied with
\(H_n=\bm Q_n\) and \(H=\bm Q_0\), since Assumption~\ref{ass:limit_design_irregular}
gives \(\bm Q_n\to_{\Pr}\bm Q_0\).

To prove the two design-specific claims, write \(r_0:=\Rank(\bm Q_0)\).
\begin{enumerate}[label=(\roman*)]
\item
Suppose the design is singular. Then
\(\Range(\bm Q_{0n})=\Range(\bm Q_0)\) for all sufficiently large \(n\).
Since \(\bm Q_{0n}\) and \(\bm Q_0\) are symmetric positive semidefinite, this is equivalent to
\(\Kernel(\bm Q_{0n})=\Kernel(\bm Q_0)\) for all sufficiently large \(n\).
Let \(\bm v_1,\dots,\bm v_{p-r_0}\) be a basis of \(\Kernel(\bm Q_0)\).
Then, for each \(j=1,\dots,p-r_0\) and all sufficiently large \(n\),
\[
0=\bm v_j'\bm Q_{0n}\bm v_j
=\E[\bm v_j'\bm Q_n\bm v_j].
\]
Because \(\bm Q_n\) is symmetric positive semidefinite almost surely,
\(\bm v_j'\bm Q_n\bm v_j\ge 0\) almost surely. Hence
\(\bm v_j'\bm Q_n\bm v_j=0\) almost surely, and therefore
\(\bm Q_n\bm v_j=\bm 0\) almost surely.
Since there are finitely many basis vectors, intersecting the corresponding probability-one events yields
\[
\Kernel(\bm Q_0)\subseteq \Kernel(\bm Q_n)
\qquad\text{almost surely for all sufficiently large }n.
\]
Hence \(\Rank(\bm Q_n)\le r_0\) almost surely for all sufficiently large \(n\).
If \(r_0=0\), this already implies \(\Pr\big(\Rank(\bm Q_n)=r_0\big)\to 1\).
Suppose now that \(r_0\ge 1\). Since \(\bm Q_n\to_{\Pr}\bm Q_0\), Weyl's inequality yields
\[
\lambda_{r_0}(\bm Q_n)\to_{\Pr}\lambda_{r_0}(\bm Q_0)>0.
\]
Therefore \(\Pr\big(\Rank(\bm Q_n)\ge r_0\big)\to 1\).
Combined with \(\Rank(\bm Q_n)\le r_0\) almost surely eventually, this gives \(\Pr\big(\Rank(\bm Q_n)=r_0\big)\to 1\).
The general equivalence established above therefore implies
\(\bm Q_n^+\to_{\Pr}\bm Q_0^+\).
\item
Suppose the design is nearly-singular. Then
\(\Range(\bm Q_{0n})\supsetneq \Range(\bm Q_0)\) for all sufficiently large \(n\),
hence \(\Rank(\bm Q_{0n})>r_0\) for all sufficiently large \(n\).
If, in addition,
\(\Pr\big(\Rank(\bm Q_n)=\Rank(\bm Q_{0n})\big)\to 1\),
then
\[
\Pr\big(\Rank(\bm Q_n)>r_0\big)\to 1.
\]
In particular, \(\Pr\big(\Rank(\bm Q_n)=r_0\big)\to 0\).
The general equivalence above therefore implies
\(\bm Q_n^+\not\to_{\Pr}\bm Q_0^+\).
\end{enumerate}
\end{proof}

\begin{proof}[Proof of Proposition~\ref{pro:cons_rank_Qn_check}]
Write \(r_0:=\Rank(\bm Q_0)\),
and let \(\sigma_{01}\ge \cdots \ge \sigma_{0p}\ge 0\)
denote the ordered eigenvalues of \(\bm Q_0\). Define \(\delta_n:=\|\bm Q_n-\bm Q_0\|_{\mathrm{op}}\).
Because \(p\) is fixed, there exists a constant \(C<\infty\) such that, for every symmetric
\(p\times p\) matrix \(\bm A\),
\[
\|\bm A\|_{\mathrm{op}}\le C\|\vech(\bm A)\|_2.
\]
Hence Assumption~\ref{ass:bounded_Qn_op} implies
\[
\|\bm Q_n-\bm Q_{0n}\|_{\mathrm{op}}=O_{\Pr}(n^{-1/2}).
\]
Moreover, Assumption~\ref{ass: near singularity} gives \(\bm Q_{0n}-\bm Q_0=\tau_n^{-1}\bm\Delta\),
so
\[
\delta_n
\le
\|\bm Q_n-\bm Q_{0n}\|_{\mathrm{op}}
+
\|\bm Q_{0n}-\bm Q_0\|_{\mathrm{op}}
=
O_{\Pr}(n^{-1/2})+O(\tau_n^{-1}).
\]
Therefore,
\[
\frac{\delta_n}{\nu_n}
=
O_{\Pr}\left(\frac{1}{\sqrt n\nu_n}\right)
+
O\left(\frac{1}{\tau_n\nu_n}\right)
\to_{\Pr}0
\]
by \eqref{eq:nu_rate_conditions}. In particular, \(\delta_n=o_{\Pr}(\nu_n)\) and \(\delta_n\to_{\Pr}0\).
Next, by Weyl's inequality,
\[
\max_{1\le j\le p}|\sigma_{nj}-\sigma_{0j}|
\le
\|\bm Q_n-\bm Q_0\|_{\mathrm{op}}
=
\delta_n.
\]
We now prove rank recovery.
If \(r_0\ge 1\), let \(\underline{\sigma}_0:=\sigma_{0r_0}>0\). Since \(\nu_n\to 0\), for all
sufficiently large \(n\), \(\nu_n<\underline{\sigma}_0/2\).
On the event \(\{\delta_n<\underline{\sigma}_0/2\}\),
\[
\min_{1\le j\le r_0}\sigma_{nj}
\ge
\min_{1\le j\le r_0}\sigma_{0j}-\delta_n
\ge
\underline{\sigma}_0-\delta_n
>
\underline{\sigma}_0/2
>
\nu_n.
\]
Hence
\[
\Pr\Big(\min_{1\le j\le r_0}\sigma_{nj}>\nu_n\Big)\to 1.
\]
If \(r_0=0\),
the preceding statement becomes vacuous.
Consequently, all nonzero eigenvalues of $\bm Q_0$ are asymptotically recovered by the thresholding procedure.

We now establish that the zero eigenvalues are also consistently identified.
If \(r_0<p\), then for every \(j>r_0\) we have \(\sigma_{0j}=0\), and therefore \(0\le \max_{j>r_0}\sigma_{nj}\le \delta_n\).
Since \(\delta_n/\nu_n\to_{\Pr}0\),
\[
\Pr\Big(\max_{j>r_0}\sigma_{nj}\le \nu_n\Big)\to 1.
\]
If \(r_0=p\), this statement is vacuous.
Combining the two displays above, we therefore obtain
\[
\Pr\big(\sigma_{nj}>\nu_n \text{ for all } j\le r_0,\ \sigma_{nj}\le \nu_n \text{ for all } j>r_0\big)\to 1.
\]
On this event, \(\check{\bm Q}_n\) keeps exactly the first \(r_0\) eigenvalues of \(\bm Q_n\)
and sets the remaining \(p-r_0\) eigenvalues to zero. Therefore,
\[
\Pr\big(\Rank(\check{\bm Q}_n)=\Rank(\bm Q_0)\big)\to 1.
\]
We next prove consistency of \(\check{\bm Q}_n\). Since \(\check{\bm Q}_n\) and \(\bm Q_n\) have
the same eigenvectors,
\[
\bm Q_n-\check{\bm Q}_n
=
\bm E_n
\diag\big(
\sigma_{n1}\mathbf 1\{\sigma_{n1}\le \nu_n\},
\dots,
\sigma_{np}\mathbf 1\{\sigma_{np}\le \nu_n\}
\big)
\bm E_n'.
\]
Hence
\[
\|\bm Q_n-\check{\bm Q}_n\|_{\mathrm{op}}
=
\max_{1\le j\le p}\sigma_{nj}\mathbf 1\{\sigma_{nj}\le \nu_n\}
\le
\nu_n.
\]
It follows that
\[
\|\check{\bm Q}_n-\bm Q_0\|_{\mathrm{op}}
\le
\|\check{\bm Q}_n-\bm Q_n\|_{\mathrm{op}}
+
\|\bm Q_n-\bm Q_0\|_{\mathrm{op}}
\le
\nu_n+\delta_n
\to_{\Pr}0.
\]
Thus \(\check{\bm Q}_n\to_{\Pr}\bm Q_0\).
Finally, Proposition~\ref{prop:random_MP_continuity}, applied to the sequence
\(\check{\bm Q}_n\), yields \(\check{\bm Q}_n^+\to_{\Pr}\bm Q_0^+\),
because we have just shown both
\[
\check{\bm Q}_n\to_{\Pr}\bm Q_0
\qquad\text{and}\qquad
\Pr\big(\Rank(\check{\bm Q}_n)=\Rank(\bm Q_0)\big)\to 1.
\]
This completes the proof.
\end{proof}

\begin{proof}[Proof of Proposition~\ref{prop:consistency_modified_Ridgeless}]
By Definition~\ref{def:modified_Ridgeless}, with the displayed equality justified by Proposition~\ref{prop:modified_Ridgeless_repr},
\[
\check{\bm\beta}_n^{\mathrm{rls}}
=
\check{\bm Q}_n^+\bm X'\bm Y/n
=
\check{\bm Q}_n^+\bm Q_n\bm\beta_0
+
\check{\bm Q}_n^+\bm X'\bm\varepsilon/n,
\]
where the second equality follows from \(\bm Y=\bm X\bm\beta_0+\bm\varepsilon\) and
\(\bm Q_n=\bm X'\bm X/n\). Proposition~\ref{pro:cons_rank_Qn_check}
gives \(\check{\bm Q}_n\to_{\Pr}\bm Q_0\) and rank consistency, hence
\(\check{\bm Q}_n^+\to_{\Pr}\bm Q_0^+\). Also,
Assumptions~\ref{ass: near singularity} and~\ref{ass:bounded_Qn_op} imply
\(\bm Q_n\to_{\Pr}\bm Q_0\). Since \(\bm X'\bm\varepsilon/n\to_{\Pr}\bm0\) by assumption, it follows that
\[
\check{\bm\beta}_n^{\mathrm{rls}}
=
\bm Q_0^+\bm Q_0\bm\beta_0+o_{\Pr}(1).
\]
Finally, Assumption~\ref{ass:orthogonality_beta0} and the model imply
\(\bm\delta_0=\bm Q_0\bm\beta_0\). Therefore, by Definition~\ref{def:Ridgeless_estimand},
\[
\bm Q_0^+\bm Q_0\bm\beta_0
=
\bm Q_0^+\bm\delta_0
=
\bm\beta_0^{\mathrm{rls}}.
\]
This proves the claim.
\end{proof}

We next prove Theorem~\ref{thm:asy_modified_Ridgeless}. To do so, we make use of the following auxiliary result.

\begin{lemma}\label{lemma: equivalent projection difference formula}
Let \(\bm Q_0\) and \(\bm Q_1\) be symmetric positive semidefinite \(p\times p\) matrices with
\(k:=\Rank(\bm Q_0)\le \Rank(\bm Q_1)\).
Let \(\bm P_0\) and \(\bm P_0^\perp\) denote the orthogonal projectors onto
\(\Kernel(\bm Q_0)\) and \(\Range(\bm Q_0)\), respectively.
Further, let
\(\bm Q_1=\sum_{i=1}^p \lambda_i \bm V_i\bm V_i'\) with
\(0\le \lambda_1\le \cdots \le \lambda_p\),
be a spectral decomposition of \(\bm Q_1\), and define
\[
\bm P_1:=\sum_{i=1}^{p-k}\bm V_i\bm V_i',
\qquad
\bm P_1^\perp:=\sum_{i=p-k+1}^{p}\bm V_i\bm V_i'.
\]
Thus \(\bm P_1\) is the orthogonal projector onto the span of the \(p-k\) eigenvectors associated with the \(p-k\) smallest eigenvalues of \(\bm Q_1\), and \(\bm P_1^\perp\) is the orthogonal projector onto the span of the \(k\) eigenvectors associated with the \(k\) largest eigenvalues.

Then
\begin{align*}
\bm P_1-\bm P_0
&=
\bm Q_0^+ \sum_{i=1}^{p-k} \lambda_i \bm V_i\bm V_i'
-\bm Q_0^+(\bm Q_1-\bm Q_0)\bm P_1
-\bm P_0\sum_{i=p-k+1}^{p}\lambda_i \bm V_i\bm V_i' \bm Q_1^+.
\end{align*}
If, in addition,
\(\Rank(\bm Q_0)=\Rank(\bm Q_1)\),
then
\begin{align*}
\bm P_1-\bm P_0
=
-\bm Q_0^+(\bm Q_1-\bm Q_0)\bm P_1
-\bm P_0(\bm Q_1-\bm Q_0)\bm Q_1^+.
\end{align*}
\end{lemma}

\begin{proof}[Proof of Lemma \ref{lemma: equivalent projection difference formula}]
Since
\(\bm P_1=\sum_{i=1}^{p-k}\bm V_i\bm V_i'\) and 
\(\bm P_1^\perp=\sum_{i=p-k+1}^{p}\bm V_i\bm V_i'\),
we have
\[
\bm Q_1\bm P_1=\sum_{i=1}^{p-k}\lambda_i\bm V_i\bm V_i'.
\]
Moreover, \(\bm Q_0\bm P_0=\bm 0\), because \(\bm P_0\) projects onto \(\Kernel(\bm Q_0)\). Therefore,
\begin{align*}
\bm Q_0(\bm P_1-\bm P_0)
&=\bm Q_0\bm P_1 =\bm Q_1\bm P_1-(\bm Q_1-\bm Q_0)\bm P_1 \\
&=\sum_{i=1}^{p-k}\lambda_i\bm V_i\bm V_i'
-(\bm Q_1-\bm Q_0)\bm P_1.
\end{align*}
Multiplying by the Moore--Penrose inverse \(\bm Q_0^+\) on the left and using
\(\bm Q_0^+\bm Q_0=\bm P_0^\perp\), we obtain
\begin{align*}
\bm P_0^\perp(\bm P_1-\bm P_0)
=
\bm Q_0^+\sum_{i=1}^{p-k}\lambda_i\bm V_i\bm V_i'
-\bm Q_0^+(\bm Q_1-\bm Q_0)\bm P_1.
\end{align*}
On the other hand,
\begin{align*}
\bm P_0(\bm P_1-\bm P_0)
=\bm P_0\bm P_1-\bm P_0 =-\bm P_0(\bm I_p-\bm P_1) =-\bm P_0\bm P_1^\perp.
\end{align*}
Since \(\Rank(\bm Q_1)\ge k\), the \(k\) largest eigenvalues
\(\lambda_{p-k+1},\ldots,\lambda_p\) are strictly positive. Hence
\[
\bm Q_1^+=\sum_{i:\lambda_i>0}\lambda_i^{-1}\bm V_i\bm V_i',
\qquad
\bm P_1^\perp
=\sum_{i=p-k+1}^{p}\bm V_i\bm V_i'
=\sum_{i=p-k+1}^{p}\lambda_i \bm V_i\bm V_i' \bm Q_1^+.
\]
It follows that
\begin{align*}
\bm P_0(\bm P_1-\bm P_0)
=
-\bm P_0\sum_{i=p-k+1}^{p}\lambda_i \bm V_i\bm V_i' \bm Q_1^+.
\end{align*}
Adding the \(\bm P_0^\perp\)- and \(\bm P_0\)-components yields
\begin{align*}
\bm P_1-\bm P_0
&=
\bm P_0^\perp(\bm P_1-\bm P_0)+\bm P_0(\bm P_1-\bm P_0) \\
&=
\bm Q_0^+\sum_{i=1}^{p-k}\lambda_i\bm V_i\bm V_i'
-\bm Q_0^+(\bm Q_1-\bm Q_0)\bm P_1
-\bm P_0\sum_{i=p-k+1}^{p}\lambda_i\bm V_i\bm V_i'\bm Q_1^+.
\end{align*}
If, in addition, \(\Rank(\bm Q_1)=\Rank(\bm Q_0)=k\), then
\[
\lambda_1=\cdots=\lambda_{p-k}=0,
\qquad
\sum_{i=p-k+1}^{p}\lambda_i\bm V_i\bm V_i'=\bm Q_1.
\]
Hence the previous display simplifies to
\begin{align*}
\bm P_1-\bm P_0
&=
-\bm Q_0^+(\bm Q_1-\bm Q_0)\bm P_1
-\bm P_0\bm Q_1\bm Q_1^+ \\
&=
-\bm Q_0^+(\bm Q_1-\bm Q_0)\bm P_1
-\bm P_0(\bm Q_1-\bm Q_0)\bm Q_1^+,
\end{align*}
because \(\bm P_0\bm Q_0=\bm 0\). This proves the claim.
\end{proof}

\begin{proof}[Proof of Theorem~\ref{thm:asy_modified_Ridgeless}]
From Definition~\eqref{eq:delta0_def}, it follows
\[
\bm\delta_{0n}
=
\E[\bm X'\bm Y/n]
=
\bm Q_{0n}\bm\beta_0+\E[\bm X'\bm\varepsilon/n].
\]
Hence Assumptions~\ref{ass: near singularity} and~\ref{ass:orthogonality_beta0} imply \(\bm\delta_0=\bm Q_0\bm\beta_0\)
and \(\bm\beta_0\in\mathcal B_0\) from Definition~\eqref{eq:B0_def}. Therefore, in
Definition~\ref{def:Ridgeless_estimand} we obtain:
\[
\bm\beta_0^{\mathrm{rls}}
=
\bm Q_0^+\bm\delta_0
=
\bm Q_0^+\bm Q_0\bm\beta_0
=
\bm P_0^\perp\bm\beta_0,
\qquad
\bm\beta_0
=
\bm\beta_0^{\mathrm{rls}}+\bm P_0\bm\beta_0.\]
From Definition~\ref{def:modified_Ridgeless} and
the identity \(\bm X'\bm Y/n=\bm Q_n\bm\beta_0+\bm X'\bm\varepsilon/n\), we next obtain
we obtain
\begin{align*}
\sqrt n\big(\check{\bm\beta}_n^{\mathrm{rls}}-\bm\beta_0^{\mathrm{rls}}\big)
&=
\sqrt n\Bigl(\check{\bm Q}_n^+\bm Q_n\bm\beta_0-\bm\beta_0^{\mathrm{rls}}\Bigr)
+
\check{\bm Q}_n^+\bm X'\bm\varepsilon/\sqrt n
\nonumber\\
&=
\sqrt n\Bigl(\check{\bm Q}_n^+\bm Q_n\bm\beta_0^{\mathrm{rls}}-\bm\beta_0^{\mathrm{rls}}\Bigr)
+
\sqrt n\check{\bm Q}_n^+\bm Q_n\bm P_0\bm\beta_0
+
\check{\bm Q}_n^+\bm X'\bm\varepsilon/\sqrt n.
\label{eq:decomp_modrls_new}
\end{align*}
Recalling that \(\check{\bm Q}_n\) is obtained by hard-thresholding the eigenvalues
of \(\bm Q_n\), while retaining the same eigenvectors,
following identities also hold:
\[
\check{\bm Q}_n^+(\bm Q_n-\check{\bm Q}_n)=\bm 0,
\qquad
\check{\bm Q}_n^+\bm Q_n
=
\check{\bm Q}_n^+\check{\bm Q}_n
=
\bm I_p-\check{\bm P}_n
=
\check{\bm P}_n^\perp.
\]
Hence, we can write:
\[
\sqrt n\big(\check{\bm\beta}_n^{\mathrm{rls}}-\bm\beta_0^{\mathrm{rls}}\big)
=
T_{1n}+T_{2n}+T_{3n},
\]
where
\[
T_{1n}:=-\sqrt n\check{\bm P}_n\bm\beta_0^{\mathrm{rls}},
\qquad
T_{2n}:=\sqrt n\check{\bm Q}_n^+\bm Q_n\bm P_0\bm\beta_0,
\qquad
T_{3n}:=\check{\bm Q}_n^+\bm X'\bm\varepsilon/\sqrt n.
\]
We next prove the convergence in distribution statements of the theorem under the various designs considered. 

\begin{enumerate}[label=(\roman*)]

\item Under a singular design, the proof of Proposition~\ref{prop:random_MP_continuity}(i),
\(\Kernel(\bm Q_0)\subseteq \Kernel(\bm Q_n)\), and thus equivalently \(\Range(\bm Q_n)\subseteq \Range(\bm Q_0)\), almost surely for all sufficiently large \(n\).
Since \(\check{\bm Q}_n\) is obtained from \(\bm Q_n\) by thresholding
eigenvalues while retaining the same eigenvectors, it also follows that \(\Range(\check{\bm Q}_n)\subseteq \Range(\bm Q_0)\) almost surely for all sufficiently large \(n\).
Furthermore, Proposition~\ref{pro:cons_rank_Qn_check} gives:
\[
\Pr\big(\Rank(\check{\bm Q}_n)=\Rank(\bm Q_0)\big)\to 1.
\]
Therefore, following identities hold with probability tending to one,
\[
\Range(\check{\bm Q}_n)=\Range(\bm Q_0),
\qquad
\check{\bm P}_n=\bm P_0,
\qquad
\check{\bm P}_n^\perp=\bm P_0^\perp.
\]
recalling that \(\bm\beta_0^{\mathrm{rls}}\in\Range(\bm Q_0)\) and
\(\bm P_0\bm\beta_0\in\Kernel(\bm Q_0)\), this implies:
\[
T_{1n}=0,
\qquad
T_{2n}=0 ,
\]
with probability tending to one.
Moreover, Proposition~\ref{pro:cons_rank_Qn_check} gives
\(\check{\bm Q}_n^+\to_{\Pr}\bm Q_0^+\), while
Assumption~\ref{ass: joint_clt_Q_score} yields
\(\bm X'\bm\varepsilon/\sqrt n\to_d \bm Z\).
Hence, by Slutsky's lemma,
\[
T_{3n}\to_d \bm Q_0^+\bm Z ,
\]
which proves part~(i) of the proposition.

\item We next suppose now that the design is nearly-singular with
\(\sqrt n/\tau_n\to c\in[0,\infty)\).
We first analyze \(T_{1n}\) and
let
\[
\bm D_n:=\check{\bm Q}_n-\bm Q_0.
\]
As in the proof of Proposition~\ref{pro:cons_rank_Qn_check},
it follows
\begin{equation}\label{eq:rate_checkQ_theorem_new}
\|\bm D_n\|_{\op}
=
O_{\Pr}\bigl(1/\min\{\sqrt n,\tau_n\}\bigr).
\end{equation}
In particular, \(\sqrt n\|\bm D_n\|_{\op}=O_{\Pr}(1)\).
On the event
\(\{\Rank(\check{\bm Q}_n)=\Rank(\bm Q_0)\}\),
which has probability tending to one by Proposition~\ref{pro:cons_rank_Qn_check},
the equal-rank part of
Lemma~\ref{lemma: equivalent projection difference formula} gives
\[
\check{\bm P}_n-\bm P_0
=
-\bm Q_0^+\bm D_n\check{\bm P}_n
-\bm P_0\bm D_n\check{\bm Q}_n^+.
\]
Hence
\[
\|\check{\bm P}_n-\bm P_0\|_{\op}
=
O_{\Pr}\bigl(\|\bm D_n\|_{\op}\bigr),
\]
using the fact that \(\check{\bm Q}_n^+\to_{\Pr}\bm Q_0^+\).
Moreover, since \(\bm P_0\bm\beta_0^{\mathrm{rls}}=\bm 0\),
\[
\check{\bm P}_n\bm\beta_0^{\mathrm{rls}}
=
(\check{\bm P}_n-\bm P_0)\bm\beta_0^{\mathrm{rls}}.
\]
Therefore, from these two last properties we can write:
\begin{align*}
(\check{\bm P}_n-\bm P_0)\bm\beta_0^{\mathrm{rls}}
&=
-\bm Q_0^+\bm D_n(\check{\bm P}_n-\bm P_0)\bm\beta_0^{\mathrm{rls}}
-\bm P_0\bm D_n\check{\bm Q}_n^+\bm\beta_0^{\mathrm{rls}} \\
&=
-\bm P_0\bm D_n\bm Q_0^+\bm\beta_0^{\mathrm{rls}}
+
o_{\Pr}\bigl(\|\bm D_n\|_{\op}\bigr),
\end{align*}
using in the last equality the fact that
\(\check{\bm Q}_n^+\to_{\Pr}\bm Q_0^+\) and
\[
\bm Q_0^+\bm D_n(\check{\bm P}_n-\bm P_0)\bm\beta_0^{\mathrm{rls}}
=
O_{\Pr}\bigl(\|\bm D_n\|_{\op}^2\bigr)
=
o_{\Pr}\bigl(\|\bm D_n\|_{\op}\bigr).
\]

Multiplying by \(\sqrt n\) and using \eqref{eq:rate_checkQ_theorem_new}, we obtain
\begin{equation}\label{eq:T1_expansion_new}
T_{1n}
=
\sqrt n\bm P_0(\check{\bm Q}_n-\bm Q_0)\bm Q_0^+\bm\beta_0^{\mathrm{rls}}
+
o_{\Pr}(1).
\end{equation}
Next, we decompose
\[
\check{\bm Q}_n-\bm Q_0
=
(\bm Q_n-\bm Q_{0n})
+
(\bm Q_{0n}-\bm Q_0)
+
(\check{\bm Q}_n-\bm Q_n).
\]
Under Assumption~\ref{ass: near singularity}, \(\bm Q_{0n}-\bm Q_0=\tau_n^{-1}\bm\Delta\),
so
\[
\sqrt n\bm P_0(\bm Q_{0n}-\bm Q_0)\bm Q_0^+\bm\beta_0^{\mathrm{rls}}
=
\frac{\sqrt n}{\tau_n}
\bm P_0\bm\Delta\bm Q_0^+\bm\beta_0^{\mathrm{rls}}
\to
c\bm P_0\bm\Delta\bm Q_0^+\bm\beta_0^{\mathrm{rls}}.
\]
Assumption~\ref{ass: joint_clt_Q_score} further implies
\[
\sqrt n\bm P_0(\bm Q_n-\bm Q_{0n})\bm Q_0^+\bm\beta_0^{\mathrm{rls}}
\to_d
\bm P_0\bm\Theta\bm Q_0^+\bm\beta_0^{\mathrm{rls}}.
\]
Finally, since
\[
\check{\bm Q}_n-\bm Q_n=(\check{\bm Q}_n-\bm Q_n)\check{\bm P}_n ,
\]
and \(\bm Q_0^+\bm\beta_0^{\mathrm{rls}}\in\Range(\bm Q_0)=\Range(\bm P_0^\perp)\), we also have:
\begin{align*}
\sqrt n\
\bigl\|
\bm P_0(\check{\bm Q}_n-\bm Q_n)\bm Q_0^+\bm\beta_0^{\mathrm{rls}}
\bigr\|_2
&\le
\sqrt n\|\check{\bm Q}_n-\bm Q_n\|_{\op}
\|\check{\bm P}_n\bm P_0^\perp\|_{\op}
\|\bm Q_0^+\bm\beta_0^{\mathrm{rls}}\|_2 \\
&=
o_{\Pr}(1).
\end{align*}
Recalling that \(\|\check{\bm Q}_n-\bm Q_n\|_{\op}
=
O_{\Pr}\bigl(1/\min\{\sqrt n,\tau_n\}\bigr)\), we also have
\[
\|\check{\bm P}_n\bm P_0^\perp\|_{\op}
=
\|(\check{\bm P}_n-\bm P_0)\bm P_0^\perp\|_{\op}
=
O_{\Pr}\bigl(1/\min\{\sqrt n,\tau_n\}\bigr),
\]
which yields, since \(\sqrt n/\tau_n=O(1)\):
\[
T_{1n}
\to_d
\bm P_0(\bm\Theta+c\bm\Delta)\bm Q_0^+\bm\beta_0^{\mathrm{rls}}.
\]
We next analyze term \(T_{2n}\), using the fact that \(\bm Q_0\bm P_0\bm\beta_0=\bm 0\),
\begin{eqnarray*}
T_{2n}
&=&
\sqrt n\check{\bm Q}_n^+(\bm Q_n-\bm Q_0)\bm P_0\bm\beta_0\\
&
=&
\check{\bm Q}_n^+\sqrt n(\bm Q_n-\bm Q_{0n})\bm P_0\bm\beta_0
+
\frac{\sqrt n}{\tau_n}\check{\bm Q}_n^+\bm\Delta\bm P_0\bm\beta_0.
\end{eqnarray*}
Since \(\check{\bm Q}_n^+\to_{\Pr}\bm Q_0^+\), Assumption~\ref{ass: joint_clt_Q_score}
and Slutsky's lemma imply
\[
T_{2n}
\to_d
\bm Q_0^+\bm\Theta\bm P_0\bm\beta_0
+
c\bm Q_0^+\bm\Delta\bm P_0\bm\beta_0
=
\bm Q_0^+(\bm\Theta+c\bm\Delta)\bm P_0\bm\beta_0.
\]
Finally,
\[
T_{3n}
=
\check{\bm Q}_n^+\bm X'\bm\varepsilon/\sqrt n
\to_d
\bm Q_0^+\bm Z\ ,
\]
by Proposition~\ref{pro:cons_rank_Qn_check},
Assumption~\ref{ass: joint_clt_Q_score}, and Slutsky's lemma.
Since \(T_{1n}\) and \(T_{2n}\) are continuous transforms of
\(\sqrt n(\bm Q_n-\bm Q_{0n})\) up to \(o_{\Pr}(1)\), while \(T_{3n}\) is a
continuous transform of \(\bm X'\bm\varepsilon/\sqrt n\) up to \(o_{\Pr}(1)\),
the joint convergence in Assumption~\ref{ass: joint_clt_Q_score} yields
\[
(T_{1n},T_{2n},T_{3n})
\to_d
\Bigl(
\bm P_0(\bm\Theta+c\bm\Delta)\bm Q_0^+\bm\beta_0^{\mathrm{rls}},
\bm Q_0^+(\bm\Theta+c\bm\Delta)\bm P_0\bm\beta_0,
\bm Q_0^+\bm Z
\Bigr).
\]
Summing the three components proves part~(ii).

\item Suppose, in addition, that
Assumption~\ref{ass: vanishing stochastic term in modified Ridgeless asymptotics}
holds.
For \(T_{1n}\), the random part is
\[
\sqrt n\bm P_0(\bm Q_n-\bm Q_{0n})\bm Q_0^+\bm\beta_0^{\mathrm{rls}}.
\]
Since \(\bm Q_0^+\bm\beta_0^{\mathrm{rls}}\in\Range(\bm Q_0)=\Range(\bm P_0^\perp)\),
Assumption~\ref{ass: vanishing stochastic term in modified Ridgeless asymptotics}
implies
\[
\sqrt n\bm P_0(\bm Q_n-\bm Q_{0n})\bm Q_0^+\bm\beta_0^{\mathrm{rls}}
=
o_{\Pr}(1).
\]
Hence, by the expansion already derived for \(T_{1n}\),
\[
T_{1n}
=
c\bm P_0\bm\Delta\bm Q_0^+\bm\beta_0^{\mathrm{rls}}
+
o_{\Pr}(1).
\]
For \(T_{2n}\), the random part is
\[
\check{\bm Q}_n^+\sqrt n(\bm Q_n-\bm Q_{0n})\bm P_0\bm\beta_0.
\]
Since \(\check{\bm Q}_n^+\to_{\Pr}\bm Q_0^+\) and
\(\sqrt n(\bm Q_n-\bm Q_{0n})\bm P_0\bm\beta_0=O_{\Pr}(1)\),
\[
\check{\bm Q}_n^+\sqrt n(\bm Q_n-\bm Q_{0n})\bm P_0\bm\beta_0
=
\bm Q_0^+\sqrt n(\bm Q_n-\bm Q_{0n})\bm P_0\bm\beta_0
+
o_{\Pr}(1).
\]
Moreover, \(\bm Q_0^+=\bm Q_0^+\bm P_0^\perp\), and since
\(\bm Q_n-\bm Q_{0n}\) is symmetric,
\[
\bigl\|\bm P_0^\perp(\bm Q_n-\bm Q_{0n})\bm P_0\bigr\|_{\op}
=
\bigl\|\bm P_0(\bm Q_n-\bm Q_{0n})\bm P_0^\perp\bigr\|_{\op}
=
o_{\Pr}(n^{-1/2}).
\]
Therefore,
\[
\bm Q_0^+\sqrt n(\bm Q_n-\bm Q_{0n})\bm P_0\bm\beta_0
=
o_{\Pr}(1),
\]
and thus
\[
T_{2n}
=
c\bm Q_0^+\bm\Delta\bm P_0\bm\beta_0
+
o_{\Pr}(1).
\]
Since
\(T_{3n}\to_d \bm Q_0^+\bm Z\)
as for part (ii),
combining the last three displays proves part~(iii).
\end{enumerate}
\end{proof}

\begin{proof}[Proof of Corollary~\ref{cor:asymptotics_irregular}]
Set \(\bm\beta_0:=\bm\beta_0^{\mathrm{rls}}\), \(r_n:=\sqrt n\) and \(\hat{\bm\beta}_n^s:=\check{\bm\beta}_n^{\mathrm{rls}}\),
and define the weight matrix
\[
\bm W_n:=\overline{\check{\bm Q}}_n
:=\check{\bm Q}_n+\bm I_p-\check{\bm Q}_n\check{\bm Q}_n^+
=\check{\bm Q}_n+\check{\bm P}_n,
\qquad
\check{\bm P}_n:=\bm I_p-\check{\bm Q}_n\check{\bm Q}_n^+.
\]
We first verify Assumption~\ref{ass: W}. Since \(\check{\bm P}_n\) is the orthogonal projector onto
\(\Kernel(\check{\bm Q}_n)\), we have
\[
\Range(\check{\bm Q}_n)\perp \Range(\check{\bm P}_n),
\qquad
\R^p=\Range(\check{\bm Q}_n)\oplus\Range(\check{\bm P}_n).
\]
Hence, for every \(\bm b\in\R^p\),
\[
\bm b'\bm W_n\bm b
=
\bm b'\check{\bm Q}_n\bm b+\|\check{\bm P}_n\bm b\|_2^2.
\]
If \(\bm b\neq \bm 0\), then either \(\check{\bm P}_n\bm b\neq \bm 0\), in which case \(\bm b'\bm W_n\bm b>0\), or
\(\check{\bm P}_n\bm b=\bm 0\), in which case \(\bm b\in\Range(\check{\bm Q}_n)\setminus\{\bm 0\}\) and therefore
\(\bm b'\check{\bm Q}_n\bm b>0\). Thus \(\bm W_n\succ\bm 0\) for every \(n\).
By Proposition~\ref{pro:cons_rank_Qn_check}, \(\check{\bm Q}_n\to_{\Pr}\bm Q_0\) and \(\check{\bm Q}_n^+\to_{\Pr}\bm Q_0^+\).
Therefore, by the continuous mapping theorem,
\[
\check{\bm Q}_n\check{\bm Q}_n^+\to_{\Pr}\bm Q_0\bm Q_0^+,
\]
and hence
\[
\bm W_n
=
\check{\bm Q}_n+\bm I_p-\check{\bm Q}_n\check{\bm Q}_n^+
\to_{\Pr}
\bm Q_0+\bm I_p-\bm Q_0\bm Q_0^+
=:\overline{\bm Q}_0.
\]
Since \(\bm Q_0\bm Q_0^+\) is the orthogonal projector onto \(\Range(\bm Q_0)\), we may write
\[
\overline{\bm Q}_0=\bm Q_0+\bm P_0,
\qquad
\bm P_0:=\bm I_p-\bm Q_0\bm Q_0^+,
\]
so \(\overline{\bm Q}_0\succ\bm 0\).
Next, by the standing assumption of the corollary, equivalently by the relevant case of Theorem~\ref{thm:asy_modified_Ridgeless},
\[
r_n(\hat{\bm\beta}_n^s-\bm\beta_0)
=
\sqrt n\big(\check{\bm\beta}_n^{\mathrm{rls}}-\bm\beta_0^{\mathrm{rls}}\big)
\to_d
\bm Q_0^+\bm Z.
\]
Thus Assumption~\ref{ass: asy distr proximal} holds with
\[
\bm\eta=\bm Q_0^+\bm Z,
\qquad
\bm W_n=\overline{\check{\bm Q}}_n,
\qquad
\bm W_0=\overline{\bm Q}_0.
\]
Therefore Theorem~\ref{prop:asymptotics:alt} applies with initial estimator
\(\check{\bm\beta}_n^{\mathrm{rls}}\), rate \(r_n=\sqrt n\), center
\(\bm\beta_0^{\mathrm{rls}}\), and weight matrices \(\overline{\check{\bm Q}}_n\).
If \(\lambda_n\sqrt n\to\lambda_0>0\) and \(q_n\to_{\Pr}\rho_{\bm\beta_0^{\mathrm{rls}}}\) in epigraph,
then Theorem~\ref{prop:asymptotics:alt}\ref{prop_distr_i} yields
\[
\sqrt n\big(\check{\bm\beta}_n^{+}-\bm\beta_0^{\mathrm{rls}}\big)
\to_d
\prox_{\lambda_0\rho_{\bm\beta_0^{\mathrm{rls}}}}^{\overline{\bm Q}_0}(\bm Q_0^+\bm Z),
\]
which is \eqref{eq:asydistr_i_irregular}. If, in addition,
\(\partial f_0(\bm\beta_0^{\mathrm{rls}})\neq\emptyset\), then
Theorem~\ref{prop:asymptotics:alt}\ref{prop_distr_i} also gives
\eqref{eq:asym_moreau_i_irregular}.
If instead \(\lambda_n\sqrt n\to 0\) and \((\lambda_n\sqrt n)q_n
\to_{\Pr}
\sigma_{N_{\dom(f_0)}(\bm\beta_0^{\mathrm{rls}})}\) in epigraph,
then Theorem~\ref{prop:asymptotics:alt}\ref{prop_distr_ii} yields
\[
\sqrt n\big(\check{\bm\beta}_n^{+}-\bm\beta_0^{\mathrm{rls}}\big)
\to_d
\prox_{\sigma_{N_{\dom(f_0)}(\bm\beta_0^{\mathrm{rls}})}}^{\overline{\bm Q}_0}(\bm Q_0^+\bm Z),
\]
which is \eqref{eq:asydistr_ii_irregular}, and the corresponding Moreau representation
\eqref{eq:asym_moreau_ii_irregular} follows from the same theorem.
This concludes the proof.
\end{proof}

\begin{proof}[Proof of Corollary~\ref{cor:oracle_prox_est}]
Let
\[
\bm\beta_0:=\bm\beta_0^{\mathrm{rls}},
\qquad
r_n:=\sqrt n,
\qquad
\mathcal S_0:=\{\bm \beta\in\R^p:\ \beta_j=0\ \text{for all}\ j\in\mathcal A^c\}.
\]
Write the adaptive penalty in weighted-\(\ell_1\) form as
\[
f_n(\bm\beta)
=
\sum_{j:w_{nj}<\infty} w_{nj}|\beta_j|
+
\sum_{j:w_{nj}=+\infty}\iota_{\{0\}}(\beta_j),
\qquad
w_{nj}:=
\begin{cases}
|\check\beta_{n,j}^{\mathrm{rls}}|^{-1}, & \check\beta_{n,j}^{\mathrm{rls}}\neq 0,\\[0.3em]
+\infty, & \check\beta_{n,j}^{\mathrm{rls}}=0.
\end{cases}
\]
Since
\[
\sqrt n\big(\check{\bm\beta}_n^{\mathrm{rls}}-\bm\beta_0^{\mathrm{rls}}\big)\to_d \bm Q_0^+\bm Z,
\]
we have \(\check{\bm\beta}_n^{\mathrm{rls}}\to_{\Pr}\bm\beta_0\). Hence
Proposition~\ref{prop:epi_table1}\ref{prop:epi_table1:alasso} applies, and the epi-limit \(f_0\) of the adaptive penalty satisfies
\[
\dom(f_0)=\mathcal S_0,
\qquad
N_{\dom(f_0)}(\bm\beta_0)=\mathcal S_0^\perp.
\]
Next we verify the scaled local epi-limit required in
Corollary~\ref{cor:asymptotics_irregular}(\ref{prop_distr_ii_irregular}).
Fix \(j\in\mathcal A^c\). Since \((\beta_0)_j=0\) and
\(\sqrt n\check\beta_{n,j}^{\mathrm{rls}}=O_{\Pr}(1)\), we have
\(\check\beta_{n,j}^{\mathrm{rls}}=O_{\Pr}(n^{-1/2})\). Therefore
\[
\frac{|\check\beta_{n,j}^{\mathrm{rls}}|}{\lambda_n r_n}
=
\frac{|\sqrt n\check\beta_{n,j}^{\mathrm{rls}}|}{\lambda_n n}
\to_{\Pr}0,
\]
because \(\lambda_n n\to\infty\). Hence
Proposition~\ref{prop:local_scaled_epi_table1}\ref{prop:local_scaled_epi_table1:alasso},
applied with \(\tilde{\bm\beta}_n=\check{\bm\beta}_n^{\mathrm{rls}}\) and \(r_n=\sqrt n\), yields
\[
(\lambda_n\sqrt n)q_n
\to_{\Pr}
\sigma_{N_{\dom(f_0)}(\bm\beta_0)}
=
\sigma_{\mathcal S_0^\perp}
\qquad\text{in epigraph.}
\]
Since also \(\lambda_n\sqrt n\to 0\),
Corollary~\ref{cor:asymptotics_irregular}(\ref{prop_distr_ii_irregular}) gives
\[
\sqrt n\big(\check{\bm\beta}_n^+-\bm\beta_0\big)
\to_d
\prox_{\sigma_{\mathcal S_0^\perp}}^{\overline{\bm Q}_0}(\bm Q_0^+\bm Z).
\]
Set \(\bm\Omega_0:=\sigma^2\bm Q_0\). Then
\(\bm Z\sim\mathcal N(\bm0,\bm\Omega_0)\), and
Assumption~\ref{ass:gaussian initial estimator} holds for the benchmark
limit \(\bm\eta=\bm Q_0^+\bm Z\), with
\[
\bm M_0=\bm Q_0,
\qquad
\bm W_0=\overline{\bm Q}_0,
\qquad
\bm\Omega_0=\sigma^2\bm Q_0.
\]
If and only if the covariance identity \eqref{eq:oracle_cov_identity} holds,
specialized here to \(\bm M_0=\bm Q_0\), \(\bm W_0=\overline{\bm Q}_0\),
and \(\bm\Omega_0=\sigma^2\bm Q_0\),
Proposition~\ref{prop:op2} yields
\[
\sqrt n(\check{\bm\beta}_n^+-\bm\beta_0^{\mathrm{rls}})_{\mathcal A}
\to_d
[(\bm Q_0)_{\mathcal A}]^+(\bm Z)_{\mathcal A},
\]
which is the limit \eqref{eq:oracle_limit_oracle_prox}.
It remains to verify Oracle selection. We apply
Corollary~\ref{cor:oracle_weighted_l1} with
\[
\bm W_n=\overline{\check{\bm Q}}_n,
\qquad
r_n=\sqrt n,
\qquad
w_{nj}:=
\begin{cases}
|\check\beta_{n,j}^{\mathrm{rls}}|^{-1}, & \check\beta_{n,j}^{\mathrm{rls}}\neq 0,\\
+\infty, & \check\beta_{n,j}^{\mathrm{rls}}=0.
\end{cases}
\]
First, Assumption~\ref{ass: W} holds because
\(\overline{\check{\bm Q}}_n\to_{\Pr}\overline{\bm Q}_0\succ\bm 0\).
Second, the localization requirement \eqref{eq:oracle_localization}
holds with \(r_n=\sqrt n\) (the benchmark localization follows from the
assumed CLT, and the proximal localization follows from
Corollary~\ref{cor:asymptotics_irregular}(ii)). Finally, if
\(j\in\mathcal A^c\), then \((\beta_0)_j=0\) and
\(\sqrt n\check\beta_{n,j}^{\mathrm{rls}}=O_{\Pr}(1)\), so:\footnote{If \(\mathcal A^c=\emptyset\), condition
\eqref{eq:weighted_l1_oracle_condition} is immediate by convention.
Suppose therefore that \(\mathcal A^c\neq\emptyset\).}
\[
\lambda_n r_n \min_{j\in\mathcal A^c} w_{nj}
\ge
\lambda_n n
\left(\max_{j\in\mathcal A^c}
|\sqrt n\check\beta_{n,j}^{\mathrm{rls}}|\right)^{-1}
\to_{\Pr}\infty.
\]
Hence condition \eqref{eq:weighted_l1_oracle_condition} holds, and
Corollary~\ref{cor:oracle_weighted_l1} yields
\[
\Pr(\check{\mathcal A}_n^+=\mathcal A)\to 1.
\]
Finally, under Assumption~\ref{ass: ridgless oracle compatibility condition},
\[
[(\bm Q_0)_{\mathcal A}]^+(\bm\delta_0)_{\mathcal A}
=
(\bm\beta_0^{\mathrm{rls}})_{\mathcal A}.
\]
Hence the limiting law above is the oracle law for the active coordinates of
the full Ridgeless estimand.
This concludes the proof.
\end{proof}

\begin{proof}[Proof of Corollary~\ref{cor:oracle_prox_est_double}]
The result follows by the same arguments as in the proof of Corollary~\ref{cor:oracle_prox_est}. Indeed, by construction, the plug-in weighting matrix \(\check{\bm W}_n\) is a consistent estimator of \(\bm W_0\succ \bm 0\), and
the limiting proximal metric induced by $\bm W_0$ coincides with the oracle design geometry on the active block, so that the covariance identity \eqref{eq:oracle_cov_identity} holds automatically:
\begin{eqnarray*}
(\bm W_0 \bm M_0^+ \bm \Omega_0 \bm M_0^+ \bm W_0 )_{\cal A}
= \sigma^2 (\bm W_0 \bm Q_0^+ \bm W_0)_{\cal A}
= \sigma^2 (\bm Q_0)_{\cal A}
\end{eqnarray*}
where the second identity holds because diagonal matrix $\bm D ({\cal A})$ has nonzero diagonal elements only at inactive coordinates with index in ${\cal A}^c$.
The asserted support consistency and asymptotic distribution then follow exactly as in Corollary~\ref{cor:oracle_prox_est}.
\end{proof}

\begin{proof}[Proof of Corollary~\ref{cor:oracle_post_Ridgeless}]
By selection consistency, with probability tending to one,
\(\check{\mathcal A}_n^+=\mathcal A\). On this event, the post-selection estimator satisfies
\[
(\tilde{\bm\beta}_n)_{\mathcal A}
=
[\check{\bm Q}_n(\mathcal A)]^+
(\bm X'\bm Y/n)_{\mathcal A},
\qquad
(\tilde{\bm\beta}_n)_{\mathcal A^c}=\bm 0.
\]
By the consistency and rank consistency of the thresholded active block, 
which follow by the same argument as in Proposition~\ref{pro:cons_rank_Qn_check} applied to the fixed active block \(\mathcal A\),
\[
[\check{\bm Q}_n(\mathcal A)]^+
\to_p
[(\bm Q_0)_{\mathcal A}]^+.
\]
Combining this with the asymptotic expansion in
Theorem~\ref{thm:asy_modified_Ridgeless} for the ridgeless refit on the active variables yields the following asymptotic distribution in the nearly-singular regime \(\sqrt n/\tau_n\to c\in[0,\infty)\):
\begin{eqnarray*}
\sqrt n\big((\tilde{\bm\beta}_n)_{\mathcal A}
-
[(\bm Q_0)_{\mathcal A}]^+(\bm\delta_0)_{\mathcal A}\big)
&\to_d &
\bm P_0({\cal A})(\bm\Theta+c\bm\Delta)_{\cal A}[(\bm Q_0)_{\mathcal A}]^+[(\bm Q_0)_{\mathcal A}]^+(\bm\delta_0)_{\mathcal A}\\
&&+
[(\bm Q_0)_{\mathcal A}]^+(\bm\Theta+c\bm\Delta)_{\cal A}\bm P_0({\cal A})(\bm\beta_0)_{\cal A}\\
&&
+
[(\bm Q_0)_{\mathcal A}]^+(\bm Z)_{\mathcal A}.
\end{eqnarray*}
Here,
\[
\bm P_0({\cal A})
:=
\bm I_{|{\cal A}|}
-
(\bm Q_0)_{\mathcal A}
[(\bm Q_0)_{\mathcal A}]^+.
\]
For $\sqrt{n}=o(\tau_n)$ one has $c=0$, so the deterministic drift component vanishes. Moreover, if Assumption~\ref{ass: vanishing stochastic term in modified Ridgeless asymptotics}
holds for the active subblock:
\[
\bm P_0 ({\cal A}) ( \bm Q_n-\bm Q_{0n})_{\cal A} \bm P_0^\perp  ({\cal A}) = o_{\Pr} ({n}^{-1/2}),
\]
then the stochastic drift component also vanishes and
\[
\sqrt n\big((\tilde{\bm\beta}_n)_{\mathcal A}
-
[(\bm Q_0)_{\mathcal A}]^+(\bm\delta_0)_{\mathcal A}\big)
\to_d
[(\bm Q_0)_{\mathcal A}]^+(\bm Z)_{\mathcal A}.
\]
Finally, Assumption~\ref{ass: ridgless oracle compatibility condition} identifies
\[
[(\bm Q_0)_{\mathcal A}]^+(\bm\delta_0)_{\mathcal A}
=
(\bm\beta_0^{\mathrm{rls}})_{\mathcal A},
\]
so that $(\tilde{\bm\beta}_n)_{\mathcal A}\to_{\Pr} (\bm\beta_0^{\mathrm{rls}})_{\mathcal A}$. Since
\[
\mathcal A
=
\{j:(\bm\beta_0^{\mathrm{rls}})_j\neq0\},
\]
the selection statement follows together with \(\Pr(\check{\mathcal A}_n^+=\mathcal A)\to1\).
This concludes the proof.
\end{proof}

\section{Proofs of supporting results}\label{sec:proofs_appendix}

\subsection{Proofs of the results in Section~\ref{app:nonunique_penGMM}}

\begin{proof}[Proof of Proposition~\ref{prop:nonunique_penGMM_app}]
We prove the two parts in turn.

\begin{enumerate}[label=(\roman*)]
\item
The map \(h\) is continuous and strictly decreasing on \(\R\), with
\[
\lim_{\beta\to-\infty} h(\beta)=+\infty,
\qquad
\lim_{\beta\to+\infty} h(\beta)=-\infty.
\]
Hence \(h\) is a bijection from \(\R\) onto \(\R\).
Therefore, for every realization of \((u_i)_{i=1}^n\), the equation
\[
h(\beta)=-\bar u_n
\]
admits a unique solution, denoted \(\hat\beta_n\).
At that point,
\[
Q_n(\hat\beta_n)=\big(h(\hat\beta_n)+\bar u_n\big)^2=0.
\]
Since \(Q_n(\beta)\ge 0\) for every \(\beta\in\R\), it follows that \(\hat\beta_n\) is the unique global minimizer of \(Q_n\).

\item
Now suppose \(n\) is even and \(\lambda=1\).
On the event \(\{\bar u_n=0\}\), we have \(Q_n(\beta)=h(\beta)^2\), so
\[
F_{n,1}(\beta)=h(\beta)^2+|\beta|
=
\begin{cases}
1-2\beta, & \beta<0,\\[0.3em]
1, & 0\le \beta\le 1,\\[0.3em]
\beta^2-\beta+1, & \beta>1.
\end{cases}
\]
Therefore \(F_{n,1}(\beta)\ge 1\) for every \(\beta\in\R\), with equality if and only if \(\beta\in[0,1]\).
Hence
\(\argmin_{\beta\in\R} F_{n,1}(\beta)=[0,1]\) on 
\(\{\bar u_n=0\}\).

If \(n\) is even, the event \(\{\bar u_n=0\}\) occurs if and only if exactly \(n/2\) of the \(u_i\)'s equal \(1\).
Thus
\[
\Pr(\bar u_n=0)=\binom{n}{n/2}2^{-n}>0.
\]
Since \(\{\bar u_n=0\}\subseteq E\), we obtain
\[
\Pr(E)\ge \Pr(\bar u_n=0)=\binom{n}{n/2}2^{-n}>0.
\]
\end{enumerate}
This proves the claim.
\end{proof}

\subsection{Proofs of the results in Section~\ref{app:sec:epigraph_benchmark}}

Before proving Propositions~\ref{prop:epi_table1}–\ref{prop:local_scaled_epi_table1} in this section, we record two auxiliary results that will be used repeatedly throughout the sequel: a deterministic epigraphical limit for weighted \(\ell_1\) penalties with possibly infinite weights, and a local epigraphical limit specifically tailored to the Adaptive Lasso.

\begin{lemma}[Weighted \(\ell_1\) epi-limits with possibly infinite weights]
\label{lem:weighted_l1_epi_det}
For each \(n\in\mathbb N\), let \(a_{nj}\in[0,\infty]\), \(j=1,\dots,p\), be deterministic,
and define
\[
g_n(\bm\beta)
:=
\sum_{j:a_{nj}<\infty} a_{nj}|\beta_j|
+
\sum_{j:a_{nj}=+\infty}\iota_{\{0\}}(\beta_j),
\qquad \bm\beta\in\R^p.
\]
Assume that, for each \(j=1,\dots,p\), \(a_{nj}\to a_j\in[0,\infty]\). Define
\[
\mathcal S
:=
\{\bm\beta\in\R^p:\ \beta_j=0 \text{ whenever } a_j=+\infty\},
\]
and
\[
g_0(\bm\beta)
:=
\sum_{j:a_j<\infty} a_j|\beta_j|+\iota_{\mathcal S}(\bm\beta).
\]
Then \(g_n,g_0\in\Gamma(\R^p)\) and \(g_n\to g_0\) in epigraph.
\end{lemma}

\begin{proof}
Each coordinate term in \(g_n\) is proper, lower semicontinuous, and convex on
\(\R\), hence \(g_n\in\Gamma(\R^p)\) for every \(n\). The same holds for \(g_0\).
We verify the sequential characterization of epi-convergence \citep[Prop.~7.2]{SM:rockafellarWets2009}.
Let \(\bm\beta_n\to\bm\beta\). If \(\bm\beta\notin\mathcal S\), choose \(j\) such that
\(a_j=+\infty\) and \(\beta_j\neq 0\). Then \(|\beta_{nj}|\ge c>0\) for all sufficiently
large \(n\), and \(a_{nj}\to+\infty\). Therefore
\[
g_n(\bm\beta_n)\ge a_{nj}|\beta_{nj}|\to+\infty=g_0(\bm\beta).
\]
If instead \(\bm\beta\in\mathcal S\), then for every \(j\) with \(a_j<\infty\),
\[
a_{nj}|\beta_{nj}|\to a_j|\beta_j|,
\]
whereas the remaining terms are nonnegative. Hence
\[
\liminf_{n\to\infty} g_n(\bm\beta_n)
\ge
\sum_{j:a_j<\infty} a_j|\beta_j|
=
g_0(\bm\beta).
\]
This proves the liminf inequality.
For the recovery sequence, fix \(\bm\beta\in\R^p\) and take the constant sequence
\(\bm\beta_n=\bm\beta\). If \(\bm\beta\in\mathcal S\), then
\[
g_n(\bm\beta)\to \sum_{j:a_j<\infty} a_j|\beta_j|=g_0(\bm\beta).
\]
If \(\bm\beta\notin\mathcal S\), then for some \(j\) one has \(a_j=+\infty\) and
\(\beta_j\neq 0\), so \(g_n(\bm\beta)\ge a_{nj}|\beta_j|\to+\infty=g_0(\bm\beta)\).
Thus
\[
\limsup_{n\to\infty} g_n(\bm\beta)\le g_0(\bm\beta)
\]
in both cases. Therefore \(g_n\to g_0\) in epigraph.
\end{proof}

\begin{lemma}[Random weighted \(\ell_1\) epi-limits]
\label{lem:weighted_l1_epi_prob}
For each \(n\in\mathbb N\), let \(a_{nj}\in[0,\infty]\), \(j=1,\dots,p\), be random,
and define
\[
g_n(\bm\beta)
:=
\sum_{j:a_{nj}<\infty} a_{nj}|\beta_j|
+
\sum_{j:a_{nj}=+\infty}\iota_{\{0\}}(\beta_j),
\qquad \bm\beta\in\R^p.
\]
Assume that, for each \(j=1,\dots,p\), \(a_{nj}\to_{\Pr} a_j\in[0,\infty]\), where
\(a_j\) is deterministic. Define \(\mathcal S\) and \(g_0\) as in
Lemma~\ref{lem:weighted_l1_epi_det}. Then \(g_n\to_{\Pr} g_0\) in epigraph.
\end{lemma}

\begin{proof}
Let \((n_k)\) be an arbitrary subsequence. Since \(p\) is fixed, there exists a
further subsequence, not relabeled, and an event \(\Omega_0\) of probability one
such that \(a_{n_kj}(\omega)\to a_j\) for every \(j=1,\dots,p\) and every
\(\omega\in\Omega_0\). By Lemma~\ref{lem:weighted_l1_epi_det}, \(g_{n_k}(\cdot,\omega)\to g_0\)
in epigraph for every \(\omega\in\Omega_0\).
Hence every subsequence admits a further subsequence converging almost surely in
epigraph to \(g_0\), which is equivalent to \(g_n\to_{\Pr} g_0\) in epigraph.
\end{proof}

\begin{lemma}[Local epi-limits for the Adaptive Lasso]
\label{lem:alasso_local_epi}
Assume that \(\bm\beta_0\in\dom(f_n^{\mathrm{al}})\) for all sufficiently large \(n\),
and let
\[
q_n^{\mathrm{al}}(\bm b)
:=
r_n\Bigl[
f_n^{\mathrm{al}}\bigl(\bm\beta_0+\bm b/r_n\bigr)
-
f_n^{\mathrm{al}}(\bm\beta_0)
\Bigr],
\qquad
r_n\to\infty.
\]
Define
\[
\ell_0(\bm b)
:=
\sum_{j\in\mathcal A_0}\frac{\sign(\beta_{0j})}{|\beta_{0j}|}b_j.
\]

\begin{enumerate}[label=(\roman*)]
\item If \(\tilde{\bm\beta}_n\to_{\Pr}\bm\beta_0\), then \(q_n^{\mathrm{al}}\to_{\Pr}\ell_0+\iota_{\mathcal S_0}\) in epigraph.

\item Let \(a_n:=\lambda_n r_n\), and assume additionally that \(a_n\to 0\) and \(|\tilde\beta_{nj}|/a_n\to_{\Pr}0\) for every \(j\notin\mathcal A_0\).
Then
\[
a_n q_n^{\mathrm{al}}
\to_{\Pr}
\iota_{\mathcal S_0}
=
\sigma_{N_{\mathcal S_0}(\bm\beta_0)}
\qquad\text{in epigraph.}
\]
\end{enumerate}
\end{lemma}

\begin{proof}
We prove both parts by the subsequence principle. Let \((n_k)\) be an arbitrary
subsequence.
For part (i), pass to a further subsequence, not relabeled, such that
\(\tilde\beta_{n_kj}\to\beta_{0j}\) almost surely for every \(j\). Fix \(\omega\) in
the corresponding probability-one event. Write
\[
q_{n_k}^{\mathrm{al}}(\bm b,\omega)=\ell_k(\bm b)+g_k(\bm b),
\]
where
\[
\ell_k(\bm b)
:=
\sum_{j\in\mathcal A_0}
\frac{r_{n_k}\bigl(|\beta_{0j}+b_j/r_{n_k}|-|\beta_{0j}|\bigr)}
{|\tilde\beta_{n_kj}(\omega)|}
\]
and
\[
g_k(\bm b)
:=
\sum_{j\notin\mathcal A_0} d_{kj}|b_j|,
\qquad
d_{kj}
:=
\begin{cases}
|\tilde\beta_{n_kj}(\omega)|^{-1}, & \tilde\beta_{n_kj}(\omega)\neq 0,\\
+\infty, & \tilde\beta_{n_kj}(\omega)=0.
\end{cases}
\]
On every compact set, sign stability on the active set gives \(\ell_k\to\ell_0\)
locally uniformly. For \(j\notin\mathcal A_0\), \(\tilde\beta_{n_kj}(\omega)\to 0\),
so \(d_{kj}\to+\infty\); hence Lemma~\ref{lem:weighted_l1_epi_det} gives
\(g_k\to\iota_{\mathcal S_0}\) in epigraph. A direct liminf/recovery-sequence
argument then yields
\[
q_{n_k}^{\mathrm{al}}(\cdot,\omega)\to \ell_0+\iota_{\mathcal S_0}
\qquad\text{in epigraph.}
\]
Since the original subsequence was arbitrary, \(q_n^{\mathrm{al}}\to_{\Pr}
\ell_0+\iota_{\mathcal S_0}\) in epigraph.
For part (ii), pass to a further subsequence, not relabeled, such that
\(\tilde\beta_{n_kj}\to\beta_{0j}\) almost surely for all \(j\), and
\(|\tilde\beta_{n_kj}|/a_{n_k}\to 0\) almost surely for \(j\notin\mathcal A_0\).
Fix \(\omega\) in the corresponding probability-one event and write
\[
a_{n_k}q_{n_k}^{\mathrm{al}}(\bm b,\omega)
=
\tilde\ell_k(\bm b)+\tilde g_k(\bm b),
\]
where
\[
\tilde\ell_k(\bm b)
:=
\sum_{j\in\mathcal A_0}
\frac{a_{n_k}r_{n_k}\bigl(|\beta_{0j}+b_j/r_{n_k}|-|\beta_{0j}|\bigr)}
{|\tilde\beta_{n_kj}(\omega)|}
\]
and
\[
\tilde g_k(\bm b)
:=
\sum_{j\notin\mathcal A_0}\tilde d_{kj}|b_j|,
\qquad
\tilde d_{kj}
:=
\begin{cases}
a_{n_k}|\tilde\beta_{n_kj}(\omega)|^{-1}, & \tilde\beta_{n_kj}(\omega)\neq 0,\\
+\infty, & \tilde\beta_{n_kj}(\omega)=0.
\end{cases}
\]
On every compact set, \(\tilde\ell_k\to 0\) locally uniformly because
\(a_{n_k}\to 0\) and \(|\tilde\beta_{n_kj}(\omega)|\to |\beta_{0j}|>0\) on the active
set. For \(j\notin\mathcal A_0\), the additional assumption implies
\(\tilde d_{kj}\to+\infty\), so Lemma~\ref{lem:weighted_l1_epi_det} yields
\(\tilde g_k\to\iota_{\mathcal S_0}\) in epigraph. The same liminf/recovery-sequence
argument gives
\[
a_{n_k}q_{n_k}^{\mathrm{al}}(\cdot,\omega)\to\iota_{\mathcal S_0}
\qquad\text{in epigraph.}
\]
Since the original subsequence was arbitrary, \(a_n q_n^{\mathrm{al}}\to_{\Pr}
\iota_{\mathcal S_0}\) in epigraph.
Finally, \(\mathcal S_0\) is a linear subspace containing \(\bm\beta_0\), so
\(N_{\mathcal S_0}(\bm\beta_0)=\mathcal S_0^\perp\). Because the support function
of a cone is the indicator of its polar,
\[
\sigma_{N_{\mathcal S_0}(\bm\beta_0)}
=
\sigma_{\mathcal S_0^\perp}
=
\iota_{(\mathcal S_0^\perp)^\circ}
=
\iota_{\mathcal S_0}.
\]
\end{proof}

Now we can proceed to prove the propositions in Section~\ref{app:sec:epigraph_benchmark}.

\begin{proof}[Proof of Proposition~\ref{prop:epi_table1}]
For items~\ref{prop:epi_table1:ridge}--\ref{prop:epi_table1:constr}, the penalty is
deterministic and does not depend on \(n\). Hence \(f_n=f_0\) for every \(n\), with
\(f_0\) given by the formula displayed in the corresponding item. Therefore
\(f_n\to f_0\) in epigraph, and since the sequence is deterministic,
\(f_n\to_{\Pr} f_0\) in epigraph.

For item~\ref{prop:epi_table1:alasso}, write
\[
f_n^{\mathrm{al}}(\bm\beta)
=
\sum_{j:\tilde\beta_{nj}\neq 0}\frac{|\beta_j|}{|\tilde\beta_{nj}|}
+
\sum_{j:\tilde\beta_{nj}=0}\iota_{\{0\}}(\beta_j).
\]
This is a weighted \(\ell_1\) penalty with random weights
\[
a_{nj}
:=
\begin{cases}
|\tilde\beta_{nj}|^{-1}, & \tilde\beta_{nj}\neq 0,\\
+\infty, & \tilde\beta_{nj}=0.
\end{cases}
\]
If \(j\in\mathcal A_0\), then \(a_{nj}\to_{\Pr} |\beta_{0j}|^{-1}\); if
\(j\notin\mathcal A_0\), then \(a_{nj}\to_{\Pr}+\infty\). Hence
Lemma~\ref{lem:weighted_l1_epi_prob} yields
\[
f_n^{\mathrm{al}}\to_{\Pr} f_0^{\mathrm{al}}
\qquad\text{in epigraph.}
\]
This proves item~\ref{prop:epi_table1:alasso}.
\end{proof}

\begin{proof}[Proof of Proposition~\ref{prop:epi_scaled_table1}]
For items~\ref{prop:epi_scaled_table1:i} and \ref{prop:epi_scaled_table1:ii},
the penalties are deterministic and sample-invariant, so
\(f_n=f_0\) for every \(n\).
If \(f_0\) is Ridge, Lasso, Group Lasso, or Elastic Net, then \(f_0\in\Gamma(\R^p)\)
and \(f_0\) is finite on \(\R^p\). Hence, for every \(\bm\beta\in\R^p\),
\[
(\lambda_n f_n)(\bm\beta)=\lambda_n f_0(\bm\beta)\to 0.
\]
Since each \(\lambda_n f_0\) and the limit \(0\) belong to \(\Gamma(\R^p)\),
\citet[Thm.~7.17]{SM:rockafellarWets2009} yields
\[
\lambda_n f_n=\lambda_n f_0\to 0=\iota_{\R^p}=\iota_{\dom(f_0)}
\qquad\text{in epigraph.}
\]

If \(f_0=\iota_C\) for some nonempty closed convex set \(C\subseteq\R^p\), then
\[
\lambda_n f_n=\lambda_n\iota_C=\iota_C=\iota_{\dom(f_0)}
\qquad\text{for every }n,
\]
and therefore \(\lambda_n f_n\to \iota_{\dom(f_0)}\) in epigraph. Since these
sequences are deterministic, the same convergences hold in probability.

For item~\ref{prop:epi_scaled_table1:iii}, \(\lambda_n f_n^{\mathrm{al}}\) is again a weighted \(\ell_1\) penalty,
now with random weights
\[
a_{nj}
:=
\begin{cases}
\lambda_n|\tilde\beta_{nj}|^{-1}, & \tilde\beta_{nj}\neq 0,\\
+\infty, & \tilde\beta_{nj}=0.
\end{cases}
\]
If \(j\in\mathcal A_0\), then \(a_{nj}\to_{\Pr}0\); if \(j\notin\mathcal A_0\), the
assumption \(|\tilde\beta_{nj}|/\lambda_n\to_{\Pr}0\) implies
\(a_{nj}\to_{\Pr}+\infty\). Hence Lemma~\ref{lem:weighted_l1_epi_prob} yields
\[
\lambda_n f_n^{\mathrm{al}}\to_{\Pr}\iota_{\mathcal S_0}
\qquad\text{in epigraph.}
\]
By Proposition~\ref{prop:epi_table1}, \(\dom(f_0^{\mathrm{al}})=\mathcal S_0\), so
the limit is \(\iota_{\dom(f_0^{\mathrm{al}})}\).
\end{proof}

\begin{proof}[Proof of Proposition~\ref{prop:local_epi_table1}]
Write \(t_n:=r_n^{-1}\downarrow 0\).
\begin{enumerate}[label=(\roman*)]
    \item Direct calculations give
    \[
q_n(\bm b)
=
\frac{\frac12\|\bm\beta_0+t_n\bm b\|_2^2-\frac12\|\bm\beta_0\|_2^2}{t_n}
=
\bm\beta_0'\bm b+\frac{t_n}{2}\|\bm b\|_2^2
\to
\bm\beta_0'\bm b
\]
for every \(\bm b\in\R^p\). Since \(q_n\) and the limit are finite convex functions
on \(\R^p\), \citet[Thm.~7.17]{SM:rockafellarWets2009} yields epi-convergence.

\item Direct calculations give \[
q_n(\bm b)
=
\sum_{j=1}^p\frac{|\beta_{0j}+t_n b_j|-|\beta_{0j}|}{t_n}.
\]
If \(\beta_{0j}\neq 0\), the sign of \(\beta_{0j}+t_n b_j\) is eventually
\(\sign(\beta_{0j})\), so the \(j\)th summand converges to
\(\sign(\beta_{0j})b_j\). If \(\beta_{0j}=0\), the \(j\)th summand is identically
\(|b_j|\). Hence
\[
q_n(\bm b)\to
\sum_{j:\beta_{0j}\neq 0}\sign(\beta_{0j})b_j
+
\sum_{j:\beta_{0j}=0}|b_j|
\]
for every \(\bm b\). The limit is finite convex, so \citet[Thm.~7.17]{SM:rockafellarWets2009}
again yields epi-convergence.

\item Direct calculations give
\[
q_n(\bm b)
=
\sum_{k=1}^K
\frac{\|\bm\beta_{0k}+t_n\bm b_k\|_2-\|\bm\beta_{0k}\|_2}{t_n}.
\]
If \(\bm\beta_{0k}\neq \bm 0\), the Euclidean norm is differentiable at
\(\bm\beta_{0k}\), so the \(k\)th term converges to
\[
\left(\frac{\bm\beta_{0k}}{\|\bm\beta_{0k}\|_2}\right)'\bm b_k.
\]
If \(\bm\beta_{0k}=\bm 0\), the \(k\)th term is identically \(\|\bm b_k\|_2\).
Therefore
\[
q_n(\bm b)\to
\sum_{k:\bm\beta_{0k}\neq \bm 0}
\left(\frac{\bm\beta_{0k}}{\|\bm\beta_{0k}\|_2}\right)'\bm b_k
+
\sum_{k:\bm\beta_{0k}=\bm 0}\|\bm b_k\|_2,
\]
and epi-convergence follows from \citet[Thm.~7.17]{SM:rockafellarWets2009}.

\item By linearity,
\[
q_n(\bm b)
=
w\frac{\|\bm\beta_0+t_n\bm b\|_1-\|\bm\beta_0\|_1}{t_n}
+
(1-w)\frac{\frac12\|\bm\beta_0+t_n\bm b\|_2^2-\frac12\|\bm\beta_0\|_2^2}{t_n}.
\]
Combining the previous two calculations gives
\[
q_n(\bm b)\to
w\left(
\sum_{j:\beta_{0j}\neq 0}\sign(\beta_{0j})b_j
+
\sum_{j:\beta_{0j}=0}|b_j|
\right)
+(1-w)\bm\beta_0'\bm b,
\]
and epi-convergence again follows from \citet[Thm.~7.17]{SM:rockafellarWets2009}.

\item Since \(f=\iota_C\),
\[
q_n(\bm b)
=
\frac{\iota_C(\bm\beta_0+t_n\bm b)-\iota_C(\bm\beta_0)}{t_n}
=
\iota_{D_n}(\bm b),
\qquad
D_n:=\{\bm b\in\R^p:\ \bm\beta_0+t_n\bm b\in C\},
\]
because \(\bm\beta_0\in C\). The sets \(D_n\) are closed, convex, and increasing.
Moreover,\footnote{For a set \(A\), \(\overline{A}\) denote its closure.}
\[
\overline{\bigcup_{n=1}^\infty D_n}=T_C(\bm\beta_0).
\]
Indeed, \(\bigcup_n D_n\subseteq T_C(\bm\beta_0)\) is immediate. Conversely, if
\(\bm b\in T_C(\bm\beta_0)\), there exist \(s_k\downarrow 0\) and \(\bm b_k\to\bm b\)
such that \(\bm\beta_0+s_k\bm b_k\in C\). Choose \(n_k\) so that \(t_{n_k}\le s_k\).
By convexity of \(C\),
\[
\bm\beta_0+t_{n_k}\bm b_k
=
\Bigl(1-\frac{t_{n_k}}{s_k}\Bigr)\bm\beta_0
+
\frac{t_{n_k}}{s_k}\bigl(\bm\beta_0+s_k\bm b_k\bigr)\in C,
\]
so \(\bm b_k\in D_{n_k}\), and therefore \(\bm b\in\overline{\bigcup_n D_n}\).
We now verify epi-convergence of \(\iota_{D_n}\) to \(\iota_{T_C(\bm\beta_0)}\).
If \(\bm b_n\to\bm b\) and \(\liminf_n \iota_{D_n}(\bm b_n)<+\infty\), then
\(\bm b_n\in D_n\) eventually, so \(\bm b\in T_C(\bm\beta_0)\). This gives the
liminf inequality. For the recovery sequence, fix \(\bm b\in T_C(\bm\beta_0)\),
choose \(\bm c_k\in\bigcup_n D_n\) such that \(\bm c_k\to\bm b\), pick increasing
\(n_k\) with \(\bm c_k\in D_{n_k}\), and define \(\bm b_n:=\bm c_k\) whenever
\(n_k\le n<n_{k+1}\). Since \(D_n\) is increasing, \(\bm b_n\in D_n\) for every \(n\),
and \(\bm b_n\to\bm b\). Hence
\[
q_n=\iota_{D_n}\to \iota_{T_C(\bm\beta_0)}
\qquad\text{in epigraph.}
\]
Because \(T_C(\bm\beta_0)\) and \(N_C(\bm\beta_0)\) are polar cones,
\[
\iota_{T_C(\bm\beta_0)}=\sigma_{N_C(\bm\beta_0)}.
\]
\end{enumerate}
Since items (i)--(v) are deterministic, all these epi-convergences also hold in
probability.
\begin{enumerate}
    \item[(vi)] This item follows directly by
Lemma~\ref{lem:alasso_local_epi}(i), since
\[
\ell_0(\bm b)+\iota_{\mathcal S_0}(\bm b)
=
\sum_{j\in\mathcal A_0}\frac{\sign(\beta_{0j})}{|\beta_{0j}|}b_j
+
\iota_{\mathcal S_0}(\bm b).
\]
\end{enumerate}
\end{proof}

\begin{proof}[Proof of Proposition~\ref{prop:local_scaled_epi_table1}]
Set \(a_n:=\lambda_n r_n\to 0\) and \(t_n:=r_n^{-1}\downarrow 0\).
For item~\ref{prop:local_scaled_epi_table1:finite}, consider first the Ridge case.
From the explicit formula above,
\[
a_n q_n(\bm b)=a_n\bm\beta_0'\bm b+\frac{a_n t_n}{2}\|\bm b\|_2^2\to 0
\qquad\text{for every }\bm b\in\R^p.
\]
For the Lasso, Group Lasso, and Elastic Net cases, the pointwise limits derived in
the proof of Proposition~\ref{prop:local_epi_table1} are finite on \(\R^p\), so \(a_n q_n(\bm b)\to 0\) for every \(\bm b\in\R^p\).
In all four cases, \(a_n q_n\) and the limit \(0\) are finite convex functions on
\(\R^p\); therefore \citet[Thm.~7.17]{SM:rockafellarWets2009} gives \(a_n q_n\to 0\) in epigraph.
Since \(\dom(f_0)=\R^p\) in these cases,
\[
\sigma_{N_{\dom(f_0)}(\bm\beta_0)}
=
\sigma_{N_{\R^p}(\bm\beta_0)}
=
\sigma_{\{\bm 0\}}
=
0.
\]

For item~\ref{prop:local_scaled_epi_table1:constr}, \(q_n=\iota_{D_n}\) with
\(D_n\) as in the proof of Proposition~\ref{prop:local_epi_table1}, hence \(a_n q_n=q_n\) for every \(n\).
The same proof shows that
\[
q_n\to \iota_{T_C(\bm\beta_0)}=\sigma_{N_C(\bm\beta_0)}
\qquad\text{in epigraph.}
\]
Since \(\dom(f_0)=C\), the limit is \(\sigma_{N_{\dom(f_0)}(\bm\beta_0)}\).

Finally, item~\ref{prop:local_scaled_epi_table1:alasso} is exactly
Lemma~\ref{lem:alasso_local_epi}(ii), which yields
\[
(\lambda_n r_n)q_n^{\mathrm{al}}\to_{\Pr}\iota_{\mathcal S_0}
\qquad\text{in epigraph.}
\]
By Proposition~\ref{prop:epi_table1}, \(\dom(f_0^{\mathrm{al}})=\mathcal S_0\), and
Lemma~\ref{lem:alasso_local_epi}(ii) identifies
\(\iota_{\mathcal S_0}=\sigma_{N_{\mathcal S_0}(\bm\beta_0)}\). Hence the limit is \(\sigma_{N_{\dom(f_0^{\mathrm{al}})}(\bm\beta_0)}\).
\end{proof}

\subsection{Proofs of the results in Section~\ref{app:oracle_selection}}

\begin{proof}[Proof of Proposition~\ref{prop:op1_iff}]
For each \(n\), consider the realized proximal criterion
\[
Q_n(\bm\beta)
:=
\frac12\|\bm\beta-\hat{\bm\beta}_n^s\|_{\bm W_n}^2+\lambda_n f_n(\bm\beta),
\qquad \bm\beta\in\R^p.
\]
Since \(\bm W_n\succ\bm 0\), the quadratic term is strongly convex. As \(f_n\) is proper, lower semicontinuous, and convex, it follows that \(Q_n\) is proper, lower semicontinuous, and strongly convex. Hence \(Q_n\) admits a unique minimizer, namely \(\hat{\bm\beta}_n\).

Moreover, by the subdifferential sum rule,
\[
\partial^e Q_n(\bm\beta)
=
\bm W_n(\bm\beta-\hat{\bm\beta}_n^s)+\lambda_n\partial^e f_n(\bm\beta).
\]
Therefore, the first-order optimality condition for convex minimization gives
\[
\bm\beta=\hat{\bm\beta}_n
\quad\Longleftrightarrow\quad
\bm 0\in \partial^e Q_n(\bm\beta)
\quad\Longleftrightarrow\quad
\bm W_n(\hat{\bm\beta}_n^s-\bm\beta)\in \lambda_n\partial^e f_n(\bm\beta).
\]
Consequently, for every \(M>0\),
\begin{equation}\label{eq:event_equiv_K_revised}
\{\mathcal K_n(M)\neq\emptyset\}
=
\left\{
\|\hat{\bm\beta}_n-\bm\beta_0\|_\infty\le M/r_n,\ 
(\hat{\bm\beta}_n)_{\mathcal A^c}\neq \bm 0
\right\}.
\end{equation}
Indeed, if \(\mathcal K_n(M)\neq\emptyset\), then any \(\bm\beta\in\mathcal K_n(M)\) satisfies the KKT condition above and therefore, by uniqueness of the minimizer, must equal \(\hat{\bm\beta}_n\). The converse is immediate.

We now prove the two implications.

\medskip
\noindent\emph{Only if part.}
Suppose Oracle property~\ref{oracle property}\ref{op1} holds, that is,
\[
\Pr(\hat{\mathcal A}_n=\mathcal A)\to 1.
\]
Then, for every fixed \(M>0\),
\[
\Pr\big(\mathcal K_n(M)\neq\emptyset\big)
\le
\Pr\big((\hat{\bm\beta}_n)_{\mathcal A^c}\neq \bm 0\big)
\le
\Pr(\hat{\mathcal A}_n\neq\mathcal A)\to 0,
\]
where the first inequality follows from \eqref{eq:event_equiv_K_revised}. Hence
\[
\limsup_{n\to\infty}\Pr\big(\mathcal K_n(M)\neq\emptyset\big)=0
\qquad\text{for every }M>0.
\]
Letting \(M\to\infty\) yields
\[
\lim_{M\to\infty}\limsup_{n\to\infty}
\Pr\big(\mathcal K_n(M)\neq\emptyset\big)=0.
\]

\medskip
\noindent\emph{If part.}
Assume now that
\[
\lim_{M\to\infty}\limsup_{n\to\infty}
\Pr\big(\mathcal K_n(M)\neq\emptyset\big)=0.
\]
Define
\[
L_n(M):=\left\{\|\hat{\bm\beta}_n-\bm\beta_0\|_\infty\le M/r_n\right\}.
\]
Since \(\hat{\bm\beta}_n-\bm\beta_0=O_{\Pr}(1/r_n)\), we have
\begin{equation}\label{eq:localization_iff_revised}
\lim_{M\to\infty}\limsup_{n\to\infty}\Pr\big(L_n(M)^c\big)=0.
\end{equation}

Next, let
\[
c_{\mathcal A}
:=
\begin{cases}
\min_{j\in\mathcal A}|\beta_{0j}|, & \mathcal A\neq\emptyset,\\[0.3em]
+\infty, & \mathcal A=\emptyset.
\end{cases}
\]
If \(\mathcal A\neq\emptyset\), then \(c_{\mathcal A}>0\). Fix \(M>0\). Since \(r_n\to\infty\), for all sufficiently large \(n\),
\[
\frac{M}{r_n}<\frac{c_{\mathcal A}}{2},
\]
where this condition is vacuous when \(\mathcal A=\emptyset\). Hence, on \(L_n(M)\), for every \(j\in\mathcal A\),
\[
|\hat\beta_{nj}|
\ge
|\beta_{0j}|-|\hat\beta_{nj}-\beta_{0j}|
\ge
c_{\mathcal A}-\frac{M}{r_n}
>
\frac{c_{\mathcal A}}{2}
>0.
\]
Thus, for all sufficiently large \(n\),
\begin{equation}\label{eq:active_preserved_local_revised}
L_n(M)\subseteq \{\mathcal A\subseteq \hat{\mathcal A}_n\}.
\end{equation}

We claim that, for all sufficiently large \(n\),
\begin{equation}\label{eq:oracle_fail_subset_revised}
\{\hat{\mathcal A}_n\neq \mathcal A\}
\subseteq
L_n(M)^c\cup\{\mathcal K_n(M)\neq\emptyset\}.
\end{equation}
Indeed, on \(L_n(M)\), relation \eqref{eq:active_preserved_local_revised} rules out false negatives. Therefore, if \(\hat{\mathcal A}_n\neq\mathcal A\) and \(L_n(M)\) occurs, then necessarily
\[
(\hat{\bm\beta}_n)_{\mathcal A^c}\neq \bm 0.
\]
Together with \(L_n(M)\), this is exactly the event \(\{\mathcal K_n(M)\neq\emptyset\}\) by \eqref{eq:event_equiv_K_revised}, proving \eqref{eq:oracle_fail_subset_revised}.

Taking probabilities in \eqref{eq:oracle_fail_subset_revised} yields
\[
\limsup_{n\to\infty}\Pr(\hat{\mathcal A}_n\neq \mathcal A)
\le
\limsup_{n\to\infty}\Pr\big(L_n(M)^c\big)
+
\limsup_{n\to\infty}\Pr\big(\mathcal K_n(M)\neq\emptyset\big).
\]
Finally, let \(M\to\infty\). The first term vanishes by \eqref{eq:localization_iff_revised}, and the second vanishes by assumption. Hence
\[
\limsup_{n\to\infty}\Pr(\hat{\mathcal A}_n\neq \mathcal A)=0,
\]
that is,
\[
\Pr(\hat{\mathcal A}_n=\mathcal A)\to 1.
\]
Thus Oracle property~\ref{oracle property}\ref{op1} holds.
\end{proof}

\subsection{Proofs of the results in Section~\ref{app:Ridgeless_differentiability}}

\begin{proof}[Proof of Proposition~\ref{prop:Ridgeless_regular_targets}]
Conditional on \(\bm x\), the log-likelihood under \(P_{\bm\beta}\) is
\[
\ell_{\bm\beta}(y,\bm x)
=
-\frac12\log(2\pi\sigma_0^2)
-\frac{1}{2\sigma_0^2}(y-\bm x'\bm\beta)^2 .
\]
Hence, for \(\bm h\in\R^p\), the score in direction \(\bm h\) at
\(P_{\bm\beta_0}\) is
\[
\left.\frac{\partial}{\partial t}\ell_{\bm\beta_0+t\bm h}(y,\bm x)\right|_{t=0}
=
\sigma_0^{-2}(y-\bm x'\bm\beta_0)\bm x'\bm h.
\]
Therefore,
\[
\E\left[
\left(
\left.\frac{\partial}{\partial t}\ell_{\bm\beta_0+t\bm h}(y,\bm x)\right|_{t=0}
\right)^2
\right]
=
\sigma_0^{-4}\E\left[\varepsilon^2(\bm x'\bm h)^2\right]
=
\sigma_0^{-2}\bm h'\bm Q_0\bm h,
\]
so the Fisher information matrix for \(\bm\beta\) is
\begin{equation}\label{eq:appendix_rls_information}
\bm I_0=\sigma_0^{-2}\bm Q_0.
\end{equation}

For part~(i), \citet[Thm.~2.1]{SM:van1991differentiable} implies that
\(\sqrt n\)-regular estimation of \(K(P_{\bm\beta})\) requires
differentiability of the target relative to the admissible path class.
In the present finite-dimensional parametric model, the specialization in
\citet[Sect.~5]{SM:van1991differentiable} reduces differentiability to the
null-space condition
\[
\Kernel(\bm I_0)\subseteq \Kernel(\dot\psi_0).
\]
By \eqref{eq:appendix_rls_information},
\(\Kernel(\bm I_0)=\Kernel(\bm Q_0)\), which proves
\eqref{eq:appendix_rls_nullspace_condition}.
If \(\bm Q_0\) is singular and \(\psi(\bm\beta)=\bm\beta\), then
\(\dot\psi_0=\bm I_p\), so \eqref{eq:appendix_rls_nullspace_condition}
fails.

For part~(ii), first note that
\[
\bm\beta-\bm P_0^\perp\bm\beta=\bm P_0\bm\beta\in\Kernel(\bm Q_0),
\qquad
\bm P_0^\perp\bm\beta\in\Range(\bm Q_0),
\]
so \(\bm P_0^\perp\bm\beta\in
(\bm\beta+\Kernel(\bm Q_0))\cap\Range(\bm Q_0)\).
If \(\bm\gamma\in
(\bm\beta+\Kernel(\bm Q_0))\cap\Range(\bm Q_0)\), then
\[
\bm\gamma-\bm P_0^\perp\bm\beta
=
(\bm\gamma-\bm\beta)+\bm P_0\bm\beta
\in\Kernel(\bm Q_0),
\]
and also \(\bm\gamma-\bm P_0^\perp\bm\beta\in\Range(\bm Q_0)\).
Hence
\[
\bm\gamma-\bm P_0^\perp\bm\beta
\in
\Kernel(\bm Q_0)\cap\Range(\bm Q_0)
=
\{\bm 0\},
\]
which proves \eqref{eq:appendix_rls_singleton}.

If \(P_{\bm\beta}=P_{\tilde{\bm\beta}}\), then
\(\bm x'(\bm\beta-\tilde{\bm\beta})=0\) almost surely.
Therefore
\[
(\bm\beta-\tilde{\bm\beta})'\bm Q_0(\bm\beta-\tilde{\bm\beta})
=
\E\left[(\bm x'(\bm\beta-\tilde{\bm\beta}))^2\right]
=
0,
\]
so \(\bm\beta-\tilde{\bm\beta}\in\Kernel(\bm Q_0)\).
Hence
\(\bm P_0^\perp\bm\beta=\bm P_0^\perp\tilde{\bm\beta}\), which shows that
\eqref{eq:appendix_rls_map} is well-defined on the quotient model.

Finally, for any path
\[
\bm\beta_t=\bm\beta_0+t\bm h+o(t),
\qquad t\to 0,
\]
we have
\[
K^{\mathrm{rls}}(P_{\bm\beta_t})
=
\bm P_0^\perp\bm\beta_t
=
\bm P_0^\perp\bm\beta_0
+
t\bm P_0^\perp\bm h
+
o(t),
\]
because \(\bm P_0^\perp\) is linear and continuous. This proves
\eqref{eq:appendix_rls_derivative}.

For part~(iii), let \(\bm\beta\in\R^p\).
Since \(L\bm\beta\in\Range(\bm Q_0)\), one has
\(\bm P_0^\perp L\bm\beta=L\bm\beta\).
Since \(L\bm\beta-\bm\beta\in\Kernel(\bm Q_0)\), one also has
\(\bm P_0^\perp(L\bm\beta-\bm\beta)=\bm 0\).
Therefore
\[
L\bm\beta
=
\bm P_0^\perp L\bm\beta
=
\bm P_0^\perp\bm\beta.
\]
Since this holds for every \(\bm\beta\in\R^p\), we obtain
\(L=\bm P_0^\perp\), proving \eqref{eq:appendix_rls_unique_projection}.
If \(\tilde{\bm\beta}_0\in\mathcal B_0\), then
\(\tilde{\bm\beta}_0-\bm\beta_0^{\mathrm{rls}}\in\Kernel(\bm Q_0)\) and
\(\bm\beta_0^{\mathrm{rls}}\in\Range(\bm Q_0)\), so
\(\bm P_0^\perp\tilde{\bm\beta}_0=\bm\beta_0^{\mathrm{rls}}\).
\end{proof}


\begin{thebibliography}{}

\bibitem [\protect \citeauthoryear {%
Andrews%
}{%
Andrews%
}{%
{\protect \APACyear {1987}}%
}]{%
andrews1987asymptotic}
\APACinsertmetastar {%
andrews1987asymptotic}%
\begin{APACrefauthors}%
Andrews, D\BPBI W.%
\end{APACrefauthors}%
\unskip\
\newblock
\APACrefYearMonthDay{1987}{}{}.
\newblock
{\BBOQ}\APACrefatitle {Asymptotic results for generalized {Wald} tests}
  {Asymptotic results for generalized {Wald} tests}.{\BBCQ}
\newblock
\APACjournalVolNumPages{Econometric Theory}{3}{3}{348--358}.
\PrintBackRefs{\CurrentBib}

\bibitem [\protect \citeauthoryear {%
Andrews%
\ \BBA {} Guggenberger%
}{%
Andrews%
\ \BBA {} Guggenberger%
}{%
{\protect \APACyear {2019}}%
}]{%
andrews2019identification}
\APACinsertmetastar {%
andrews2019identification}%
\begin{APACrefauthors}%
Andrews, D\BPBI W.%
\BCBT {}\ \BBA {} Guggenberger, P.%
\end{APACrefauthors}%
\unskip\
\newblock
\APACrefYearMonthDay{2019}{}{}.
\newblock
{\BBOQ}\APACrefatitle {Identification-and singularity-robust inference for
  moment condition models} {Identification-and singularity-robust inference for
  moment condition models}.{\BBCQ}
\newblock
\APACjournalVolNumPages{Quantitative Economics}{10}{4}{1703--1746}.
\PrintBackRefs{\CurrentBib}

\bibitem [\protect \citeauthoryear {%
Aswani%
, Bickel%
\BCBL {}\ \BBA {} Tomlin%
}{%
Aswani%
\ \protect \BOthers {.}}{%
{\protect \APACyear {2011}}%
}]{%
aswani2011regression}
\APACinsertmetastar {%
aswani2011regression}%
\begin{APACrefauthors}%
Aswani, A.%
, Bickel, P.%
\BCBL {}\ \BBA {} Tomlin, C.%
\end{APACrefauthors}%
\unskip\
\newblock
\APACrefYearMonthDay{2011}{}{}.
\newblock
{\BBOQ}\APACrefatitle {Regression on manifolds: Estimation of the exterior
  derivative} {Regression on manifolds: Estimation of the exterior
  derivative}.{\BBCQ}
\newblock
\APACjournalVolNumPages{The Annals of Statistics}{39}{1}{48--81}.
\PrintBackRefs{\CurrentBib}

\bibitem [\protect \citeauthoryear {%
Bauschke%
\ \BBA {} Combettes%
}{%
Bauschke%
\ \BBA {} Combettes%
}{%
{\protect \APACyear {2016}}%
}]{%
bauschke2016convex}
\APACinsertmetastar {%
bauschke2016convex}%
\begin{APACrefauthors}%
Bauschke, H\BPBI H.%
\BCBT {}\ \BBA {} Combettes, P\BPBI L.%
\end{APACrefauthors}%
\unskip\
\newblock
\APACrefYear{2016}.
\newblock
\APACrefbtitle {Convex analysis and monotone operator theory in {Hilbert}
  spaces} {Convex analysis and monotone operator theory in {Hilbert} spaces}\
  (\PrintOrdinal{2nd}\ \BEd).
\newblock
\APACaddressPublisher{}{Springer}.
\PrintBackRefs{\CurrentBib}

\bibitem [\protect \citeauthoryear {%
Caner%
}{%
Caner%
}{%
{\protect \APACyear {2008}}%
}]{%
caner2008nearly}
\APACinsertmetastar {%
caner2008nearly}%
\begin{APACrefauthors}%
Caner, M.%
\end{APACrefauthors}%
\unskip\
\newblock
\APACrefYearMonthDay{2008}{}{}.
\newblock
{\BBOQ}\APACrefatitle {Nearly-singular design in {GMM} and generalized
  empirical likelihood estimators} {Nearly-singular design in {GMM} and
  generalized empirical likelihood estimators}.{\BBCQ}
\newblock
\APACjournalVolNumPages{Journal of Econometrics}{144}{2}{511--523}.
\PrintBackRefs{\CurrentBib}

\bibitem [\protect \citeauthoryear {%
Caner%
}{%
Caner%
}{%
{\protect \APACyear {2009}}%
}]{%
caner2009lasso}
\APACinsertmetastar {%
caner2009lasso}%
\begin{APACrefauthors}%
Caner, M.%
\end{APACrefauthors}%
\unskip\
\newblock
\APACrefYearMonthDay{2009}{}{}.
\newblock
{\BBOQ}\APACrefatitle {Lasso-type {GMM} estimator} {Lasso-type {GMM}
  estimator}.{\BBCQ}
\newblock
\APACjournalVolNumPages{Econometric Theory}{25}{1}{270--290}.
\PrintBackRefs{\CurrentBib}

\bibitem [\protect \citeauthoryear {%
Charles%
\ \BBA {} Hurst%
}{%
Charles%
\ \BBA {} Hurst%
}{%
{\protect \APACyear {2003}}%
}]{%
charles2003correlation}
\APACinsertmetastar {%
charles2003correlation}%
\begin{APACrefauthors}%
Charles, K\BPBI K.%
\BCBT {}\ \BBA {} Hurst, E.%
\end{APACrefauthors}%
\unskip\
\newblock
\APACrefYearMonthDay{2003}{}{}.
\newblock
{\BBOQ}\APACrefatitle {The correlation of wealth across generations} {The
  correlation of wealth across generations}.{\BBCQ}
\newblock
\APACjournalVolNumPages{Journal of Political Economy}{111}{6}{1155--1182}.
\PrintBackRefs{\CurrentBib}

\bibitem [\protect \citeauthoryear {%
Fan%
\ \BBA {} Li%
}{%
Fan%
\ \BBA {} Li%
}{%
{\protect \APACyear {2001}}%
}]{%
fan2001variable}
\APACinsertmetastar {%
fan2001variable}%
\begin{APACrefauthors}%
Fan, J.%
\BCBT {}\ \BBA {} Li, R.%
\end{APACrefauthors}%
\unskip\
\newblock
\APACrefYearMonthDay{2001}{}{}.
\newblock
{\BBOQ}\APACrefatitle {Variable selection via nonconcave penalized likelihood
  and its oracle properties} {Variable selection via nonconcave penalized
  likelihood and its oracle properties}.{\BBCQ}
\newblock
\APACjournalVolNumPages{Journal of the American Statistical
  Association}{96}{456}{1348--1360}.
\PrintBackRefs{\CurrentBib}

\bibitem [\protect \citeauthoryear {%
Frank%
\ \BBA {} Friedman%
}{%
Frank%
\ \BBA {} Friedman%
}{%
{\protect \APACyear {1993}}%
}]{%
frank1993statistical}
\APACinsertmetastar {%
frank1993statistical}%
\begin{APACrefauthors}%
Frank, I\BPBI E.%
\BCBT {}\ \BBA {} Friedman, J\BPBI H.%
\end{APACrefauthors}%
\unskip\
\newblock
\APACrefYearMonthDay{1993}{}{}.
\newblock
{\BBOQ}\APACrefatitle {A statistical view of some chemometrics regression
  tools} {A statistical view of some chemometrics regression tools}.{\BBCQ}
\newblock
\APACjournalVolNumPages{Technometrics}{35}{2}{109--135}.
\PrintBackRefs{\CurrentBib}

\bibitem [\protect \citeauthoryear {%
Gabaix%
\ \BBA {} Ibragimov%
}{%
Gabaix%
\ \BBA {} Ibragimov%
}{%
{\protect \APACyear {2011}}%
}]{%
gabaix2011rank}
\APACinsertmetastar {%
gabaix2011rank}%
\begin{APACrefauthors}%
Gabaix, X.%
\BCBT {}\ \BBA {} Ibragimov, R.%
\end{APACrefauthors}%
\unskip\
\newblock
\APACrefYearMonthDay{2011}{}{}.
\newblock
{\BBOQ}\APACrefatitle {Rank-1/2: a simple way to improve the {OLS} estimation
  of tail exponents} {Rank-1/2: a simple way to improve the {OLS} estimation of
  tail exponents}.{\BBCQ}
\newblock
\APACjournalVolNumPages{Journal of Business \& Economic
  Statistics}{29}{1}{24--39}.
\PrintBackRefs{\CurrentBib}

\bibitem [\protect \citeauthoryear {%
Green%
}{%
Green%
}{%
{\protect \APACyear {1987}}%
}]{%
green1987penalized}
\APACinsertmetastar {%
green1987penalized}%
\begin{APACrefauthors}%
Green, P\BPBI J.%
\end{APACrefauthors}%
\unskip\
\newblock
\APACrefYearMonthDay{1987}{}{}.
\newblock
{\BBOQ}\APACrefatitle {Penalized likelihood for general semi-parametric
  regression models} {Penalized likelihood for general semi-parametric
  regression models}.{\BBCQ}
\newblock
\APACjournalVolNumPages{International Statistical Review/Revue Internationale
  de Statistique}{}{}{245--259}.
\PrintBackRefs{\CurrentBib}

\bibitem [\protect \citeauthoryear {%
Hansen%
}{%
Hansen%
}{%
{\protect \APACyear {1982}}%
}]{%
hansen1982large}
\APACinsertmetastar {%
hansen1982large}%
\begin{APACrefauthors}%
Hansen, L\BPBI P.%
\end{APACrefauthors}%
\unskip\
\newblock
\APACrefYearMonthDay{1982}{}{}.
\newblock
{\BBOQ}\APACrefatitle {Large sample properties of generalized method of moments
  estimators} {Large sample properties of generalized method of moments
  estimators}.{\BBCQ}
\newblock
\APACjournalVolNumPages{Econometrica}{50}{4}{1029--1054}.
\newblock
\begin{APACrefDOI} \doi{10.2307/1912775} \end{APACrefDOI}
\PrintBackRefs{\CurrentBib}

\bibitem [\protect \citeauthoryear {%
Hastie%
, Montanari%
, Rosset%
\BCBL {}\ \BBA {} Tibshirani%
}{%
Hastie%
\ \protect \BOthers {.}}{%
{\protect \APACyear {2022}}%
}]{%
hastie2022surprises}
\APACinsertmetastar {%
hastie2022surprises}%
\begin{APACrefauthors}%
Hastie, T.%
, Montanari, A.%
, Rosset, S.%
\BCBL {}\ \BBA {} Tibshirani, R\BPBI J.%
\end{APACrefauthors}%
\unskip\
\newblock
\APACrefYearMonthDay{2022}{}{}.
\newblock
{\BBOQ}\APACrefatitle {Surprises in high-dimensional ridgeless least squares
  interpolation} {Surprises in high-dimensional ridgeless least squares
  interpolation}.{\BBCQ}
\newblock
\APACjournalVolNumPages{Annals of statistics}{50}{2}{949}.
\PrintBackRefs{\CurrentBib}

\bibitem [\protect \citeauthoryear {%
Hoerl%
\ \BBA {} Kennard%
}{%
Hoerl%
\ \BBA {} Kennard%
}{%
{\protect \APACyear {1970}}%
}]{%
hoerl1970ridge}
\APACinsertmetastar {%
hoerl1970ridge}%
\begin{APACrefauthors}%
Hoerl, A\BPBI E.%
\BCBT {}\ \BBA {} Kennard, R\BPBI W.%
\end{APACrefauthors}%
\unskip\
\newblock
\APACrefYearMonthDay{1970}{}{}.
\newblock
{\BBOQ}\APACrefatitle {Ridge regression: applications to nonorthogonal
  problems} {Ridge regression: applications to nonorthogonal problems}.{\BBCQ}
\newblock
\APACjournalVolNumPages{Technometrics}{12}{1}{69--82}.
\PrintBackRefs{\CurrentBib}

\bibitem [\protect \citeauthoryear {%
Knight%
}{%
Knight%
}{%
{\protect \APACyear {2008}}%
}]{%
knight2008shrinkage}
\APACinsertmetastar {%
knight2008shrinkage}%
\begin{APACrefauthors}%
Knight, K.%
\end{APACrefauthors}%
\unskip\
\newblock
\APACrefYearMonthDay{2008}{}{}.
\newblock
{\BBOQ}\APACrefatitle {Shrinkage estimation for nearly singular designs}
  {Shrinkage estimation for nearly singular designs}.{\BBCQ}
\newblock
\APACjournalVolNumPages{Econometric Theory}{24}{2}{323--337}.
\PrintBackRefs{\CurrentBib}

\bibitem [\protect \citeauthoryear {%
Knight%
\ \BBA {} Fu%
}{%
Knight%
\ \BBA {} Fu%
}{%
{\protect \APACyear {2000}}%
}]{%
knight2000asymptotics}
\APACinsertmetastar {%
knight2000asymptotics}%
\begin{APACrefauthors}%
Knight, K.%
\BCBT {}\ \BBA {} Fu, W.%
\end{APACrefauthors}%
\unskip\
\newblock
\APACrefYearMonthDay{2000}{}{}.
\newblock
{\BBOQ}\APACrefatitle {Asymptotics for lasso-type estimators} {Asymptotics for
  lasso-type estimators}.{\BBCQ}
\newblock
\APACjournalVolNumPages{The Annals of Statistics}{28}{5}{1356--1378}.
\PrintBackRefs{\CurrentBib}

\bibitem [\protect \citeauthoryear {%
Liao%
}{%
Liao%
}{%
{\protect \APACyear {2013}}%
}]{%
liao2013adaptive}
\APACinsertmetastar {%
liao2013adaptive}%
\begin{APACrefauthors}%
Liao, Z.%
\end{APACrefauthors}%
\unskip\
\newblock
\APACrefYearMonthDay{2013}{}{}.
\newblock
{\BBOQ}\APACrefatitle {Adaptive {GMM} shrinkage estimation with consistent
  moment selection} {Adaptive {GMM} shrinkage estimation with consistent moment
  selection}.{\BBCQ}
\newblock
\APACjournalVolNumPages{Econometric Theory}{29}{5}{857--904}.
\PrintBackRefs{\CurrentBib}

\bibitem [\protect \citeauthoryear {%
L{\"u}tkepohl%
\ \BBA {} Burda%
}{%
L{\"u}tkepohl%
\ \BBA {} Burda%
}{%
{\protect \APACyear {1997}}%
}]{%
lutkepohl1997modified}
\APACinsertmetastar {%
lutkepohl1997modified}%
\begin{APACrefauthors}%
L{\"u}tkepohl, H.%
\BCBT {}\ \BBA {} Burda, M\BPBI M.%
\end{APACrefauthors}%
\unskip\
\newblock
\APACrefYearMonthDay{1997}{}{}.
\newblock
{\BBOQ}\APACrefatitle {Modified {Wald} tests under nonregular conditions}
  {Modified {Wald} tests under nonregular conditions}.{\BBCQ}
\newblock
\APACjournalVolNumPages{Journal of Econometrics}{78}{2}{315--332}.
\PrintBackRefs{\CurrentBib}

\bibitem [\protect \citeauthoryear {%
Phillips%
}{%
Phillips%
}{%
{\protect \APACyear {2001}}%
}]{%
phillips2001regression}
\APACinsertmetastar {%
phillips2001regression}%
\begin{APACrefauthors}%
Phillips, P\BPBI C.%
\end{APACrefauthors}%
\unskip\
\newblock
\APACrefYearMonthDay{2001}{}{}.
\newblock
\APACrefbtitle {Regression with slowly varying regressors} {Regression with
  slowly varying regressors}\ \APACbVolEdTR {}{Cowles Foundation Discussion
  Paper\ \BNUM\ 1310}.
\newblock
\APACaddressInstitution{}{Cowles Foundation for Research in Economics, Yale
  University}.
\newblock
\begin{APACrefURL}
  \url{https://elischolar.library.yale.edu/cowles-discussion-paper-series/1569/}
  \end{APACrefURL}
\PrintBackRefs{\CurrentBib}

\bibitem [\protect \citeauthoryear {%
Phillips%
}{%
Phillips%
}{%
{\protect \APACyear {2016}}%
}]{%
phillips2016inference}
\APACinsertmetastar {%
phillips2016inference}%
\begin{APACrefauthors}%
Phillips, P\BPBI C.%
\end{APACrefauthors}%
\unskip\
\newblock
\APACrefYearMonthDay{2016}{}{}.
\newblock
{\BBOQ}\APACrefatitle {Inference in near-singular regression} {Inference in
  near-singular regression}.{\BBCQ}
\newblock
\APACjournalVolNumPages{Advances in Econometrics}{36}{}{461--486}.
\newblock
\begin{APACrefDOI} \doi{10.1108/S0731-905320160000036022} \end{APACrefDOI}
\PrintBackRefs{\CurrentBib}

\bibitem [\protect \citeauthoryear {%
Polson%
, Scott%
\BCBL {}\ \BBA {} Willard%
}{%
Polson%
\ \protect \BOthers {.}}{%
{\protect \APACyear {2015}}%
}]{%
polson2015proximal}
\APACinsertmetastar {%
polson2015proximal}%
\begin{APACrefauthors}%
Polson, N\BPBI G.%
, Scott, J\BPBI G.%
\BCBL {}\ \BBA {} Willard, B\BPBI T.%
\end{APACrefauthors}%
\unskip\
\newblock
\APACrefYearMonthDay{2015}{}{}.
\newblock
{\BBOQ}\APACrefatitle {Proximal algorithms in statistics and machine learning}
  {Proximal algorithms in statistics and machine learning}.{\BBCQ}
\newblock
\APACjournalVolNumPages{Statistical Science}{30}{4}{559--581}.
\PrintBackRefs{\CurrentBib}

\bibitem [\protect \citeauthoryear {%
Puri%
, Russell%
\BCBL {}\ \BBA {} Mathew%
}{%
Puri%
\ \protect \BOthers {.}}{%
{\protect \APACyear {1984}}%
}]{%
Madanetal1984}
\APACinsertmetastar {%
Madanetal1984}%
\begin{APACrefauthors}%
Puri, M\BPBI L.%
, Russell, C\BPBI T.%
\BCBL {}\ \BBA {} Mathew, T.%
\end{APACrefauthors}%
\unskip\
\newblock
\APACrefYearMonthDay{1984}{}{}.
\newblock
{\BBOQ}\APACrefatitle {Convergence of generalized inverses with applications to
  asymptotic hypothesis testing} {Convergence of generalized inverses with
  applications to asymptotic hypothesis testing}.{\BBCQ}
\newblock
\APACjournalVolNumPages{Sankhyā: The Indian Journal of Statistics, Series A
  (1961-2002)}{46}{2}{277--286}.
\PrintBackRefs{\CurrentBib}

\bibitem [\protect \citeauthoryear {%
Rockafellar%
\ \BBA {} Wets%
}{%
Rockafellar%
\ \BBA {} Wets%
}{%
{\protect \APACyear {2009}}%
}]{%
rockafellarWets2009}
\APACinsertmetastar {%
rockafellarWets2009}%
\begin{APACrefauthors}%
Rockafellar, R\BPBI T.%
\BCBT {}\ \BBA {} Wets, R\BPBI J\BHBI B.%
\end{APACrefauthors}%
\unskip\
\newblock
\APACrefYear{2009}.
\newblock
\APACrefbtitle {Variational analysis} {Variational analysis}\ (\BVOL~317).
\newblock
\APACaddressPublisher{}{Springer}.
\PrintBackRefs{\CurrentBib}

\bibitem [\protect \citeauthoryear {%
Salinetti%
\ \BBA {} Wets%
}{%
Salinetti%
\ \BBA {} Wets%
}{%
{\protect \APACyear {1981}}%
}]{%
salinetti1981convergence}
\APACinsertmetastar {%
salinetti1981convergence}%
\begin{APACrefauthors}%
Salinetti, G.%
\BCBT {}\ \BBA {} Wets, R\BPBI J\BHBI B.%
\end{APACrefauthors}%
\unskip\
\newblock
\APACrefYearMonthDay{1981}{}{}.
\newblock
{\BBOQ}\APACrefatitle {On the convergence of closed-valued measurable
  multifunctions} {On the convergence of closed-valued measurable
  multifunctions}.{\BBCQ}
\newblock
\APACjournalVolNumPages{Transactions of the American Mathematical
  Society}{266}{1}{275--289}.
\PrintBackRefs{\CurrentBib}

\bibitem [\protect \citeauthoryear {%
Salinetti%
\ \BBA {} Wets%
}{%
Salinetti%
\ \BBA {} Wets%
}{%
{\protect \APACyear {1986}}%
}]{%
salinetti1986convergence}
\APACinsertmetastar {%
salinetti1986convergence}%
\begin{APACrefauthors}%
Salinetti, G.%
\BCBT {}\ \BBA {} Wets, R\BPBI J\BHBI B.%
\end{APACrefauthors}%
\unskip\
\newblock
\APACrefYearMonthDay{1986}{}{}.
\newblock
{\BBOQ}\APACrefatitle {On the convergence in distribution of measurable
  multifunctions (random sets) normal integrands, stochastic processes and
  stochastic infima} {On the convergence in distribution of measurable
  multifunctions (random sets) normal integrands, stochastic processes and
  stochastic infima}.{\BBCQ}
\newblock
\APACjournalVolNumPages{Mathematics of Operations Research}{11}{3}{385--419}.
\PrintBackRefs{\CurrentBib}

\bibitem [\protect \citeauthoryear {%
Solon%
}{%
Solon%
}{%
{\protect \APACyear {1992}}%
}]{%
solon1992intergenerational}
\APACinsertmetastar {%
solon1992intergenerational}%
\begin{APACrefauthors}%
Solon, G.%
\end{APACrefauthors}%
\unskip\
\newblock
\APACrefYearMonthDay{1992}{}{}.
\newblock
{\BBOQ}\APACrefatitle {Intergenerational income mobility in the {United
  States}} {Intergenerational income mobility in the {United States}}.{\BBCQ}
\newblock
\APACjournalVolNumPages{The American Economic Review}{82}{3}{393--408}.
\PrintBackRefs{\CurrentBib}

\bibitem [\protect \citeauthoryear {%
Tibshirani%
}{%
Tibshirani%
}{%
{\protect \APACyear {1996}}%
}]{%
tibshirani1996regression}
\APACinsertmetastar {%
tibshirani1996regression}%
\begin{APACrefauthors}%
Tibshirani, R.%
\end{APACrefauthors}%
\unskip\
\newblock
\APACrefYearMonthDay{1996}{}{}.
\newblock
{\BBOQ}\APACrefatitle {Regression shrinkage and selection via the lasso}
  {Regression shrinkage and selection via the lasso}.{\BBCQ}
\newblock
\APACjournalVolNumPages{Journal of the Royal Statistical Society: Series B
  (Methodological)}{58}{1}{267--288}.
\PrintBackRefs{\CurrentBib}

\bibitem [\protect \citeauthoryear {%
van~der Vaart%
}{%
van~der Vaart%
}{%
{\protect \APACyear {1991}}%
}]{%
van1991differentiable}
\APACinsertmetastar {%
van1991differentiable}%
\begin{APACrefauthors}%
van~der Vaart, A.%
\end{APACrefauthors}%
\unskip\
\newblock
\APACrefYearMonthDay{1991}{}{}.
\newblock
{\BBOQ}\APACrefatitle {On differentiable functionals} {On differentiable
  functionals}.{\BBCQ}
\newblock
\APACjournalVolNumPages{The Annals of Statistics}{19}{1}{178--204}.
\PrintBackRefs{\CurrentBib}

\bibitem [\protect \citeauthoryear {%
Yuan%
\ \BBA {} Lin%
}{%
Yuan%
\ \BBA {} Lin%
}{%
{\protect \APACyear {2006}}%
}]{%
yuan2006model}
\APACinsertmetastar {%
yuan2006model}%
\begin{APACrefauthors}%
Yuan, M.%
\BCBT {}\ \BBA {} Lin, Y.%
\end{APACrefauthors}%
\unskip\
\newblock
\APACrefYearMonthDay{2006}{}{}.
\newblock
{\BBOQ}\APACrefatitle {Model selection and estimation in regression with
  grouped variables} {Model selection and estimation in regression with grouped
  variables}.{\BBCQ}
\newblock
\APACjournalVolNumPages{Journal of the Royal Statistical Society: Series B
  (Statistical Methodology)}{68}{1}{49--67}.
\PrintBackRefs{\CurrentBib}

\end{thebibliography}

\begin{thebibliography}{}

\bibitem [\protect \citeauthoryear {%
Andrews%
}{%
Andrews%
}{%
{\protect \APACyear {1987}}%
}]{%
SM:andrews1987asymptotic}
\APACinsertmetastar {%
SM:andrews1987asymptotic}%
\begin{APACrefauthors}%
Andrews, D\BPBI W.%
\end{APACrefauthors}%
\unskip\
\newblock
\APACrefYearMonthDay{1987}{}{}.
\newblock
{\BBOQ}\APACrefatitle {Asymptotic results for generalized {Wald} tests}
  {Asymptotic results for generalized {Wald} tests}.{\BBCQ}
\newblock
\APACjournalVolNumPages{Econometric Theory}{3}{3}{348--358}.
\PrintBackRefs{\CurrentBib}

\bibitem [\protect \citeauthoryear {%
Attouch%
}{%
Attouch%
}{%
{\protect \APACyear {1984}}%
}]{%
SM:attouch1984variational}
\APACinsertmetastar {%
SM:attouch1984variational}%
\begin{APACrefauthors}%
Attouch, H.%
\end{APACrefauthors}%
\unskip\
\newblock
\APACrefYear{1984}.
\newblock
\APACrefbtitle {Variational convergence for functions and operators}
  {Variational convergence for functions and operators}\ (\BVOL~1).
\newblock
\APACaddressPublisher{}{Pitman Advanced Publishing Program}.
\PrintBackRefs{\CurrentBib}

\bibitem [\protect \citeauthoryear {%
Bauschke%
\ \BBA {} Combettes%
}{%
Bauschke%
\ \BBA {} Combettes%
}{%
{\protect \APACyear {2016}}%
}]{%
SM:bauschke2016convex}
\APACinsertmetastar {%
SM:bauschke2016convex}%
\begin{APACrefauthors}%
Bauschke, H\BPBI H.%
\BCBT {}\ \BBA {} Combettes, P\BPBI L.%
\end{APACrefauthors}%
\unskip\
\newblock
\APACrefYear{2016}.
\newblock
\APACrefbtitle {Convex analysis and monotone operator theory in {Hilbert}
  spaces} {Convex analysis and monotone operator theory in {Hilbert} spaces}\
  (\PrintOrdinal{2nd}\ \BEd).
\newblock
\APACaddressPublisher{}{Springer}.
\PrintBackRefs{\CurrentBib}

\bibitem [\protect \citeauthoryear {%
Billingsley%
}{%
Billingsley%
}{%
{\protect \APACyear {2013}}%
}]{%
SM:billingsley2013convergence}
\APACinsertmetastar {%
SM:billingsley2013convergence}%
\begin{APACrefauthors}%
Billingsley, P.%
\end{APACrefauthors}%
\unskip\
\newblock
\APACrefYear{2013}.
\newblock
\APACrefbtitle {Convergence of probability measures} {Convergence of
  probability measures}.
\newblock
\APACaddressPublisher{}{Wiley}.
\PrintBackRefs{\CurrentBib}

\bibitem [\protect \citeauthoryear {%
Hiriart-Urruty%
\ \BBA {} Lemar{\'e}chal%
}{%
Hiriart-Urruty%
\ \BBA {} Lemar{\'e}chal%
}{%
{\protect \APACyear {2004}}%
}]{%
SM:hiriart2004fundamentals}
\APACinsertmetastar {%
SM:hiriart2004fundamentals}%
\begin{APACrefauthors}%
Hiriart-Urruty, J\BHBI B.%
\BCBT {}\ \BBA {} Lemar{\'e}chal, C.%
\end{APACrefauthors}%
\unskip\
\newblock
\APACrefYear{2004}.
\newblock
\APACrefbtitle {Fundamentals of convex analysis} {Fundamentals of convex
  analysis}.
\newblock
\APACaddressPublisher{}{Springer}.
\PrintBackRefs{\CurrentBib}

\bibitem [\protect \citeauthoryear {%
Knight%
}{%
Knight%
}{%
{\protect \APACyear {2008}}%
}]{%
SM:knight2008shrinkage}
\APACinsertmetastar {%
SM:knight2008shrinkage}%
\begin{APACrefauthors}%
Knight, K.%
\end{APACrefauthors}%
\unskip\
\newblock
\APACrefYearMonthDay{2008}{}{}.
\newblock
{\BBOQ}\APACrefatitle {Shrinkage estimation for nearly singular designs}
  {Shrinkage estimation for nearly singular designs}.{\BBCQ}
\newblock
\APACjournalVolNumPages{Econometric Theory}{24}{2}{323--337}.
\PrintBackRefs{\CurrentBib}

\bibitem [\protect \citeauthoryear {%
Knight%
\ \BBA {} Fu%
}{%
Knight%
\ \BBA {} Fu%
}{%
{\protect \APACyear {2000}}%
}]{%
SM:knight2000asymptotics}
\APACinsertmetastar {%
SM:knight2000asymptotics}%
\begin{APACrefauthors}%
Knight, K.%
\BCBT {}\ \BBA {} Fu, W.%
\end{APACrefauthors}%
\unskip\
\newblock
\APACrefYearMonthDay{2000}{}{}.
\newblock
{\BBOQ}\APACrefatitle {Asymptotics for lasso-type estimators} {Asymptotics for
  lasso-type estimators}.{\BBCQ}
\newblock
\APACjournalVolNumPages{The Annals of Statistics}{28}{5}{1356--1378}.
\PrintBackRefs{\CurrentBib}

\bibitem [\protect \citeauthoryear {%
Lachout%
}{%
Lachout%
}{%
{\protect \APACyear {2006}}%
}]{%
SM:lachout2006epi}
\APACinsertmetastar {%
SM:lachout2006epi}%
\begin{APACrefauthors}%
Lachout, P.%
\end{APACrefauthors}%
\unskip\
\newblock
\APACrefYearMonthDay{2006}{}{}.
\newblock
{\BBOQ}\APACrefatitle {Epi-convergence almost surely, in probability and in
  distribution} {Epi-convergence almost surely, in probability and in
  distribution}.{\BBCQ}
\newblock
\APACjournalVolNumPages{Annals of Operations Research}{142}{1}{187--214}.
\PrintBackRefs{\CurrentBib}

\bibitem [\protect \citeauthoryear {%
Pflug%
}{%
Pflug%
}{%
{\protect \APACyear {1995}}%
}]{%
SM:pflug1995asymptotic}
\APACinsertmetastar {%
SM:pflug1995asymptotic}%
\begin{APACrefauthors}%
Pflug, G\BPBI C.%
\end{APACrefauthors}%
\unskip\
\newblock
\APACrefYearMonthDay{1995}{}{}.
\newblock
{\BBOQ}\APACrefatitle {Asymptotic stochastic programs} {Asymptotic stochastic
  programs}.{\BBCQ}
\newblock
\APACjournalVolNumPages{Mathematics of Operations Research}{20}{4}{769--789}.
\PrintBackRefs{\CurrentBib}

\bibitem [\protect \citeauthoryear {%
Rockafellar%
\ \BBA {} Wets%
}{%
Rockafellar%
\ \BBA {} Wets%
}{%
{\protect \APACyear {2009}}%
}]{%
SM:rockafellarWets2009}
\APACinsertmetastar {%
SM:rockafellarWets2009}%
\begin{APACrefauthors}%
Rockafellar, R\BPBI T.%
\BCBT {}\ \BBA {} Wets, R\BPBI J\BHBI B.%
\end{APACrefauthors}%
\unskip\
\newblock
\APACrefYear{2009}.
\newblock
\APACrefbtitle {Variational analysis} {Variational analysis}\ (\BVOL~317).
\newblock
\APACaddressPublisher{}{Springer}.
\PrintBackRefs{\CurrentBib}

\bibitem [\protect \citeauthoryear {%
Salinetti%
\ \BBA {} Wets%
}{%
Salinetti%
\ \BBA {} Wets%
}{%
{\protect \APACyear {1981}}%
}]{%
SM:salinetti1981convergence}
\APACinsertmetastar {%
SM:salinetti1981convergence}%
\begin{APACrefauthors}%
Salinetti, G.%
\BCBT {}\ \BBA {} Wets, R\BPBI J\BHBI B.%
\end{APACrefauthors}%
\unskip\
\newblock
\APACrefYearMonthDay{1981}{}{}.
\newblock
{\BBOQ}\APACrefatitle {On the convergence of closed-valued measurable
  multifunctions} {On the convergence of closed-valued measurable
  multifunctions}.{\BBCQ}
\newblock
\APACjournalVolNumPages{Transactions of the American Mathematical
  Society}{266}{1}{275--289}.
\PrintBackRefs{\CurrentBib}

\bibitem [\protect \citeauthoryear {%
Salinetti%
\ \BBA {} Wets%
}{%
Salinetti%
\ \BBA {} Wets%
}{%
{\protect \APACyear {1986}}%
}]{%
SM:salinetti1986convergence}
\APACinsertmetastar {%
SM:salinetti1986convergence}%
\begin{APACrefauthors}%
Salinetti, G.%
\BCBT {}\ \BBA {} Wets, R\BPBI J\BHBI B.%
\end{APACrefauthors}%
\unskip\
\newblock
\APACrefYearMonthDay{1986}{}{}.
\newblock
{\BBOQ}\APACrefatitle {On the convergence in distribution of measurable
  multifunctions (random sets) normal integrands, stochastic processes and
  stochastic infima} {On the convergence in distribution of measurable
  multifunctions (random sets) normal integrands, stochastic processes and
  stochastic infima}.{\BBCQ}
\newblock
\APACjournalVolNumPages{Mathematics of Operations Research}{11}{3}{385--419}.
\PrintBackRefs{\CurrentBib}

\bibitem [\protect \citeauthoryear {%
Tibshirani%
}{%
Tibshirani%
}{%
{\protect \APACyear {1996}}%
}]{%
SM:tibshirani1996regression}
\APACinsertmetastar {%
SM:tibshirani1996regression}%
\begin{APACrefauthors}%
Tibshirani, R.%
\end{APACrefauthors}%
\unskip\
\newblock
\APACrefYearMonthDay{1996}{}{}.
\newblock
{\BBOQ}\APACrefatitle {Regression shrinkage and selection via the lasso}
  {Regression shrinkage and selection via the lasso}.{\BBCQ}
\newblock
\APACjournalVolNumPages{Journal of the Royal Statistical Society: Series B
  (Methodological)}{58}{1}{267--288}.
\PrintBackRefs{\CurrentBib}

\bibitem [\protect \citeauthoryear {%
van~der Vaart%
}{%
van~der Vaart%
}{%
{\protect \APACyear {1991}}%
}]{%
SM:van1991differentiable}
\APACinsertmetastar {%
SM:van1991differentiable}%
\begin{APACrefauthors}%
van~der Vaart, A.%
\end{APACrefauthors}%
\unskip\
\newblock
\APACrefYearMonthDay{1991}{}{}.
\newblock
{\BBOQ}\APACrefatitle {On differentiable functionals} {On differentiable
  functionals}.{\BBCQ}
\newblock
\APACjournalVolNumPages{The Annals of Statistics}{19}{1}{178--204}.
\PrintBackRefs{\CurrentBib}

\bibitem [\protect \citeauthoryear {%
Vogel%
\ \BBA {} Lachout%
}{%
Vogel%
\ \BBA {} Lachout%
}{%
{\protect \APACyear {2003}}%
}]{%
SM:vogel2003continuous}
\APACinsertmetastar {%
SM:vogel2003continuous}%
\begin{APACrefauthors}%
Vogel, S.%
\BCBT {}\ \BBA {} Lachout, P.%
\end{APACrefauthors}%
\unskip\
\newblock
\APACrefYearMonthDay{2003}{}{}.
\newblock
{\BBOQ}\APACrefatitle {On continuous convergence and epi-convergence of random
  functions. Part {I}: Theory and relations} {On continuous convergence and
  epi-convergence of random functions. part {I}: Theory and relations}.{\BBCQ}
\newblock
\APACjournalVolNumPages{Kybernetika}{39}{1}{75--98}.
\PrintBackRefs{\CurrentBib}

\end{thebibliography}
\end{document}